\documentclass[10pt,leqno]{article}
\usepackage[utf8]{inputenc}
\usepackage{tikz}

\usepackage[official]{eurosym}
\usepackage{fullpage}
\usepackage{amsmath,amsthm}
\usepackage{amssymb}
\usepackage{epic,eepic}
\usepackage{cases}

\def\Mfor{\mathbb{M}[\vec{U}]}

\def\min{{\rm min}}
\def\max{{\rm max}}
\def\sup{{\rm sup}}
\def\otp{{\rm otp}}
\def\Lim{{\rm Lim}}

\def\succ{{\rm Succ}}

\def\llvdash{{\|\hskip-2pt \raise 3pt\hbox{\vrule
height 0.25pt width 0.4cm}}}

\def\l{{\langle}}
\def\r{{\rangle}}

\def\calP{\mathcal P}

\def\oa{{\overline A^{\,\lower 7pt_{\hbox{$\scriptstyle\bet$}}
\hbox{$\scriptstyle 0\tau$}}}}

\def\bet{\beta}

\def\llvdash{{\|\hskip-2pt \raise 3pt\hbox{\vrule height 0.25pt
width 0.4cm}}}

\newtheorem{theorem}{Theorem}[section]
\newtheorem{lemma}[theorem]{Lemma}
\newtheorem{corollary}[theorem]{Corollary}
\newtheorem{proposition}[theorem]{Proposition}
\newtheorem*{claim*}{Claim}
\newtheorem{definition}[theorem]{Definition}

\newtheorem{question}[theorem]{Question}
\newtheorem{conjecture}[theorem]{Conjecture}
\theoremstyle{remark}
\newtheorem{remark}[theorem]{Remark}

\newtheorem{example}[theorem]{Example}

\newcommand{\pr}{\medskip\noindent\textit{Proof}. }

\newcommand{\lusim}[1]{\smash{\underset{\raisebox{1.2pt}[0cm][0cm]{$\sim$}}
{{#1}}}}
\def\Lev{{\rm Lev}}

\def\otp{{\rm otp}}

\def\crit{{\rm crit}}

\def\llvdash{{\|\hskip-2pt \raise 3pt\hbox{\vrule
height 0.25pt width 0.15cm}}}

\def\Vdashbks{\hbox{$\Vdash\!\!\!\!{\raise2pt\hbox
{$\scriptscriptstyle\backslash$}}$}}

\title{Intermediate models of Magidor-Radin forcing-Part II}
\author{ Tom Benhamou and Moti Gitik\footnote{ The work of the second author was partially supported by ISF grant No.1216/18.}}
\date{\today}
\begin{document}

\maketitle
\begin{abstract}
    We continue the work done in \cite{PrikryCaseGitikKanKoe},\cite{TomMoti},\cite{partOne}. We prove that for every set of ordinals $A$ in a Magidor-Radin generic extension using a coherent sequence such that $o^{\vec{U}}(\kappa)<\kappa^+$, there is $C'\subseteq C_G$, such that $V[A]=V[C']$.
    Also we prove that the supremum of a fresh set in a Prikry, tree Prikry, Magidor, Radin-Magidor and Radin forcing, changes cofinality to $\omega$.
\end{abstract}
\section{Introduction}

A basic fact about the Cohen and Random forcings is that every subforcing of the Cohen (Random) forcing is equivalent to it.
Kanovey, Koepke and the second author showed in \cite{PrikryCaseGitikKanKoe} that the same is true for the standard Prikry forcing.
The result was generalized  to the Magidor forcing in \cite{TomMoti}. This was pushed further to versions of the Magidor-Radin forcing with $o^{\vec{U}}(\kappa)< \kappa$, in \cite{partOne}. 
The result for $o^{\vec{U}}(\kappa)<\kappa$, splits into two parts. The first is to prove that for every $V$-generic filter $G$, for the Magidor-Radin forcing, and any set of ordinals $A\in V[G]$, there is a subsequence of the generic club $C\subseteq C_G$ such that $V[A]=V[C]$. Thus, in order to analyse the intermediate models of $V[G]$, it suffices to study models of the form $V[C]$, where $C\subseteq C_G$. The second part is to show that each model of the form $V[C]$ is a $V$-generic extension for a Magidor-Radin-like forcing.

The main purpose of the present paper is to study  sets in generic extension of the version of Magidor-Radin forcing for $o^{\vec{U}}(\kappa)<\kappa^+$. It turns out that the first statement holds and every set  in the extension is equivalent to a subsequence of a generic Magidor-Radin sequence.
There are considerable additional difficulties here and new ideas are used to overcome them.
However, we do not give here a classification for models of the form $V[C]$.

The major difference between the case $o^{\vec{U}}(\kappa)<\kappa$ and $o^{\vec{U}}(\kappa)\geq\kappa$, is that we cannot split $\Mfor$ to the part below $o^{\vec{U}}(\kappa)$ and above it. As proven in \cite{partOne}, this decomposition provided the ability to run over all possible extension types. In terms of $C_G$ this means that we cannot split $C_G$ below $\kappa$ in a way that will determine what are the measures used in the construction of $C_G$. The classical example for such a sequence is
 $$C_G(0),C_G(C_G(0)),C_G(C_G(C_G(0))),...$$
 in which every element in the sequence is taken from a measure which depends on the previous element in the sequence. This example suggests that some sort of tree construction is needed in order to refer to such sequences in the ground model.
 
 In context of \cite{partOne} and \cite{TomMoti}, we are working by induction on $\kappa$. Formally we prove the following inductive step:
 \begin{theorem}\label{MainResaultParttwo}
 Let $\vec{U}$ be a coherent sequence with maximal measurable $\kappa$, such that $o^{\vec{U}}(\kappa)<\kappa^+$. Assume the inductive hypothesis:
 \begin{center}
 $(IH)$ \ \ \ For every $\delta<\kappa$, any coherent sequence $\vec{W}$ with maximal measurable $\delta$ and any set  of ordinals
 
 \ \ \ \ \ $A\in V[H]$ for  $H\subseteq\mathbb{M}[\vec{W}]$, there is  $C\subseteq C_H$, such that $V[A]=V[C]$.
 \end{center}
 Then for every $V$-generic filter $G\subseteq\Mfor$ and any set of ordinals $A\in V[G]$, there is $C\subseteq C_G$ such that $V[A]=V[C]$.
 \end{theorem}
 As a corollary of this, we obtain the main result of this paper:
 \begin{theorem}\label{MainResaultPartwo}
 Let $\vec{U}$ be a coherent sequence such that $o^{\vec{U}}(\kappa)<\kappa^+$. Then for every $V$-generic filter $G\subseteq\Mfor$, such that $\forall\alpha\in C_G. o^{\vec{U}}(\alpha)<\alpha^+$ and every set of ordinals $A\in V[G]$, there is $C\subseteq C_G$ such that $V[A]=V[C]$.
 \end{theorem}
Since every intermediate $ZFC$ model $V\subseteq M\subseteq V[G]$ is of the form $M=V[A]$ for some set of ordinals $A$ (See for example \cite[Corollary 15.42, Lemma 15.43]{Jech2003}), we conclude that every such $M$ is of the form $M=V[C]$ for some $C\subseteq C_G$. 
In this paper, the models $V[A]$ considered are always $ZFC$ models as $A$ would be a set of ordinals or can be coded as a set of ordinals using a function in $V$.\footnote{ For example if $A\subseteq V$ or if $A$ is a \textit{sequence} of sets of ordinals.} In any case, $V[A]$ \cite[Lemma 15.43]{Jech2003} is the minimal $ZFC$ model which contains both $V$ and $A$ as an element. 

 Distinguishing from the case where $o^{\vec{U}}(\kappa)<\kappa$, we do not have a classification of what are exactly the subforcings which generate the models $V[C']$. Let us give some examples of subforcings of $\Mfor$ in the case of $o^{\vec{U}}(\kappa)=\kappa$.
 \begin{example}\label{canonicalsequence}
  Let $G$ be a generic and let $C_G$ be the generic club added by $\Mfor$, consider the increasing continuous enumeration of $C_G$, $\langle C_G(i)\mid i<\kappa\rangle$. Assume that $C_G(0)>0$, and consider again the sequence $\langle \kappa_n\mid n<\omega\rangle$ which is defined as follows:
     $$\kappa_0=C_G(0), \ \kappa_{n+1}=C_G(\kappa_n).$$
     Consider the following tree of measures:
     $$\vec{W}=\langle W_{\vec{\alpha}}\mid \vec{\alpha}\in [\kappa]^{<\omega}\rangle$$
     where $W_{\vec{\alpha}}=U(\kappa,\max(\vec{\alpha}))$. Note here that since $o^{\vec{U}}(\kappa)=\kappa$, this is well defined. It is not hard to check the Mathias criterion for the tree-Prikry forcing with $\vec{W}$, given in \cite{TomTreePrikry}, to conclude that $\langle \kappa_n\mid n<\omega\rangle$ is a tree-Prikry generic sequence with respect to $\vec{W}$. Note that, since the sequence of measures $\langle U(\kappa,i)\mid i<\kappa\rangle$ is a discrete family of normal measures, this tree-Prikry forcing falls under the framework of \cite{MinimalPrikry} and therefore the model $V[\langle\kappa_n\mid n<\omega\rangle]$ is minimal above $V$. This phenomenon does not occur in generic extensions of $\Mfor$ with $o^{\vec{U}}(\kappa)<\kappa$.
    
\end{example}
\begin{example}
 The previous example can be made more complex. Let $f:[\kappa]^{<\omega}\rightarrow\kappa$ be any function. Then $\langle \alpha_n\mid n<\omega\rangle$ is defined as follows:
     $\alpha_0=C_G(\l\r)$ and $\alpha_{n+1}$ is obtained by applying $f$ to some finite $\vec{C}_n\in[C_G]^{<\omega}$ i.e. $\alpha_{n+1}=C_G(f(\vec{C}_n))$.
\end{example}

Another theorem proven in section $6$ determines the cofinality of the supremum of a fresh set in
 Prikry, Magidor, Magidor-Radin and Radin extensions.
\begin{theorem}\label{freshfresh}
Assume that $\mathbb{P}$ is either Prikry, tree Prikry, Magidor, Magidor-Radin or Radin forcing.
Let $G\subseteq\mathbb{P}$ be $V$-generic. If $A\in V[G]$ is a fresh set of ordinals with respect to $V$, then $cf^{V[G]}(\sup(A))=\omega$.
\end{theorem}

The paper is organized as follows:
\begin{itemize}
    \item Section $2$: Subsections $2.1,2.2$ consist of basic definitions and properties of the forcing. Then $2.3$ provides several general definitions and previous results. In subsection $2.4$ we develop the theory of fat trees.
    \item Section $3$: We deal with the case of  sets with cardinality less than $\kappa$.
    \item Section $4$: The proof for subsets of $\kappa$ is presented.
    \item Section $5$: In $5.1$ an argument for general sets is given. In $5.2$, we prove some general results above the quotient forcing of several Prikry type forcing.
    \item Section $6$: Devoted to the proof of \ref{freshfresh}. 
    \item Section $7$: Presents further research directions and open questions related to this paper.
\end{itemize}

\section{Preliminaries}
 Most of the basic definitions are identical to \cite{partOne} and \cite{Gitik2010}.

\subsection{Magidor forcing}
 
 Let $\vec{U}=\langle U(\alpha,\beta)\mid \alpha\leq \kappa \ ,\beta<o^{\vec{U}}(\alpha)\rangle$ be a coherent sequence. For every $\alpha\leq\kappa$, denote $$\cap\vec{U}(\alpha)=\underset{i<o^{\vec{U}}(\alpha)}{\bigcap}U(\alpha,i).$$
\begin{definition}\label{Magidor-conditions}
$\mathbb{M}[\vec{U}]$ consists of elements $p$ of the form
$p=\langle t_1,...,t_n,\langle\kappa,B\rangle\rangle$.
 For every  $1\leq i\leq n $, $t_i$ is either an ordinal
 $\kappa_i$ if $ o^{\vec{U}}(\kappa_i)=0$
 or a pair $\langle\kappa_i,B_i\rangle$  if \ $o^{\vec{U}}(\kappa_i)>0$.
\begin{enumerate}
\item $B\in\cap\vec{U}(\kappa)$, \  $\min(B)>\kappa_n$.
    \item  For every  $1\leq i\leq n$.
    \begin{enumerate}
    \item $\langle\kappa_1,...,\kappa_n\rangle\in [\kappa]^{<\omega}$ (increasing finite sequence below $\kappa$).
    \item $B_i\in \cap\vec{U}(\kappa_i)$.
    \item  $\min(B_i)>\kappa_{i-1}$ \ $(i>1)$.
    \end{enumerate}
\end{enumerate}
 Moreover, denote $t_{n+1}=\l\kappa,B\r$.
\end{definition}
\begin{definition}\label{Magidor-order}
 For $p=\langle t_1,t_2,...,t_n,\langle\kappa,B\rangle\rangle,q=\langle s_1,...,s_m,\langle\kappa,C\rangle\rangle\in \Mfor$ , define  $p \leq q$ ($q$ extends $p$) iff:
\begin{enumerate}
    \item $n \leq m$.
    \item $B \supseteq C$.
    \item $\exists 1 \leq i_1 <...<i_n \leq m$ such that for every $1 \leq j \leq m$:
    \begin{enumerate}
        \item If $\exists 1\leq r\leq n$ such that $i_r=j$ then $\kappa(t_r)=\kappa( s_{i_r})$ and $C(s_{i_r})\subseteq B(t_r)$.
        \item Otherwise $\exists \ 1 \leq r \leq n+1$ such that $ i_{r-1}<j<i_{r}$ then 
        \begin{enumerate}
        \item $\kappa(s_j) \in B(t_r)$.
        \item $B(s_j)\subseteq B(t_r)\cap \kappa(s_j)$.
        \end{enumerate}
    \end{enumerate}
\end{enumerate}
We also use ``p directly extends q", $q \leq^{*} p$ if:
\begin{enumerate}
    \item $q \leq p$.
    \item $n=m$.
\end{enumerate}
\end{definition}
Let us add some notation, for a pair $t=\langle \alpha, X\rangle$ we denote $\kappa(t)=\alpha,\ B(t)=X$. If $t=\alpha$ is an ordinal then $\kappa(t)=\alpha$,  $B(t)=\emptyset$, and $\cap \vec{U}(\alpha)=P(\alpha)$ (the power set of $\alpha$).

For a condition $p=\langle t_1,...,t_n,\langle \kappa,B\rangle\rangle\in\Mfor$ we denote
$n=l(p)$, $p_i=t_i$, $B_i(p)=B(t_i)$ and $\kappa_i(p)=\kappa(t_i)$ for any $1\leq i\leq l(p)$, $t_{l(p)+1}=\langle\kappa,B\rangle$, $t_0=0$. Also denote
$$\kappa(p)=\{\kappa_i(p)\mid i\leq l(p)\}\text{ and }B(p)=\bigcup_{i\leq l(p)+1}B_i(p).$$
\begin{remark}\label{changes with respect to part one}
In \cite{TomMoti},\cite{partOne} we had another requirement in Definition \ref{Magidor-order}, that given a condition $p$, if we would like to add an ordinal $\alpha$ to the sequence in the interval $(\kappa_{i-1}(p),\kappa_i(p))$ then we needed to make sure that $o^{\vec{U}}(\alpha)<o^{\vec{U}}(\kappa_i(p))$. 
This condition is not essential as any condition $p$ can be directly extended to a condition in the set
$$D=\{q\in \Mfor\mid\forall i\leq l(q)+1. \forall\alpha\in B_i(q). o^{\vec{U}}(\alpha)<o^{\vec{U}}(\kappa_i(q))\}.$$
The order defined in \ref{Magidor-order} on elements of $D$ automatically satisfies the extra requirement. 

For this reason we will point out along this section some points where this assumption changes properties of $\Mfor$. The major one, is in Propositions \ref{indc},\ref{IndCG}.
\end{remark}
\begin{definition}\label{end extension}
  Let $p\in\Mfor$. For every $ i\leq l(p)+1$, $\alpha\in B_{i}(p)$ with $o^{\vec{U}}(\alpha)>0$, and $B\in \cap\vec{U}(\alpha)$, define 
$$p^{\frown}\l\alpha,B\r=\langle p_1,...,p_{i-1},\langle\alpha,B_{i}(p)\cap B\rangle,\langle\kappa_{i}(p),B_{i}(p)\setminus(\alpha+1)\rangle,p_{i+1},...,p_{l(p)+1}\rangle.$$
Also $p^{\smallfrown}\l\alpha\r=p^{\smallfrown}\l\alpha,\alpha\r$. 
If $o^{\vec{U}}(\alpha)=0$, define
$$p^{\frown}\langle\alpha\rangle=\langle p_1,...,p_{i-1},\alpha,\langle\kappa_{i}(p),B_{i}(p)\setminus(\alpha+1)\rangle,...,p_{l(p)+1}\rangle.$$
For $\langle\alpha_1,...,\alpha_n\rangle\in[\kappa]^{<\omega}$ and $\l B_1,...,B_n\r$, where $B_i\in\cap\vec{U}(\alpha_i)$, define recursively,
$$p^{\smallfrown}\l\l\alpha_1,...,\alpha_n\r,\l B_1,...,B_n\r\r=(p^{\smallfrown}\l\l\alpha_1,...,\alpha_{n-1}\r,\l B_1,...,B_{n-1}\r\r)^{\smallfrown}\l \alpha_n,B_n\r$$
and
$$p^{\frown}\langle\alpha_1,...,\alpha_n\rangle=(p^{\frown}\langle\alpha_1,...,\alpha_{n-1}\rangle)^{\frown}\langle\alpha_n\rangle.$$
\end{definition}
For $\vec{\alpha}=\l \alpha_1,...,\alpha_n\r$, denote $|\vec{\alpha}|=n$ and $\vec{\alpha}(i)=\alpha_i$. If $I\subseteq\{1,...,n\}$ then $\vec{\alpha}\restriction I=\l \vec{\alpha}(i_1),...,\vec{\alpha}(i_k)\r$ where $\{i_1,i_2,...,i_k\}$ is the increasing enumeration of $I$. For $Y\subseteq\omega$, $\vec{\alpha}\restriction Y=\vec{\alpha}\restriction (Y\cap\{1,...,n\})$. We will usually identify $\vec{\alpha}$ with the set $\{\alpha_1,...,\alpha_n\}$. Also for two sequences $\vec{\alpha},\vec{\beta}$, we denote their concatenation by $\vec{\alpha}{}^{\smallfrown}\vec{\beta}$.

Note that if we add a pair of the form $\langle \alpha , B\cap\alpha\rangle$ then in $B\cap\alpha$ there might be many ordinals  which are irrelevant to the forcing and cannot be added. Namely, ordinals $\beta$ such that $B\cap\beta\notin\cap\vec{U}(\beta)$. Note that we no longer have to require $o^{\vec{U}}(\beta)\geq o^{\vec{U}}(\alpha)$.
We can avoid such ordinals by shrinking the large sets.
\begin{proposition}\label{BetterSet}  Let $\alpha\leq\kappa$, and $A\in\cap\vec{U}(\alpha)$. Then there exists $A^*\subseteq A$ such that:
 \begin{enumerate}
  \item $A^*\in\cap\vec{U}(\alpha)$.
  \item For every  $x\in A^*$,  $A^*\cap x\in\cap\vec{U}(x)$.
\end{enumerate}
\end{proposition}
\pr  For any $j<o^{\vec{U}}(\alpha)$,
$$Ult(V,U(\alpha,j))\models A=j_{U(\alpha,j)}(A)\cap\alpha\in \underset{i<j}{\bigcap}U(\alpha,i).$$
Coherency of the sequence implies that $A':=\{\alpha<\kappa\mid A\cap\alpha\in\cap\vec{U}(\alpha)\}\in U(\alpha,j)$, this is for every $j<o^{\vec{U}}(\alpha)$. 

 Define inductively $A^{(0)}=A$, $A^{(n+1)}=(A^{(n)})'$. By definition, $\forall\alpha\in A^{(n+1)}_j$,  $A^{(n)}\cap\alpha\in\cap\vec{U}(\alpha)$. Define $A^*=\underset{n<\omega}{\bigcap}A^{(n)}\in\cap\vec{U}(\kappa)$, this set
has the required property.
$\blacksquare$

The conditions $p^{\smallfrown}\vec{\alpha}$ and $p^{\smallfrown}\l\vec{\alpha},\vec{B}\r$ are minimal extensions of $p$ in a sense given in the following proposition. The proof of the proposition is a direct verification of \ref{Magidor-conditions},\ref{Magidor-order}.
\begin{proposition}\label{frown extension}
Let $p\in\Mfor$ and $\vec{\alpha}\in [\kappa]^n$. Suppose that $\vec{\alpha}$ decomposes according to the ordinals of $p$ as $\vec{\alpha}= \vec{\alpha}_1{}^{\smallfrown}...^{\smallfrown}\vec{\alpha}_{l(p)+1}\in\prod_{i=1}^{l(p)+1}[B_i(p)]^{l_i}$. Let $\vec{C}= \vec{C}_1{}^{\smallfrown}...^{\smallfrown}\vec{C}_{l(p)+1}$, be a sequence of sets such that $|\vec{C}_i|=l_i$ (in particular $|\vec{C}|=n$)and for each $i\leq n$, $\vec{C}(i)\subseteq \vec{\alpha}(i)$.
\begin{enumerate}
    \item $p^{\smallfrown}\l \vec{\alpha},\vec{C}\r\in \Mfor$ if and only if $\forall i\leq n. \ \exists j\leq l(p)$ such that $\vec{\alpha}(i)\in B_j(p)$ and $ B_{j}(p)\cap \vec{C}(i)\in\cap\vec{U}(\vec{\alpha}(i))$.
    \item  Suppose that $p^{\frown}\l\vec{\alpha},\vec{C}\r\in\Mfor$, then for any extension $q$ of $p$, if:
    \begin{enumerate}
    \item $\kappa(p)\cup\vec{\alpha}\subseteq \kappa(q)$. 
    \item For all $j\leq l(q)$, if $\kappa_{i-1}(p)<\kappa_j(q)\leq \kappa_i(p)$, for some $i\leq l(p)+1$, then:
    \begin{enumerate}
        \item $B_j(q)\subseteq \bigcup_{1\leq j\leq l_i}\vec{C}_i(j)\cup \big(B_i(p)\setminus \max(\vec{\alpha}_i+1)\big)$\footnote{Note that by Definition \ref{end extension}, the set $\bigcup_{1\leq i\leq l(p)+1}\Big(\bigcup_{1\leq j\leq l_i}\vec{C}_i(j)\cup \big(B_i(p)\setminus \max(\vec{\alpha}_i+1)\big)\Big)$ is exactly $B(p^{\smallfrown}\l\vec{\alpha},\vec{C}\r$.}.
    \item  If $\kappa_j(q)<\kappa_i(p)$, then $$\kappa_j(q)\in \bigcup_{1\leq j\leq l_i}\vec{C}_i(j)\cup \big(B_i(p)\setminus \max(\vec{\alpha}_i+1)\big).$$
    \end{enumerate}
    \end{enumerate}
    Then $p^{\frown}\l\vec{\alpha},\vec{B}\r\leq q$.
\end{enumerate}
\end{proposition}
The previous proposition also provides a criterion for $p^{\smallfrown}\vec{\alpha}\in \Mfor$, and establishes the minimality of the extension $p^{\smallfrown}\vec{\alpha}$. Namely, 
$\text{for every }p\leq q,\text{ if } \kappa(q)\cup\vec{\alpha}\subseteq \kappa(q) \text{ then }p^{\smallfrown}\vec{\alpha}\leq q$.
\begin{definition}
Let $p\in \Mfor$, $\alpha<\kappa$ and let $i\leq l(p)$ be such that $\alpha\in[\kappa_{i}(p),\kappa_{i+1}(p))$
 $$p\restriction\alpha= \langle p_1,...,p_i\rangle \ and \ p\restriction(\alpha,\kappa]=\langle p_{i+1},...,p_{l(p)+1}\rangle.$$
 Also, for $\lambda$ with $o^{\vec{U}}(\lambda)>0$ define
 $$\Mfor\restriction\lambda=\Big\{p\restriction\lambda\mid p\in\Mfor, \lambda \text{ appears in } p\Big\}, \ \ \ 
 \Mfor\restriction(\lambda,\kappa]=\{p\restriction (\lambda,\kappa]\mid p\in\Mfor, \lambda \text{ appears in } p\}.$$
\end{definition}
Note that $\Mfor\restriction\lambda$ is just Magidor forcing on $\lambda$ and $\Mfor\restriction (\lambda,\kappa]$ is a subset of $\Mfor$ which generates a Magidor club in the interval $(\lambda,\kappa]$. 
\begin{remark}\label{Magidor forcing in extension}
Let $\lambda<\kappa$ which $o^{\vec{U}}(\lambda)>0$ and let $p\in\Mfor$ be such that $\lambda$ appears in $p$. Let $H\subseteq \Mfor\restriction\lambda$ with $p\restriction\lambda\in H$, then in $V[H]$ we added new (bounded) subsets of $\kappa$, hence $\vec{U}$ is no longer a sequence of ultrafilters. However, for the relevant interval $(\lambda,\kappa]$, $\vec{U}\restriction(\lambda,\kappa]$ generates a coherent sequence of ultrafilters $\vec{W}$ and formally we force with $\mathbb{M}[\vec{W}]$. Note that the ground model forcing $\Mfor\restriction(\lambda,\kappa]$ is dense in $\mathbb{M}[\vec{W}]$, hence we can simply force with $\Mfor\restriction(\lambda,\kappa]$ over $V[H]$ to complete to a generic extension of $\Mfor$.
\end{remark}
The following propositions can be found in \cite{partOne}:
\begin{proposition}\label{dec} Let $p\in\Mfor$ and $\langle\lambda,B\rangle$ a pair in $p$. Then
$$\Mfor/p\simeq \Big(\Mfor\restriction \lambda\Big)/\Big(p\restriction\lambda\Big)\times\Big(\Mfor\restriction(\lambda,\kappa]\Big)/\Big(p\restriction(\lambda,\kappa]\Big).$$
\end{proposition}
\begin{proposition}
Let $p\in\Mfor$ and $\langle\lambda,B\rangle$ be a pair in $p$. Then the order $\leq^*$ in the forcing $\Big(\Mfor\restriction(\lambda,\kappa]\Big)/\Big(p\restriction(\lambda,\kappa]\Big)$ is $\delta$-directed where $\delta=\min\{\nu>\lambda\mid o^{\vec{U}}(\nu)>0\}$. Meaning that for every $X\subseteq \Mfor\restriction (\lambda,\kappa]$ such that $|X|<\delta$ and for every $q\in X, \ p\leq^* q$, there is an  $\leq^*$-upper bound for $X$.
\end{proposition}
 \begin{lemma}
$\Mfor$ satisfies $\kappa^+$-cc.
\end{lemma}
The following lemma is the well known Prikry condition:
\begin{lemma}\label{prikrycondition}
$\Mfor$ satisfies the Prikry condition i.e. for any statement in the forcing language $\sigma$ and any $p\in\Mfor$ there is $p\leq^*p^*$ such that $p^*||\sigma$ i.e. either $p^*\Vdash\sigma$ or $p\Vdash\neg\sigma$.
\end{lemma} 

The next lemma can be found in \cite{ChangeCofinality} and the proof in \cite{partOne}:
\begin{lemma}\label{MagLemma}
  Let $G\subseteq \Mfor$ be generic and
suppose that $A\in V[G]$ is such that $A\subseteq V_\alpha$. Let $p\in G$ and $\l\lambda,B\r$ a pair in $p$ such that $\alpha<\lambda$, then $A\in V[G\restriction\lambda]$.
\end{lemma}

\begin{corollary}
$\Mfor$ preserves all cardinals.
\end{corollary}
 \begin{definition}
  Let $G\subseteq \Mfor$ be generic, define the \textit{Magidor club}
  $$C_{G}=\{ \nu \mid \exists \ A\exists p\in G \ s.t. \ \langle \nu,A\rangle\in p\}.$$

\end{definition}
  We will abuse notation by sometimes considering $C_G$ as the canonical enumeration of the set $C_G$. The set $C_{G}$ is closed and unbounded in $\kappa$, therefore, the order type of $C_{G}$ determines the cofinality of $\kappa$ in $V[G]$. The next propositions can be found in \cite{Gitik2010}.
 \begin{proposition}\label{decprop}
 Let $G\subseteq\Mfor$ be generic. Then $G$ can be reconstructed from $C_{G}$ as follows
     $$ G=\{p\in\Mfor\mid (\kappa(p)\subseteq C_{G}) \wedge (C_{G}\setminus\kappa(p)\subseteq B(p))\}.$$
      In particular $V[G]=V[C_{G}]$.
\end{proposition}
\begin{proposition}\label{genericproperties}
Let $G\subseteq\Mfor$ be generic. 
\begin{enumerate}
    \item $C_G$ is a club at $\kappa$.
    \item For every $\delta\in C_G$, $o^{\vec{U}}(\delta)>0$ iff $\delta\in Lim(C_G)$\footnote{The set of limit points of $X\subseteq\kappa$ is $Lim(X):=\{\alpha\mid \sup(\alpha\cap X)=\alpha\}\subseteq\kappa+1$.}.
    \item For every $\delta\in Lim(C_G)$, and every $A\in \cap\vec{U}(\delta)$, there is $\xi<\delta$ such that $C_G\cap(\xi,\delta)\subseteq A$.
    
    \item If $\l \delta_i\mid i<\theta\r$ is an increasing sequence of elements of $C_G$, let $\delta^*=\sup_{i<\theta}\delta_i$, then $o^{\vec{U}}(\delta^*)\geq\limsup_{i<\theta}o^{\vec{U}}(\delta_i)+1$.\footnote{ For a sequence of ordinals $\l \rho_j\mid j<\gamma\r$, $\limsup_{j<\gamma}\rho_j=\min(\sup_{i<j<\gamma}\rho_j\mid i<\gamma)$.}
    \item Let $\delta\in Lim(C_G)$ and let $A$ be a positive set, $A\in (\cap\vec{U}(\delta))^+$, i.e. $
\delta\setminus A\notin \cap\vec{U}(\delta)$. \footnote{Equivalently,if there is some $i<o^{\vec{U}}(\delta)$ such that $A\in U(\delta,i)$.} Then, $\sup(A\cap C_G)=\delta$.
\item If $A\subseteq V_\alpha$, then $A\in V[C_G\cap\lambda]$, where $\lambda=\max(Lim(C_G)\cap\alpha+1)$.
\item For every $V$-regular cardinal $\alpha$, if $cf^{V[G]}(\alpha)<\alpha$ then $\alpha\in Lim(C_G)$.
\end{enumerate}
\end{proposition}
\pr 
The proof of $(1),(2),(3),(5),(6),(7)$ can be found in \cite{Gitik2010} and does not use the extra property of \ref{Magidor-order} (see Remark \ref{changes with respect to part one}).

To see $(4)$, use the closure of $C_G$, to find $q\in G$ such that $\delta^*$ appears in $q$. Clearly, $A:=\{\alpha<\delta^*\mid o^{\vec{U}}(\alpha)<o^{\vec{U}}(\delta^*)\}\in\cap\vec{U}(\delta^*)$, thus by $(3)$, there is $\xi<\delta^*$  such that $C_G\cap(\xi,\delta^*)\subseteq A$. Let $i<\theta$ be such that for every $j>i$, $o^{\vec{U}}(\delta_j)<o^{\vec{U}}(\delta^*)$.
By definition of $limsup$, $$\limsup_{j<\theta}o^{\vec{U}}(\delta_j)+1\leq \sup_{i<j<\theta}o^{\vec{U}}(\delta_j)+1\leq o^{\vec{U}}(\delta^*).$$
$\blacksquare$

\begin{proposition}\label{indc}
 Let $G\subseteq \Mfor$ be a $V$-generic filter and $C_{G}$ the corresponding Magidor sequence. Let $p\in G$, then for every $i\leq l(p)+1$ 
 \begin{enumerate}
     \item If  $o^{\vec{U}}(\kappa_i(p))\leq \kappa_i(p)$, and $\forall\alpha\in B_i(p)$, $o^{\vec{U}}(\alpha)<o^{\vec{U}}(\kappa_i(p))$, then
    $$\otp( [\kappa_{i-1}(p),\kappa_i(p))\cap C_{G} )=\omega^{o^{\vec{U}}(\kappa_i(p))}.$$
    \item If $o^{\vec{U}}(\kappa_i(p))\geq \kappa_i(p)$, then 
    $$\otp( [\kappa_{i-1}(p),\kappa_i(p))\cap C_{G} )=\kappa_i(p).$$
    
 \end{enumerate}
 \end{proposition}
\pr  The same as in \cite{partOne}, replacing the usage of Definition \ref{Magidor-order} with the assumption that $\forall\alpha\in B_i(p)$, $o^{\vec{U}}(\alpha)<o^{\vec{U}}(\kappa_i(p))$.$\blacksquare$

  Proposition \ref{indc} suggests a connection between the index in $C_G$ of ordinals appearing in $p$ 
  and Cantor normal form.

  \begin{definition}
 Let $p\in G$. For each $i\leq l(p)$ define
   $$\gamma_i(p)=\sum_{j=1}^{i}\omega^{o^{\vec{U}}(\kappa_j(p))}.$$ 
   \end{definition}
 \begin{corollary}\label{IndCG}
 Let G be $\Mfor$-generic and $C_{G}$ the corresponding Magidor sequence. Let $p\in G$, such that for every $1\leq i\leq l(p)$, and every $\alpha\in B_i(p)$, $o^{\vec{U}}(\alpha)<o^{\vec{U}}(\kappa_i(p))$, then
 $$p\Vdash \lusim{C}_G(\gamma_i(p))=\kappa(t_i).$$
 \end{corollary}

For more details and basic properties of Magidor forcing see \cite{ChangeCofinality},\cite{Gitik2010}, \cite{TomMoti} or \cite{partOne}. 
\subsection{Magidor forcing with $o(\kappa)<\kappa^+$}

When we assume $o^{\vec{U}}(\kappa)<\kappa$, the measure $U(\kappa,\xi)$ concentrates on measurables $\alpha$ with $o^{\vec{U}}(\alpha)=\xi$, which is a canonical discrete family for those measures. In our more general situation, $ o^{\vec{U}}(\kappa)<\kappa^+$, we can still separate the measures but the decomposition is not canonical. 
More precisely, for every $\alpha\leq\kappa$, we would like to have sets which witness the fact that the sequence of ultrafilters $\l U(\alpha,\beta)\mid \beta<o^{\vec{U}}(\alpha)\r$ is discrete.
\begin{proposition}\label{decompositionHighOrder}
Assume $o^{\vec{U}}(\alpha)<\alpha^+$, then there are pairwise disjoint sets $\langle X^{(\alpha)}_i\mid i<o^{\vec{U}}(\alpha)\rangle$ such that $X^{(\alpha)}_i\in U(\alpha,i)$. \end{proposition}
\pr By assumption, $|o^{\vec{U}}(\alpha)|\leq\alpha$. Enumerate the measures $$\{U(\alpha,i)\mid i<o^{\vec{U}}(\alpha)\}=\{W_j\mid j<\rho\}$$ where $\rho\leq\alpha$. For every $i\neq j$ below $\rho$, find $Y_{i,j}\in W_i\setminus W_j$. By normality $Y_i=\Delta_{j<\rho} Y_{i,j}\in W_i$. Also, for $j\neq i$, $Y_i\notin W_j$ since $Y_i\subseteq Y_{i,j}\cup j\notin W_j$. Set $Z_i=Y_i\setminus(\cup_{j<i}Y_j)$, then $Z_i\in W_i$  and $\l Z_i\mid i<\rho\r$ are pairwise disjoint. Finally, define $X^{(\alpha)}_\xi=Z_i$ where $\xi<o^{\vec{U}}(\kappa)$ is such that $W_i=U(\alpha,\xi)$.$\blacksquare$
\begin{definition}
Let $\alpha\leq\kappa$.
\begin{enumerate}
    \item For $o^{\vec{U}}(\alpha)\leq\alpha$ define for every $i<o^{\vec{U}}(\alpha)$ $$X^{(\alpha)}_i=\{x<\alpha\mid i=o^{\vec{U}}(x)\}\in U(\alpha,i).$$
    \item For $\alpha< o^{\vec{U}}(\alpha)<\alpha^+$ fix a decomposition of $\alpha$, $\l X^{(\alpha)}_i\mid i<o^{\vec{U}}(\alpha)\r$ guaranteed by the previous proposition such that $X^{(\alpha)}_i\in U(\alpha,i)$.
    \item
  For $\beta<\alpha$ denote by $o^{(\alpha)}(\beta)=\xi$ the unique $\xi<o^{\vec{U}}(\alpha)$ such that $\beta\in X^{(\alpha)}_\xi$. Also let $o^{(\alpha)}(\alpha)=o^{\vec{U}}(\alpha)$.
  \end{enumerate}
\end{definition}
Note that if $o^{\vec{U}}(\alpha)\leq\alpha$ then $o^{(\alpha)}(\beta)=o^{\vec{U}}(\beta)$.
\begin{proposition}\label{increasing order at limits}
For every $V$-generic $G\subseteq\Mfor$ and for every $\kappa_0\in Lim(C_G)$ (Recall that $\kappa\in Lim(C_G)$) such that $o^{\vec{U}}(\kappa_0)<\kappa_0^+$, there is $\xi<\kappa_0$ such that for every $\alpha\in Lim(C_G)\cap(\xi,\kappa_0]$ $$o^{(\kappa_0)}(\alpha)\geq limsup(o^{(\kappa_0)}(\beta)+1\mid \beta\in C_G\cap\alpha).$$
In other words, there is $\xi_\alpha<\alpha$ such that for every $\beta\in C_G\cap(\xi_\alpha,\alpha)$, $o^{(\kappa_0)}(\beta)<o^{(\kappa_0)}(\alpha)$.
\end{proposition}
\pr If $o^{\vec{U}}(\kappa_0)<\kappa_0$, then $o^{(\kappa_0)}(\alpha)=o^{\vec{U}}(\alpha)$ and the proposition follows from \ref{genericproperties}(4). Also if $\alpha=\kappa_0$, then clearly for every $\beta<\kappa_0$, $o^{(\kappa_0)}(\beta)<o^{\vec{U}}(\kappa_0)$ by definition.
Assume that $\kappa_0\leq o^{\vec{U}}(\kappa_0)<\kappa_0^+$ and let $\pi:\kappa_0\longleftrightarrow o^{\vec{U}}(\kappa_0)$ be a bijection. For every $\rho<o^{\vec{U}}(\kappa_0)$ denote by
$$E_\rho=\pi^{-1''}\rho\subseteq \kappa_0$$
and for every $\alpha<\kappa_0$ define $Y_\alpha=X^{(\kappa_0)}_{\pi(\alpha)}$. 
In $M_{U(\kappa_0,\rho)}$, define $$j_{U(\kappa_0,\rho)}(\l Y_\alpha\mid \alpha<\kappa_0\r)=\l Y'_\alpha\mid \alpha<j_{U(\kappa_0,\rho)}(\kappa_0)\r.$$ Since $crit(j_{U(\kappa_0,\rho)})=\kappa_0$, for $\alpha<\kappa_0$, $Y'_\alpha=j_{U(\kappa_0,\rho)}(Y_\alpha)$. Moreover, $j_{U(\kappa_0,\rho)}(E_{\rho})\cap\kappa_0=E_{\rho}$ and $j_{U(\kappa_0,\rho)}(Y_\alpha)\cap\kappa_0= Y_\alpha$. Hence $$\cup_{\alpha\in j_{U(\kappa_0,\rho)}(E_{\rho})\cap\kappa_0} Y'_\alpha\cap\kappa_0=\cup_{\alpha\in E_\rho}j_{U(\kappa_0,\rho)}(Y_\alpha)\cap\kappa_0=\cup_{\alpha\in E_\rho}Y_\alpha=\cup_{\xi<\rho}X^{(\kappa_0)}_\xi.$$
By coherency, $\cap_{\xi<\rho}U(\kappa_0,\xi)=\cap j_{U(\kappa_0,\rho)}(\vec{U})(\kappa_0)$, thus
$$(*) \ \ \ \ \ M_{U(\kappa_0,\rho)}\models\cup_{\alpha\in j_{U(\kappa_0,\rho)}(E_{\rho})\cap\kappa_0} Y'_\alpha\cap\kappa_0\in \cap j_{U(\kappa_0,\rho)}(\vec{U})(\kappa_0).$$
%%Also, for every $\alpha\in j_{U(\kappa,\rho)}(E_{\rho})\cap\kappa}$, $Y'_\alpha\cap\kappa=Y_\alpha=X^{(\kappa_0)}_{\pi(\alpha)}\in U(\kappa,\pi(\alpha))$, hence $$ \ \ \ \ \ (**) \ \ \ \ \  M_{U(\kappa,\rho)}\models \forall \alpha\in j_{U(\kappa,\rho)}(E_{\rho})\cap\kappa}. Y'_\alpha\cap\kappa\in (\cap j_{U(\kappa,\rho)}(\vec{U}) )^+$$
Reflecting $(*)$ we get
$$X'_\rho=\Big\{\beta\in X^{(\kappa_0)}_\rho \mid \cup_{\alpha\in E_{\rho}\cap\beta} Y_\alpha\cap\beta\in \cap \vec{U}(\beta) \Big\}\in U(\kappa_0,\rho).$$
Now let $\xi<\kappa_0$ be such that $C_G\cap (\xi,\kappa_0)\subseteq \cup_{\rho<o^{\vec{U}}(\kappa_0)}X'_{\rho}$, and let $\alpha\in Lim(C_G)\cap(\xi,\kappa_0)$. Denote  $o^{(\kappa_0)}(\alpha)=\rho$, and since $X^{(\kappa_0)}_i$ are pairwise disjoint,  $\alpha\in X'_{\rho}$. By definition of $X'_{\rho}$, $$\cup_{i\in E_{\rho}\cap\alpha} Y_i\cap\alpha\in \cap \vec{U}(\alpha)\text{ and } \forall i\in E_\rho\cap\alpha. Y_i\cap \beta\in(\cap \vec{U}(\alpha))^+.$$ By \ref{genericproperties}(3) there is $\xi_\alpha<\alpha$ such that $C_G\cap (\xi_\alpha,\alpha)\subseteq \cup_{i\in E_{\rho}\cap\alpha} Y_i\cap\alpha$. In particular, for every $\beta\in C_G\cap(\xi_\alpha,\alpha)$, there is $i\in E_{\rho}\cap\alpha$ such that $\beta\in Y_i=X^{(\kappa_0)}_{\pi(i)}$. Since $i\in E_\rho$, $\pi(i)<\rho$ so $o^{(\kappa_0)}(\beta)<\rho=o^{(\kappa_0)}(\alpha)$, hence $limsup_{\beta\in C_G\cap\alpha}(o^{(\kappa_0)}(\beta)+1)\leq o^{(\kappa_0)}(\alpha)$. $\blacksquare$

 \begin{corollary}\label{The Wintessing sequence}
For every $V$-generic $G\subseteq\Mfor$ and for every $\kappa_0\in Lim(C_G)$ with $o^{\vec{U}}(\kappa_0)<\kappa_0^+$ there is $\eta<\kappa_0$ such that for every $\alpha\in Lim(C_G)\cap(\eta,\kappa_0]$ the following hold:
\begin{enumerate}
    \item If $o^{(\kappa_0)}(\alpha)=\beta+1$ is a successor ordinal, then there is $\xi<\alpha$ such that $\otp(C_G\cap X^{(\kappa_0)}_{\beta}\cap(\xi,\alpha))=\omega$, hence $cf^{V[G]}(\alpha)=\omega$.
    \item If $cf^V(o^{(\kappa_0)}(\alpha))=\lambda<\kappa_0$, then $\lambda<\alpha$ and let $\l \rho_i\mid i<\lambda\r$ be cofinal in $o^{(\kappa_0)}(\alpha)$, then there is $\xi<\alpha$ such that the sequence $x_i=\min(C_G\cap X^{(\kappa)}_{\rho_i}\setminus\xi)$ is increasing and unbounded in $\alpha$, hence $cf^{V[G]}(\alpha)=cf^{V[G]}(\lambda)$.
    \item Assume that $cf^V(o^{\kappa_0}(\alpha))=\kappa$, and let $\l \rho_i\mid i<\kappa\r$ be cofinal in $o^{\kappa_0}(\alpha)$, then there is $\xi<\alpha$ such that the sequence $x_0=\min(C_G\cap(\xi,\kappa_0))$ and $x_{n+1}=\min( C_G\cap X^{(\kappa)}_{\rho_{x_n}})$ is increasing and unbounded in $\alpha$, hence $cf^{V[G]}(\alpha)=\omega$.
\end{enumerate}
\end{corollary}

\pr 
For each successor $\rho=\beta+1$ 
 consider the set
 $$S_\rho=\{\alpha\in X^{(\kappa_0)}_\rho\mid X^{(\kappa_0)}_\beta\cap\alpha\in(\cap\vec{U}(\alpha))^+\}.$$
 Since $j_{U(\kappa_0,\rho)}(X^{(\kappa_0)}_\beta)\cap \kappa_0=X^{(\kappa_0)}_\beta\in U(\kappa_0,\beta)$, the coherency implies that $j_{U(\kappa_0,\rho)}(X^{(\kappa_0)}_\beta)\cap \kappa_0\in (\cap j_{U(\kappa_0,\rho)}(\vec{U})(\kappa_0))^+$. By elementarity, $\kappa_0\in j_{U(\kappa_0,\rho)}(S_\rho)$, hence $S_{\rho}\in U(\kappa_0,\rho)$.
 
 For $\rho$ such that $cf^{V}(\rho)=:\lambda<\kappa_0$, fix a cofinal sequence $\l \rho_i\mid i<\lambda\r\in V$.
  Consider the set
$$S'_{\rho}=\{\alpha\in X^{(\kappa_0)}_\rho\mid\forall i<\lambda. X^{(\kappa_0)}_{\rho_i}\cap\alpha\in(\cap\vec{U}(\alpha))^+\}.$$
 Also $S'_\rho\in U(\kappa_0,\rho)$.
 Indeed, since $\lambda<\kappa_0$ $$j_{U(\kappa_0,\rho)}(\l X^{(\kappa_0)}_{\rho_i}\mid i<\lambda\r)=\l j_{U(\kappa_0,\rho)}(X^{(\kappa_0)}_{\rho_i})\mid i<\lambda\r.$$ As before, for every $i<\lambda$ it follows that $j_{U(\kappa_0,\rho)}(X^{(\kappa_0)}_{\rho_i})\cap \kappa_0\in (\cap j_{U(\kappa_0,\rho)}(\vec{U})(\kappa_0))^+$, thus $S'_\rho\in U(\kappa_0,\rho)$. We shrink $S'_{\rho}$ a bit more, consider
 $$S_{\rho}=\Big\{\alpha\in S'_{\rho}\mid \{\beta<\alpha\mid \forall  i<\lambda. \beta\in X^{(\kappa_0)}_{\rho_i}\rightarrow\forall j<i.X^{(\kappa_0)}_{\rho_j}\cap\beta\in (\cap\vec{U}(\beta))^+\}\in\cap\vec{U}(\alpha)\Big\}.$$
 To see that $S_\rho\in U(\kappa_0,\rho)$, 
 for every $i<\lambda$ consider the set $$E_{\rho_i}=\{\beta\in X^{(\kappa_0)}_{\rho_i}\mid\forall j<i.X^{(\kappa_0)}_{\rho_j}\cap\beta\in(\cap\vec{U}(\beta))^+\}.$$
 In $M_{U(\kappa_0,\rho_i)}$, for every  $j<i$, $j_{U(\kappa_0,\rho_i)}(X^{(\kappa_0)}_{\rho_j})\cap\kappa_0\in(\cap j_{U(\kappa_0,\rho_i)}(\vec{U})(\kappa_0))^+$ it follows that $E_{\rho_i}\in U(\kappa_0,\rho_i)$. For $y\in \rho\setminus\{\rho_i\mid i<\lambda\}$, set $E_{y}=X^{(\kappa_0)}_y$. Then $E:=\cup_{y<\rho}E_y\in\cap_{y<\rho}U(\kappa_0,y)$. The set $E$ has the property that for every $\beta\in E$, if $\beta\in X^{(\kappa_0)}_{\rho_i}$ for some $i<\lambda$, then $\beta\in E_{\rho_i}$ and therefore $\forall j<i.X^{(\kappa_0)}_{\rho_j}\cap\beta\in(\cap\vec{U}(\beta))^+$.
 
 In $M_{U(\kappa_0,\rho)}$, by coherency $o^{j_{U(\kappa_0,\rho)}(\vec{U})}(\kappa_0)=\rho$ and for every $\beta<\kappa_0$, $\cap j_{U(\kappa_0,\rho)}(\vec{U})(\beta)=\cap\vec{U}(\beta)$. Also  $E\in M_{U(\kappa_0,\rho)}$ (by $\kappa_0$-closure) and $E\in \cap j_{U(\kappa,\rho)}(\vec{U})(\kappa_0)$. Denote $X'_i=j_{U(\kappa_0,\rho)}(X^{(\kappa_0)}_{\rho_i})$, then for every $\beta\leq\kappa_0$, $X'_i\cap\beta= X^{(\kappa_0)}_{\rho_i}\cap\beta$. It follows that $$M_{U(\kappa_0,\rho)}\models \{\beta<\kappa_0\mid \forall  i<\lambda. \beta\in X'_{i}\rightarrow\forall j<i.X'_j\cap\beta\in (\cap j_{U(\kappa_0,\rho)}(\vec{U})(\beta))^+\}\in\cap j_{U(\kappa_0,\rho)}(\vec{U})(\kappa_0) .$$
 Reflecting this, we get that $S_\rho\in U(\kappa_0,\rho)$.
 
 If $cf^V(\rho)=\kappa_0$, fix a continuous cofinal sequence $\l \rho_i\mid i<\kappa_0\r\in V$, consider
 $$S'_\rho=\{\alpha\in X^{(\kappa_0)}_\rho\mid \forall i<\alpha. X^{(\kappa_0)}_{\rho_i}\cap\alpha\in(\cap\vec{U}(\alpha))^+\}.$$
 Then as before $S'_\rho\in U(\kappa_0,\rho)$. Next, consider
 $$S_\rho=\Big\{\alpha\in S'_{\rho}\mid \{\beta<\alpha\mid \exists\zeta<\beta.\cup_{i<\rho_{\zeta}}X^{(\kappa_0)}_i\cap\beta\in\cap\vec{U}(\beta)\}\in \cap\vec{U}(\alpha)\Big\}.$$
 To see that $S_{\rho}\in U(\kappa_0,\rho)$, let $\xi<\rho$, find $\zeta<\kappa_0$ such that $\rho_\zeta>\xi$. Denote
 $$j_{U(\kappa_0,\xi)}(\l X^{(\kappa_0)}_i\mid i<o^{\vec{U}}(\kappa_0)\r)=\langle X'_i\mid i<o^{j(\vec{U})}(j_{U(\kappa_0,\xi)}(\kappa_0))\r, \ j_{U(\kappa_0,\xi)}(\l \rho_i\mid i<\kappa_0\r)=\l \rho'_i\mid i<j_{U(\kappa_0,\xi)}(\kappa_0)\r.$$
  then $\rho_\zeta\leq j_{U(\kappa_0,\xi)}(\rho_\zeta)=\rho'_{\zeta}$. If follows that   $\cup_{i<\rho'_{\zeta}}X'_i\cap\kappa_0\in\cap_{i<\xi}U(\kappa_0,i)=\cap j_{U(\kappa_0,\xi)}(\vec{U})(\kappa_0)$. To see this, note that for every $y<\xi$, $j_{U(\kappa_0,\xi)}(y)<j_{U(\kappa_0,\xi)}(\rho_{\zeta})=\rho'_{\zeta}$, hence $$X^{(\kappa_0)}_y=j_{U(\kappa_0,\xi)}(X^{(\kappa_0)}_y)\cap\kappa_0=X'_{j_{U(\kappa_0,\xi)}(y)}\cap\kappa_0\subseteq \cup_{i<\rho'_{\zeta}}X'_i\cap\kappa_0.$$
   This means that in $M_{U(\kappa_0,\xi)}$,
 $$\exists\zeta<\kappa_0. \ \cup_{i<\rho'_\zeta}X'_i\cap\kappa_0\in\cap j_{U(\kappa_0,\xi)}(\vec{U})(\kappa_0).$$
 Reflecting this, we get that for every $\xi<\rho$,
 $$\{\beta<\kappa_0\mid \exists\zeta<\beta.\cup_{i<\rho_\zeta}X^{(\kappa_0)}_i\cap\beta\in\cap\vec{U}(\beta)\}\in U(\kappa_0,\xi).$$
 Now in $M_{U(\kappa_0,\rho)}$ using coherency it follows that
 $$\{\beta<\kappa_0\mid \exists\zeta<\beta.\cup_{i<\rho_\zeta}X^{(\kappa_0)}_i\cap\beta\in\cap\vec{U}(\beta)\}\in\cap_{\xi<\rho}U(\kappa_0,\xi)=\cap j_{U(\kappa_0,\rho)}(\vec{U})(\kappa_0).$$
 Finally, reflect this to conclude that $S_\rho\in U(\kappa_0,\rho)$.
 
 By \ref{genericproperties}(3) there is $\eta'$ such that $C_G\cap(\eta',\kappa_0)\subseteq\cup_{\rho<o^{\vec{U}}(\kappa_0)}S_\rho$, define $\eta=\max\{\eta',\xi\}<\kappa_0$ where $\xi$ is from Proposition \ref{increasing order at limits}. Let $\alpha\in Lim(C_G)\cap(\eta,\kappa_0)$, then $\alpha\in S_{o^{(\kappa_0)}(\alpha)}$. Since $\alpha>\xi$, there is $\xi\leq \xi_\alpha<\alpha$ such that for every $\nu\in C_G\cap(\xi_\alpha,\alpha)$, $o^{(\kappa_0)}(\nu)<o^{(\kappa_0)}(\alpha)$.
 
 If $o^{(\kappa_0)}(\alpha)=\beta+1$ then $X^{(\kappa_0)}_\beta\cap\alpha\in(\cap\vec{U}(\alpha))^+$ hence by \ref{genericproperties}(5), $\sup(X^{(\kappa_0)}_{\beta}\cap\alpha\cap C_G)=\alpha$.  Let us argue that $\otp(X^{(\kappa_0)}_\beta\cap C_G\cap(\xi_\alpha,\alpha))=\omega$. Just otherwise denote by $\mu$ the $\omega$-th element of $X^{(\kappa_0)}_\beta\cap C_G\cap(\xi_\alpha,\alpha)$, then $\mu<\alpha$. Since $\mu>\xi$, Proposition \ref{increasing order at limits} implies that $o^{\kappa_0}(\mu)\geq \beta+1$. On the other hand, $\mu>\xi_\alpha$, thus $o^{(\kappa_0)}(\mu)<o^{(\kappa_0)}(\alpha)$, contradiction.
 
 If $cf^V(o^{(\kappa_0)}(\alpha)):=\lambda<\kappa_0$, then by definition of $S_{o^{(\kappa_0)}(\alpha)}$, $\forall i<\lambda. X^{(\kappa_0)}_{\rho_i}\cap\alpha\in(\cap\vec{U}(\alpha))^+$. hence by \ref{genericproperties}(5), for every $ i<\lambda$, $\sup(X^{(\kappa_0)}_{\rho_i}\cap\alpha\cap C_G)=\alpha$, thus the sequence of $x_i$'s defined in the proposition starting above any $\xi<\alpha$ is well defined. The second property of $S_{o^{(\kappa_0)}(\alpha)}$ is that $$Y:=\{\beta<\alpha\mid \forall  i<\lambda. \beta\in X^{(\kappa_0)}_{\rho_i}\rightarrow\forall j<i.X^{(\kappa_0)}_j\cap\beta\in (\cap\vec{U}(\beta))^+\}\in\cap\vec{U}(\alpha).$$ By \ref{genericproperties}(3) there is $\xi\leq \zeta_\alpha<\alpha$ such that $C_G\cap(\zeta_\alpha,\alpha)\subseteq Y$. Start the definition of $x_i$'s above $\zeta_\alpha$. To see it is increasing, note that $x_i\in C_G\cap (\zeta_\alpha,\alpha)\cap X^{(\kappa_0)}_{\rho_i}$ so by definition of $Y$, $\forall j<i$, $X^{(\kappa_0)}_{\rho_j}\cap x_i\in (\cap\vec{U}(x_i))^+$, again by \ref{genericproperties}(5), for every $j<i$ $\sup(X^{(\kappa_0)}_{\rho_j}\cap x_i\cap C_G)=x_i$ and therefore by minimality of $x_j$ it follows that for $j<i$, $x_j<x_i$. To see that the sequence of $x_i$'s is unbounded, notice that otherwise its limit point would be some $\zeta\in (\zeta_\alpha,\alpha)$. Since the $x_i$'s are increasing and by Proposition \ref{increasing order at limits}, $$o^{(\kappa_0)}(\zeta)\geq limsup_{i<\lambda} o^{(\kappa_0)}(x_i)+1 =limsup_{i<\lambda}\rho_i+1= o^{(\kappa_0)}(\alpha).$$
 contradicting the choice of $\xi_\alpha$.
 
 Finally, if $cf^V(o^{(\kappa_0)}(\alpha))=\kappa_0$, then $\forall i<\alpha. X^{(\kappa_0)}_{\rho_i}\cap\alpha\in(\cap\vec{U}(\alpha))^+$ hence by \ref{genericproperties}(5), $\forall i<\alpha.\sup(X^{(\kappa_0)}_{\rho_i}\cap\alpha\cap C_G)=\alpha$. If 
 the limit $x^*$ of the $x_n$'s defined in the proposition would be less than $\alpha$, then by the definition of $S_{o^{(\kappa_0)}(\alpha)}$ there is $\zeta<x^*$ such that $\cup_{i<\rho_{\zeta}}X^{(\kappa_0)}_i\cap x^*\in\cap \vec{U}(x^*)$. 
 To see the contradiction, on one hand there is $\sigma< x^*$ such that $C_G\cap (\sigma,x^*)\subseteq \cup_{i<\rho_{\zeta}}X^{(\kappa_0)}_i\cap\zeta$ so there is $N<\omega$ such that $\forall n\geq N$, $x_n\in \cup_{i<\rho_{\zeta}}X^{(\kappa_0)}_i\cap\zeta$.
 On the other find $N\leq n<\omega$  such that $x_n>\zeta$, then $o^{(\kappa_0)}(x_{n+1})=\rho_{x_n}>\rho_{\zeta}$, which implies $x_{n+1}\notin \cup_{i<\rho_{\zeta}}X^{(\kappa_0)}_i\cap\zeta$.$\blacksquare$
\begin{corollary}\label{VGenericCardinals}
Let $G\subseteq \Mfor$ be $V$-generic. Assume that $o^{\vec{U}}(\kappa)<\kappa^+$, then for every $V$-regular cardinal $\alpha$, $cf^{V[G]}(\alpha)<\alpha$ iff $\alpha\in C_G\cup \{\kappa\}$ and $0<o^{\vec{U}}(\alpha)<\alpha^+$.
\end{corollary}
\subsection{Other preliminaries}
In the last part of the proof we will need to analyze the quotient forcing. Let us recall some basic facts about it:
\begin{definition}
Let $\mathbb{P},\mathbb{Q}$ be forcing notions. A function $\tau:\mathbb{P}\rightarrow\mathbb{Q}$ is a projection iff $\tau$ is order preserving, $Im(\tau)$ is dense, and $$\forall p\in \mathbb{P}.\forall q\geq\tau(p).\exists p'\geq p.\pi(p')\geq q$$
\end{definition}
\begin{definition}\label{definition of quotient}
Let $\mathbb{P},\mathbb{Q}\in V$ be forcing notions,  $\tau:\mathbb{P}\rightarrow\mathbb{Q}$ be any projection and let $H\subseteq\mathbb{Q}$ be $V$-generic. Define \textit{the quotient forcing}
$\mathbb{P}/H=\tau^{-1''}H$.
Also if $G\subseteq \mathbb{P}$ is a $V$-generic filter, \textit{the projection of $G$} is the filter 
$$\tau_*(G):=\{q\in\mathbb{Q}\mid\exists p\in G. q\leq_{\mathbb{Q}}\tau(p)\}$$
\end{definition}
\begin{proposition}\label{properties of quotient} Let $\tau:\mathbb{P}\rightarrow\mathbb{Q}$ be a projection, then:
\begin{enumerate}
\item If $G\subseteq\mathbb{P}$ is $V$-generic then $\tau_*(G)$ is $V$-generic filter for $\mathbb{Q}$
    \item If $G\subseteq\mathbb{P}$ is $V$-generic then $G\subseteq\mathbb{P}/\tau_*(G)$ is $V[\tau_*(G)]$-generic filter.
    \item If $H\subseteq\mathbb{Q}$ is $V$-generic and $G\subseteq\mathbb{P}/H$ is $V[H]$-generic, then $\tau_*(G)=H$ and $G\subseteq\mathbb{P}$ is $V$-generic.
\end{enumerate}
\end{proposition}
\begin{definition}\label{definition of equivalent subalgebra}
 Let $\mathbb{P}$ be a forcing notion and $\lusim{D}$ be a $\mathbb{P}$-name for a subset of $\kappa$. Define $\mathbb{P}_{\lusim{D}}$, the complete subalgebra of regular open cuts $\l RO(\mathbb{P}),\leq_B\r$ \footnote{$RO(\Mfor)$ is the set of all regular open cuts of $\Mfor$(see for example \cite[Thm. 14.10]{Jech2003}), as usual we identify $\Mfor$ as a dense subset of $RO(\Mfor)$. The order $\leq_B$ is in the standard definition of Boolean algebras orders i.e. $p\leq_B q$ means $p\Vdash q\in \hat{G}$.} generated by the set $X=\{||\alpha\in \lusim{D}||\mid \alpha<\kappa\}$.\end{definition}
 \begin{definition}\label{definition of projection}
Define the function $\pi:\mathbb{P}\rightarrow \mathbb{P}_{\lusim{D}}$ by 
$\pi(p)=\inf\{b\in \mathbb{P}_{\lusim{D}}\mid p\leq_B b\}$.
\end{definition}
It not hard to check that $\pi$ is a projection.
Let $G$ be $V$-generic for $\mathbb{P}$ and $D\subseteq \kappa$ the interpretation of $\lusim{D}$ under $G$ i.e. $\lusim{D}_G=D$. Denote by $H=\pi_*(G)$ the $V$-generic filter for $\mathbb{P}_{\lusim{D}}$ induced by $G$, then $V[D]=V[H]$ (see for example \cite[Lemma 15.42]{Jech2003}). In fact $$D=\{\alpha<\kappa\mid ||\alpha\in \lusim{D}||\in X\cap H\}$$
As for the other direction, any generic filter $H$ is definable and uniquely determined (see \cite[Lemma 15.40]{Jech2003}) by the set
$$X\cap H=\{||\alpha\in\lusim{D}||\mid \alpha\in D\}$$
We sometimes abuse notation by defining $\mathbb{P}/D=\mathbb{P}/\pi_*(G)$.
It is important to note that $\mathbb{P}/D$ depends on the choice of the name $\lusim{D}$.

\begin{definition}\label{Index}
 Let $X,X'$ be sets of ordinals such that $X'\subseteq X\subseteq On$. Let $\alpha=otp(X,\in)$ be the order type of $X$ and $\phi:\alpha\rightarrow X$ be the order isomorphism witnessing it. The indices of $X'$ in $X$ are $$Ind(X',X)=\phi^{-1''}X'=\{\beta<\alpha\mid \phi(\beta)\in X'\}.$$
\end{definition}
\begin{definition}
  We denote $X\subseteq^* Y$ if $X\setminus Y$ is finite. Also define $X=^*Y$ if $X\subseteq^* Y\wedge Y\subseteq^* X$, equivalently, if $X\bigtriangleup Y$ is finite.
\end{definition}
Notice that the $X\subseteq^*Y$ sometimes denotes inclusion modulo \textit{bounded}, however in this paper, $X\subseteq^* Y$ means inclusion modulo \textit{finite}.
In the next theorem, we will need the Erd\"{o}s-Rado theorem \cite{ErdesRado}, which is stated here for the convenience of the reader (For the proof see \cite[Theorem 7.3]{kanamori1994} or \cite{YairEskew}). 
\begin{theorem}\label{ER}
If $\theta$ is a regular cardinal
then for every $\rho<\theta$, $$(2^{<\theta})^+\rightarrow(\theta+1)^2_{\rho}$$
i.e. for every $f:[(2^{<\theta})^+]^2\rightarrow \rho$ there is $H\subseteq (2^{<\theta})^+$ such that $\otp(H)=\theta+1$ and $f\restriction [H]^2$ is constant.
\end{theorem}

\begin{theorem}\label{Modfinitestab}
Let $\aleph_0<\lambda$ be a strong limit cardinal, and $\mu>\lambda$ be regular. Let $\langle D_\alpha\mid \alpha<\mu\rangle$ be any $\subseteq^*$-increasing sequence of subsets of $\lambda$. Then the sequence $=^*$-stabilizes i.e. there is $\alpha^*<\mu$ such that for every $\alpha^*\leq \alpha<\mu$, $D_\alpha=^*D_{\alpha^*}$. 
\end{theorem}
\begin{remark}
The theorem fails for $\lambda=\aleph_0$. Let us construct a counter example:

Define $\l D_i\mid i<\omega_1\r$  a sequence of subsets of $\omega$ by induction, such that:
\begin{enumerate}
    \item $\l D_i\mid i<\omega_1\r$ is $\subseteq^*$-increasing.
    \item For all $i<j<\omega_1$, $|D_j\setminus D_i|=\aleph_0$.
    \item For every $i<\omega_1$, $|\omega\setminus D_i|=\aleph_0$.
\end{enumerate}Let $D_0=\emptyset$. Assume that for $\alpha<\omega_1$, $\l D_i\mid i<\alpha\r$ is $\subseteq^*$-increasing, and let us define $D_\alpha$. If $\alpha=\beta+1$, then by $3$,  $|\omega\setminus D_\beta|=\aleph_0$, Let $\omega\setminus D_{\beta}=X\uplus Y$ where $|X|=|Y|=\aleph_0$. Define $D_\alpha=D_{\beta}\cup X$. If $\alpha$ is limit, then $cf(\alpha)=\omega$, let $\l \alpha_n\mid n<\omega\r$ be increasing and cofinal in $\alpha$ and denote $E_n=D_{\alpha_n}$. 
We construct natural numbers $x_n,y_n$. By $3$, $|\omega\setminus E_0|=\omega$, let $x_0,y_0\in \omega\setminus E_0$ be distinct. Assume that $x_k, y_k$ are defined for every $k\leq n$, then $Z=\omega\setminus((\cup_{m\leq n+1}E_{m})\cup\{x_k,y_k\mid k\leq n\})$ is infinite. Indeed, for each $m\leq n$, $E_m\subset^* E_{n+1}$ hence $R_m:=E_m\setminus E_{n+1}$ is finite. It follows that $R=\cup_{m\leq n}R_m$ a is finite and that $\cup_{m< n+1}E_{m}=E_{n+1}\cup R$. Apply $3$ to $E_{n+1}$, to see that $Z=\omega\setminus((\cup_{m\leq n+1}E_{m})\cup\{x_k,y_k\mid k\leq n\})$ is infinite, and  pick $x_{n+1},y_{n+1}\in Z$ distinct.
Clearly $$|\{x_n\mid n<\omega\}|=|\{y_n\mid n<\omega\}|=\aleph_0\text{ and }\{x_n\mid n<\omega\}\cap\{y_n\mid n<\omega\}=\emptyset.$$ Let $D_{x,\alpha}=\omega\setminus\{x_n\mid n<\omega\}$ and $D_{y,\alpha}=\omega\setminus\{y_n\mid n<\omega\}$. We claim that for every $n<\omega$, $E_n\subseteq^* D_{x,\alpha},D_{y,\alpha}$. By symmetry it suffices to show it for $D_{x,\alpha}$. If $r\in E_n\setminus D_{x,\alpha}$, then there is $m$ such that $r=x_m$, since for every $m\geq n$, $x_m\notin E_n$, it follows that $m<n$. Thus $E_n\setminus D_{x,\alpha}\subseteq \{x_m\mid m<n\}$, implying $E_n\subseteq^* D_{x,\alpha}$. Let us argue that either for every $n<\omega$, $|D_{x,\alpha}\setminus E_n|=\omega$, or for every $n<\omega$, $|D_{y,\alpha}\setminus E_n|=\omega$. Assume otherwise, so there is $n<\omega$ such that $D_{x,\alpha}=^*E_n$  and there is $k<\omega$ such that $D_{y,\alpha}=^*E_k$. For every $n\leq m<\omega$, $$D_{x,\alpha}=^*E_n\subseteq^* E_m\subseteq^*D_{x,\alpha}.$$ Hence $E_m=^*D_{x,\alpha}$. In the same way we see that for every $k\leq m<\omega$, $E_m=^*D_{y,\alpha}$. Let $m>\max\{n,k\}$. Then $D_{y,\alpha}=^*E_m=^*D_{x,\alpha}$, contradiction.

Without loss of generality, assume that for every $n<\omega$, $|D_{x,\alpha}\setminus E_n|=\omega$. Define $D_\alpha=D_{x,\alpha}$. Let us prove $(1),(2),(3)$. To see $(1)$, for each $\beta<\alpha$ find $n<\omega$ such that $\beta<\alpha_n$, then $D_\beta\subseteq^* D_{\alpha_n}=E_n\subseteq^* D_\alpha$. Also $D_{\alpha}\setminus D_{\alpha_n}\subseteq (D_{\alpha}\setminus D_{\beta})\cup( D_{\beta}\setminus D_{\alpha_n})$. Since $|D_{\alpha}\setminus D_{\alpha_n}|=\omega$ and $|D_{\beta}\setminus D_{\alpha_n}|<\omega$ it follows that $|D_{\alpha}\setminus D_{\beta}|=\omega$, so $(2)$ holds. Finally, $(3)$ follows since $\{x_n\mid n<\omega\}\subseteq\omega\setminus D_\alpha$.
 
\end{remark}
\textit{Proof of \ref{Modfinitestab}.}
Toward a contradiction, assume that the theorem fails, then by regularity of $\mu$, there is $Y\subseteq \mu$ such that $|Y|=\mu$ and for every $\alpha,\beta\in Y$, if $\alpha<\beta$ then $D_\alpha\subseteq^* D_\beta$ and $|D_\beta\setminus D_\alpha|\geq\omega$. For every $\xi<\kappa$, find $E_{\xi}\subseteq \xi$ such that the set
$$X_\xi:=\{\nu<\mu\mid D_\nu\cap \xi=E_\xi\}$$
is unbounded in $\mu$, set $\alpha_\xi:=\min(X_\xi)$. Since $D_\alpha$ is $\subseteq^*$-increasing,  for every $\alpha_\xi\leq\alpha<\mu$, 
$D_\alpha\cap \xi=^* E_\xi$. To see this, find $\beta\in X_\xi$ such that $\alpha_\xi\leq\alpha\leq\beta$, then
$D_{\alpha_\xi}\subseteq^* D_{\alpha}\subseteq^* D_{\beta}$
Hence
$$E_\xi=D_{\alpha_\xi}\cap \xi\subseteq^* D_{\alpha}\cap \xi\subseteq^* D_{\beta}\cap \xi=E_\xi.$$
Set $\alpha^*=\sup\{\alpha_i\mid i<\lambda\}$, by regularity, $\alpha^*<\mu$.
It follows that
$$(*) \ \ \text{For every }\delta<\lambda\text{ and every }\alpha^*\leq \beta_1<\beta_2<\mu. \ D_{\beta_1}\cap\delta=^*E_{\delta}=^* D_{\beta_2}\cap\delta.$$
and that
$$(**) \ \ \ \ \ \text{For every }\alpha^*\leq \beta_1<\beta_2<\mu. \ \ \  |D_{\beta_1}\Delta D_{\beta_2}|\leq\omega. \ \ \ \  \ \ \ \ \ \ \ \ \ $$
To see $(**)$, assume otherwise, then there are $\beta_1,\beta_2$ such that $|D_{\beta_1}\Delta D_{\beta_2}|\geq\omega_1$.  Thus there is $\delta<\lambda$ such that $|D_{\beta_1}\cap\delta\Delta D_{\beta_2}\cap\delta|\geq\aleph_0$
contradicting $(*)$.

Also $cf(\lambda)=\aleph_0$, since for any distinct $\beta_1,\beta_2\in Y\setminus\alpha^*$, $|D_{\beta_1}\Delta D_{\beta_2}|\geq\aleph_0$, and by $(**)$, $|D_{\beta_1}\Delta D_{\beta_2}|\leq\aleph_0$ so by Cantor–Bernstein $|D_{\beta_1}\Delta D_{\beta_2}|=\aleph_0$. Since $\beta_1,\beta_2>\alpha^*$, $D_{\beta_1}\Delta D_{\beta_2}$ cannot be bounded, hence   $cf(\lambda)=\aleph_0$. 

Let $\chi:=(2^{<\aleph_1})^+=(2^{\aleph_0})^+$ . Since $\lambda>\aleph_0$ is strong limit $\chi<\lambda<\mu$. Fix any $X\subseteq Y\setminus\alpha^*$ such that $|X|=\chi$.
Define a partition $f:[X]^2\rightarrow \omega$:

Let $\l \eta_n\mid n<\omega\r$ be cofinal in $\lambda$. 
For any $i<j$ in $X$, $D_i\subseteq^* D_j$, hence there is $n_{i,j}<\omega$ such that $(D_{i}\setminus \eta_{n_{i,j}})\subseteq (D_{j}\setminus \eta_{n_{i,j}})$. Simply pick some  $\eta_{n_{i,j}}$ above all the finitely many elements in $D_{i}\setminus D_{j}$.  Then set
$$f(i,j)=n_{i,j}.$$
Apply the Erd\"{o}s-Rado theorem and find $I\subseteq X$ such that $\otp(I)=\omega_1+1$ which is homogeneous with color $n^*<\omega$.  This means that for any $i<j$ in $I$, $D_{i}\setminus\eta_{n^*}\subseteq D_{j}\setminus \eta_{n^*}$. Recall that $i,j\in Y\setminus\alpha^*$, then by $(*)$, it follows also that $(D_j\setminus \eta_{n^*})\setminus (D_i\setminus \eta_{n^*})$ is infinite.

Let $\langle i_\rho\mid \rho<\omega_1+1\rangle$ be the increasing enumeration of $I$. We will prove that $|D_{i_{\omega_1}}\setminus D_{i_0}|\geq\omega_1$, and since $i_0,i_{\omega_1}\geq\alpha^*$, this is a contradiction to $(**)$. 

Indeed, for every $r<\omega_1$, pick any $\delta_r$ from the infinite set $(D_{i_{r+1}}\setminus\eta_{n^*})\setminus (D_{i_r}\setminus \eta_{n^*})$.
Since the sequence $\l D_{i_r}\setminus \eta_{n^*}\mid r\leq\omega+1\r$ is $\subseteq$-increasing, for every $\beta\leq r<\alpha\leq\omega_1$, $\delta_r\in D_{i_{\alpha}}\setminus D_{i_\beta}$. 

In particular, for every $r<\omega_1$, $\delta_r\in D_{i_{\omega_1}}\setminus D_{i_0}$ so the map $r\mapsto \delta_r$ is well defined from $\omega_1$ to $D_{i_{\omega_1}}\setminus D_{i_0}$. Also if $r_1<r_2<\omega_1$, then $\delta_{r_2}\notin D_{i_{r_1}+1}$ and $\delta_{r_1}\in D_{i_{r_1+1}}$ so $\delta_{r_1}\neq \delta_{r_2}$. Thus we found an injection of $\omega_1$ to $D_{i_{\omega_1}}\setminus D_{i_0}$,  contradicting $(**)$.$\blacksquare$

\subsection{Fat Trees}
In case $o^{\vec{U}}(\kappa)$ is for example $\omega_1$, the strong Prikry property for $\Mfor$ ensures that given $p\in\Mfor$ and a dense open set $D\subseteq\Mfor$, there is a choice of measures $U(\kappa_1,i_1),...,U(\kappa_n,i_n)$ where $\kappa_1\leq...\leq\kappa_n\leq\kappa$ and a direct extension $p\leq^*p^*$ such that for every choice $\vec{\alpha}\in A_1\times...\times A_n $ from the typical sets associated to $U(\kappa_1,i_1),...,U(\kappa_n,i_n)$,  $p^{*\smallfrown}\l\alpha_1,...,\alpha_n\r\in D$. This means that in the ground model we can determine measures which are necessary to enter $D$. 

For higher order of $\kappa$ this is no longer the case. For example, assume that $o^{\vec{U}}(\kappa)=\kappa$ and consider the first element of $C_G$ i.e. $C_G(0)$. Since $\otp(C_G)=\kappa$, consider $C_G(C_G(0))$. Let $\lusim{x}$ be such that $\Vdash_{\Mfor} \lusim{x}=C_{\lusim{G}}(C_{\lusim{G}}(0))$. Consider any condition of the form $p=\l \l\kappa,A\r\r$. There is no choice of measures in the ground model and no direct extension of $p$ which determines $\lusim{x}$. Instead, we can construct a tree $T$ with two levels. The first level is simply all the ordinals which can be $C_G(0)$, namely $\Lev_1(T)=\{\alpha\in A\mid o^{\vec{U}}(\alpha)=0\}\in U(\kappa,0)$. 
Now any extension of the form $p^{\smallfrown}\alpha$ for $\alpha\in A$ forces that $C_G(0)=\alpha$, so to determine $\lusim{x}$ we only need to pick some ordinal in the set $\{\beta\in A\setminus\alpha+1\mid o^{\vec{U}}(\beta)=\alpha\}\in U(\kappa,\alpha)$. Hence we define $\succ_{T}(\l\alpha\r)=\{\beta\in A\setminus\alpha+1\mid o^{\vec{U}}(\beta)=\alpha\}$. Since the measure used in the second level is different for every choice of $\alpha$, we cannot find a single measure that will turn this tree into a product. 

This section is devoted to the study of some combinatorial aspects of such trees.
\begin{definition}\label{fat-tree}
Let $\vec{U}$ be a coherent sequence of normal measures and $\theta_1\leq...\leq\theta_n$ be measurables with $o^{\vec{U}}(\theta_i)>0$. A $\vec{U}-fat \ tree$ on $\theta_1\leq...\leq\theta_n$ is a tree $\langle T, \leq_{T}\rangle$ such that 
\begin{enumerate}
\item $T\subseteq\prod_{i=1}^n\theta_i$ and $\langle \ \rangle\in T$.
\item $\leq_{T}$ is end-extension i.e. $t\leq_{T}s \Leftrightarrow t=s\cap\max(t)+1$.
\item $T$ is downward closed with respect to end-extension.
\item For any $t\in T$ one of the following holds:
\begin{enumerate}
\item $|t|=n$.
\item $|t|<n$ and there is $\beta< o^{\vec{U}}(\theta_{|t|+1})$ such that $\{\alpha\mid t^{\frown}\langle\alpha\rangle\in T\}\in U(\theta_{|t|+1},\beta)$.
     
\end{enumerate} 
\end{enumerate}
\end{definition}
Some usual notations of trees:
\begin{enumerate}
\item $\succ_T(t)=\{\alpha\mid t^{\smallfrown}\langle\alpha\rangle\in T\}$.
\item For each $t\in T$ with $|t|<n$, choose $\xi(t)$ such that $\succ_T(t)\in U(\theta_{|t|+1},\xi(t))$, and define $U^{(T)}_t=U(\theta_{|t|+1},\xi(t))$  (We drop the superscript $(T)$ when there is no risk of confusion).
\item Note that if the measures in $\vec{U}$ can be separated i.e. there are $\langle X(\alpha,\beta)\mid \langle\alpha,\beta\rangle\in Dom(\vec{U})\rangle$ such that $X_i\in U_i\wedge \forall j\neq i X_i\notin U_j$, then we can intersect each set of the form $\succ_T(t)$ with appropriate $X_i$ and then $\xi(t)$ has a unique choice. 
\item $ht(t)=\otp(s\in T\mid s<_T t)$.
\item $\Lev_i(T)=\{t\in T\mid ht(t)=i\}$.
\item The height of a tree is $ht(T)=\max(\{n<\omega\mid \Lev_n(T)\neq\emptyset\})$. 
\item We will assume that if $\theta_i<\theta_{i+1}$ then for every $t\in \Lev_i(T)$, $\min(\succ_T(t))>\theta_i$.
\item For $t\in T$ the tree above $t$ is $T/t=\{s\in T\mid t\leq_Ts\}$. We identify $T/t$ with the $\vec{U}$-fat tree $\{s\setminus t\mid s\in T/t\}$.
\item The set of all maximal branches of $T$ is denoted by $mb(T)=\Lev_{ht(T)}(T)$. In general, we identify maximal branches of the tree with points at the top level. Note that $mb(T)$ completely determines $T$.
\item Let $J\subseteq \{0,1,...,ht(T)\}$ then $T\restriction J=\{t\restriction J\mid t\in T\}$.
\end{enumerate}

For every $\vec{U}$-fat tree $T$ in $\theta_1\leq...\leq\theta_n$ of height $n$,  define \textit{the iteration associated to $T$},  $\l j^{(T)}_{m,k},M_k\mid 0\leq m\leq k\leq n\r$, usually we drop the superscript  $T$.  Let $V=M_0$,
$$j_1=j_{0,1}:=j_{U^{(T)}_{\l\r}}:V\rightarrow Ult(V,U^{(T)}_{\l\r})\simeq M_{U^{(T)}_{\l\r}}:=M_1$$ then $crit(j_{1})=\theta_1\in j_{1}(\succ_T(\l\r))=\succ_{j_{1}(T)}(\l\r)$. Thus $\l\theta_1\r\in \Lev_1(j_1(T))$.

Assume that $\l j_{m',m},M_m\mid 0\leq m\leq m'\leq k\r$ is defined for some $k<n$, for every $1\leq i\leq k<n$, denote $\kappa_i:=crit(j_{i-1,i})=j_{i-1}(\theta_i)$ and assume $\l \kappa_1,...,\kappa_k\r\in \Lev_k(j_k(T))$. Let $$j_{k,k+1}:= j_{U^{(j_k(T))}_{\l\kappa_1,...,\kappa_k\r}}:M_k\rightarrow Ult(M_k,U^{(j_k(T))}_{\l\kappa_1,...,\kappa_k\r})\simeq M_{k+1}$$
$j_{i,k+1}=j_{k,k+1}\circ j_{i,k}$ and $j_{k+1}=j_{0,k+1}$. Note that $\succ_{j_k(T)}(\l\kappa_1,...,.\kappa_k\r)\in U^{(j_k(T))}_{\l\kappa_1,...,\kappa_k\r}$ which is a normal measure on $j_{k}(\theta_k+1)$. Thus  $$\kappa_{k+1}:=j_{k}(\theta_{k+1})=\crit(j_{k,k+1})\in j_{k,k+1}(\succ_{j_k(T)}(\l\kappa_1,...,.\kappa_k\r))=\succ_{j_{k+1}(T)}(\l\kappa_1,...,\kappa_k\r).$$ Therefore, $\l\kappa_1,...,.\kappa_k,\kappa_{k+1}\r\in \Lev_{k+1}(j_{k+1}(T))$. We denote $j_T=j_n$ and $M_T=M_n$.

More generally, a \textit{tree iteration} of $\vec{U}$-measures is a finite iteration  $\l j_{m,k},M_k\mid 0\leq m\leq k\leq n\r$ of $V$ such that for some measurable cardinals $\theta_1\leq...\leq \theta_n$, for every $0\leq m<n$, there is a normal measure $W_{m+1}\in j_m(\vec{U})$ on $j_m(\theta_{m+1})$ such that $$j_{m,m+1}=j_{W_m}: M_m\rightarrow Ult(M_m,W_m)\simeq M_{m+1}.$$
Denote $\kappa_m=j_{m-1}(\theta_m)$
and derive an  ultrafilter $U$ on $\prod_{i=1}^n\theta_i$ by the formula:
$$X\in U\longleftrightarrow \l \kappa_1, \kappa_2,...,\kappa_n\r\in j_n(X).$$
Let us verify some standard properties of such an iteration:
\begin{proposition}\label{Los for trees}
Let $\l j_{m,k},M_k\mid 0\leq m\leq k\leq n\r$ be a tree iteration of $\vec{U}$-measures. Then:
\begin{enumerate}
\item $U$ is a $\theta_1$-complete ultrafilter on $\prod_{i=1}^n\theta_i$. 
    \item For any formula $\Phi(y_1,...,y_{m})$ and any $f_1,...,f_m:\prod_{i=1}^n\theta_i\rightarrow V$,
$$M_n\models \Phi(j_n(f_1)(\kappa_1,...,\kappa_n),...,j_n(f_m)(\kappa_1,...,\kappa_n))\Leftrightarrow \{\vec{\alpha}\in \prod_{i=1}^n\theta_i\mid \Phi(f_1(\vec{\alpha}),...,f_m(\vec{\alpha}))\}\in U.$$ 
    \item Let $j_U:V\rightarrow Ult(V,U)\simeq M_U$ be the elementary embedding associated to $U$, then $M_U=M_n$ and $j_U=j_n$.
\item For every $R\in U$ there is a $\vec{U}$-fat tree $S$ such that $mb(S)\subseteq R$, $mb(S)\in U$. Moreover, if $j_{i-1}(f_i)(\kappa_1,...,\kappa_{i-1})=W_i$ (the ultrafilter used in $j_{i,i-1}$), then for every $s\in \Lev_{i-1}(S)$, $\succ_{S}(s)\in f_i(s)$.
\end{enumerate}
\end{proposition}
\pr $(1)$ is a standard consequence of the critical point of the iteration being $\theta_1$.

For $(2)$, by elementarity of $j_n$,
$$j_n(\{\vec{\alpha}\in \prod_{i=1}^n\theta_i\mid \Phi(f_1(\vec{\alpha}),...,f_m(\vec{\alpha}))\})=\{\vec{\alpha}\in \prod_{i=1}^nj_n(\theta_i)\mid M_n\models\Phi(j_n(f_1)(\vec{\alpha}),...,j_n(f_m)(\vec{\alpha}))\}$$
Note that $\kappa_i=j_{i-1}(\theta_i)=crit(j_{i,i+1})$, thus $\kappa_i<j_{i,i+1}(j_{i-1}(\theta_i))\leq j_n(\theta_i)$. By definition of $U$, 
$$M_n\models \Phi(j_n(f_1)(\kappa_1,...,\kappa_n),...,j_n(f_m)(\kappa_1,...,\kappa_n))\leftrightarrow$$ $$\leftrightarrow \l\kappa_1,...,\kappa_n\r\in j_n(\{\vec{\alpha}\in \prod_{i=1}^n\theta_i\mid \Phi(f_1(\vec{\alpha}),...,f_m(\vec{\alpha}))\})\leftrightarrow$$ $$\leftrightarrow\{\vec{\alpha}\in \prod_{i=1}^n\theta_i\mid \Phi(f_1(\vec{\alpha}),...,f_m(\vec{\alpha}))\}\in U$$

For $(3)$, it suffices to prove $M_U\simeq M_n$ via an isomorphism $k:M_U\rightarrow M_n$ such that $k\circ j_U=j_n$.
Define $k([f]_U)=j_n(f)(\kappa_1,...,\kappa_n)$. By $(2)$, $k$ is well defined and elementary embedding. Moreover, by elementarity of $j_n$, if $c_x$ is the constant function with value $x$ then $j_n(c_x)$ is constant with value $j_n(x)$. Thus,
$$k(j_U(x))=k([c_x]_U)=j_n(c_x)(\kappa_1,...,\kappa_n)=j_n(x)$$
To see $j_U$ is onto, let $x\in M_n$, since $M_n$ is the ultrapower of $M_{n-1}$ by $W_n$, there is $f_{n-1}\in M_{n-1}$,  $f_{n-1}:j_{n-1}(\theta_n)\rightarrow M_{n-1}$, such that $j_{n,n-1}(f_{n-1})(\kappa_n)=x$. Inductively, assume that $x=j_{n,i}(f_i)(\kappa_{i+1},...,\kappa_n)$, where $f_i:\prod_{k=i+1}^nj_i(\theta_k)\rightarrow M_i$. Since $M_i$ is the ultrapower of $W_{i-1}$, there is $g_{i-1}: j_{i-1}(\theta_i)\rightarrow M_{i-1}$ such that $j_{i,i-1}(g_{i-1})(\kappa_i)=f_i$. By elementarity, for every $\alpha<j_{i-1}(\theta_i)$, $g_{i-1}(\alpha):\prod_{k=i+1}^nj_{i-1}(\theta_k)\rightarrow M_{i-1}$. Define $$f_{i-1}:\prod_{k=i}^nj_{i-1}(\theta_k)\rightarrow M_{i-1}\text{ by } f_{i-1}(\alpha_i,...,\alpha_n)=g(\alpha_i)(\alpha_{i+1},...,\alpha_n)$$ Since $\kappa_i=crit(j_{i,i-1})<crit(j_{n,i})=\kappa_{i+1}$, $j_{n,i}(\kappa_i)=\kappa_i$ and $$j_{n,i-1}(f_{i-1})(\kappa_i,...,\kappa_n)=j_{n,i}(j_{i,i-1}(g_{i-1})(\kappa_i))(\kappa_{i+1},...,\kappa_n)=j_{n,i}(f_{i})(\kappa_{i+1},...,\kappa_n)=x$$
We conclude that there is $f_0:\prod_{i=1}^n\theta_i:\rightarrow V$ such that $k([f]_U)=j_n(f_0)(\kappa_1,...,\kappa_n)=x$.

To see $(4)$, let $W_{i}\in M_{i-1}$ be the ultrafilter used in $j_{i,i-1}$. Apply $(3)$, and fix for every $1\leq i\leq n$  $f_i:\prod_{k=1}^{i-1}\theta_k\rightarrow V$ such that $j_{i-1}(f_i)(\kappa_1,...,\kappa_{i-1})=W_i$. We prove $(4)$ by induction on the length of the iteration $n$. For $n=1$ we can take $S$ such that $\Lev_1(S)=R$, also  $\succ_{S}(\l\r)\in W_1=j_0(f_1)(\l\r)=f_1(\l\r)$. Assume this holds for iterations of length $i-1$. Let $R\in U$, where $U$ is derived from an iteration of length $i$. Since $R\in U$, by definition $\l\kappa_1,...,\kappa_i\r\in j_i(R)$. It follows that $\kappa_i\in j_{i,i-1}(\{\alpha<\kappa_i\mid \l\kappa_1,...,\kappa_{i-1}\r^{\smallfrown}\alpha\in R\})$. Since $j_{i,i-1}$ is the ultrapower by $W_i$,
$$(\star) \ \ Z:=\{\alpha<\kappa_i\mid \l\kappa_1,...,\kappa_{i-1}\r^{\smallfrown}\alpha\in j_{i-1}(R)\}\in W_i=j_{i-1}(f_i)(\l\kappa_1,...,\kappa_{i-1}\r)$$
Let $R'=\{\vec{\alpha}\mid \{\alpha<\theta_i\mid \vec{\alpha}^{\smallfrown}\alpha\in R\}\in f_i(\vec{\alpha})\}$, then by $(\star)$, $\l\kappa_1,...,\kappa_{i-1}\r\in j_{i-1}(R')$. Apply induction to $R'$ and $j_{i-1}$ to find $S'$, $\vec{U}$-fat tree such that $mb(S')\subseteq R'$, $\l\kappa_1,...,\kappa_{i-1}\r\in j_{i-1}(mb(S'))$ and for every $s\in \Lev_{k-1}(S')$, $\succ_{S'}(s)\in f_k(s)$. Define $$S\restriction\{1,...,i-1\}=S'\text{ and for every }s\in mb(S'), \  \succ_{S}(s)=\{\alpha<\theta_i\mid s^{\smallfrown}\alpha\in R\}$$ Clearly $mb(S)\subseteq R$ and by definition of $R'$, for every $ s\in \Lev_{i-1}(S)$, $\succ_{S}(s)\in f_i(s)$, which is a $\vec{U}$-measure over $\theta_i$. Together with the induction hypothesis, we conclude that $S$ is a $\vec{U}$-fat tree on $\theta_1\leq...\leq\theta_i$. Finally, $\l\kappa_1,...,\kappa_{i-1}\r\in j_{i-1}(mb(S'))=j_{i-1}(\Lev_{i-1}(S))$, and by elementarity, $$\succ_{j_{i-1}(S)}(\l\kappa_1,...,\kappa_{i-1}\r)=\{\alpha<\kappa_i\mid \l\kappa_1,...,\kappa_{i-1}\r^{\smallfrown}\alpha\in j_{i-1}(R)\}\in W_i$$ Hence $\kappa_i\in j_{i,i-1}(\succ_{j_{i-1}(S)}(\l\kappa_1,...,\kappa_{i-1}\r))=\succ_{j_{i}(S)}(\l\kappa_1,...,\kappa_{i-1}\r)$. It follows that $\l\kappa_1,...,\kappa_i\r\in j_i(mb(S))$ as wanted.$\blacksquare$

If $T$ is a $\vec{U}$-tree then by definition, the iteration of $T$ is a tree iteration of $\vec{U}$-measures. We denote by $U_T$ the ultrafilter derived from $j_T:V\rightarrow M_n$.
\begin{proposition}\label{the drived ultrafilter of a tree}
Let $T$ be a $\vec{U}$-fat tree on $\theta_1\leq...\leq\theta_n$. Then:
\begin{enumerate}
    \item $mb(T)\in U_T$.
    \item If $S\subseteq T$ is such that 
    \begin{enumerate}
        \item $ht(S)=ht(T)=n$.
        \item $\succ_S(\l\r)\in U^{(T)}_{\l\r}$.
        \item For every $\alpha\in \succ_S(\l\r)$, $mb(S/{\l\alpha\r})\in U_{T/{\l\alpha\r}}$ 
    \end{enumerate} 
    Then $mb(S)\in U_T$.
    \item If $S\subseteq T$ is such that 
    \begin{enumerate}
        \item $ht(S)=ht(T)=n$.
        \item $mb(S\restriction\{1,...,n-1\})\in U_{T\restriction\{1,...,n-1\}}$.
        \item For every $s\in \Lev_{n-1}(S)$, $\succ_{S}(s)\in U^{(T)}_{\l s\r}$ \footnote{If $\l U_i\mid i<\lambda\r$ is a sequence of $\lambda$-complete ultrafilters over a set $B$ and $U$ is a $\lambda$-complete ultrafilter over $\lambda$, then $U-lim_{i<\lambda}U_i$ is a $\lambda$-complete ultrafilter over $\lambda\times B$, defined by: $$ U-lim_{i<\lambda}U_i:=\Big\{ X\subseteq \lambda\times B\mid \{i<\lambda\mid \{b\in B\mid (i,b)\in X\}\in U_i\}\in U\Big\}.$$ We can inductively conclude that $$U_T=U^{(T)}_{\l\r}-lim_{\alpha_1<\theta_1} U^{(T)}_{\l\alpha\r}-lim_{\alpha_1<\theta_2}U^{(T)}_{\l\alpha_1,\alpha_2\r}-... -lim_{\alpha_{n-1}<\theta_{n-1}}U^{(T)}_{\l\alpha_1,...,\alpha_{n-2},\alpha_{n-1}\r}.$$}.

    \end{enumerate} 
    Then $mb(S)\in U_T$.
    \item If $S\subseteq T$ is such that 
    \begin{enumerate}
        \item $ht(S)=ht(T)=n$.
        \item For every $s\in S\setminus mb(S)$, $\succ_{S}(s)\in U^{(T)}_{\l s\r}$.
    \end{enumerate} 
    Then $mb(S)\in U_T$.
    \item  If $S$ is a $\vec{U}$-fat tree, and $mb(S)\in U_T$, then there is a choice of measures $U^{(S)}_s$ such that $j^{(S)}_n=j^{(T)}_n$ and in particular, $U_S=U_T$.
\end{enumerate}
\end{proposition}

\pr For $(1)$, by definition of $j_T$, we have that $\l\kappa_1,...,\kappa_n\r\in mb(j_n(T))=j_n(mb(T))$, hence by definition of $U_T$, $mb(T)\in \vec{U}(T)$.
For $(2)$, note that in $M_1$ we have the tree $j_1(T)/{\l\kappa_1\r}$. 
By $(b),(c)$ it follow that in $M_1$, $ mb(j_1(S)/{\l\kappa_1\r})\in U_{j_1(T)/{\l\kappa_1\r}}$.
By definition, the iteration defined inside $M_1$ of $j_1(T)_{\l\kappa_1\r}$ is simply the iteration $j_T$ starting from the second step inside $M_1$, namely, $\l j_{m,k}\mid 1\leq k\leq m\leq n\r$. Hence  $$\l\kappa_2,...,\kappa_n\r\in j_{n,1}(mb(j_1(S)/{\l\kappa_1\r}))=mb(j_n(S)/{\l\kappa_1\r}).$$ It follows that $\l\kappa_1,...,\kappa_n\r\in mb(j_n(S))$ and by definition $mb(S)\in U_T$.

As for $(3)$, note that $j_{T\restriction\{1,...,n-1\}}$ is by definition the first $n-1$ steps of the iteration of $j_T$. By $(b)$, $mb(S\restriction\{1,...,n-1\})\in U_{T\restriction\{1,...,n-1\}}$, thus $\l\kappa_1,...,\kappa_{n-1}\r\in \Lev_{n-1}(j_{n-1}(S))$. By $(c)$, and elementarity of $j_{n-1}$, it follows that $\succ_{j_{n-1}(S)}(\l\kappa_1,...,\kappa_{n-1}\r)\in U^{(j_{n-1}(T))}_{\l\kappa_1,...,\kappa_{n-1}\r}$, hence $\kappa_n\in j_{n,n-1}(\succ_{j_{n-1}(S)}(\l\kappa_1,...,\kappa_{n-1}\r))=\succ_{j_{n}(S)}(\l\kappa_1,...,\kappa_{n-1}\r)$. In other words, $\l\kappa_1,...,\kappa_n\r\in j_n(mb(S))$ and by definition $mb(S)\in U_T$.

For $(4)$, by induction on $i\leq n$ let us argue that $\Lev_i(S)=mb(S\restriction \{1,...,i\})\in U_{T\restriction\{1,...,i\}}$. If $i=1$ then  $\Lev_1(S)=\succ_{S}(\l\r)\in U^{(T)}_{\l\r}$. Assume that
$mb(S\restriction \{1,...,i-1\})\in\vec{U}_{T\restriction\{1,...,i-1\}}$. By $(b)$, for every $s\in \Lev_{i-1}(S)$, $\succ_{S}(s)\in U^{(T)}_{s}$, now apply $(3)$ to $S\restriction\{1,...,i\}$ and $T\restriction\{1,...,i\}$ to conclude that $mb(S\restriction \{1,...,i\})\in U_{T\restriction\{1,...,i\}}$. 

To see $(5)$, again argue by induction  on $i$ that $j^{(T)}_{i}=j^{(S)}_{i}$. Since $mb(S)\in U_T$,  $\l\kappa_1,...,\kappa_n\r\in mb(j_n(S))$, hence $\kappa_1\in \Lev_1(j_n(S))$. Since $crit(j_{1,n})=\kappa_2$, $\kappa_1\in \Lev_1(j_1(S))$, and therefore $\Lev_1(S)\in U^{(T)}_{\l\r}$, choose $U^{(S)}_{\l\r}=U^{(T)}_{\l\r}$ which implies that $j^{(T)}_{0,1}=j^{(S)}_{0,1}$. Assume that $j^{(T)}_i=j^{(S)}_i=j_i$. Since $\kappa_{i+1}\in \succ_{j^{(T)}_n(S)}(\l\kappa_1,...,\kappa_i\r)$ then $\kappa_{i+1}\in j^{(T)}_{i+1,i}(\succ_{j_i(S)}(\l\kappa_1,...,\kappa_i\r))$ thus $$(*)\ \ \ \ \succ_{j_i(S)}(\l\kappa_1,...,\kappa_i\r)\in U^{(j_i(T))}_{\l\kappa_1,...,\kappa_i\r}.$$ Back in $V$, for every $s\in \Lev_i(S)$, if $\succ_S(s)\in U^{(T)}_s$, let $U^{(S)}_s=U^{(T)}_s$, otherwise, we pick a random ultrafilter. Then by $(*)$, and elementarity $U^{(j_i(S))}_{\l\kappa_1,...,\kappa_i\r}=U^{(j_1(T))}_{\l\kappa_1,...,\kappa_i\r}$ hence $j^{(T)}_{i+1}=j^{(S)}_{i+1}$.
$\blacksquare$

The following lemma is a generalization of a combinatorial property that was proven in \cite{TomMoti} for product of measures. It can be stated for more general trees, however, let us restrict the attention to our needs.

\begin{lemma}\label{StabTree}
Let $\vec{U}$ be a sequence of normal measures and let $T$ be a $\vec{U}$-fat tree on $\theta_1\leq\theta_2\leq...\leq\theta_n$. For  $f:mb(T)\rightarrow \theta_1$ regressive i.e. $f(t)<\min(t)$ there is a $\vec{U}$-fat tree $T'\subseteq T$ such that $mb(T')\in U_T$ and  $f\restriction mb(T')=const$.
\end{lemma}
\pr By induction on the height of the tree. If $ht(T)=1$ it is the case of one normal measure, namely $U_{\langle\rangle}$, which is well known. Assume the lemma holds for $n$ and fix $T,f$  such that $ht(T)=n+1$. For $\vec{\alpha}\in \Lev_{n}(T)$ consider $\succ_T(\vec{\alpha})\in U^{(T)}_{\vec{\alpha}}$. Define $f_{\vec{\alpha}}:\succ_T(\vec{\alpha})\rightarrow\theta_1$ by $f_{\vec{\alpha}}(\beta)=f(\vec{\alpha}^{\frown}\beta)$. Then there exist $H_{\vec{\alpha}}\in U_{\vec{\alpha}}$ homogeneous for $f_{\vec{\alpha}}$ with color $c_{\vec{\alpha}}<\min(\vec{\alpha})$. 
Consider the regressive function $$g:mb(T\restriction \{1,...,n\})\rightarrow\theta_1
\  \ \ g(\vec{\alpha})=
c_{\vec{\alpha}}.$$
Since $ht(T\restriction \{1,...,n\})=n$ we can apply the induction hypothesis to $g$, so let $T'\subseteq T\restriction \{1,...n\}$ be such that $mb(T')\in U_{T\restriction\{1,...,n\}}$ be a homogeneous $\vec{U}$-fat tree with color $c^*$. Extend $T'$ by adjoining $H_{\vec{\alpha}}$ as the successors of $\vec{\alpha}\in mb(T')$, denote the resulting tree by $T^*$. Note that by the induction, $T^*\subseteq T$ is a $\vec{U}$-fat tree with $ht(T^*)=n+1$, and by \ref{the drived ultrafilter of a tree}(3) $mb(T^*)\in U_T$ and $f\restriction mb(T^*)$ is constantly $c^*$. 
 $\blacksquare$

In what come next, we will generalize (Corollary \ref{corollary for important coordinates}) a well known combinatorial property of normal measures (Corollary \ref{important coordinates for single measure}), which is a consequence of \textit{Weak compactness} of normal measures:
\begin{proposition}[folklore]
Let $U$ be a normal ultrafilter over $\kappa$, and $f:[A]^2\rightarrow \{0,1\}$ such that $A\in U$. Then there is $A'\subseteq A$ such that $A'\in U$ and $f\restriction [A']^2$ is constant.$\blacksquare$
\end{proposition}
\begin{corollary}\label{important coordinates for single measure}
 Let $U$ be a normal ultrafilter over $\kappa$, let $X$ be an arbitrary set, and $f:A\rightarrow X$ any function such that $A\in U$. Then there is $A'\subseteq A$ such that $A'\in U$ and $f\restriction A'$ is either constant or $1-1$.
\end{corollary}
\pr Define $g:[A]^2\rightarrow\{0,1\}$ by
$$g(\alpha,\beta)=1\leftrightarrow f(\alpha)=f(\beta).$$
By weak compactness, there is $A'\subseteq A$, $A'\in U$, and $c\in\{0,1\}$ such that for every $\alpha,\beta\in A'$, $\alpha<\beta$, $g(\alpha,\beta)=c$. If $c=1$, then $f\restriction A'$ is constant and if $c=0$ then $f\restriction A'$ is $1-1$.$\blacksquare$ 

In this argument  we compare $f(\alpha),f(\beta)$ for distinct $\alpha,\beta$. It is always the case that $\alpha<\beta\vee \beta<\alpha$ hence we can think about this comparison as a function defined on $[A]^2$ which is a set in $U\times U$. One problem to generalize this argument to $\vec{U}$-fat trees is the following: Although for a given function $f:mb(T)\rightarrow X$ , a $\vec{U}$-fat tree $T$, and distinct pair $t,t'\in mb(T)$ we can identify this pair as a branch of some $\vec{U}$-fat tree $S$, $S$ might vary for different $t,t'$.

 For example, if $t=\langle \alpha_1,\alpha_2,\alpha_3\rangle$ and $t'=\langle\alpha_1',\alpha_2',\alpha_3'\rangle$ the following is a possible such interweaving:
$$\alpha_1<\alpha_1'=\alpha_2<\alpha_2'<\alpha_3'<\alpha_3$$
then we can think of  $t,t'$ as a single branch from a tree $S$ of height $5$ such that any branch $s=\l s_1,s_2,s_3,s_4,s_5\r\in mb(S)$ decomposes back to $t=\l s_1,s_2,s_5\r$ ans $t'=\l s_2,s_3,s_4\r$.  However there can be different interweaving of $t,t'$ for which we need a different tree.

Generally, if $t:=\l\alpha_1,...,\alpha_n\r, \ t':=\l\alpha'_1,...,\alpha'_n\r\in mb(T)$, the set $\{\alpha_1,...,\alpha_n\}\cup \{\alpha'_1,...,\alpha'_n\}$ naturally orders in one of finitely many ways and induces an interweaving of $t,t'$:
\begin{definition}\label{interweaving}
   $p$ is an \textit{interweaving of $T$}, if it is a pair of order embedding $\langle g,g'\rangle$ where $g,g':ht(T)\rightarrow \{1,...,k\}$ so that $Im(g)\cup Im(g')=\{1,...,k\}$. Denote $A_p=Im(g), A_p'=Im(g')$ and $k=|p|$.
\end{definition} 
Let $T$ be a $\vec{U}$-fat tree on $\theta_1\leq...\leq \theta_n$. For every interweaving $p=\l g,g'\r$, define the iteration associated with $p$, $j_p=j^{(T)}_p$:

The length of the iteration is $|p|$. Let $M_0=V$ and $j_0=Id$. Assume that we are at the $m$th step of the iteration and denote the critical points $\kappa_1,...,\kappa_m$. Also assume inductively that $$\langle k_i\mid i\in A_p\cap \{1,...,m\}\rangle\in j_m(T),\langle k_i\mid i\in A'_p\cap \{1,...,m\}\rangle\in j_m(T).$$ If $m+1\in A_p\setminus A_p'$, let $r\leq ht(T)$ be such that $g(r)=m+1$. Then $\langle k_i\mid i\in A_p\cap \{1,...,m\}\rangle\in \Lev_{r-1}(j_m(T))$ and the ultrafilter $\vec{U}^{(j_m(T))}_{\langle \kappa_i\mid i\in A_p\cap \{1,...,m\}\rangle}$ which is an ultrafilter over $j_m(\theta_{r})$ is defined in $M_m$ and
for every $i\in A_p\cap\{1,...,m\}$,  $\kappa_i<j_m(\theta_r)$. 

If there is $i\in A'_p\cap\{1,...,m\}$ such that $\kappa_i\geq j_m(\theta_r)$, then  declare that the iteration is undefined. Otherwise, perform the ultrapower of $M_m$ by $\vec{U}^{(j_m(T))}_{\langle \kappa_i\mid i\in A_p\cap \{1,...,m\}\rangle}$. It follows that $\kappa_{m+1}:=crit(j_{m,m+1})=j_m(\theta_r)$ and $$\langle \kappa_i\mid i\in A_p\cap \{1,...,m+1\}\rangle=\langle \kappa_i\mid i\in A_p\cap \{1,...m\}\rangle^{\smallfrown}\kappa_{m+1}\in j_{m+1}(T).$$ If $m+1\in A'_p\setminus A_p$ we perform the symmetric procedure.
If $m+1\in A_p\cap A'_p$, let $r,r'\leq ht(T)$ be such that $m+1=g(r)=g'(r')$
 there are two possibilities, either $$\vec{U}^{(j_m(T))}_{\langle\kappa_i\mid i\in A_p\cap \{1,...,m\}\rangle}\neq \vec{U}^{(j_m(T))}_{\langle\kappa_j\mid j\in A'_p\cap \{1,...,m\}\rangle}.$$
 In this case, declare that the iteration is undefined. Otherwise $$\vec{U}^{(j_m(T))}_{\langle\kappa_i\mid i\in A_p\cap \{1,...,m\}\rangle}= \vec{U}^{(j_m(T))}_{\langle\kappa_j\mid j\in A'_p\cap \{1,...,m\}\rangle}$$
 then $j_m(\theta_r)=j_m(\theta_{r'})$ and perform the ultrapower with this measure. Thus for every $i\leq m$, $$\kappa_{m+1}:=crit(j_{m,m+1})=j_m(\theta_r)=j_m(\theta_{r'})>\kappa_i$$  and $$\langle \kappa_i\mid i\in A_p\cap \{1,...,m+1\}\rangle\in j_{m+1}(T), \ \langle \kappa_i\mid i\in A'_p\cap \{1,...,m+1\}\rangle\in j_{m+1}(T).$$
 In any case we denote $\theta(m)=\theta_r$ so that $\kappa_m=j_{m-1}(\theta(m))$, by construction, if $m=g(r)$ then $\theta(m)=\theta_r$ and if $m=g'(r')$ then $\theta(m)=\theta_{r'}$. If $j_p$ is defined then
 $$\theta(1)<j_1(\theta(2))<...<j_{|p|-1}(\theta(|p|))$$
 and since $j_{m-1}(\theta(m))=crit(j_{m,m-1})$,  $\theta(1)\leq\theta(2)\leq...\leq\theta(|p|)$. It follows that $j_p$ is a tree iteration of $\vec{U}$-measures.
 
 \begin{proposition}\label{TreeUlt}
Let $T$ be a  $\vec{U}$-fat tree and fix an interweaving $p=\l g,g'\r$ such that $j_p$ is defined. Then
\begin{enumerate}
\item There is a $\vec{U}$-fat tree, $S_p$, with $ht(S_p)=|p|$ and for every $s\in mb(S_p)$, $s\restriction A_p,s\restriction A'_p\in mb(T)$ interweave as $p$.  Moreover, for every $r\in \Lev_m(S_p)$, if $m\in A_p$ then $U^{(S_p)}_r=U^{(T)}_{r\restriction A_p\cap \{1,...,m\}}$ and if $m\in A_p'$ then $U^{(s_p)}_r=U^{(T)}_{r\restriction A'_p\cap \{1,...,m\}}$.
\item We can shrink $T$ to $R$ such that $mb(R)\in \vec{U}_T$ and if $t,t'\in mb(R)$ interweave as  $p$ then $t\cup t'\in S_{p}$
\item  If $g'(1)<g(1)$, then we can shrink $T$ to $R$ such that $mb(R)\in \vec{U}_T$ and for every $t\in mb(R)$ and $\alpha\in \succ_R(\langle\rangle)\cap\min(t)$ there is $t'\in mb(T)$ such that $t,t'$ interweave as $p$ and $\min(t')=\alpha$.

\end{enumerate}
\end{proposition}

\pr
For $(1)$, if the iteration $j_p$ is defined, then in particular for every $m$,  $j_{m,m+1}$ is the ultrapower by $U^{(j_{m}(T))}_{\l \kappa_i\mid i\in A_p\cap\{1,...,m\}\r}$ or by $U^{(j_{m}(T))}_{\l \kappa_i\mid i\in A'_p\cap\{1,...,m\}\r}$ which is a measure over $j_m(\theta_{r_{m+1}})$ for some $r_{m+1}\leq ht(T)$. Since $j_p$ is defined, we can derive the ultrafilter $U_p$ from $j_p$ over $\prod_{i=1}^{|p|}\theta(i)$. In $M_{|p|}$ we have that
$$\l \kappa_1,...,\kappa_{|p|}\r\restriction A_p, \l \kappa_1,...,\kappa_{|p|}\r\restriction A'_p\in mb(j_p(T))\text{ interweave as }p.$$
Then by \ref{Los for trees}(2), $R=\{\vec{\alpha}\in\prod_{i=1}^{|p|}\theta(i)\mid \vec{\alpha}\restriction A_p,\vec{\alpha}\restriction A'_p\in mb(T)\text{ interweave as } p\}\in\vec{U}_p$. By construction of $j_p$, if $m\in A_p$ then the function $f_m(t)=U^{(T)}_{t\restriction A_p\cap\{1,...,m-1\}}$ satisfies that the measure $j_{m-1}(f_m)(\l\kappa_1,...,\kappa_{m-1}\r)$ is the one applied at the $m$-th step of the iteration. If $m\in A_p'$ define a similar function $f'_m$ depending on $t\restriction A_p'\cap\{1,...,m-1\}$.  By \ref{Los for trees}(4), there is a $\vec{U}$-fat tree $S_p$ such that $mb(S_p)\subseteq R$ and $mb(S_p)\in U_p$. Then any $s\in mb(S_p)$ is in $R$ and therefore $s\restriction A_p,\ s\restriction A_p'$ interweave as $p$. Moreover, for every $r\in \Lev_{m-1}(S_p)$, $U^{(S_p)}_r=f_m(r)=U^{(T)}_{r\restriction A_p\cap\{1,...,m-1\}}$ or $U^{(S_p)}_r=f'_m(r)=U^{(T)}_{r\restriction A'_p\cap\{1,...,m-1\}}$.

To see $(2)$, for every $\vec{\alpha}\in \Lev_{m+1}(S_p)$ define $t(\vec{\alpha})\in T$ to be $\vec{\alpha}\restriction A_p\cap \{1,...,m\}$ and $t'(\vec{\alpha})=\vec{\alpha}\restriction A_p'\cap \{1,...,m\}$. From $(1)$ it follows that if $m+1\in A_p$ then $\succ_{S_p}(\vec{\alpha})\in U^{(T)}_{t(\vec{\alpha})}$ and similarly for $m+1\in A_p'$. Define $R$ inductively, the levels of $S_p$ which correspond to the first level are $g(1)$ and $g'(1)$ are the successors of nodes at levels $g(1)-1$ and $g'(1)-1$. Note that at least one of $g(1),g'(1)$ must be $1$. Also note that for every $\vec{\alpha}\in \Lev_{g(1)-1}(S_p)$, $t(\vec{\alpha})=\l\r$ and that for every $\vec{\beta}\in \Lev_{g'(1)-1}(S_p)$, $t'(\vec{\beta})=\l\r$. Define $$B_{\l\r}=\Delta_{\vec{\alpha}\in \Lev_{g(1)-1}(S_p)}\succ_{S_p}(\vec{\alpha}), C_{\l\r}=\Delta_{\vec{\alpha}\in \Lev_{g'(1)-1}(S_p)}\succ_{S_p}(\vec{\alpha})\in U^{(T)}_{\langle\rangle}.$$
Let $\succ_{R}(\langle\rangle)=B_{\l\r}\cap C_{\l\r}\in U^{(T)}_{\l\r}$. Moreover, at least one of $B_{\l\r},C_{\l\r}$ is simply $\succ_{S_p}(\l\r)$. 

Assume $r\in \Lev_{m}(R)$ is defined, the levels of $S_p$ which correspond to the $m$th level are $g(m),g'(m)$ (Which might be the same level), thus for every $\vec{\alpha}\in \Lev_{g(m)}(S_p)$, $t(\vec{\alpha})\in \Lev_m(T)$ and for every $\vec{\beta}\in \Lev_{g'(m)}(S_p)$, $t'(\vec{\beta})\in \Lev_m(T)$. Define
$$\succ_R(r)=\underset{\vec{\alpha}\in \Lev_{g(m)}(S_p), t(\vec{\alpha})=r}{\Delta}\succ_{S_p}(\vec{\alpha})\cap
\underset{\vec{\alpha}\in \Lev_{g(m)}(S_p), t'(\vec{\alpha})=r}{\Delta}\succ_{S_p}(\vec{\alpha})\in U^{(T)}_r.$$
By \ref{the drived ultrafilter of a tree}(4) $mb(R)\in U_T$. If $t,t'\in mb(R)$ interweave as $p$, we prove inductively that $(t\cup t')\restriction\{1,...,k\}\in \Lev_k(S_p)$. Clearly $(t\cup t')\restriction \{1\}=\l\alpha\r\in \Lev_1(S_p)$, as $\alpha\in B_{\l\r}\cap C_{\l\r}\subseteq\succ_{S_p}(\l\r)$. Assume that $(t\cup t')\restriction\{1,...,k\}\in \Lev_k(S_p)$, if $k+1\in A_p$, let $r$ be such that $g(r)=k+1$, then $(t\cup t')(k+1)=t(r)>(t\cup t')(k)$. Also $t((t\cup t')\restriction\{1,...,k\})=t\restriction \{1,...,r-1\}$. By definition of diagonal intersection and $R$ it follows that $$t(r)\in \succ_{R}(t\restriction \{1,...,r-1\})\subseteq\succ_{S_p}((t\cup t')\restriction\{1,...,k\}))$$ hence $(t\cup t')\restriction\{1,...,k+1\}\in \Lev_{k+1}(S_p)$. The case where $k+1\in A_p'$ is similar.

To see $(3)$, suppose that $g'(1)<g(1)$. Define a sequence inductively, let  $\vec{\eta}_1=\langle\beta_1,...,\beta_{g(1)-1}\rangle\in S_p$. Then by $(1)$, $\succ_{S_p}(\vec{\eta}_1)\in U^{(T)}_{\langle\rangle}$,  thus by definition of $j_T$, 
$$\kappa_1\in j_1(\succ_{S_p}(\vec{\eta}_1))=\succ_{j_1(S_p)}(\vec{\eta}_1).$$ Consider $\vec{\eta}_1^{\frown}\langle \kappa_1\rangle\in \Lev_{g(1)}(j_1(S_p))$, pick any $\vec{\eta}_2$ such that $\vec{\eta}_1^{\frown}\langle \kappa_1\rangle^{\frown}\vec{\eta}_2\in \Lev_{g(2)-1}(j_1(S_p))$,
then $$\succ_{j_1(S_p)}(\vec{\eta}_1^{\frown}\langle \kappa_1\rangle^{\frown}\vec{\eta}_2)\in j_1(\vec{U})^{(j_1(T))}_{\l\kappa_1\r}\text{ thus }\kappa_2\in \succ_{j_2(S_p)}(\vec{\eta}_1^{\frown}\langle \kappa_1\rangle^{\frown}\vec{\eta}_2)$$ continuing in this fashion we end up with a witness for the statement 
$$M_n\models \exists t\in mb(j_n(T))\ s.t. \ \langle\kappa_1,...,\kappa_n\rangle,t\text{ interweave as } p.$$ Since $\beta_1\in \succ_{S_p}(\langle\rangle)=\succ_T(\langle\rangle)=\succ_{j_n(T)}(\l\r)\cap\kappa_1$ was arbitrary, it follows that 
$$M_n\models \forall\beta\in \succ_{j_n(T)}(\langle\rangle)\cap\kappa_1\exists t\in mb(j_n(T)) \  s.t. \  min(t)=\beta\wedge \langle\kappa_1,...,\kappa_n\rangle,t \text{  interweave as }p.$$  By \ref{Los for trees}(2) $$\{s\in mb(T)\mid \forall\beta\in\succ_{T}(\l\r)\cap s_1\exists t\in mb(T).min(t)=\beta\wedge s,t\text{ interweave as }p\}\in U_T.$$ By \ref{Los for trees}(4) we can find $R$ as wanted.
$\blacksquare$

 \begin{proposition}\label{undefined interweaving}
 Let $T$ be a  $\vec{U}$-fat tree, and let $p=\l g,g'\r$ be an interweaving. If $j_p$ is undefined then there is $T'\subseteq T$ such that $mb(T')\in U_T$ and every $t,t'\in mb(T')$ do not interweave as $p$.
 \end{proposition}
 \pr Let $m$ be the step of the iteration where we declared that $j_p$ is undefined. By definition, there are two cases to consider:
 
 \textbf{Case 1: Assume that $m+1\in A_p\setminus A'_p$ and there is $i\in A'_p\cap \{1,...,m\}$ such that $j_{i-1}(\theta(i))\geq j_m(\theta(m+1))$.} Then $\theta(i)>\theta(m+1)$, otherwise, $\theta(m+1)\geq \theta(i)$ hence $$j_{i-1}(\theta(m+1))\geq j_{i-1}(\theta(i))\geq j_m(\theta(m+1))\geq j_{i-1}(\theta(m+1))$$ hence $j_{i-1}(\theta(m+1))=j_{i-1}(\theta(i))$ and $\theta(m+1)=\theta(i)$. But $j_{i-1}(\theta(i))=crit(j_{i,i-1})$ $j_{i,i-1}(j_{i-1}(\theta(m+1))>j_{i-1}(\theta(m+1))$ hence $j_m(\theta(m+1))\geq j_i(\theta(m+1))>j_{i-1}(\theta(m+1))=j_{i-1}(\theta(i))$, contradiction, thus $\theta(i)>\theta(m+1)$. Let $r_1,r_2\leq ht(T)$ be such that $g(r_1)=m+1$ and $g'(r_2)=i$. Then $\theta_{r_1}=\theta(m+1)$ and $\theta_{r_2}=\theta(i)$. The tree $T'$ is obtained from $T$ by shrinking $\succ_T(t)$ for each $t\in Lev_{r_2-1}(T)$ such that $\min(\succ_{T'}(t))>\theta(m+1)$. To see that $T'$ is as wanted, assume that $s,s'\in mb(T')$ interweave as $p$. Then $$s'(r_2)=(s\cup s')(i)<(s\cup s')(m+1)=s(r_1).$$
 On the other hand, $s'(r_2)\in \succ_{T'}(s'\restriction\{1,...,r_2-1\})$, hence $s(r_1)<\theta(m+1)<s'(r_2)$, contradiction.
 
 \textbf{Case 2: Assume that} $m+1\in A_p\cap A'_p$ \textbf{and}  $U^{(j_m(T))}_{\l\kappa_i\mid i\in A_p\cap\{1,...,m\}\r}\neq U^{(j_m(T))}_{\l\kappa_i\mid i\in A'_p\cap\{1,...,m\}\r}$. 
 These are measures over $j_m(\theta(m+1)),j_m(\theta'(m+1))$ respectively. If $\theta(m+1)\neq \theta'(m+1)$, then for example $\theta(m+1)<\theta'(m+1)$ and we can shrink $T$ as in case $1$ to eliminate such an interweaving, hence assume  $\theta(m+1)=\theta'(m+1)$. Consider the first $m$ steps of the iteration $j_p$, let $A_p\cap\{1,...,m\}=\{g(1),...,g(k)\}, \ A_p'\cap\{1,...,m\}=\{g'(1),...,g'(k')\}$, then $\theta_{k+1}=\theta(m+1)=\theta'(m+1)=\theta_{k'+1}$.  Similar to  \ref{TreeUlt}(1), $$M_m\models \l\kappa_1,...\kappa_m\r\restriction \{g(1),...,g(k)\}\in \Lev_k(j_m(T)),\ \l\kappa_1,...\kappa_m\r\restriction \{g'(1),...,g'(k')\}\in \Lev_{k'}(j_m(T
)).$$
Moreover,
$$M_m\models U^{(j_m(T))}_{\l\kappa_1,...\kappa_m\r\restriction \{g(1),...,g(k)\}}\neq U^{(j_m(T))}_{\l\kappa_1,...\kappa_m\r\restriction \{g'(1),...,g'(k')\}}$$ 
  since the iteration up to $m$ is defined we can find a $\vec{U}$-fat tree $S$ such that:
  \begin{enumerate}
      \item $\l\kappa_1,...,\kappa_m\r\in j_m(mb(S))$.
      \item For every $s\in mb(S)$, $s\restriction \{g(1),...,g(k)\}\in \Lev_k(T),\ s\restriction \{g'(1),...,g'(k')\}\in \Lev_{k'}(T)$.
      \item $U^{(T)}_{s\restriction \{g(1),...,g(k)\}}\neq U^{(T)}_{s\restriction \{g'(1),...,g'(k')\}}$.
      \item $s\restriction \{g(1),...,g(k)\},\ s\restriction \{g'(1),...,g'(k')\}$ interweave as in $\l g\restriction\{1,...,k\},g'\restriction\{1,...,k'\}\r$.
      \item For every $s\in \Lev_r(S)$, let $t(s):=s\restriction A_p\cap\{1,...,r\}$ and  $t'(s):=s\restriction A'_p\cap\{1,...,r\}$. Then $U^{(S)}_s$ is either $U^{(T)}_{t(s)}$ if $r+1\in A_p$ or $U^{(T)}_{t'(s)}$ if $r+1\in A'_p$.
  \end{enumerate}  
Since $T$ mentions at most $|T\cap [\theta_{k+1}]^{<\omega}|\leq \theta_{k+1}$ measures on $\theta_{k+1}$ we can use the normality of the measures to separate them. Namely, for every $r\in T$ such that $U^{(T)}_r$ is a measure on $\theta_{k+1}$, find  $X_r\in U^{(T)}_r$ such that if $U^{(T)}_r\neq U^{(T)}_{r'}$ then $X_r\cap X_{r'}=\emptyset$. Now we shrink the tree $T$ similar to \ref{TreeUlt}(2), from $(5)$ it follows that if $j\in A_p$ then for every $\vec{\alpha}\in \Lev_{j-1}(S)$, $\succ_{S}(\vec{\alpha})\in U^{(T)}_{t(\vec{\alpha})}$ and similarly for $j\in A_p'$. Define $R\subseteq T$ inductively,  $$B_{\l\r}=\Delta_{\vec{\alpha}\in \Lev_{g(1)-1}(S)}\succ_{S}(\vec{\alpha}), C_{\l\r}=\Delta_{\vec{\alpha}\in \Lev_{g'(1)-1}(S)}\succ_{S}(\vec{\alpha})\in U^{(T)}_{\langle\rangle}.$$
Let $\succ_{R}(\langle\rangle)=B_{\l\r}\cap C_{\l\r}\in U^{(T)}_{\l\r}$. As before, at least one of $B_{\l\r},C_{\l\r}$ is simply $\succ_{S}(\l\r)\subseteq\succ_T(\l\r)$. 
Given $r\in \Lev_{j}(R)$, define
$$B_{r}=\begin{cases}\underset{\vec{\alpha}\in \Lev_{g(j)}(S),\ t(\vec{\alpha})=r}{\Delta}\succ_{S}(\vec{\alpha}) & j<k\\ \succ_T(r)\cap X_r & j=k\\ \succ_T(r) & j>k\end{cases}$$
$$C_{r}=\begin{cases}\underset{\vec{\alpha}\in \Lev_{g'(j)}(S),\ t'(\vec{\alpha})=r}{\Delta}\succ_{S}(\vec{\alpha}) & j<k'\\ \succ_T(r)\cap X_r & j=k'\\ \succ_T(r) & j>k'\end{cases}.$$
Then $B_r,C_r\in U^{(T)}_r$ and let $\succ_{R}(r)=B_r\cap C_r$.
So by \ref{the drived ultrafilter of a tree}(4) $mb(R)\in U_T$. Let us argue that $R$ is as wanted. Toward a contradiction, assume that $t,t'\in mb(R)$ interweave as $p$, in particular $t\restriction\{1,...,k\},t'\restriction\{1,...,k'\}$ interweave as $\l g\restriction\{1,...,k\},g'\restriction\{1,...,k'\}\r$, as in the proof of \ref{TreeUlt}(2), we conclude that $s=(t\restriction\{1,...,k\})\cup(t'\restriction\{1,...,k'\})\in mb(S)$ and by $(3)$ $U^{(T)}_{s\restriction \{g(1),...,g(k)\}}\neq U^{(T)}_{s\restriction \{g'(1),...,g'(k')\}}$. However, $s\restriction \{g(1),...,g(k)\}=t\restriction\{1,...,k\}$, $s\restriction \{g'(1),...,g'(k')\}=t'\restriction\{1,...,k'\}$ hence $\succ_{R}(t\restriction\{1,...,k\})\subseteq X_{t\restriction\{1,...,k\}}$ is disjoint from  $\succ_{R}(t\restriction\{1,...,k'\})\subseteq X_{t\restriction\{1,...,k'\}}$. On the other hand, $p$ impose that $t(k'+1)=t(k+1)$ is a member of the intersection, contradiction.$\blacksquare$
 
To illustrate the second problem of generalizing weak compactness, consider for example the function $f:mb(T)\rightarrow \kappa$, $f(\alpha,\beta)=\alpha$. No matter how we shrink $T$ to $S$, $f\restriction mb(S)$ will be neither constant nor $1-1$. However, we can ignore the coordinate $\beta$ and obtain a $1-1$ function. Generally, we will argue that $f$ might depend on some of the levels of the tree and the other levels can be ignored. Let us formulate this precisely:
\begin{definition}\label{the induced function}
Let $T$ be a tree of height $n$. For every $I\subseteq\{1,...,n\}$ define an equivalence relation $\sim_I$ on $mb(T)$ by $t\sim_I t'\leftrightarrow t\restriction I=t'\restriction I$.
For $f:mb(T)\rightarrow X$, \textit{the induced function} denoted by $f_I:mb(T\restriction I)\rightarrow X$ is the relation $\{\l t\restriction I,f(t)\r\mid t\in mb(T)\}$. 
\end{definition}
Clearly $f_I$ is a well defined function if and only if $f$ is constant on equivalence classes of $\sim_I$. For example, if $I=\emptyset$ and $f_\emptyset$ is well defined then $f$ is constant. 
\begin{definition}\label{ Definition of important coordinates}
  Let $T$ be a $\vec{U}$-fat tree of height $n$, and let $f:mb(T)\rightarrow B$ be any function.
  \begin{enumerate}
      \item A coordinate $i\in\{1,...,n\}$ is called an \textit{important coordinate} for $f$ if $\forall t_1,t_2\in mb(T)$, $t_1(i)\neq t_2(i)$ implies $f(t_1)\neq f(t_2)$.
      \item The \textit{set of important coordinates for $f$} is the set $$I(T,f)=\{i\in\{1,...,n\}\mid i\text{ is an important coordinate}\}.$$ We say that $I(T,f)$ is \textit{complete} if $f_{I(T,f)}$ is well defined i.e. $\forall t,t'\in mb(T). t\sim_{I(T,f)} t'$ implies $f(t)=f(t')$. Also we say that $I(T,f)$ is \textit{consistent} if for every $\vec{U}$-fat tree $S\subseteq T$ such that $mb(S)\in U_T$, $I(S,f\restriction mb(S))\subseteq I(T,f)$. 
  \end{enumerate}
\end{definition}
\begin{remark}\label{Remark important coordinates}
\begin{enumerate}
    \item The structure of the tree $T$, imposes some dependency between the levels of the tree which are not related to the function. For example, assume that $o^{\vec{U}}(\kappa)=\kappa$ and that $\l X^{(\kappa)}_i\mid i<\kappa\r$ is a discrete family for $\l U(\kappa,i)\mid i<\kappa\r$. Let $T$ be the tree of height $2$ such that:
    $\succ_T(\l\r)=X^{(\kappa)}_0$ and for every $\alpha\in\succ_T(\l\r)$, $\succ_T(\l\alpha\r)=X^{(\kappa)}_\alpha$. 
    Define the function $f:mb(T)\rightarrow \kappa$ by $f(\l\alpha,\beta\r)=\beta$. Clearly, we see that the function $f$ depends only on the second coordinate i.e. for every $\l\alpha,\beta\r,\l\gamma,\delta\r\in mb(T)$, $f(\l\alpha,\beta\r)=f(\l\gamma,\delta\r)\leftrightarrow \beta=\delta$ and $f_{\{2\}}$ is well defined. However, the structure of the tree is such that if $\alpha\neq \gamma$ then $X^{(\kappa)}_\alpha\cap X^{(\kappa)}_\gamma=\emptyset$ and $\beta\neq \delta$, which imposes that $1$ is important. Note that in this case, by definition, $I(T,f)=\{1,2\}$.
    \item If $S\subseteq T$  then $I(T,f)\subseteq I(S,f\restriction mb(S))$. Hence if $I(T,f)$ is complete then also $I(S,f\restriction mb(S))$ is complete, and if $I(T,F)$ is consistent, then $I(T,f)=I(S,f\restriction mb(S))$ and also $I(S,f\restriction mb(S))$ is consistent. 
\end{enumerate}
\end{remark}
\begin{lemma}\label{important coordinates}
Let $T$ be a $\vec{U}$-fat tree on $\theta_1\leq...\leq\theta_n$ and $f:mb(T)\rightarrow B$ where B is any set.  Then there is a $\vec{U}$-fat tree $T'\subseteq T$, with $mb(T')\in U_T$ and $I\subseteq\{1,...,ht(T)\}$ such that for any $t,t'\in mb(T')$ 
$$t\restriction I=t'\restriction I \Leftrightarrow f(t)=f(t').$$
\end{lemma}
Before proving the lemma, let us state as a corollary the generalization we desired:
\begin{corollary}\label{corollary for important coordinates}
 Let $T$ be a $\vec{U}$-fat tree on $\theta_1\leq...\leq\theta_n$ and $f:mb(T)\rightarrow B$ where B is any set. Then there is a $\vec{U}$-fat tree $T'\subseteq T$, with $mb(T')\in U_T$ such that the set of important coordinates $T^*:=I(T',f\restriction mb(T'))$ is complete and consistent. In particular $(f\restriction mb(T'))_{I^*}$ is well defined and $1-1$. 
\end{corollary}
\textit{Proof of corollary \ref{corollary for important coordinates}}. Let $I\subseteq\{1,...,n\}$ be as guaranteed by \ref{important coordinates}, then $I\subseteq I(T',f\restriction mb(T'))$. Indeed, every $i\in I$ is important, since if $t_1,t_2\in mb(T')$, $t_1(i)\neq t_2(i)$ then $t_1\restriction I\neq t_2\restriction I$, thus $f(t_1)\neq f(t_2)$. 

Therefore, $f_{I(T',f\restriction mb(T'))}$ is well defined, since for every $t_1,t_2\in mb(T')$, $t_1\restriction I(T',f\restriction mb(T'))=t_2\restriction I(T',f\restriction mb(T'))$ implies that $t_1\restriction I=t_2\restriction I$, hence $f(t_1)=f(t_2)$. We conclude that $I(T',f\restriction mb(T'))$ is complete.

To ensure consistency, we shrink $T'$ even more. For every $i\notin I(T',f\restriction mb(T'))$ if there is $R\subseteq T'$, $mb(R)\in U_T$ such that $i\in I(R,f\restriction mb(R))$, pick any such $R$ and denote it by $R_i$, otherwise let $R_i=T'$. Define $X^*=\cap_{i\notin I(T',f\restriction mb(T'))}mb(R_i)$. Clearly $mb(X^*)\in U_T$. By \ref{Los for trees}(4) there is a $\vec{U}$-fat tree $T^*$ such that $mb(T^*)\subseteq X^*$ and $mb(T^*)\in U_T$. It follows that $T^*\subseteq R_i\subseteq T'$ for every $i$. By \ref{Remark important coordinates}, $ I(T',f\restriction mb(T'))\subseteq I(T^*,f\restriction mb(T^*))$ and therefore $I(T^*,f\restriction mb(T^*))$ is also complete. To see it if consistent, let $S\subseteq T^*$, $mb(S)\in U_T$, and let $i\in I(S,f\restriction mb(S))$, then $S\subseteq T$, so by definition of $R_i$, $i\in I(R_i,f\restriction mb(R_i))$. Since $T^*\subseteq R_i$, then $I(R_i,f\restriction mb(R_i))\subseteq I(T^*,f\restriction mb(T^*))$. $\blacksquare$

\textit{Proof of Lemma \ref{important coordinates}.} Again we go by induction on $ht(T)$. For $ht(T)=1$ it is well known. Assume $ht(T)=n+1$ and fix $\alpha\in \Lev_1(T)$ and consider the function $$f_{\alpha}:mb(T/\langle\alpha\rangle)\rightarrow B \ \ f_\alpha(\vec{\beta})=f(\alpha^{\frown}\vec{\beta}).$$ By the induction hypothesis there is $T'_\alpha\subseteq T/\langle\alpha\rangle$ such that $mb(T'_\alpha)\in U_{T/{\l\alpha\r}}$ and $I_\alpha\subseteq\{2,...,n+1\}$ such that 
$$(\star) \ \ \forall t_1,t_2\in mb(T'_\alpha). t_1\restriction I_\alpha=t_2\restriction I_\alpha\leftrightarrow f_\alpha(t_1)=f_\alpha(t_2).$$
Find $H\in U^{(T)}_{\langle\rangle}$ and $I'\subseteq\{2,...,n\}$ such that $I_\alpha=I'$ for $\alpha\in H$. Let $S$ be the tree with $\Lev_1(S)=H$ and for every $\alpha\in H$, $S/{\l\alpha\r}=T'_\alpha$, then by \ref{the drived ultrafilter of a tree}(2), $mb(S)\in U_T$. It follows that for every $t,s\in mb(S)$, if $t\restriction \{1\}\cup I'=s\restriction \{1\}\cup I'$ then $$f(t)=f_{t(1)}(t\restriction\{2,...,n\})=f_{s(1)}(s\restriction \{2,...,n\})=f(s).$$
If the implication $f(t)=f(t')\rightarrow t \restriction\{1\}\cup I'=t'\restriction \{1\}\cup I'$ holds for every $t,t'\in mb(S)$, then we can take $I=I'\cup\{1\}$ and we are done. However there can still be a counter example i.e. $t,t'\in mb(S)$, such that $$t\restriction I'\cup\{1\}\neq t'\restriction I'\cup\{1\}\wedge f(t)=f(t').$$ Our strategy will be to go over all possible interweaving of counter examples  and shrink the tree $S$ to eliminate them. We will see that if we fail to do so, then we can take $I=I'$. 
Note that if $t(1)=t'(1)$ then by the construction of $S$, $t,t'$ cannot be a counter example, hence a counter example is one with  $t(1)\neq t'(1)$.

Fix any interweaving $p=\l g,g'\r$ with $g(1)\neq g'(1)$, and consider the iteration, $j_p$.
 If this iteration is undefined then by \ref{undefined interweaving} we can shrink $S$ such that we have eliminated this of kind interweaving. If the iteration is defined, compare  $j_p(f)(\langle \kappa_i\mid i\in A_p\rangle),j_p(f)(\langle\kappa_{j}\mid j\in A'_p\rangle)$. Suppose the interweaving is such that for some $i\in I'$, $g(i)\neq g'(i)$ we claim that $$(\star\star) \ \  j_p(f)(\langle \kappa_i\mid i\in A_p\rangle)\neq j_p(f)(\langle\kappa_{j}\mid j\in A'_p\rangle).$$
 Otherwise by \ref{TreeUlt}(1) find $\vec{U}$-fat tree $S_p$ such that $(\star\star)$ holds for maximal branches of $S_p$. Let $i$ be maximal such that $g(i)\neq g'(i)$, without loss of generality, suppose that $g'(i)<g(i)$. Note that  $q=g(i)\in A_p\setminus A_p'$, otherwise, if $q\in A'_p$, then for some $j>i$, $g'(j)=g(i)$ and therefore $g(j)>g(i)=g'(j)$ hence $g(j)\neq g'(j)$ contradicting the maximality of $i$.
 We construct recursively $t,r\in S_p$, so pick any element in $s\in \Lev_{q-1}(S_p)$, set $$t\restriction \{1,...,.q-1\}=s=r\restriction \{1,...,q-1\}.$$ Pick $t(q)<r(q)\in \succ_{S_p}(t)$, since $q\notin A'_p$, then $t\restriction A'_p\cap\{1,...,q\}= r\restriction A'_p\cap\{1,...,q\}$. Assume that $t\restriction\{1,...,k\},r\restriction\{1,...,k\}\in \Lev_k(S_p)$ are defined such that $$t\restriction A'_p\cap \{1,...,k\}=r\restriction A'_p\cap \{1,...,k\}.$$ If $k+1\in A'_p$ then $U^{(S_p)}_{t\restriction\{1,...,k\}}=U^{(S_p)}_{r\restriction\{1,...,k\}}$, as it depends only on $t\restriction A'_p\cap \{1,...,k\}$. Thus we can choose $$t(k+1)=r(k+1)\in \succ_{S_p}(t\restriction\{1,...,k\})\cap \succ_{S_p}(r\restriction\{1,...,k\}).$$ 
 If $k+1\in A_p\setminus A'_p$, pick $t(k+1)\in \succ_{S_p}(t\restriction\{1,...,k\})$ and $r(k+1)\in \succ_{S_p}(r\restriction\{1,...,k\})$ randomly. Note that in any case $t\restriction A_p'\cap\{1,...,k+1\}=r\restriction\{1,...,k+1\}$ Eventually we obtain $t,r\in mb(S_p)$ with $t\restriction A'_p=r\restriction A'_p=\vec{\alpha}'$ and $\min(t)=\min(r)=\min(s)$.
 Hence $t\restriction A_p,r\restriction A_p,\vec{\alpha}'\in mb(S)$, note that both $t\restriction A_p,\vec{\alpha}'$ and $r\restriction A_p,\vec{\alpha}'$ interweave as $p$. Consequently, 
 $$f(t\restriction A_p)=f(\vec{\alpha}')=f(r\restriction A_p).$$ This means we found a counter example with the same first coordinate which is a contradiction, concluding that $j_p(f)(\langle \kappa_i\mid i\in A_p\rangle)\neq j_p(f)(\langle\kappa_{j}\mid j\in A'_p\rangle)$. By \ref{TreeUlt}(1) and \ref{TreeUlt}(2) we can shrink $S$ so that for every $t,t'$ which interweaves as $p$, $f(t)\neq f(t')$, in other words, we have eliminated all counter examples which interweave as $p$. Next, consider $p$ for which  $g(i)=g'(i)$ for every $i\in I'$. If $$j_p(f)(\langle \kappa_i\mid i\in A_p\rangle)= j_p(f)(\langle\kappa_{j}\mid j\in A'_p\rangle)$$ then we can shrink $S$ so that whenever $t,t'\in mb(S)$ interweave as $p$, $f(t)=f(t')$. By \ref{TreeUlt}(3) we can shrink $S$ further to $S^*$ so that for every $t\in mb(S^*)$ and $\alpha<\min(t)$ there is $s\in mb(S)$ so that $\min(s)=\alpha\wedge t,s$ interweave as $p$. We claim that we can drop $1$ i.e. $I'=I$ is the set desired. To see this, assume that $t,t'\in mb(S^*)$. Without loss of generality, assume that $\min(t')=\alpha<\min(t)$, by the construction of $S^*$, there is $t''\in mb(S)$ such that 
 \begin{enumerate}
     \item $t,t''$ interweave as $p$.
     \item $t\restriction I= t''\restriction I$.
     \item $\min(t')=\alpha=\min(t'')$.
 \end{enumerate} 
Hence $$f(t)=f(t')\Leftrightarrow^{(1)} f(t'')=f(t')\Leftrightarrow^{(3)} t''\restriction I=t'\restriction I\Leftrightarrow^{(2)} t\restriction I=t'\restriction I.$$
 Finally if $j_p(f)(\langle \kappa_i\mid i\in A_p\rangle)\neq j_p(f)(\langle\kappa_{j}\mid j\in A'_p\rangle)$ then we shrink $S$ and eliminate counter examples which interweave as $p$. Obviously, if we went through all possible interweaving of all counter examples and eliminated them, then $I=I'\cup\{1\}$ will be as desired.
 $\blacksquare$

\begin{lemma}\label{function separation}
Let $T$ and $S$ be $\vec{U}$-fat trees on $\kappa_1\leq...\leq\kappa_n$, $\theta_1\leq...\leq\theta_m$ respectively. Suppose $F:mb(T)\rightarrow \kappa$ and $G:mb(S)\rightarrow \kappa$ are any functions such that $I:=I(T,F), \ J:=I(S,G)$ are complete and consistent. Then there exists $\vec{U}$-fat subtrees $T^*,S^*$ with $mb(T^*)\in U_T$ and $ mb(S^*)\in U_S$ such that one of the following holds:
\begin{enumerate}
 \item $mb(T^*)\restriction I=mb(S^*)\restriction J$\footnote{ Denote $mb(T)\restriction I=\{t\restriction I\mid t\in mb(T)\}$.} and $(F\restriction mb(T^*))_{I}=(G\restriction mb(S^*))_{J}$.
 \item $Im(F\restriction mb(T^*))\cap Im(G\restriction mb(S^*))=\emptyset$.
\end{enumerate}
\end{lemma}

\pr The argument is similar to product of measures version  in \cite{partOne}. Fix $F,G$, we proceed by induction on $\langle ht(T),ht(S)\rangle=:\l n,m\r$. Let us first deal with some trivial cases:

If $I=J=\emptyset$ i.e. $F,G$ are constantly $d_F,d_G$, respectively. Either $d_F\neq d_G$ and $(2)$ holds, or $d_F=d_G$ and $(1)$ holds. If $I=\emptyset$ and $j_0\in J\neq\emptyset$, then $F$ is constantly $d_F$.
If $d_F\notin Im(G)$ then $(2)$ holds, otherwise, there is $\vec{\beta}\in mb(S)$ such that $G(\vec{\beta})=d_F$,  remove $\vec{\beta}(j_0)$
from $\Lev_{j_0}(S)$ i.e. define:\begin{enumerate}
    \item $S^*\restriction \{1,...,j_0-1\}:=S\restriction  \{1,...,j_0-1\}$.
    \item For every $t\in \Lev_{j_0-1}(S)$, define $\succ_{S^*}(t):=\succ_{S}(t)\setminus\{\vec{\beta}(j_0)\}$.
    \item For every $t\in \Lev_{j_0}(S^*)$, $S^*/t:=S/t$.
\end{enumerate}   
By \ref{the drived ultrafilter of a tree}, $mb(S^*)\in U_S$.
If $\vec{\beta}'\in mb(S^*)$, then $G(\vec{\beta}')\neq d_F$, just otherwise, $\vec{\beta}'\restriction J=\vec{\beta}\restriction J$
and in particular $\vec{\beta}(j_0)=\vec{\beta}'(j_0)$, contradiction, then again $(2)$ holds. Similarly, if $J=\emptyset$ and $I\neq \emptyset$ then we can prove $(2)$. This argument includes the case that one of the trees is $\{\l\r\}$ in which case the functions are constantly $f(\l\r)$ or $g(\l\r)$. Thus we can assume that $n,m\geq 1$.
Without loss of generality, assume that $\theta_1\leq\kappa_1$.

For every $\beta\in \succ_S(\l\r)$,
consider the function\footnote{Note that if $m=1$ then $S/\l\beta\r=\{\l\r\}$ and $G_{\beta}$ is constant.} $$G_{\beta}: mb(S/\l\beta\r)\rightarrow \kappa, \ G_\beta(\vec{\beta})=G(\beta^{\smallfrown}\vec{\beta}).$$
Then for every $\beta\in \succ_S(\l\r)$, $I(S/\l\beta\r,G_\beta)\supseteq J\setminus \{1\}$. Shrink $\succ_S(\l\r)$ to stabilize $I(S/\l\beta\r,G_\beta)=J^*$. Then $J^*=J\setminus \{1\}$, since if we let $S^*$ be the tree obtained from $S$ by shrinking $\succ_S(\l\r)$, and $S^*/\l\beta\r=S/\l\beta\r$, then by \ref{the drived ultrafilter of a tree}(4) $mb(S^*)\in U_S$. By coherency $I(S^*,G\restriction mb(S^*))\subseteq J$. So if $j\in J^*$ then it follows by definition of important coordinate that $j\in I(S^*,G)$, hence $j\in J$. It follows now that for every $\beta$,  $I(S/\l\beta\r,G_\beta)$ is complete. For consistency, the argument given in corollary \ref{important coordinates} applies by shrinking $S/\l\beta\r$ if necessary. To ease notation we keep denoting the shrinked tree by $S$. 
Apply induction to $F$ and $G_\beta$, $I,J^*$,  to find
$T^\beta\subseteq T,\ S^{\beta}\subseteq S/\l\beta\r$ for which $mb(T^\beta)\in U_T,\ mb(S^{\beta})\in U_{S/\l\beta\r}$ such that one of the following holds:
\begin{enumerate}
    \item  $mb(T^\beta)\restriction I=mb(S^\beta)\restriction J^*$ and $(F\restriction mb(T^\beta))_I=(G_{\beta}\restriction mb(S^\beta))_{J^*}$.
    \item $Im(F\restriction T^\beta)\cap Im(G_{\beta}\restriction mb(S^\beta))=\emptyset$.
\end{enumerate}

Denote by $i_\beta\in\{1,2\}$ the relevant case. There is $H\subseteq \succ_S(\langle\rangle)$, $H\in U^{(S)}_{\langle\rangle}$ and $i^*\in\{1,2\}$ such that for every $\beta\in H$, $i_\beta=i^*$.
Let $S^*$ be the tree such that $\succ_{S^*}(\langle\rangle)=H$ and for every $\beta\in H$, $S^*/\l\beta\r=S^\beta\in \vec{U}_{S/\l\beta\r}$. By \ref{the drived ultrafilter of a tree}(2), $S^*\subseteq S$ and $mb(S^*)\in U_S$.

If $i^*=1$, let $T^*=\cup_{\beta\in H}T^\beta\subseteq T$ then $mb(T^*)\in U_T$.
Argue that $1\notin J$ and therefore $J^*=J$. Indeed, fix some $\beta_1<\beta_2\in H$, Pick some $t\in mb(T^{\beta_1})\cap mb(T^{\beta_2})$ (this is possible since they are both in $U_T$) then
$$t\restriction I\in ( mb(T^{\beta_1})\restriction I)\cap( mb(T^{\beta_2})\restriction I).$$ Since for every $\beta\in H$, $mb(T^\beta)\restriction I=mb(S^\beta)\restriction J^*$ there are $s_1\in mb(S^{\beta_1})$ and $s_2\in mb(S^{\beta_2})$
such that $s_1\restriction J^*=t\restriction I=s_2\restriction J^*$. Hence $\beta_1^{\smallfrown}s_1,\beta_2^{\smallfrown}s_2\in mb(S)$ and
$$G(\beta_1^{\smallfrown}s_1)=G_{\beta_1}(s_1)=(G_{\beta_1})_{J^*}(s_1\restriction J^*)=F_I(t\restriction I)=(G_{\beta_2})_{J^*}(s_2\restriction J^*)=G_{\beta_2}(s_2)=G(\beta_2^{\smallfrown}s_2)$$
we found two maximal branches  $x,y\in mb(S)$ which differ on $\{1\}$ such that $G(x)=G(y)$, by the definition of important coordinates it follows that $1\notin J$. Moreover, $mb(T^*)\restriction I=mb(S^*)\restriction J$ and that $(F\restriction mb(T^*))_I=(G\restriction mb(S^*))_J$, namely, $(1)$ holds. To see this, 
$$mb(T^*)\restriction I=\cup_{\beta\in H} mb(T^\beta)\restriction I=\cup_{\beta\in H} mb(S^\beta)\restriction J^*=\cup_{\beta\in \succ_{S^*}(\l\r)}mb(S^*_{\l\beta\r})\restriction J=mb(S^*)\restriction J.$$
Also if $\rho\in mb(T^*)\restriction I=mb(S^*)\restriction J$, there is $\beta\in H$ such that $\rho\in mb(T^{\beta})\restriction I=mb(S^\beta)\restriction J$,
hence $$(G\restriction mb(S^*))_J(\rho)=(G_{\beta}\restriction mb(S^{\beta}))_J(\rho)=(F\restriction mb(T^{\beta}))_I(\rho)=(F\restriction mb(T^*))_I(\rho).$$
Assume $i^*=2$. 

We repeat the same process, consider now $F_\alpha$ for every $\alpha\in \succ_{T}(\l\r)$, we can shrink $T$ so that $I\setminus\{1\}=I(T/{\l\alpha\r},F_\alpha)$ is complete and consistent. Apply induction to $F_\alpha,G$,
such that for every $\alpha$, we have $j_\alpha\in\{1,2\}$ which correspond to $i_\beta$. We shrink $\succ_{T}(\l\r)$ to some $W$ and stabilize $j_\alpha$. If $j^*=1$ then $1\notin I$, and we can find $S^*\subseteq S$, $T^*\subseteq T$ such that $mb(S^*)\in U_S$ and $mb(T^*)\in U_T$ such that $$mb(S^*)\restriction J= mb(T^*)\restriction I\text{ and }(F\restriction mb(T^*))_I=(G\restriction mb(S^*))_J$$ so $(1)$ holds. Assume that $j^*=2$.

\textbf{Case 1:
Assume $\theta_1<\kappa_1$}. shrink $\succ_{T}(\l\r)$ so that $\min(\succ_{T}(\l\r))>\theta_1$.
Since $U_T$ is $\kappa_1$-complete and $|H|=\theta_1$, $\cap_{\beta\in H}mb(T^{\beta})\in U_T$. By \ref{TreeUlt}(4) there is a $\vec{U}$-fat tree $T^*$ such that $mb(T^*)\in U_T$ and $mb(T^*)\subseteq \cap_{\beta\in H}mb(T^{\beta})$ in particular $T^*\subseteq T$.
It follows that
$$(\star) \ \ \ \forall t\in mb(T^*)\forall s\in mb(S^*).  F(t)\neq G(s).$$
To see this, note that $s(1)\in \succ_{S^*}(\l\r)=H$, $t\in mb(T^{s(1)})$ and $s\restriction\{2,...,n\}\in mb(S^{s(1)})$. Since $i^*=2$, $Im(F\restriction mb(T^{s(1)}))\cap Im(G_\beta\restriction mb(S^{s(1)}))=\emptyset$, hence $F(t)\neq G_{s(1)}(s\restriction\{2,...,n\})=G(s)$.

\textbf{Case 2: Assume that $\theta_1=\kappa_1$}. Shrink the trees $T$ and $S$ in the following way:
$\succ_{T'}(\l\r)=\Delta_{\beta\in H} \succ_{T^{\beta}}(\l\r)\in U^{(T)}_{\l\r}, \ \succ_{S'}(\l\r)=\Delta_{\alpha\in W}\succ_{S^{\alpha}}(\l\r)\in U^{(S)}_{\l\r}$.
Also for every $\alpha\in \succ_{T'}(\l\r)$, find a $\vec{U}$-fat tree $T'/\l\alpha\r$ such that $mb(T'/\l\alpha\r)\subseteq \cap_{\beta\in H\cap\alpha} mb(T^{\beta}/\l\alpha\r)$. In the same fashion for every $\beta\in \succ_{S'}(\l\r)$, find $S'/\l\beta\r$ such that $mb(S'/\l\beta\r)\subseteq \cap_{\alpha\in W\cap\beta} mb(S^{\alpha}/\l\beta\r)$. Then we claim the following:
  $$(\star\star) \ \ \ \forall t\in mb(T')\forall s\in mb(S'). t(1)\neq s(1)\rightarrow F(t)\neq G(s).$$
To see this, assume for example that $s(1)<t(1)$ (the case $t(1)<s(1)$ is symmetric), note that $s(1)\in \succ_{S^*}(\l\r)=H$, and by the definition of diagonal intersection, $t(1)\in \succ_{T^{s(1)}}(\l\r)$. Also, $t\restriction\{2,...,n\}\in mb(T^{s(1)}/\l t(1)\r)$ and therefore $t\in T^{s(1)}$. Clearly, $s\restriction\{2,...,n\}\in mb(S'/\l s(1)\r)=mb(S^{s(1)})$. Since $i^*=2$, $Im(F\restriction mb(T^{s(1)}))\cap Im(G_{s(1)}\restriction mb(S^{s(1)}))=\emptyset$, hence $F(t)\neq G_{s(1)}(s\restriction\{2,...,n\})=G(s)$.

So we are left with the situation that $s=\min(s)=\min(t)$. If $U^{(S)}_{\l\r}\neq U^{(T)}_{\l\r}$ we can shrink $\succ_{T^*}(\l\r),\succ_{S^*}(\l\r)$ so that they are disjoint, avoid this situation and conclude $(2)$.
If $U^{(T)}_{\langle\rangle}=U^{(S)}_{\langle\rangle}$, let $A=\succ_{T^*}(\langle\rangle)\cap \succ_{S^*}(\langle\rangle)$. For every $\alpha\in A$, apply the induction hypothesis to the functions $F_\alpha,G_\alpha$, $I\setminus\{1\},J\setminus\{1\}$ we obtain $T^\alpha\subseteq T/\l\alpha\r$ and $S^\alpha\subseteq S/\l\alpha\r$ such that $(1)$ or $(2)$ holds. We denote the relevant case by $r_\alpha$. Again, shrink $A$ to $A^*$ and find $r^*\in\{1,2\}$ so that for every $\alpha\in A^*$, $r_\alpha=r^*$. Define $\succ_{T^*}(\l\r)=\succ_{S^*}(\l\r)=A^*$ and for every $\alpha\in A^*$, $T^*/\l\alpha\r=T^{\alpha}$ and $S^*/\l\alpha\r=S^{\alpha}$. Clearly $T^*\subseteq T$, $S^*\subseteq S$  and $mb(T^*)\in \vec{U}_T,\ mb(S^*)\in\vec{U}_S$.

If $r^*=2$,  For every $\alpha^{\smallfrown}t\in mb(T^*),\ \alpha^{\smallfrown}s\in mb(S^*)$, we have that $r_{\alpha}=2$, then $F(\alpha^{\smallfrown}t)=F_\alpha(t)\in Im(F_\alpha\restriction mb(T^{\alpha}))$ and $G(\alpha^{\smallfrown}s)=G_\alpha(s)\in Im(G_\alpha\restriction mb(S^{\alpha}))$. By $r_{\alpha}=2$, $G(\alpha,s)\neq F(\alpha,t)$ and we have eliminated the possibility of $F(t)=G(s)$ where $\min(s)=\min(t)$, we conclude that $(2)$ holds. 

Finally, assume $r^*=1$, namely that for $I\setminus\{1\}= I^*\subseteq\{2,...,ht(T)\},J\setminus\{1\}= J^*\subseteq\{2,...,ht(S)\}$, and every $\alpha\in A^*$ $$mb(T^\alpha)\restriction I^*=mb(S^\alpha)\restriction J^*\ \wedge \ (F_\alpha\restriction mb(T^\alpha))_{I^*}=(G_{\alpha}\restriction mb(S^\alpha))_{J^*}.$$ It follows that
$$(\triangle) \ \ \  mb(T^*)\restriction I^*\cup\{1\}=\cup_{\alpha\in A^*}\{\alpha\}\times mb(T^\alpha)\restriction I^*=\cup_{\alpha\in A^*}\{\alpha\}\times mb(S^\alpha)\restriction J^*=mb(S^*)\restriction J^*\cup\{1\}.$$
Moreover, for every $\l\alpha\r^{\smallfrown}\rho\in mb(T^*)\restriction I^*\cup\{1\}$, $$(\triangle\triangle) \ \ \ (F\restriction_{ mb(T^*)})_{I^*\cup\{1\}}(\alpha,\rho)=(F_\alpha\restriction mb(T^{\alpha}))_{I^*}(\rho)=(G_\alpha\restriction mb(S^{\alpha}))_{J^*}(\rho)=(G\restriction mb(S^*))_{J^*\cup\{1\}}(\alpha,\rho).$$
If $1\notin I$ then $1$ is not an important coordinate for $F\restriction mb(T^*)$ and by definition this means that there are $t_1,t_2\in mb(T^*)$ such that $t_1(1)\neq t_2(1)$ and $F(t_1)=F(t_2)$. Then $$t_1\restriction I\in mb(T^{t_1(1)})\restriction I=mb(S^{(t_1(1)})\restriction J^*$$ 
$$t_2\restriction I\in mb(T^{t_2(1)})\restriction I=mb(S^{(t_2(1)})\restriction J^*.$$
So there are $s_1,s_2\in mb(S^*)$ such that $s_1(1)=t_1(1), s_2(1)=t_2(1)$ and $s_1\restriction J^*=t_1\restriction I,s_2\restriction J^*=t_2\restriction I$.
It follows that
$$G(s_1)=G_{s_1(1)}(s_1\restriction J^*)=F_{t_1(1)}(t_1\restriction I)=F(t_1)\neq F(t_2)=F_{t_2(1)}(t_2\restriction I)=G_{s_2(1)}(s_2\restriction J^*)=G(s_2).$$
So $1$ is not important for $G\restriction mb(S^*)$, hence $1\notin J$. In a similar way, we conclude that If $1\notin J$ then $1\notin I$.
In either case, from $(\triangle),(\triangle\triangle)$ we conclude that $(1)$ holds.
$\blacksquare$

%Let us assume that $ht(T)=ht(T')=2$, Then $f,f'$ are binary functions. We would like to compare $f(x,y)$ and $f'(x',y')$. As in the proof of \ref{important coordinates}, for every possible $x,y,x',y'$ there are only 13 many possible ways to arrange them (assuming that $x<y$ and $x'<y'$) For example:
%$x'<y'<x<y,\ x<y=x'<y',\ x<x'<y'<y,\ x'<x<y=y' ,\ x=x'<y=y'.$
%For every such arrangement $\sigma$, define a tree $T_\sigma$. Let us build an example tree for $\sigma=x<y=x'<y'$. The first level is $succ_T(\langle\rangle)$. Next consider 
\section{The proof for short sequences}
Let us return to $\Mfor$ and use the combinatorial tools developed in the last section.
\begin{definition}
 Let $p\in\Mfor$ be a condition.
 A tree of extensions of $p$ is a $\vec{U}$-fat tree $T$  on $\theta_1\leq...\leq\theta_n$, such that for every $1\leq i\leq n$, $\theta_i\in \kappa(p)$ and each $t\in T$ is a legal extension of $p$ i.e. $p^{\frown}t\in\Mfor$. Denote by $\xi(t),\kappa(t)$ the ordinals such that $\succ_{T}(t)\in U(\kappa(t),\xi(t))$.
 \end{definition}

 If $T$ is a tree of extensions of $p$ and $T'\subseteq T$ is a $\vec{U}$-fat tree such that $mb(T')\in U_T$ then $T'$ is also a tree of extensions of $p$.

Let $p^{\smallfrown}\vec{\alpha}\in\Mfor $, and for every $r\leq |\vec{\alpha}|=:n$ let $B_r\in \cap\vec{U}(\vec{\alpha}(r))$. Define $$p^{\smallfrown}\l \vec{\alpha},\vec{B}^{\vec{\alpha}}\r:=p^{\smallfrown}\l \vec{\alpha}(1),B_1\cap \vec{\alpha}(1)\r^{\smallfrown}...^{\smallfrown}\l \vec{\alpha}(n),B_{n}\cap( \vec{\alpha}(n-1),\vec{\alpha}(n))\r.$$
\begin{proposition}\label{amalgamate}
Let $T$ be a $\vec{U}$-fat tree of extensions of $p$, and let for every $t\in mb(T)$, $p_t\geq^* p^{\frown}t$ be a condition. Then there are $p^*,T^*$ and $B^s$ for $s\in T^*\setminus mb(T^*)$ such that:
\begin{enumerate}
    \item $p\leq^* p^*$.
    \item $T^*\subseteq T$ is a $\vec{U}$-fat tree of extensions for $p^*$ with $mb(T^*)\in U_T$.
    \item $B^{s}\in \cap_{\xi<\xi(s)}U(\kappa(s),\xi)$.
    \item For every $t\in mb(T^*)$
$$p_t\leq^*p^{*\smallfrown}\l t,\vec{B}^{t}\r:=p^{*\smallfrown}\l t(1), B^{\l\r}\cap t(1)\r^\smallfrown...^{\smallfrown}\l t(n),B^{t\restriction\{1,...,n-1\}}\cap t(n)\r.$$
\end{enumerate}    
\end{proposition}
\pr Assume that $T$ is on $\kappa_{j_1}(p)\leq...\leq \kappa_{j_n}(p)$, and let us proceed by induction on $ht(T)$.
If $ht(T)=1$, then for every $\alpha\in \succ_{T}(\l\r)\in U(\kappa_{j_1}(p),\xi(\l\r))$ denote $$p^{\smallfrown}\alpha\leq^* p_\alpha=\l p_\alpha\restriction \kappa_{j_1-1}(p),\l\alpha, B_\alpha\r, \l\kappa_{j_1}(p), C_\alpha\r,  p_\alpha\restriction(\kappa_{j_1}(p),\kappa]\r.$$  The order $\leq^*$ is more than $\kappa_{j_1}(p)$-closed in $\Mfor\restriction(\kappa_{j_1}(p),\kappa]$, so we can find $p^*_>\in \Mfor\restriction(\kappa_{j_1}(p),\kappa]$ such that $p_\alpha\restriction(\kappa_{j_1}(p),\kappa]\leq p^*_>$ for every $\alpha\in \succ_T(\l\r)$. For the lower part, shrink $\succ_T(\l\r)$ to $H\in U(\kappa_{j_1}(p),\xi(\l\r))$ and find $p^*_<\in \Mfor\restriction \kappa_{j_1-1}(p)$ such that for every $\alpha\in H$, $p^*_<=p_\alpha\restriction\kappa_{j_1-1}(p)$. 
Next, by normality $$C:=\Delta_{\alpha<\kappa_{j_1}(p)}C_\alpha\in \cap \vec{U}(\kappa_{j_1}(p)).$$  Use \ref{BetterSet} to find $C^*\subseteq C$ such that for every $\alpha\in C^*$, $C^*\cap\alpha\in\vec{U}(\alpha)$. As for the $B_\alpha$'s, for every $\alpha\in H$,
$B_\alpha\in \cap \vec{U}(\alpha)$. Use ineffability and shrink $H$ to $H'\in U(\kappa_{j_1}(p),\xi(\l\r))$ and find a single set $X$ such that for every $\alpha\in H'$, $X\cap\alpha=B_\alpha$, it follows that,
$B^{\l\r}:=C^*\cap X\in\cap_{j<\xi(\l\r)}U(\kappa_{j_1}(p),j)$. Set $\succ_{T^*}(\l\r)=H'\cap C^*$ and let $$p\leq^*\l p^*_{<},\l\kappa_{j_1}(p),C^*\r,p^*_>\r=:p^*.$$  
To see that $p^*,B^{\l\r},T^*$ is as wanted, let $\alpha\in \succ_{T^*}(\l\r)$. Since $\alpha\in H'$,  $B^{\l\r}\cap \alpha=B_\alpha\cap C^*\subseteq B_\alpha$. Since $\alpha\in H$,  $p_\alpha\restriction\kappa_{j_1-1}=p^*_<$ and since $\alpha\in \succ_{T}(\l\r)$,  $p_\alpha\restriction(\kappa_{j_1},\kappa]\leq^*p^*_>$. Finally note that $$B_{j_1}(p^*)\setminus\alpha+1=C^*\setminus\alpha+1\subseteq C_\alpha.$$
Thus $p_\alpha\leq^*p^{*\smallfrown}\l \alpha, B^{\l\r}\cap\alpha\r$.
Assume that $n=ht(T)>1$, then for every $t\in T\setminus mb(T)$, and for every $\alpha\in \succ_{T}(t)$, we are given some condition $p^{\smallfrown}t^{\smallfrown}\alpha\leq^*p_{t^{\smallfrown}\alpha}$. Apply the case $ht(T)=1$ to $p^{\smallfrown}t$ and $\succ_{T}(t)$ to find $p^{\smallfrown}t\leq^* p^*_t$, $\succ_{T^*}(t)$ and a set $B^{t}\in\cap_{\xi<\xi(t)}U(\kappa(t),\xi)$ such that for for every $\alpha\in \succ_{T^*}(t)$, $p_{t^{\smallfrown}\alpha}\leq^* p_t^{*\smallfrown}\l\alpha, B^t\cap\alpha\r$. 
Apply the induction hypothesis to $p, T\setminus mb(T)$, to find $p\leq^*p^*$, $T^*\subseteq T\setminus mb(T)$ and sets $B^s$ such that for every $t\in mb(T^*)$, $p^*_t\leq^* p^{*\smallfrown}\l t,\vec{B}^t\r$. Hence for every $\alpha\in \succ_{T^*}(t)$, $$p_t\leq^*p_t^{*\smallfrown}\l \alpha,B^t\cap\alpha\r\leq^*p^{*\smallfrown}\l t,\vec{B}^t\r^{\smallfrown}\l\alpha,B^t\cap\alpha\r=p^{*\smallfrown}\l t^{\smallfrown}\alpha,\vec{B}^{t^{\smallfrown}\alpha}\r.$$ It follows that $p^*$, $T^*$ and $B^{t}$ are as wanted.
$\blacksquare$ 

The following lemma is the strong Prikry property for $\Mfor$.
\begin{lemma}\label{strongPrkryProperty}
Let $D\subseteq \Mfor$ be dense open, and let $p\in\Mfor$ be any condition, then there is $p\leq^* p^*$ and a tree of extensions of $p^*$, $T$ and sets $B^s\in \cap_{\xi<\xi(s)}U(\kappa(s),\xi)$ for every $s\in T\setminus mb(T)$ such that for every  $t\in mb(T)$, $p^{*\smallfrown}\l t,\vec{B}^t\r\in D$. 
\end{lemma}
\pr
 Let $r\leq l(p)+1$, $\vec{\alpha}\in[\kappa_r(p)]^{<\omega}$, such that $p^{\frown}\vec{\alpha}\in \Mfor$ is a condition. Set
 $$A^0_r(\vec{\alpha})=\{\alpha\in B_r(p)\setminus(\max(\vec{\alpha})+1)\mid \exists q\geq^* p^{\frown}\vec{\alpha}^{\smallfrown}\l\alpha\r. \ q\in D\}, \ \ A^1_r(\vec{\alpha})=B_r(p)\setminus A^0_r(\vec{\alpha}).$$
For every $i< o^{\vec{U}}(\kappa_r(p))$, only one of $A^0_r(\vec{\alpha}),A^1_r(\vec{\alpha})$ is in $U(\kappa_r(p),i)$. Denote it by $A_{r,i}(\vec{\alpha})$ and let $C_{r,i}(\vec{\alpha})\in\{0,1\}$ be such that $A_{r,i}(\vec{\alpha})=A_r^{C_{r,i}(\vec{\alpha})}(\vec{\alpha})$. Define 
 $$A_{r,i}=\underset{\vec{\alpha}\in[\kappa_r(p)]^{<\omega}}{\Delta}A_{r,i}(\vec{\alpha})\cap B_{r}(p)\in U(\kappa_r(p),i).$$
so far $A_{r,i}$ has the property that for $\vec{\alpha}\in[\kappa_r(p)]^{<\omega}$ if $\exists\alpha\in A_{r,i}$ and $p^{\frown}\vec{\alpha}^{\smallfrown}\l\alpha\r\leq^* q\in D$ then for every $\alpha \in  A_{r,i}$ there is $p^{\smallfrown}\vec{\alpha}^{\smallfrown}\l \alpha\r\leq^*q\in D$. 

For every $\langle \alpha_1,...,\alpha_{n-1}\r\in[\kappa_r(p)]^{n-1}$, define $D_{r,i}^{(1)}(\alpha_1,...,\alpha_{n-1},*):A_{r,i}\rightarrow\{0,1\}$ by
 $$D_{r,i}^{(1)}(\alpha_1,...,\alpha_{n-1},\alpha)=0 \Leftrightarrow \exists r\leq s\leq l(p)+1\exists j<o^{\vec{U}}(\kappa_s(p)) \ C_{s,j}(\alpha_1,...,\alpha_{n-1},\alpha)=0.$$
 Find a homogeneous set for $D^{(1)}_{r,i}$, $A^{(1)}_{r,i}(\alpha_1,...,\alpha_{n-1})\in U(\kappa_r(p),i)$ with color $C^{(1)}_{r,i}(\alpha_1,...,\alpha_{n-1})$. Define
 $$A^{(1)}_{r,i}=\underset{\vec{\alpha}\in[\kappa_r(p)]^{n-1}}{\Delta} A^{(1)}_{r,i}(\vec{\alpha})\cap B_r(p)\in U(\kappa_r(p),i).$$ In similar fashion, define recursively for $k\leq n$
$$D_{r,i}^{(k)}(\alpha_1,...,\alpha_{n-k},\alpha)=0\Leftrightarrow \exists r\leq s\leq l(p)+1\exists j<o^{\vec{U}}(\kappa) \ C_{s,j}^{(k-1)}(\alpha_1,...,\alpha_{n-k},\alpha)=0.$$
find homogeneous $A_{r,i}^{(k)}(\alpha_1,...,\alpha_{n-k})\in U(\kappa_r(p),i)$ with color $C_{r,i}^{(k)}(\alpha_1,...,\alpha_{n-k})$ and let
$$A^{(k)}_{r,i}=\underset{\vec{\alpha}\in[\kappa_r(p)]^{n-k}}{\Delta} A^{(k)}_{r,i}(\vec{\alpha})\cap B_r(p)\in U(\kappa_r(p),i).$$

Eventually, set  
$$A_{r,i,n}=\underset{k\leq n}{\bigcap}A^{(k)}_i, \ A_{r,i}=\underset{n<\omega}{\bigcap}A_{r,i,n}\in U(\kappa_r(p),i)\text{ and } A_r=\underset{i<o^{\vec{U}}(\kappa_r(p))}{\bigcup}A_{r,i}.$$
 Let $p\leq^*p_1$, where $p_1$ is obtained from $p$ by shrinking $B_r(p)$ to the set obtained from \ref{BetterSet} to $A_r$ such that for every $\alpha\in B_r(p_1)$, $\alpha\cap B_r(p_1)\in\cap\vec{U}(\alpha)$. By density, there exists $p'\geq p_1$ such that $p'\in D$. 
 There is $\langle\vec{\alpha},\alpha\rangle\in [B(p^*)]^{<\omega}$ such that $p_1^{\frown}\langle\vec{\alpha},\alpha\rangle\leq^*p'$.
Find $s_1\leq...\leq s_n\leq r$, $i_j\leq o^{\vec{U}}(\kappa_{s_j}(p))$ and $k<o^{\vec{U}}(\kappa_{r}(p))$ such that $\alpha\in A_{r,k}$ and  $\vec{\alpha}=\langle\alpha_1,...,\alpha_{n-1}\rangle\in \prod^{n-1}_{j=1}A_{s_j,i_j}$.
 It follows that $A_{r,k}(\vec{\alpha})=A^0_{r,k}(\vec{\alpha})$. Hence,
$$C_{r,k}(\vec{\alpha})=0\Rightarrow D_{s_n,i_n}^{(1)}(\alpha_1,...,\alpha_{n})=0\Rightarrow C^{(1)}_{s_n,i_n}(\alpha_1,...,\alpha_{n-1})=0\Rightarrow D^{(2)}_{s_{n-1},i_{n-1}}(\alpha_1,...,\alpha_{n-1})=0\Rightarrow$$ $$ C^{(2)}_{s_{n-1},i_{n-1}}(\alpha_1,...,\alpha_{n-2})=0\Rightarrow...\Rightarrow D^{(n)}_{s_1,i_1}(\alpha_1)=0\Rightarrow C^{(n)}_{s_1,i_1}(\langle\rangle)=0.$$
Define the tree $T'$: Let $s(\l\r)=s_1$, $\xi(\l\r)=i_1$ and define  $$\succ_{T'}(\langle\rangle)=A_{s(\l\r),\xi(\l\r)}\cap B_{s(\l\r)}(p_1)\in U(\kappa_{s(\l\r)}(p),\xi(\l\r)).$$ Since $A_{s_1,i_1}\subseteq A^{(n)}_{s_1,i_1}(\l\r)$ is homogeneous, $D^{(n)}_{i_1}(x)=0$ for every $x\in A_{s_1,i_1}$. Hence, there are $\kappa_{s(x)}(r)$ and $\xi(x)$ such that $D^{(n-1)}_{s(x),\xi(x)}(x,*)$ takes the color $0$ on $A_{s(x),\xi(x)}$. Let $$\succ_{T'}(\l\alpha\r)=A_{s(\alpha),\xi(\alpha)}\cap B_{s(\alpha)}(p_1).$$ Recursively, define the other levels in a similar fashion. By \ref{frown extension}, for every $t\in mb(T')$, $p_1\leq p_1^{\smallfrown}t\in \Mfor$. Consider the function $t\in mb(T')\mapsto \l s(t\restriction 0),s(t\restriction 1),...,s(t\restriction n)\r$, then by \ref{StabTree}, we can find a $\vec{U}$-fat tree $T''\subseteq T'$, $mb(T'')\in U_{T'}$ such that $\l s(t\restriction 0),s(t\restriction 1),...,s(t\restriction n)\r$ is stabilized for $t\in mb(T'')$. 

By the construction of the tree $T''$, for every $t\in mb(T'')$ there is $p_1^{\frown}t\leq^* p_t$ such that
$p_t\in D$. By Proposition \ref{amalgamate} we can amalgamate all those $p_t$'s and find a single $p\leq^* p^*$, shrink $T''$ to $T^*$ and find $B^s$ for $s\in T^*\setminus mb(T^*)$ such that for every $t\in mb(T^*)$,  $p_t\leq^*p^{*\frown}\l t, \vec{B}^t\r$. Since $D$ is open then $p^{*\frown}\l t, \vec{B}^t\r\in D$.
$\blacksquare$

 \begin{proposition}\label{densetree}
 Let $p\in\Mfor$ be a condition, $T$ a $\vec{U}$-fat tree of extensions of $p$, and sets $B^{s}\in\cap_{\xi<\xi(s)}U(\kappa(s),\xi)$ for every $s\in T\setminus mb(T)$ such that for every $t\in mb(T)$, $p\leq p^{\smallfrown}\l t,\vec{B}^t\r\in \Mfor$. Then there are $p\leq^*p^*$, a tree $T^*\subseteq T$ of extensions of $p^*$, $mb(T^*)\in U_T$ and sets $A^s\subseteq B^s$, $A^s\in \cap_{\xi<\xi(s)}U(\kappa(s),\xi)$  such that
 $$D_{T^*,\vec{A}}:=\{p^{*    \frown}\l t,\vec{A}^t\r \mid  t\in mb(T)\}$$
 is pre-dense above $p^*$. In particular, for any generic $G$ with $p^*\in G$, $G\cap D_{T^*}\neq\emptyset$.
\end{proposition}
\pr
Assume that $T$ is on $\kappa_{j_1}(p)\leq...\leq\kappa_{j_n}(p)$ and again we argue by induction on $ht(T)$.
Assume that $ht(T)=1$, use \ref{BetterSet} to find $A_<\subseteq B^{\l\r}\cap B_{j_1}(p)$ such that $A_<\in\cap_{\xi<\xi(\l\r)} U(\kappa_{j_1}(p),\xi)$ and for every $\alpha\in A_<$, $\alpha\cap A_<\in \cap\vec{U}(\alpha)$. Consider the sets
$$A_{\xi(\l\r)}=\succ_T(\l\r)\cap B_{j_1}(p)\cap\{\alpha<\kappa_{j_1}(p)\mid A_{<}\cap\alpha\in\cap\vec{U}(\alpha)\}\in U(\kappa_{j_1}(p),\xi(\l\r))$$
$$A_>=B_{j_1}(p)\cap\{\alpha<\kappa_{j_1}(p)\mid \exists A_{\xi(\l\r)}\cap\alpha\in (\cap\vec{U}(\alpha))^+\}\in \bigcap_{\xi(\l\r)<\xi<o^{\vec{U}}(\kappa_{j_1}(p))}U(\kappa_{j_1}(p),\xi).$$
Let $p\leq^*p^*$ be the condition obtained from $p$ by shrinking $B_{j_1}(p)$ to $$B_{j_1}(p^*):=A_<\cup A_{\xi(\l\r)}\cup A_>$$ let $A^{\l\r}:=A_<$ and
shrink $\succ_{T}(\l\r)$ to $\succ_{T^*}(\l\r):=A_{\xi(\l\r)}$. Clearly, $T^*$ is a tree of extensions for $p^*$ as for every $\alpha\in \succ_{T^*}(\l\r)$, $A_<\cap \alpha\in \cap\vec{U}(\alpha)$ and $A_<\cap\alpha\subseteq B_{j_1}(p^*)\cap\alpha$. To see that $p^*,T^*,A^{\l\r}$ are as wanted, let $p^*\leq q$. Let $\vec{\alpha}$ be such that $p^{*\smallfrown}\vec{\alpha}\leq^* q$. Without loss of generality, assume that $\vec{\alpha}\in[(\kappa_{j_1-1}(p),\kappa_{j_1}(p))]^n$ and let $X_i$ denote the sets of the pairs $\l\vec{\alpha}(i),X_i\r$ and $\l\kappa_{j_1}(p),X\r$ appearing in $q$.

If $\vec{\alpha}\in[ A_<]^n$, since $X\in \cap\vec{U}(\kappa_{j_1}(p))$, then  $$X^*:=X\cap \succ_{T^*}(\l\r)\cap\{\alpha \mid  \alpha\cap X\in \cap\vec{U}(\alpha)\}\in U(\kappa_{j_1}(p)),\xi(\l\r)).$$ In particular $X^*$ is unbounded and we can find $\alpha\in X^*\setminus\max(\vec{\alpha})+1$.  It follows that $p^{*\smallfrown}\l\alpha, A^{\l\r}\cap\alpha\r\in D_{T^*,\vec{A}}$. We claim that $q,p^{*\smallfrown}\l\alpha, A^{\l\r}\cap\alpha\r\leq q'$, where
$$q'=p^{*\smallfrown}\l\vec{\alpha}(1),X_1\cap A_<\r^{\smallfrown}...^{\smallfrown}\l\vec{\alpha}(n),X_{n}\cap A_<\r^{\smallfrown}\l\alpha, X\cap A_{<}\cap\alpha\r.$$
 Indeed, for every $\beta\in A_<$, $\beta\cap A_<\in\cap\vec{U}(\beta)$. In particular for every $i$, $\vec{\alpha}(i)\cap A_<\in \cap\vec{U}(\vec{\alpha}(i))$, thus  $X_i\cap A_<\in\cap\vec{U}(\vec{\alpha}(i))$. Also by definition of $X^*$, $\alpha\cap X\in \cap\vec{U}(\alpha)$ and by definition of $\succ_{T^*}(\l\r)$, $A_<\cap\alpha\in\cap \vec{U}(\alpha)$. By \ref{frown extension}, $q\leq q'$ and $p^{*\smallfrown}\l\alpha,A^{\l\r}\cap\alpha\r\leq q'$.

If there is $j\leq n$ such that $\vec{\alpha}(j)\notin A_<$, let $r$ be the minimal such $j$. Since $\vec{\alpha}(r)\in B_{j_1}(p)$, there are two cases here, either $\vec{\alpha}(r)\in A_{\xi(\l\r)}$ or $\vec{\alpha}(r)\in A_>$.
If $\vec{\alpha}(r)\in A_{\xi(\l\r)}=\succ_{T^*}(\l\r)$, then $p^{*\smallfrown}\l\vec{\alpha}(r),A^{\l\r}\cap\alpha\r\in D_{T^*,\vec{A}}$ and we claim that $p^{*\smallfrown}\l\vec{\alpha}(r),A^{\l\r}\cap\vec{\alpha}(r)\r,q\leq q'$ where
$$q'=p^{*\smallfrown}\l\vec{\alpha}(1),X_1\cap A_<\r^{\smallfrown}...^{\smallfrown}\l\vec{\alpha}(r),A_<\cap X_r\r^{\smallfrown}\l \vec{\alpha}(r+1),X_{r+1}\r^{\smallfrown}...^{\smallfrown}\l\vec{\alpha}(n),X_{n}\r.$$
 By minimality of $r$, $\vec{\alpha}(i)\in A_<$ for every $i<r$ and the same argument as before justifies that, $X_i\cap A_<\in \cap\vec{U}(\vec{\alpha}(i))$. Since $\vec{\alpha}(r)\in A_{\xi(\l\r)}$, by definition we have that $A_<\cap\vec{\alpha}(r)\in\cap\vec{U}(\vec{\alpha}(r))$, hence $X_r\cap A_<\in \cap\vec{U}(\vec{\alpha}(r))$, then again we use \ref{frown extension}. 
Finally, if $\vec{\alpha}(r)\in A_>$, then $A_{\xi(\l\r)}\cap\vec{\alpha}(r)\in (\cap\vec{U}(\vec{\alpha}(r)))^+$. In particular $$X^*:=A_{\xi(\l\r)}\cap X_r\cap\{\alpha\mid \alpha\cap X_r\in\cap\vec{U}(\alpha)\}\in (\cap\vec{U}(\vec{\alpha}(r)))^+$$ hence there is $\alpha\in X_r\cap A_{\l\xi(\l\r)}\setminus \vec{\alpha}(r-1)+1$. 
This time, the witness for the compatibility of $p^{*\smallfrown}\l\alpha, A^{\l\r}\cap\alpha\r,q$ will be $$q'=p^{*\smallfrown}\l\vec{\alpha}(1),X_1\cap A_<\r^{\smallfrown}...^{\smallfrown}\l\vec{\alpha}(r-1),A_<\cap X_{r-1}\r^{\smallfrown}\l\alpha , X_r\cap A_<\cap\alpha\r^{\smallfrown}\l \vec{\alpha}(r),X_{r}\setminus\alpha\r^{\smallfrown}...^{\smallfrown}\l\vec{\alpha}(n),X_{n}\r.$$
This concludes the case $ht(T)=1$. Let $T$ be such that $n=ht(T)>1$, for every $s\in T\setminus mb(T)$, apply the case $n=1$ to $\succ_{T}(s)$ and the condition $p^{\smallfrown}s$ to find $$p^{\smallfrown}s\leq p^*_s,\text{ a set } A^{s}\subseteq B^{s},\text{ and } \succ_{T^*}(s)\subseteq \succ_T(s), \  \succ_{T^*}(s)\in U(\kappa_{j_n}(p),\xi(s))$$
such that $\{p_s^{*\smallfrown}\l\alpha,  A^s\cap\alpha\mid \alpha\in\succ{T^*}(s)\}$ is pre-dense above $p^*_s$. Apply \ref{amalgamate}, and find a condition $p\leq^* p_1$, $T_1\subseteq T\setminus mb(T)$, $mb(T_1)\in U_{T\setminus mb(T)}$ and sets
$C^s\in\cap_{\xi<\xi(s)}U(\kappa(s),\xi)$ such that for every $t\in mb(T_1)$, $p^*_t\leq^*p_1^{\smallfrown}\l t,\vec{C}^t\r$. Now apply induction hypothesis to $p_1$, $T_1$ and the sets $B^{s}\cap C^{s}$, find $p_1\leq^* p^*$ and $T^*\restriction\{1,...,n-1\}\subseteq T_1$  and sets $A^s$ such that $\{p^{*\smallfrown}\l s,\vec{A}^s\r\mid s\in mb(T^*\restriction\{1,...,n-1\})\}$ is pre-dense. Let us prove that above $p^*$,
$\{p^{*\smallfrown} \l t,\vec{A}^t\r\mid t\in mb(T^*)\}$ is pre-dense above $p^*$. Let $p^*\leq q$, then there is $s\in mb(T^*\restriction\{1,...,n-1\})$ such that $p^{*\smallfrown}\l s,\vec{A}^s\r$ and $q$ are compatible via some $q'$. Since $A^s\subseteq C^s$, it follows that
$$p^{*}_s\leq^* p_1^{\smallfrown}\l s,\vec{C}^s\r\leq^* p^{*\smallfrown}\l s,\vec{A}^s\r\leq q'.$$ Therefore, there is $\alpha\in \succ_{T^*}(s)$ such that $p_s^{*\smallfrown}\l \alpha,A^s\cap\alpha\r,q'$ are compatible via $q''$. It follows that $p_s^{*\smallfrown}\l \alpha,A^s\cap\alpha\r\leq q''$ and also $p^{*\smallfrown} \l s,\vec{A}^s\r\leq q'\leq q''$. So $\l \alpha,A^s\cap\alpha\r$ can be added to
$p^{*\smallfrown}\l s,\vec{A}^s\r$ and $p^{*\smallfrown}\l s^{\smallfrown}\alpha,\vec{A}^{s^{\smallfrown}\alpha}\r=p^{*\smallfrown}\l s,\vec{A}^s\r^{\smallfrown}\l\alpha, A^s\cap\alpha\r\leq q''$. 
We conclude that $q''$ is a witness for the compatibility of $q$ and $p^{*\smallfrown}\l s^{\smallfrown}\alpha,\vec{A}^{s^{\smallfrown}\alpha}\r$.$\blacksquare$

We will often have two conditions $p\leq^* p^*$ and a tree of extensions $T$ of $p$ as in \ref{densetree}, so there are sets $B^t$ such that $D_{T,\vec{B}}$ is pre-dense above $p$. We would like to remove some of the branches in $T$ to get a tree of extensions of $p^*$, $T^*\subseteq T$, such that $D_{T^*,\vec{B}}$ is pre-dense above $p^*$. $T^*$ can simply be defined as:
$$T^*=\{t\in mb(T)\mid p^{*\smallfrown}\l t,\vec{B}^t\r\in mb(T)\}$$
It is not hard to check that $T^*$ is a $\vec{U}$-fat tree and $mb(T^*)\in U_T$. To see that $D_{T^*,\vec{B}}$ is pre-dense above $p^*$, let $p^*\leq q$ then there is $t\in mb(T)$ such that $p^{\smallfrown}\l t,\vec{B}^t\r,q$ are compatible via a condition $q''$. Since $t$ appears in $q''$ and $p^*\leq q\leq q''$, it follows by \ref{frown extension} that $t\in mb(T^*)$ and $p^{*\smallfrown}\l t,\vec{B}^t\r,q\leq q''$.
\begin{corollary}\label{Rad-prop2}
 Let $p\in\Mfor$ and $\langle\lambda,B\rangle$ in the stem of $p$. Consider the decomposition, $p=\langle q,r\rangle$, where $q\in\Mfor\restriction\lambda\wedge r\in \Mfor\restriction(\lambda,\kappa]$, and $\kappa$ is the maximal measurable in $\vec{U}$. Let $\lusim{x}$ be a $\Mfor$-name for an ordinal. Then there is $r\leq^*r^*\in \Mfor\restriction(\lambda,\kappa]$ such that for any $q\leq q'\in\Mfor\restriction \lambda$ if there exist $r^*\leq r'\in \Mfor\restriction(\lambda,\kappa]$ such that
$$\langle q',r'\rangle \ || \lusim{x}$$
then there is a tree of extensions of $r^*$, $T/{q'}$, and sets $B^{t,q'}$ such that $D_{T_{q'},\vec{B}^{t,q'}}$ is pre-dense above $r^*$ and
 $$\forall t\in mb(T/{q'}). \  \l q',r^{*\frown}\l t,\vec{B}^t\r\r \ || \lusim{x}.$$
\end{corollary}
\pr 
For every $q\in \Mfor\restriction\lambda$, let
$$D_q=\Big\{p'\in\Mfor\restriction(\lambda,\kappa]\ \big| \  (\langle q,p'\rangle ||\lusim{x})\vee (\forall p''\geq p'. \l q,p''\r \text{ does not decide }\lusim{x})\Big\}.$$
Clearly, $D_q\subseteq \Mfor\restriction(\lambda,\kappa]$ is dense open, hence by the strong Prikry property, there is $r\leq^* r_q$, a tree of extensions $T'_q$ and sets $A^{s,q}$ for $s\in T'_q\setminus mb(T'_q)$ such that for every $t\in mb(T'_q)$, $r_q^{\smallfrown}\l t,\vec{A}^{t,q}\r\in D_q$.
For each $t\in mb(T'_q)$ one of the following holds:
\begin{enumerate}
    \item $ \l q,r_q^{\smallfrown}\l t,\vec{A}^{t,q}\r\r ||\lusim{x}$.
    \item $\forall p''\geq r_q^{\smallfrown}\l t,\vec{A}^{t,q}\r. \l q,p''\r \text{ does not decide }\lusim{x}$.
\end{enumerate}
Denote by $i_{t}\in\{1,2\}$ the case which holds. This defines a function $g:mb(T'_q)\rightarrow \{1,2\}$. Apply \ref{StabTree}, shrink $T'_q$ to $T''_q$ and find $i^*\in\{1,2\}$ such that for every $t\in mb(T''_q)$ $i_{t}=i^*$. Finally, apply \ref{densetree}, extend $r_q$ to $r^{*}_q$ shrink $T''_q$ to $T^*_q$ and find sets $B^{s,q}\subseteq A^{s,q}$ so that $$D_{T^*_q,\vec{B}^q}=\{r^{*\smallfrown}_q\l s,\vec{B}^{s,q}\r\mid s\in mb(T^*_q)\}$$ is pre-dense above $r^*_q$. There is sufficient $\leq^*$-closure in $\Mfor\restriction(\lambda,\kappa]$ to find a single $r^*$ such that $r_q\leq^* r^*$ for every $q\in\Mfor\restriction\lambda$.
Let us prove that $r^*$ is as wanted. We can shrink the trees $T^*_q$ to $T_q$ as in the discussion before \ref{Rad-prop2}, to be extension trees of $r^*$ such that $D_{T^*_q,\vec{B}^q}$ is pre-dense.  To see that $r^*,T_q,B^{s,q}$ are as wanted, let $q'\geq q$ and assume that there is $r'\geq r^*$ such that $\l q',r'\r ||\lusim{x}$. Since the set $\{r^{*\smallfrown}\l t,\vec{B}^{t,q}\r\mid 
t\in mb(T_{q'})\}$ is pre-dense above $r^*$, there is $t\in mb(T_{q'})$ such that $r^{*\smallfrown}\l t,\vec{B}^{t,q}\r,r'$ are compatible. In particular, there is $r''\geq r^{*\smallfrown}\l t,\vec{B}^{t,q}\r$ such that $\l q',r''\r ||\lusim{x}$, indicating that $i^*=i_{t}=1$.
Hence for every $s\in mb(T^*_{q'})$, $i_{s}=1$, thus
$$\l q',r^{*\smallfrown}\l s,\vec{B}^{s,q}\r\r ||\lusim{x}.$$
$\blacksquare$

The next lemma is the first step toward Theorem \ref{MainResaultParttwo}. Recall the inductive hypothesis $(IH)$: for every $\mu<\kappa$ and every coherent sequence $\vec{W}$ with maximal element $\mu$, every $V$-generic $G_\mu\subseteq \mathbb{M}[\vec{W}]$ and a set of ordinals $X\in V[G]$ there is $C'\subseteq C_{G_\mu}$ such that $V[X]=V[C']$.
\begin{lemma}\label{very short}
  Let $G\subseteq\Mfor$ be $V$-generic filter and assume $(IH)$. Let $A\in V[G]$ be a set of ordinals such that $|A|<\kappa$, where $\kappa$ is the maximal measurable in $\vec{U}$. Then there exists $C'\subseteq C_G$, $|C'|\leq |A|$, such that $V[A]=V[C']$.
\end{lemma}
\pr Let $A=\langle a_i \mid i<\lambda\rangle$ where $\lambda=|A|<\kappa$ be an enumeration of $A$. In $V$, pick a sequence of $\Mfor$-names for $A$, $\langle \lusim{a}_i\mid i<\lambda\rangle$.
We proceed by a density argument, let $p\in\Mfor\restriction(\lambda,\kappa]$
be any condition, using Lemma \ref{Rad-prop2}, find a $\leq^*$-increasing sequence $\langle p_i \mid i<\lambda\rangle$ above $p$ and maximal antichains
$Z_i\subseteq\Mfor\restriction \lambda$ such that for every $q\in Z_i$ there is a $\vec{U}$-fat tree $T_{q,i}$ and sets $B^{s,q}_i$ such that any extension of $p_i$ from $mb(T_{q,i})$ together with $q$ and the sets $B^{s,q}_i$ decides
$\lusim{a}_i$, and the set $$D_{T_{q,i},\vec{B}^q_i}:=\{p_i^{\smallfrown}\l t,\vec{B}^{t,q}_i\r\mid t\in mb(T_{q,i})\}.$$ is pre-dense above $p_i$. The forcing $\Mfor\restriction(\lambda,\kappa]$ has sufficiently $\leq^*$-closure to find $p'$ such that for every $i<\lambda, \ p_i\leq^* p'$.
Define the function $F_{q,i}:mb(T_{q,i})\rightarrow On$ by:
$$F_{q,i}(t)=\gamma \ \ \ \Leftrightarrow \ \ \  \l q,p^{'\smallfrown}\l t, \vec{B}^{t,q}_i\r\r\Vdash \lusim{a}_i=\check{\gamma}.$$
By Lemma \ref{important coordinates}, we can find $T'_{q,i}\subseteq T_{q,i}$, $mb(T'_{q,i})\in U_{T_{q,i}}$ such that $I_{q,i}:=I(T'_{q,i},F_{q_i}\restriction mb(T'_{q,i}))$ is complete and consistent. For any $q,q'\in Z_i$ apply Lemma \ref{function separation} to the functions $F_{q,i},F_{q',i}$ and shrink $T'_{q,i},T'_{q',i}$ to $T^{q,q'}_{q,i},T^{q,q'}_{q',i}$, $mb(T^{q,q'}_{q,i})\in U_{T_{q,i}},mb(T^{q,q'}_{q',i})\in U_{T_{q',i}}$ so that either 
\begin{enumerate}
    \item $mb(T^{q,q'}_{q,i})\restriction I_{q,i}=mb(T^{q,q'}_{q',i})\restriction I_{q',i}$ and $(F_{q,i}\restriction mb(T^{q,q'}_{q,i}))_{I_{q,i}}=(F_{q',i}\restriction mb(T^{q,q'}_{q',i}))_{I_{q',i}}$.
    \item $Im(F_{q,i}\restriction mb(T^{q,q'}_{q,i}))\cap Im(F_{q',i}\restriction mb(T^{q,q'}_{q',i}))=\emptyset$.
\end{enumerate}   The ultrafilter $U_{T_{q,i}}$ is sufficiently closed to ensure that $X^*_q=\cap_{q'\in \Mfor\restriction\lambda}mb(T^{q,q'}_{q,i})\in U_{T_{q,i}}$ and by \ref{Los for trees} there is a $\vec{U}$-fat tree $T'_{q,i}\subseteq T_{q,i}$ such that $mb(T'_{q,i})\subseteq X^*_q,$ and $mb(T'_{q,i})\in U_{T_{q,i}}$. By \ref{densetree}, there is $p'\leq^*p^*_q$, $T^*_{q,i}$ and $A^{s,q}_i\subseteq B^{s,q}_i$ such that $D_{T^*_{q,i},\vec{A}^{q}_i}$ is pre-dense above $p^*_q$. Since $|\Mfor\restriction\lambda|$ is small enough there is a single $p^*\in \Mfor\restriction(\lambda,\kappa]$ such that $p^*_{q}\leq^* p^*$ for every $q\in \Mfor\restriction\lambda$. Restrict the trees to this condition $p^*$ as in the discussion before \ref{Rad-prop2}, so that $D_{T^*_q,\vec{B}^q_i}$ are pre-dense above $p^*$. we abuse  notation here by keeping the same notation after the restriction.

Denote $G=G_<\times G_>$ 
 so that $G_<\subseteq \Mfor\restriction\lambda$ is $V$-generic and $G_>\subseteq\Mfor\restriction(\lambda,\kappa]$ is $V[G_<]$-generic. By density, find $p^*\in G_>$ as above. For every $i<\lambda$, since $Z_i$ is a maximal antichain, there is $q_i$ such that $G_<\cap Z_i=\{q_i\}$. Since $D_{T^*_{q_i,i},\vec{A}^{q_i}_i}$ is pre-dense above $p^*$, find $t_i\in mb(T^*_{q_i,q})$ such that $p^{*\frown}\l t_i,\vec{A}^{t_i,q_i}_i\r\in G_>$, define $C_i=t_i\restriction I_{q_i,i}$ and let $C'=\underset{i<\lambda}{\bigcup}C_i\subseteq C_{G_>}$. Clearly $|C'|\leq \lambda=|A|$. Let us prove that $\langle C_i\mid i<\lambda\rangle\in V[A]$. Indeed, define in $V[A]$ the sets $$M_i=\{q\in Z_i\mid a_i\in Im(F_{q,i})\}$$  then, for any $q,q'\in M_i$, $a_i\in Im(F_{q,i})\cap Im(F_{q',i})\neq\emptyset$. Hence $(1)$ must hold for $F_{q,i},F_{q',i}$ i.e.
 $$mb(T^*_{q,i})\restriction I_{q,i}=mb(T^*_{q',i})\restriction I_{q',i}\wedge(F_{q,i}\restriction mb(T^*_{q,i}))_{I_{q,i}}=(F_{q',i}\restriction mb(T^*_{q',i}))_{I_{q',i}}.$$
This means that no matter how we pick $q_i'\in M_i$, we will end up with the same function $(F_{q'_i,i}\restriction mb(T^*_{q'_i,i}))_{I_{q_i',i}}$ and the same important values $mb(T^*_{q'_i,i})\restriction I_{q_i,i}$. In $V[A]$, choose any $q'_i\in M_i$, let $D_i'\in F_{q_i',i}^{-1''}\{a_i\}\cap mb(T^*_{q_i',i})$ and $C_i'=D'_i\restriction I_{q_i',i}$. Since $q_i,q'_i\in M_i$ we have $ C_i=C_i'$, hence $\langle C_i\mid i<\lambda\rangle\in V[A]$. In order the reconstruct $A$ from the union $C'$ we still have to code some information from the part of $G_<$, namely, $\{q'_i\mid i<\lambda\},\langle Ind(C_i,C')\mid i<\lambda\rangle\in V[A]$. These sets can be coded as a subset of ordinals below $(2^{\lambda})^+$, by \ref{genericproperties}(6) $$\{q'_i\mid i<\lambda\},\langle Ind(C_i,C')\mid i<\lambda\rangle\in V[G_<].$$ By the induction hypothesis   applied to $G_<$, we can find $C''\subseteq C_{G_<}$ such that $$V[\l q'_i\mid i<\lambda\r,\langle Ind(C_i,C')\mid i<\lambda\rangle]=V[C''].$$ 
Also $|C''|\leq |C_{G_<}|\leq\lambda$ hence $C:=C'\uplus C''$ is of cardinality at most $\lambda$. Note that $C',C''\in V[C]$ as $C''=C\cap\lambda, \ C'=C\setminus\lambda$. Finally, all the information  about the function $F_{q,i}$ needed to restore $A$ is coded in $C',C''$. Namely, 
$A=\{(F_{q'_i,i})_{I_{q'
_i,i}}(C'\restriction Ind(C_i,C'))\mid i<\lambda\}$. Hence $V[A]=V[C]$.
$\blacksquare$
\begin{corollary}\label{specialTree}
Suppose that $p\in\Mfor$ and $\lusim{x}$ is a name such that $p\Vdash \lusim{x}\in \lusim{C}_G$. Then there is $p^*\geq^* p$ such that either $p^* || \lusim{x}$ or there is a $\vec{U}$-fat tree, $T$ and sets $A^{s}$ such that $\forall t\in mb(T)$ $p^{\frown}\l t,\vec{A}^{t}\r\Vdash \lusim{x}=\max(t)$. Moreover, in the latter case, let $i\leq l(p)+1$ be such that $mb(T)$ splits on $\kappa_i(p)$ and assume that $o^{\vec{U}}(\kappa_i(p))<\kappa_i(p)^+$, then for every $t\in \Lev_{ht(T)-1}(T)$,  $$p^{*\frown}\langle t,\vec{A}^{t}\rangle || o^{(\kappa_i(p))}(\lusim{x}).$$
In other words, there is $\gamma<o^{\vec{U}}(\kappa_i(p))$ such that
$$p^{*\frown}\langle t,\vec{A}^{t}\rangle\Vdash\lusim{x}\in X^{(\kappa_i(p))}_\gamma.$$ 
\end{corollary}
\pr Assume that there is no $p^*\geq^* p$ which decides $\lusim{x}$. By \ref{Rad-prop2}
find $T$ with minimal $ht(T)$ such that there is $p^*\geq p$, sets $B^{s}$ and for every $t\in mb(T)$, $p^{*\frown}\l t,B^{s}\r || \lusim{x}$. Assume that $\kappa(p)=\{\nu_1,...,\nu_n\}$ are the ordinals appearing in $p$, denote by $x_t$ the forced value and shrink $T$ so that the function $$f(t)=\begin{cases}
i & x_t=\nu_i\\
n+1 & x_t\notin\{\nu_1,...,\nu_n\}
\end{cases}$$ 
is constant. If $f$ would be constantly some $i\leq n$ then by Proposition \ref{densetree} there is $p\leq^* p'$, $T'\subseteq T$ and sets $A^s\subseteq B^s$  such that $\{p^{'\smallfrown}\l t,\vec{A}^t\r\mid t\in mb(T')\}$ is pre-dense above $p'$, it follows that $p'\Vdash \lusim{x}=\nu_i$, contradiction. So we may assume that $x_t\notin \{\nu_1,...,\nu_n\}$. Keep shrinking $T$ so that there is a unique $i\leq ht(T)$, such that $x_t\in [t(i),t(i+1))$ (where $t(ht(T)+1)=\kappa$). If $i<ht(T)$ then for every $t\in \Lev_i(T)$, the function $g_t:mb(T/t)\rightarrow \kappa$, defined by $g_t(s)=x_{t^{\smallfrown}s}$ is regressive and therefore by \ref{StabTree} can be stabilized on some $S_t\subseteq T/t$, $mb(S_t)\in U_{T/t}$ so that for every $t\in S_t$, $x_{t^{\smallfrown}s}=y_t$, depending only on $t$. As in the situation that $f$ was constant, for every $t\in \Lev_i(T)$ we can find $p^{*\smallfrown}t\leq^*p_t$ such that $p_t\Vdash \lusim{x}=x_t$. By \ref{amalgamate}, there is $T^*\subseteq T\restriction\{1,...,i\}$, $p^{*}\leq^*p^{**}$ and sets $Z^s\subseteq A^s$ such that for every $t\in \Lev_i(t)$, $p
_t\leq^*p^{**\smallfrown}\l t,\vec{Z}^t\r$, this contradicts the minimality of $ht(T)$. Hence it must be that for every $t\in mb(T)$, $x_t\geq t(ht(T))=\max(t)$. It is impossible that $x_t>\max(t)$, otherwise, $$x_t\notin \{\nu_1,...,\nu_n\}\cup t$$ and we can remove from the large sets of the condition $p^{*\smallfrown}\l t,\vec{A}^t\r$ the single ordinal $x_t$ and obtain a condition $q$ such that $q\Vdash \lusim{x}=x_t\notin C_{\lusim{G}}$, but $p\leq q$, then $q\Vdash \lusim{x}\in C_{\lusim{G}}$, contradiction. We conclude that $\forall t\in mb(T).x_t=\max(t)$. Which is what we desired.

For the second part, assume that for $i\leq l(p)+1$,  $mb(T)$ splits on $\kappa_i(p)$ and that $o^{\vec{U}}(\kappa_i(p))<\kappa_i(p)^+$. It follows that the measures in $\vec{U}(\kappa_i(p))$ are separated by the sets $X^{(\kappa_i(p))}_{\gamma}$. For every $t\in \Lev_{ht(T)-1}(T)$, shrink $\succ_T(t)\in U(\kappa_i(p),\xi(t))$ to $\succ_{T^*}(t)=\succ_{T}(t)\cap X^{(\kappa_i(p))}_{\xi(t)}$. It follows that for every $t\in \Lev_{ht(T)-1}(T)$, and for every $\beta\in \succ_{T^*}(t)$, $\beta\in X^{(\kappa_i(p))}_{\xi(t)}$. Since $p^{\frown}\l t,\vec{A}^t\r \Vdash \lusim{x}\in \succ_{T^*}(t)$, we conclude that
$p^{\frown}\l t,\vec{A}^t\r \Vdash  o^{(\kappa_i(p))}(\lusim{x})=\xi(t)$
$\blacksquare$

The following lemma is analogous to a lemma proven in \cite{TomTreePrikry} for Prikry forcing.
\begin{lemma}\label{MagidorHousdorf}
Let $G\subseteq\Mfor$ be $V$-generic and let $\delta\leq \kappa$ be a limit point of $C_G$. Then for every set of ordinals $D\in V[C_G]$  such that $$|D|<\delta\wedge C_G\cap D=\emptyset$$ there is $X\in\bigcap\vec{U}(\delta)$ such that $X\cap D=\emptyset$.
\end{lemma}
\pr
Let $\lambda:=|D|$, note that $D\in V[C_G\cap\delta]$ and since $C_G\cap\delta$ is $V$-generic for $\Mfor\restriction\delta$, we can assume without loss of generality that $\delta=\kappa$. We start with a single $\Mfor$-name of an ordinal $\lusim{x}$  and $p\in G$ such that $p\Vdash \lusim{x}\notin \lusim{C}_G$. Assume that $p=\langle q_0,r\r$, is a decomposition of $p$ such that $\max(\kappa(q_0))\geq\lambda$. Then by \ref{Rad-prop2} there is $r\leq^*r^*$ and  a maximal antichain $Z\subseteq\mathbb{M}[\vec{U}]\restriction\max(q_0)$ above $q_0$, such that for every $q\in Z$ there is a tree $T_q$ and sets $A^{s,q}$ for which the set $\{r^{*\smallfrown}\l t,\vec{A}^{t,q}\r\mid t\in mb(T_{q})\}$ is pre-dense above $r^*$ and for every $t\in mb(T_q)$, $$\l q,r^{*\smallfrown}\l t,\vec{A}^{s,q}\r\r \Vdash \lusim{x}=f_q(t).$$   
Since $p\Vdash \lusim{x}\notin C_{\lusim{G}}$, for every $\vec{b}\in mb(T_q)$, $f_q(\vec{b})\notin\vec{b}$ hence it falls in one of the intervals 
$$(0,\vec{b}(1)),(\vec{b}(1),\vec{b}(2)),...,(\vec{b}(ht(T_q)),\kappa)$$
let $n_{\vec{b}}$ be the index of this interval. Apply \ref{StabTree} to find a tree $T'_q\subseteq T_q$, $mb(T'_q)\in U_{T_q}$ on which the value $n_{\vec{b}}$ is constantly $n^*_q$. Since for every $t\in \Lev_{n^*_q}(T_q)$, the function $s\mapsto f_q(t^{\smallfrown} s)$ defined in $mb((T'_q)/t)$, is regressive, apply \ref{StabTree}, obtain a tree $(T^*_q)_t\subseteq (T'_q)/t$ on which the value is constant. Let  $T^*_q\restriction \{1,...,n^*_q\}=T'_q\restriction \{1,...,n^*_q\}$ and for every $t\in \Lev_{n^*_q}(T^*_q)$, $(T^*_q)/t=(T^*_q)_t$ is defined as above.  Then on $T^*_q$, $f_q(t)$ depends only on $t\restriction\{1,...,n^*_q\}$ and $f_q(t\restriction\{1,...,n^*\})>t(n^*_q)$.
Extend $r^*_0\leq^*r^*_q$, shrink $S_q\subseteq T^*_q$ to a tree of extensions of $r^*_q$ and find $B^{q,s}\subseteq A^{q,s}$ such that for every $s\in \Lev_{n^*}(S^*_q)$, $D_{S_q,\vec{B}^{q}}$ is pre-dense above $r^*_q$ and $\l q,r_1^{*\smallfrown}\l s,\vec{B}^{q,s}\r\r\Vdash \max(s)<\lusim{x}=f_q(s)$. Finally find a single $r^*$ such that $r^*_q\leq^*r^*$, shrink the trees and sets to this condition and denote $S_q\restriction \{1,...,n^*_q\}=S^*_q$.
Apply \ref{BetterSet} and let $A_{\lusim{x}}=\{\alpha\in B_{l(r)+1}(r^*)\mid \alpha\cap B_{l(r)+1}(r^*)\in\cap\vec{U}(\alpha)\}\in \cap\vec{U}(\kappa)$.
It must be that for every $q\in Z$ and every $s\in mb(S^*_q)$, $f_q(s)\notin A_{\lusim{x}}\setminus \max(s)$,
otherwise, add the ordinal $f_q(s)$ and obtain the condition $$\langle q, r^{*\smallfrown}\langle s,\vec{B}^{q,s}_<\r ^{\smallfrown}\l f_q(s)\r\r\Vdash\lusim{x}=f_q(s)\in C_G$$
contradiction. Since $f_q(s)>\max(s)$, we conclude that $f_q(s)\notin A_{\lusim{x}}$. 
we claim that $$p\leq^*\langle q_0,r^*\rangle\Vdash \lusim{x}\notin A_{\lusim{x}}.$$
Otherwise, there is $q\in Z$, $s\in mb(S^*_q)$ and $p'$ such that $$\langle q,r^{*\frown}\l s, B^{s,q}\r\leq p'\Vdash \lusim{x}\in A_{\lusim{x}}$$
but also $p'\Vdash\lusim{x}=f_q(s)$ so $f_q(s)\in A_{\lusim{x}}$ which is a contradiction. 
Now the lemma follows easily, let $\{d_i\mid i<\lambda<
\kappa\}\in V[C_G]$ be some set of ordinals such that $$C_G\cap\{d_i\mid i<\lambda\}=\emptyset$$ then we can take names $\{\lusim{d}_i\mid i<\lambda\}$ and some $p=\l q_0,r_0\r$ forcing $\forall i<\lambda. \lusim{d}_i\notin \lusim{C}_G$, as before we can define the sets $A_{\lusim{d}_i}\in\cap\vec{U}(\kappa)$ and for $i<\lambda$ find a $\leq^*$-increasing  sequence $\langle q_0, r_i\rangle$, find $p^*$ which bounds all of them and $A^*=\underset{i<\lambda}{\bigcap}A_{\lusim{d}_i}\in \cap\vec{U}(\kappa)$, then $p^*$ forces that $\forall i<\lambda$ $\lusim{d}_i\notin A^*$. By density argument we can find such $p^*$ in $G$.$\blacksquare$
\section{The proof for subsets of $\kappa$}
 Let $A\in V[G]$, we do not assume that $A\subseteq\kappa$, since some of the results will be applied for other type of sets.  Define $$\kappa^*:=\max\{\alpha\in Lim(C_G)\mid o^{\vec{U}}(\alpha)\geq\alpha^+\}.$$ 
If $o^{\vec{U}}(\kappa)<\kappa^+$, then $\{\beta<\kappa\mid o^{\vec{U}}(\beta)<\beta^+\}\in\cap\vec{U}(\kappa)$, it follows that $\kappa^*<\kappa$ is well defined. Moreover, for every $\alpha\in C_G\setminus\kappa^*+1$, $o^{\vec{U}}(\alpha)<\alpha^+$ and thus $o^{(\alpha)}$ is defined.

\begin{definition}
 Let $A\in V[G]$ be any set of ordinals. In $V[A]$, consider the crucial set
$$X_A=\{\nu\mid \nu \text{ is }V-\text{regular  and } \nu>  cf^{V[A]}(\nu)\}$$
Denote $\overline{X}_A=X_A\cup Lim(X_A)\subseteq \kappa\cup\{\kappa\}$.
\end{definition}
\begin{proposition}\label{propertiesofXA}
\begin{enumerate}
    \item $\overline{X}_A\subseteq Lim(C_G)$.
    \item $X_A\in V[A]$.
    \item If $o^{\vec{U}}(\kappa)<\kappa^+$, $X_A\setminus\kappa^*+1$ is closed i.e. for every $\kappa^*<\alpha\leq\kappa$, if $\sup(X_A\cap\alpha)=\alpha$ then $\alpha\in X_A$.
    \item If $C\subseteq^* C_G$ and $C\in V[A]$, then $Lim(C)\subseteq \overline{X}_A$
\end{enumerate}
\end{proposition} 
\pr For every $\alpha\in X_A$, $cf^{V[G]}(\alpha)\leq cf^{V[A]}(\alpha)<\alpha$, and $\alpha$ is $V$-regular, it follows by \ref{genericproperties}(7) that $X_A\subseteq Lim(C_G)$, and since $Lim(C_G)$ is closed, then $\overline{X}_A\subseteq Lim(C_G)$.

$(2)$ is trivial as the definition of $X_A$ occurs in $V[A]$. As for $(3)$,  by induction on $\alpha\in\Lim(X_A\setminus\kappa^*)$. Suppose $\alpha=\sup(X_A\cap\alpha)$, then by induction, $X_A\cap(\kappa^*,\alpha)$ is a club at $\alpha$ and  by $(1)$, $\alpha\in Lim(C_G)\setminus\kappa^*$. Define in $V[A]$, $$o_A(\alpha)=limsup_{\gamma\in X_A\cap \alpha}o^{(\alpha)}(\gamma)+1.$$ By definition of $o^{(\alpha)}$, $o_A(\alpha)\leq o^{\vec{U}}(\alpha)<\alpha^+$, hence $cf^{V}(o_A(\alpha))\leq\alpha$. By the definition of $limsup$, $o_A(\alpha)$ satisfies two properties:
\begin{enumerate}
    \item For every $\nu<\alpha$ and every $j<o_A(\alpha)$ there is $j\leq j'<o_A(\alpha)$ such that $X_A\cap X^{(\alpha)}_{j'}\cap(\nu,\alpha)\neq\emptyset$.
    \item There is some $\xi_\alpha<\alpha$ such that for every $\nu\in X_A\cap(\xi_\alpha,\alpha)$,  $o^{(\alpha)}(\nu)<o_A(\alpha)$.
\end{enumerate}
We split into cases:

If $o_A(\alpha)=\beta+1$, then by property $(1)$ $\sup (X_A\cap X^{(\alpha)}_\beta\cap(\xi_\alpha,\alpha))=\alpha$. Let us argue that $\otp( X_A\cap X^{(\alpha)}_\beta\cap(\xi_\alpha,\alpha))=\omega$, this is enough to conclude $ cf^{V[A]}(\alpha)=\omega$, hence $\alpha\in X_A$. In the interval $(\xi_\alpha,\alpha)$ it is impossible to have a limit point $\zeta$ of $X^{(\alpha)}_\beta\cap X_A$. Otherwise, by induction $\zeta\in X_A$ and by \ref{increasing order at limits}, $o^{(\alpha)}(\zeta)\geq\beta+1$ contradicting property $(2)$.

If $\lambda:=cf^V(o_A(\alpha))<\alpha$, let $\l \lambda_i\mid i<\lambda\r\in V$ be increasing and cofinal in $o_A(\alpha)$.
Define inductively $\l x_i\mid i<\lambda\r$, first, $x_0=\min (X_A\cap(\xi_\alpha,\alpha))<\alpha$. At successor step, $i+1$, $x_i\in X_A\cap(\xi_\alpha,\alpha)$ is defined, by property $(1)$, there is $$\lambda_{i+1}\leq j'<\alpha\text{ and }x_{i+1}\in X_A\cap(x_i,\alpha)\cap X^{(\alpha)}_{j'}.$$ At limit step $\delta<\lambda$, if $\sup(x_i\mid i<\delta)$ is unbounded in $\alpha$, then clearly $\alpha$ changes cofinality in $V[A]$. Otherwise, let $y_\delta=\sup_{i<\delta}x_i<\alpha$ and there is some $x_{\delta}\in X_A\cap X^{(\alpha)}_{\lambda_\delta}\cap(y_{\delta},\alpha)$. Assume that $\l x_i\mid i<\lambda\r$ is defined, if $x^*=\sup_{i<\lambda}x_i\in(\xi,\alpha)<\alpha$, then by induction hypothesis $x^*\in X_A\cap(\xi_\alpha,\alpha)$ and $$o^{(\alpha)}(x^*)\geq limsup_{i<\lambda} o^{(\alpha)}(x_i)+1=limsup_{i<\lambda}\lambda_i+1=o_A(\alpha)$$ contradicting  property $(2)$.

Finally, if $cf^{V}(o_A(\alpha))=\alpha$,
we take $\l \alpha_i\mid i<\alpha\r\in V$ cofinal continuous sequence in $o_A(\alpha)$ which witnesses this. 
Let $Z:=\{\beta<\alpha\mid o^{(\alpha)}(\beta)<\alpha_{\beta}\}$ let us argue that $Z\in\cap_{i<o_A(\alpha)}U(\alpha,i)$. Let $i<o_A(\alpha)$, denote 
$$j_{U(\alpha,i)}(\l \alpha_\xi\mid \xi<\alpha\r)=\l \alpha'_\xi\mid \xi<j_{U(\alpha,i)}(\alpha)\r, \ \ j_{U(\alpha,i)}(\l X^{(\alpha)}_\xi\mid \xi<o^{\vec{U}}(\alpha)\r)=\l X'_\xi\mid \xi< j_{U(\alpha,i)}(o^{\vec{U}}(\alpha))\r.$$
Since $X^{(\alpha)}_i\in U(\alpha,i)$ it follows that $\alpha\in j_{U(\alpha,i)}(X^{(\alpha)}_i)=X'_{j_{U(\alpha,i)}(i)}$ which by definition implies that $$(\star) \ \ \ o^{(j_{U(\alpha,i)}(\alpha))}(\alpha)=j_{U(\alpha,i)}(i).$$ Also, since $i< o_A(\alpha)$, then $$j_{U(\alpha,i)}(i)<\cup j_{U(\alpha,i)}''[o_A(\alpha)]=\cup_{\xi<\alpha}j_{U(\alpha,i)}(\alpha_\xi)=\cup_{\xi<\alpha}\alpha'_\xi.$$
By elementarity, the sequence $\l\alpha'_\xi\mid \xi<j_{U(\alpha,i)}(\alpha)\r$ is also continuous, hence $$(\star\star) \ \ j_{U(\alpha,i)}(i)<\cup_{z<\alpha}\alpha'_z=\alpha'_\alpha.$$
We conclude from $(\star),(\star\star)$ that
$$o^{j_{U(\alpha,i)}(\alpha)}(\alpha)=j_{U(\alpha,i)}(i)<\alpha'_\alpha.$$
Hence $\alpha\in j_{U(\alpha,i)}(Z)$ so $Z\in U(\alpha,i)$ as wanted. 

Consider the set $Z_*:=Z\uplus(\cup_{o_A(\alpha)\leq j<o^{\vec{U}}(\alpha)}X^{(\alpha)}_j)$. Then $Z_*\in\cap\vec{U}(\alpha)$ and by \ref{genericproperties}(3), there is $\eta<\alpha$ such that $C_G\cap (\eta,\alpha)\subseteq Z_*$. In particular $X_A\cap(\eta,\alpha)\subseteq Z_*$. By property $(2)$, if $\rho\in X_A\cap(\max\{\eta,\xi_\alpha\},\alpha)$, then $o^{(\alpha)}(\rho)<o_A(\alpha)$, hence $\rho\in Z$ hence $X_A\cap (\max\{\eta,\xi_\alpha\},\alpha)\subseteq Z$. By definition of $Z$, for every $\rho\in (\max\{\eta,\xi_\alpha\},\alpha)\cap X_A$, $o^{(\alpha)}(\rho)<\alpha_\rho$. Now to see that $cf^{V[A]}(\alpha)=\omega$, define $x_0=\min(X_A\cap(\max\{\eta,\xi_\alpha\},\alpha))$, recursively assume that $x_n<\alpha$ is defined. Then by property $(1)$, there is $x_{n}'\geq x_n$ and some $x_{n+1}\in X_A\cap X^{(\alpha)}_{\alpha_{x_n'}}\cap(x_n,\alpha)$. To see that $\l x_n\mid n<\omega\r$ is unbounded in $\alpha$, assume otherwise, then  $x^*=\sup_{n<\omega}x_n<\alpha$ and by induction $x^*\in X_A\cap(\max\{\eta,\xi_\alpha\},\alpha)$ hence $x^*\in Z$.  By Proposition  \ref{increasing order at limits} $$o^{(\alpha)}(x^*)\geq limsup_{n<\omega} o^{(\alpha)}(x_n)=limsup_{n<\omega} \alpha_{x'_n}\geq \alpha_{x^*}$$ contradiction the definition of $Z$.

To see $(4)$, if $C\setminus C_G$ is finite then clearly $Lim(C)\subseteq Lim(C_G)$ and every $\delta\in Lim(C_G)$ is $V$-regular.  Let $\delta\in Lim(C)$, it suffices to prove that $X_A$ is unbounded in $\delta$.  Fix any $\rho<\delta$, and let $\rho'=\min(Lim(C\setminus\rho+1))$, then $\rho<\rho'\leq\delta$, and also by minimality $\otp(C\cap(\rho,\rho'))=\omega$.  Since $C\in V[A]$, it follows that $cf^{V[A]}(\rho')=\omega$  and since $\rho'\in Lim(C)$, it is $V$-regular. By definition it follows that $\rho'\in X_A\cap(\rho,\delta]$.
$\blacksquare$

It is possible that $X_A$ below $\kappa^*$ is not closed:
\begin{example}
If there is $\alpha\in C_G$ such that $o^{\vec{U}}(\alpha)=\alpha^+$, then $\alpha$ stays regular in $V[G]$. Set $A=C_G$, then $X_A\cap\alpha$ will be unbounded in $\alpha$, but $\alpha\notin X_A$.
\end{example}

There are trivial examples for $A$ in which the set $X_A$ is bounded. However the following definition filters this situation.
\begin{definition}
  Let $A\subseteq On$, we say that $A$ \textit{stabilizes} if there is $\beta<\kappa$ such that $\forall\alpha<\sup(A), \  A\cap\alpha\in V[G\restriction\beta]$
\end{definition}
This definition is more general than the notion of fresh set:
\begin{definition}\label{freshSetDef}
Let $M\subseteq M'$ be two $ZFC$ models. A set of ordinals $X\in M'\setminus M$ is \textit{Fresh with respect to $M$} if $\forall\alpha<\sup(X). X\cap\alpha\in M$.
\end{definition}
\begin{proposition}\label{nonstabcofinality}
Suppose that $A\in V[G]$ such that $A$ does not stabilize. Assume that $\forall\beta<\sup(A)$ there is $C_\beta\subseteq C_G$ such that $V[C_\beta]=V[A\cap\beta]$. Then:
 \begin{enumerate}
    \item  If $A\subseteq\kappa$, then $X_A\cap\kappa$ is unbounded in $\kappa$. 
     \item If $o^{\vec{U}}(\kappa)<\kappa^+$, then $cf^{V[A]}(\kappa)<\kappa$.
 \end{enumerate}
\end{proposition}
\pr The following argument works for both $(1),(2)$, we try to prove that $X_A$ is unbounded. Let $\kappa^*\leq\delta<\kappa$, take some $\beta<\sup(A)$ such that $A\cap\beta\notin V[G\restriction\delta]$ which exists by our assumption that $A$ does not stabilize. By assumption, there exists $C_{\beta}\subseteq C_G$ such that $$V[C_\beta]=V[A\cap\beta]\subseteq V[A].$$ It is impossible that $C_\beta\setminus (C_G\cap\delta)$ is finite, otherwise $$A\cap\beta\in V[C_\beta]\subseteq V[G\restriction\delta]$$ which contradicts the choice of $\beta$. Let $\gamma_\delta$ be the first limit point of $C_\beta$ above $\delta$. By minimality, $\otp(C_\beta\cap(\delta,\gamma_\delta))=\omega$, hence $cf^{V[A]}(\gamma_\delta)=\omega$ and $\gamma_\delta\in X_A\setminus\delta$. 

To see $(1)$, if $A\subseteq\kappa$, then necessarily $\gamma_\delta<\kappa$ for every $\delta$, this is since $\gamma_\delta\in Lim(C_\beta)$, and $\beta<\sup(A)\leq\kappa$, so $V[C_\beta]=V[A\cap\beta]\subseteq V[C_G\cap\beta]$. This implies that $\gamma_\delta\leq\beta$, otherwise, in $V[C_G\cap\beta]$ the cofinality of some measurable above $\beta$ changes, which contradicts $\beta^+$-cc of $\Mfor\restriction\beta$.
To see $(2)$, if some $\gamma_\delta=\kappa$, then $\kappa\in X_A$ and $cf^{V[A]}(\kappa)<\kappa$. Otherwise, $\gamma_\delta<\kappa$, and we conclude that $X_A$ is unbounded in $\kappa$. By the assumption $o^{\vec{U}}(\kappa)<\kappa^+$, thus by \ref{propertiesofXA}(3), $X_A\setminus\kappa^*$ is closed, and $\kappa$ is a limit point of this set, so $\kappa\in X_A$.
$\blacksquare$
\begin{corollary}\label{nonstabcofinalitysubsetofkappa}
Assume $(IH)$ and suppose that $A\in V[G]$ such that $A\subseteq\kappa$ does not stabilize and $o^{\vec{U}}(\kappa)<\kappa^+$.  Then $X_A\cap(\kappa^*,\kappa)$ is a club at $\kappa$ and $cf^{V[A]}(\kappa)<\kappa$.
\end{corollary}
\pr Since $A\subseteq\kappa$, then by \ref{very short}, for every $\beta<\kappa$, there is $C_\beta$ such that $V[A\cap\beta]=V[C_\beta]$ so we can apply \ref{nonstabcofinality}(1), \ref{nonstabcofinality}(2) and \ref{propertiesofXA}(3) applies to conclude that $X_A\cap(\kappa^*,\kappa)$ is a club and $cf^{V[A]}(\kappa)<\kappa$.$\blacksquare$

Note that it is possible that $cf^{V[G]}(\kappa)< cf^{V[A]}(\kappa)<\kappa$, however $cf^{V[A]}(\kappa)$ must be some member of the generic club that will eventually change its cofinality to $cf^{V[G]}(\kappa)$. 

\begin{example}
Assume that $o^{\vec{U}}(\kappa)=\kappa$, then $cf^{V[G]}(\kappa)=\omega$.
Using the enumeration $C_G=\langle C_G(i)\mid i<\kappa\rangle$ and the canonical sequence $\alpha_n$ that was defined in example \ref{canonicalsequence}, we can define in $V[G]$ the set
$$A=\bigcup_{n<\omega} \{C_G(\alpha_n)+\alpha\mid \alpha<C_G(n)\}$$
then $A$ does not stabilize. Moreover, we cannot construct the sequence $\langle \alpha_n\mid n<\omega\rangle$ or any other $\omega$-sequence unbounded in $\kappa$ inside $V[A]$ since  $A$ is generic for the forcing $\mathbb{M}[\vec{U}\restriction(C_G(\omega),\kappa]]$ which does not change the cofinality of $\kappa$ to $\omega$. For this kind of examples the case $o^{\vec{U}}(\kappa)<\kappa$ suffices.

\end{example}

The following definition will allow us to refer to subsets of $C_G$ in $V[A]$.
\begin{definition}
 Let $A\in V[G]$ be any set. A set $D\in V[A]$ is a \textit{Mathias set} if 
 \begin{enumerate}
     \item $Lim(D)\subseteq \overline{X}_A$.
     \item For every $\delta\in Lim(D)$,  every $Y\in \bigcap \vec{U}(\delta)$ there is $\xi<\delta$ such that $D\cap(\xi,\delta)\subseteq Y$.
 \end{enumerate}
\end{definition}

\begin{lemma}\label{FiniteNoise}
For every $D\in V[A]$, $D$ is Mathias if and only if $D\subseteq^* C_G$ i.e. $D\setminus C_G$ is finite.
\end{lemma}
\pr 
If $D\setminus C_G$ is finite then by \ref{propertiesofXA}(4), $Lim(D)\subseteq \overline{X}_A$. For the second condition of a Mathias set, simply use \ref{genericproperties}(3). 

In the other direction, assume that $D$ is a Mathias set. Toward a contradiction, assume $|D\setminus C_G|\geq\omega$, and let $\delta\leq sup(D)$ be minimal such that $|D\cap\delta\setminus C_G|\geq\omega$ then $\delta\in Lim(D)\subseteq \overline{X}_A\subseteq Lim(C_G)$. By minimality, $\{d_n\mid n<\omega\}=D\cap\delta\setminus C_G$ is unbounded in $\delta$. By \ref{MagidorHousdorf} there is $Y\in \bigcap \vec{U}(\delta)$ such that $Y\cap\{d_n\mid n<\omega\}=\emptyset$ contradicting condition $(2)$ of the Mathias set $D$.$\blacksquare$

\begin{proposition}\label{CodingBoundedInformation}
Let $A\in V[G]$ and $\lambda<\kappa$, let $\lambda_0:=\max(Lim(C_G)\cap\lambda+1)$ and assume $(IH)$.  Then there is a Mathias set $F_\lambda\subseteq \lambda_0$ such that $V[F_\lambda]=V[A]\cap V[C_G\cap\lambda]$.
\end{proposition}
\pr Consider in $V[A]$ the sets $$B:=\{D\subseteq\lambda\mid D\text{ is a Mathias set}\}.$$
Then $|B|\leq 2^\lambda$, enumerate $B=\l D_i\mid i<2^\lambda\r$,  let $E=\{\l i, d\r\mid i<2^{\lambda}, d\in D_i\}\subseteq 2^{\lambda}\times\lambda$, clearly $V[B]=V[E]$  and  $E\subseteq V_{2^{\lambda}}$. Also, since elements of $Lim(C_G)$ are strong limits in $V[C_G]$, $$\max (Lim(C_G)\cap 2^\lambda+1)=\max(Lim(C_G)\cap\lambda+1)=\lambda_0.$$ By Proposition \ref{genericproperties}(6), $E\in V[C_G\cap\lambda_0]$ and by induction hypothesis there is $F_\lambda\subseteq C_G\cap\lambda_0$ such that $V[F_\lambda]=V[E]$. Since $E\in V[A]$, also $F_\lambda\in V[A]$, and since $F_\lambda\subseteq C_G\cap\lambda_0$, $F_\lambda\in V[C_G\cap\lambda_0]$  so $V[F_\lambda]\subseteq V[A]\cap V[C_G\cap\lambda_0]$. For the other direction, if $X\in V[A]\cap V[C_G\cap\lambda_0]$, then by induction there is $C\subseteq C_G\cap\lambda_0$ such that $V[X]=V[C]$, and also $C\in V[A]$. Then $C\subseteq\lambda$ is a Mathias set, hence $C\in B$, and therefore, $C\in V[B]=V[F_\lambda]$. $\blacksquare$

%\subsection{Stabilization of Subsets of $\kappa$}

%Let us start by proving a lemma which will help us code the information we need into one sequence.
The following lemma will be crucial to pack information given by two sets $D,C\subseteq C_G$ into a single set $E\subseteq C_G$.
\begin{proposition}\label{countableinformation}
Assume that $o^{\vec{U}}(\kappa)<\kappa^+$ and $(IH)$. Let $D,E\in V[A]$ be Mathias sets such that  $\lambda:=|D|<\kappa$. Denote $\theta=\max\{\lambda,\kappa^*\}$, Then there is $F\in V[A]$ such that:
\begin{enumerate}
    \item $F$ is a Mathias set. $F\cap\theta= F_\theta$.
    \item $(D\cup E)\setminus\theta\subseteq F\subseteq \sup(D\cup E)$.
    \item $D,E\in V[F]$.
\end{enumerate}
\end{proposition}
\begin{remark}\label{problemunion}
Note that simply taking the union $D\cup E$ will not suffice for the proposition:

For example, assume that $o^{\vec{U}}(\kappa)=\delta$ and $o^{\vec{U}}(\delta)=1$, and pick any generic $G$ with the condition $\l \l \delta, \{\alpha<\delta\mid o^{\vec{U}}(\alpha)=0\}\r,\l\kappa,\{\delta<\alpha<\kappa\mid o^{\vec{U}}(\alpha)<\delta\}\r\r\in G$. Then $G$ is generic such that $\otp(C_G)=C_G(\omega)=\delta$. Let $$D=\{C_G(C_G(n))\mid n<\omega\}\text{ and }E=\{C_G(\alpha)\mid \omega\leq\alpha<C_G(\omega)\}\setminus D$$ Then $D\cup E=\{C_G(\alpha)\mid \omega\leq \alpha< C_G(\omega)\}$, hence in $V[D\cup E]$, $C_G(\omega)$ is still measurable. On the other hand, from $D$, we can reconstruct $\l C_G(n)\mid n<\omega\r$ as $o^{\vec{U}}(C_G(C_G(n)))=C_G(n)$. So it is impossible that $D\in V[D\cup E]$.
\end{remark}
\textit{Proof of \ref{countableinformation}.} 
 Fix $\Mfor$-names $\lusim{E},\langle \lusim{d}_i\mid i<\lambda\rangle$ for the elements of $E\setminus\theta$ and $D\setminus\theta$ respectively. Split the forcing at $\theta$ and find $\l q',r'\r\in G$ such that 
 $$(1)\ \ \ \l q',r'\r\Vdash \lusim{E},\{ \lusim{d}_i\mid i<\lambda\}\subseteq\lusim{C}_G\setminus\theta\text{ and } \forall\alpha\in \lusim{C}_G\setminus\kappa^*.o^{\vec{U}}(\alpha)<\alpha^+.$$
The idea is that for every $\delta\in D\setminus \kappa(r')$, there is $i\leq l(r')+1$ such that $\delta\in (\kappa_{i-1}(r),\kappa_i(r))$.  Then $\delta$ is definable from $D\cup E$ and two other parameters: $$\gamma(\delta):=o^{(\kappa_i(r'))}(\delta)\text{ and }\beta(\delta):=\sup(x\in (D\cup E)\cap \delta\mid \gamma(x)\geq\gamma(\delta)).$$
 Indeed, 
 $$\delta=\min(y\in (D\cup E)\setminus\beta(\delta)\mid \gamma(y)=\gamma(\delta)).$$
 Then $\beta(\delta)$ is a member of $E\cup D$ below $\delta$\footnote{Actually, since we can always shrink the large set of $\delta$ to filter from a final segment of $E\cup D$ ordinals $\rho$ with $\gamma(\rho)\geq\gamma(\delta)$, it will follow that $\beta(\delta)$ is strictly below $\delta$}. As for $\gamma(\delta)$, we use \ref{specialTree}, there is a $\vec{U}$-fat tree $T$ deciding $\delta$ to be the top most ordinal in a maximal branche of $T$, and  $\gamma(\delta)$ will be decided by the lower part of the branch, and hence below $\delta$, and therefore by finitely many elements of $C_G$ below $\delta$. After adding these finitely many element to $E$, we repeat this process on the added points. This process should stabilize after $\omega$ many steps, since we are creating a decreasing sequence of ordinals.

Formally, proceed by a density argument, let $r'\leq r\in \Mfor\restriction(\theta,\kappa]$.
 Define recursively for every $k<\omega$:
 $r\leq^* r^*_k$, maximal anti chains $\l Z^{(k)}_{i,j}\mid i<\lambda,j<\omega\r$, $\Mfor$-names  $$\l \lusim{\delta}^{(k)}_{i,j}\mid i<\lambda,j<\omega\r\text{ and }\l T^{(k)}_{q,i,j},I^{(k)}_{q,i,j},F^{(k)}_{q,i,j},\vec{A}^{(k)}_{q,i,j}\mid i<\lambda,j<\omega,q\in Z^{(k)}_{i,j}\r.$$ First for every $j<\omega$ and $i<\lambda$, let $\lusim{\delta}^{(0)}_{i,j}=\lusim{d}_i$. Assume  $r\leq^*r^*_k$ and $\lusim{\delta}^{(k)}_{i,j}$ are defined such that such that for all $i<\lambda,j<\omega$,
 $\l q',r^*_k\r\Vdash \lusim{\delta}^{(k)}_{i,j}\in \lusim{C}_G\setminus\theta$.
 
Fix $i<\lambda,j<\omega$, use \ref{specialTree} to find $r^*_k\leq^*r_{i,j}$ and a maximal antichain $Z^{(k)}_{i,j}\subseteq \Mfor\restriction\theta$ above $q'$, such that for every $q\in Z^{(k)}_{i,j}$, either $\l q,r_{i,j}\r ||\lusim{\delta}^{(k)}_{i,j}$, or there is a $\vec{U}$-fat tree of extensions of $r_{i,j}$ $T^{(k)}_{q,i,j}$ and sets $A^t_{q,i,j}$ such that $D_{T^{(k)}_{q,i,j},\vec{A}^{(k)}_{q,i,j}}$ is pre-dense above $r_{i,j}$, and for every  $t^{\smallfrown}\alpha\in mb(T^{(k)}_{q,i,j})$, $$(2) \ \ \ \l q,r_{i,j}^{\smallfrown}\l t^{\smallfrown}\alpha,\vec{A}^{t^{\smallfrown}\alpha}_{q,i,j}\r\r\Vdash\lusim{\delta}^{(k)}_{i,j}=\alpha, \ \ \ \l q,r_{i,j}^{\smallfrown}\l t,\vec{A}^t_{q,i,j}\r\r||\ o^{(\lusim{\kappa}^{(k)}_{i,j})}(\lusim{\delta}^{(k)}_{i,j})$$
where $\lusim{\kappa}^{(k)}_{i,j}$ is an $\Mfor$-name for the unique $\kappa_y(r)$, $y\leq l(r)+1$ such that $\lusim{\delta}^{(k)}_{i,j}\in (\kappa_{y-1}(r),\kappa_{y}(r))$. Note that  $\lusim{\kappa}^{(k)}_{i,j}$ is also a $\Mfor$-name for the measurable on which $mb(T^{(k)}_{q^*,i,j})$ splits, for the unique $q^*$ in $Z^{(k)}_{i,j}\cap G\restriction \theta$.     %Since 
%, there is a $\Mfor$-name,  $\lusim{\gamma}^{(k)}_{q,i,j}$ such that $\l q,r_{i,j}\r\Vdash o^{\kappa^{(k)}_{q,i,j}}(\lusim{\delta}^{(k)}_{i,j})=\lusim{\gamma}^{(k)}_{q,i,j}$ (in case that $\l q,r_{i,j}\r||\lusim{\delta}^{(k)}_{i,j}$, then take a canonical name for the value decided for $o^{\vec{U}}(\lusim{\delta}^{(k)}_{i,j})$).
Let $F^{(k)}_{q,i,j}:\Lev_{ht(T^{(k)}_{q,i,j})-1}(T^{(k)}_{q,i,j})\rightarrow \kappa$ be the function defined by
$$(3) \ \ \ F^{(k)}_{q,i,j}(s)=\gamma\leftrightarrow \l q,r_{i,j}^{\smallfrown}\l s,\vec{A}^{s}_{q,i,j}\r\r\Vdash o^{(\lusim{\kappa}^{(k)}_{i,j})}(\lusim{\delta}^{(k)}_{i,j})=\gamma.$$
This notation works in case that $\l q,r_{i,j}\r||\lusim{\delta}^{(k)}_{i,j}$ by taking the tree of height $0$ and $F^{(k)}_{q,i,j}(\l\r)$ is the decided value for $o^{\vec{U}}(\lusim{\delta}^{(k)}_{i,j})$.
Shrink $T^{(k)}_{q,i,j}$, and find a complete and consistent set of important coordinates $I^{(k)}_{q,i,j}$. Also as in \ref{very short}, we shrink the trees even more so that for every $q_1,q_2\in Z^{(k)}_{i,j}$ one of the following holds:
\begin{enumerate}
    \item[(4.1)] $Im(F^{(k)}_{q_1,i,j})\cap Im(F^{(k)}_{q_2,i,j})=\emptyset$.
    \item[(4.2)] $T^{(k)}_{q_1,i,j}\restriction I^{(k)}_{q_1,i,j}=T^{(k)}_{q_2,i,j}\restriction I^{(k)}_{q_2,i,j}$ and $(F^{(k)}_{q_1,i,j})_{I^{(k)}_{q_1,i,j}}=(F^{(k)}_{q_2,i,j})_{I^{(k)}_{q_2,i,j}}$.
\end{enumerate}
Note that for every $V$-generic filter $H\subseteq \Mfor$ such that $\l q,r_{i,j}\r\in H$, there is $t\in mb(T^{(k)}_{q,i,j})$ such that $\l q,r_{i,j}^{\smallfrown}\l t,\vec{A}^t_{q,i,j}\r\r\in H$, and if $t_1^{\smallfrown}\alpha_1,t_2^{\smallfrown}\alpha_2\in mb(T^{(k)}_{q,i,j})$ are two such branches, then by $(2)$ $\alpha_1=(\lusim{\delta}^{(k)}_{i,j})_H=\alpha_2$ and in particular $F^{(k)}_{q,i,j}(t_1)=F^{(k)}_{q,i,j}(t_2)$ which implies that $t_1\restriction I^{(k)}_{q,i,j}=t_2\restriction I^{(k)}_{q,i,j}$, thus $t_1\restriction I^{(k)}_{q,i,j}$  is unique.  Let  $\lusim{\vec{\alpha}}^{(k)}_{q,i,j}$ be a $\Mfor$-name such that
$$(5) \ \ \ \l q,r_{i,j}\r\Vdash \forall t\in mb(T^{(k)}_{q,i,j}). \ \l q,r_{i,j}^{\smallfrown}\l t,\vec{A}^{t}_{q,i,j}\r\r\in \lusim{G}\rightarrow \lusim{\vec{\alpha}}^{(k)}_{q,i,j}=t\restriction I^{(k)}_{q,i,j}.$$ Note that if $q_1,q_2\in Z^{(k)}_{i,j}$ are such that $(4.2)$ holds, then both $\l q_1, r_{i,j}\r,\l q_2,r_{i,j}\r$ force that $\lusim{\vec{\alpha}}^{(k)}_{q_1,i,j}=\lusim{\vec{\alpha}}^{(k)}_{q_2,i,j}$. Moreover, it is forced by $\l q,r_{i,j}\r$ that $|\lusim{\vec{\alpha}}^{(k)}_{q,i,j}|=|I^{(k)}_{q,i,j}|$ so we assume that $\lusim{\vec{\alpha}}^{(k)}_{q,i,j}=\l \lusim{\vec{\alpha}}^{(k)}_{q,i,j}(w)\mid w\leq |I^{(k)}_{q,i,j}|\r$. 
 Next, let $\lusim{\beta}^{(k)}_{q,i,j}$ be a $\Mfor$-name such that $$(6) \ \ \l q,r_{i,j}\r\Vdash\lusim{\beta}_{q,i,j}^{(k)}=\sup(\{x\in(\lusim{D}\cup\lusim{E})\cap\lusim{\delta}^{(k)}_{i,j}\mid o^{(\lusim{\kappa}^{(k)}_{i,j})}(\lusim{\delta}^{(k)}_{i,j})\leq o^{(\lusim{\kappa}^{(k)}_{i,j})}(x)\}\cup\{\theta\}).$$
 By definition of $\lusim{\beta}^{(k)}_{q,i,j}$ and since we split the forcing at $\theta$, the trees $T^{(k)}_{q,i,j}$ are extension trees of $r_{i,j}$ and for every $w$, $$(7)\ \ \ \l q,r_{i,j}\r \Vdash \lusim{\vec{\alpha}}^{(k)}_{q,i,j}(w),\lusim{\beta}^{(k)}_{q,i,j}\in C_{\lusim{G}}\cap[\theta,\lusim{\delta}^{(k)}_{i,j})$$ 
 just otherwise, there is a generic $H$ with $\l q,r_{i,j}\r \in H$ and $(\lusim{\beta}^{(k)}_{i,j})_H=(\lusim{\delta}^{(k)}_{i,j})_H$. However, by \ref{increasing order at limits}, $$o^{((\lusim{\kappa}^{(k)}_{i,j})_H)}((\lusim{\beta}^{(k)}_{i,j})_H)>o^{((\lusim{\kappa}^{(k)}_{i,j})_H)}((\lusim{\delta}^{(k)}_{q,i,j})_H)$$ contradiction.
  By $\leq^*$-closure of $\Mfor\restriction(\theta,\kappa]$, find a single $r^*_{k+1}$ such that $r_{i,j}\leq^* r^*_{k+1}$ for every $i,j$.
We conclude that for every $q\in Z^{(k)}_{i,j} $, we have defined $T^{(k)}_{q,i,j},F^{(k)}_{q,i,j},I^{(k)}_{q,i,j}$ and names $\l \lusim{\vec{\alpha}}^{(k)}_{q,i,j}(w)\mid w\leq |I^{(k)}_{q,i,j}|\r,\lusim{\beta}^{(k)}_{q,i,j}$.
We would like to turn these names to be independent of $q\in Z^{(k)}_{i,j}$.  For $\lusim{\beta}^{(k)}_{q,i,j}$  it is easy to find $\Mfor$-names $\lusim{\beta}^{(k)}_{i,j}$ such that for every $q\in Z^{(k)}_{i,j}$, $\l q,r^*_{k+1}\r\Vdash  \lusim{\beta}^{(k)}_{q,i,j}=\lusim{\beta}^{(k)}_{i,j}$. As for $\l \lusim{\vec{\alpha}}^{(k)}_{q,i,j}(w)\mid w\leq |I^{(k)}_{q,i,j}|\r$, the length $|I^{(k)}_{q,i,j}|$ might depend on $q$, so we define $\lusim{\vec{\alpha}}^{(k)}_{q,i,j}(w)=\theta$ if $|I^{(k)}_{q,i,j}|<w<\omega$, and we can find names $\lusim{\vec{\alpha}}^{(k)}_{i,j}(w)$ independent of $q$.
With these new names, in $(6),(7)$ we can replace $\l q,r_{i,j}\r$ by $\l q',r^*_{k+1}\r$.
Enumerate the names $$\{\lusim{\vec{\alpha}}^{(k)}_{i,j}(w),\lusim{\beta}^{(k)}_{i,j}\mid j,w<\omega\}=\{\lusim{\delta}^{(k+1)}_{i,s}\mid s<\omega\}.$$
 This concludes the inductive definition. Use $\sigma$-closure to find $r^*_n\leq^* r_\omega$, and shrink all the trees to be extension trees of $r_{\omega}$ such that for every $i<\lambda,\ k,j<\omega$ and $q\in Z^{(k)}_{i,j}$, $D_{T^{(k)}_{q,i,j},\vec{A}^{(k)}_{q,i,j}}$ is pre-dense above $r_\omega$. By density there is such $r_\omega\in G$. Define
$$\langle (\lusim{\delta}^{(k)}_{i,j})_G\mid k,j<\omega, i<\lambda\rangle.$$ By $(7)$,
$\l q',r_\omega\r\Vdash \lusim{\delta}^{(k)}_{i,j}\in\lusim{C}_G\setminus\theta$, thus $(\lusim{\delta}^{(k)}_{i,j})_G\in C_G\setminus\theta$. 

\begin{claim*}
$\langle (\lusim{\delta}^{(k)}_{i,j})_G\mid k,j<\omega, i<\lambda\rangle\in V[A]$.
\end{claim*} 
\textit{Proof of claim:} Work inside $V[A]$, recall that $D,E\in V[A]$, therefore 
$\langle (\lusim{\delta}^{(0)}_{i,j})_G\mid i<\lambda,j<\omega\rangle$ is in $V[A]$. Assume we have successfully defined
$\langle (\lusim{\delta}^{(k)}_{i,j})_G\mid i<\lambda,j<\omega\rangle$, let us define inside $V[A]$ from this sequence the sequence $\langle (\lusim{\delta}^{(k+1)}_{i,j})_G\mid i<\lambda,j<\omega\rangle$.
First, in $V[G]$, for each $i<\lambda,j<\omega$,  let 
$Z^{(k)}_{i,j}\cap G\restriction\theta=\{q^G_{i,j}\}$ 
and let $t_{i,j}\in mb(T^{(k)}_{q^G_{i,j},i,j})$ 
such that $\l q^G_{i,j},r_\omega^{\smallfrown}\l t_{i,j},\vec{A}^{t_{i,j}}_{q^G_{i,j},i,j}\r\in G$. Let 
 $y\leq l(r_\omega)+1$ be such that $(\lusim{\kappa}^{(k)}_{i,j})_G=\kappa_y(r_\omega)$, which is definable in $V[A]$ using $r_\omega, (\lusim{\delta}^{(k)}_{i,j})_G$, as the unique $y\leq l(r_\omega)+1$ such that $(\lusim{\delta}^{(k)}_{i,j})_G\in(\kappa_{y-1}(r_\omega),\kappa_y(r_\omega))$. By $(3)$,
 $$\l q^G_{i,j},r_\omega^{\smallfrown}\l t_{i,j}\setminus\{\max(t_{i,j})\},\vec{A}^{t_{i,j}\setminus\{\max(t_{i,j})\}}_{q^G_{i,j},i,j}\r\Vdash o^{(\kappa_y(r_\omega))}(\lusim{\delta}^{(k)}_{i,j})=F^{(k)}_{q^G_{i,j},i,j}(t_{i,j}\setminus\{\max(t_{i,j})\})$$
  hence it must be that $F^{(k)}_{q^G_{i,j},i,j}(t_{i,j}\setminus\{\max(t_{i,j})\})=o^{(\kappa_y(r_\omega))}((\lusim{\delta}^{(k)}_{i,j})_G)$.
Although the sequence $\l q^{G}_{i,j}\mid i<\lambda,j<\omega\r$ might not be in $V[A]$, we can do something similar to \ref{very short}. Back in $V[A]$, $o^{(\kappa_y(r_\omega))}((\lusim{\delta}^{(k)}_{i,j})_G)$ is definable since in $V$ we have the decomposition $$\l X^{(\kappa_y(r_\omega))}_\gamma\mid \gamma<o^{\vec{U}}(\kappa_y(r_\omega))\r$$ and  $o^{(\kappa_y(r_\omega))}((\lusim{\delta}^{(k)}_{i,j})_G)$ is the unique $\gamma_{i,j}<o^{\vec{U}}(\kappa_y(r_\omega))$ such that $(\lusim{\delta}^{(k)}_{i,j})_{G}\in X^{(\kappa_y(r_\omega))}_{\gamma_{i,j}}$. Let
$$M^{(k)}_{i,j}=\{q\in Z^{(k)}_{i,j}\mid o^{(\kappa_y(r_\omega))}((\lusim{\delta}^{(k)}_{i,j})_G)\in Im(F^{(k)}_{q,i,j})\}.$$
Notice that $q^G_{i,j}\in M^{(k)}_{i,j}$, as witnessed by $t_{i,j}\setminus\{\max(t_{i,j})\}$, hence $Im(F^{(k)}_{q,i,j})\cap Im( F^{(k)}_{q^G_{i,j},i,j})\neq\emptyset$ for any $q\in M^{(k)}_{i,j}$ and we conclude that  $(4.2)$ must hold. Choose in $V[A]$ any $q^{(k)}_{i,j}\in M^{(k)}_{i,j}$ and any $s^{(k)}_{i,j}\in mb(T^{(k)}_{q^{(k)}_{i,j},i,j})$ such that $F_{q^{(k)}_{i,j},i,j}^{(k)}(s^{(k)}_{i,j})=o^{(\kappa_y(r_\omega))}((\lusim{\delta}^{(k)}_{i,j})_G)$.  By $(5)$, $(\lusim{\vec{\alpha}}_{i,j})_G=(t_{i,j})\restriction I^{(k)}_{q^G_{i,j},i,j}$ and since $$(F^{(k)}_{q^{(k)}_{i,j},i,j})_{I^{(k)}_{q^{(k)}_{i,j},i,j}}(s_{i,j}\restriction I^{(k)}_{q^{(k)}_{i,j},i,j})=o^{(\kappa_y(r_\omega))}((\lusim{\delta}^{(k)}_{i,j})_G)=(F^{(k)}_{q^G_{i,j},i,j})_{I^{(k)}_{q^G_{i,j},i,j}}(t_{i,j}\restriction I^{(k)}_{q^{G}_{i,j},i,j})$$ it follows that  $t_{i,j}\restriction I^{(k)}_{q^{G}_{i,j},i,j}=(\lusim{\vec{\alpha}}_{i,j})_G=s_{i,j}\restriction I^{(k)}_{q^{(k)}_{i,j},i,j}$. Hence
$\l(\lusim{\vec{\alpha}}^{(k)}_{i,j}(w))_G\mid w<\omega\r$ is definable in $V[A]$.
Also, by $(6)$, $(\lusim{\beta}^{(k)}_{i,j})_G$ is definable from $(\lusim{\delta}^{(k)}_{i,j})_G$, $\kappa_{y}(r_{\omega})$ and $D\cup E$ which are all available in $V[A]$. By definition of the sequence  $\langle (\lusim{\delta}^{(k+1)}_{i,j})_G\mid i<\lambda,j<\omega\rangle$ it is definable in $V[A]$. So we conclude that $\langle (\lusim{\delta}^{(k)}_{i,j})_G\mid k,j<\omega, i<\lambda\rangle\in V[A]$. $\blacksquare_{\text{Claim}}$

We keep the notation of $q^{(k)}_{i,j}$ from the proof of the claim, use Proposition \ref{CodingBoundedInformation} to find $F_{\theta}$ such that
$$V[F_{\theta}]=V[A]\cap V[C_G\cap\theta].$$
Define $$F_*=\{(\lusim{\delta}^{(k)}_{i,j})_G\mid k,j<\omega, i<\lambda\}, \ \ F^*=(E\cup F_*)\setminus \theta \uplus F_\theta\in V[A].$$ 
Clearly, $F^*$ is a Mathias set and $F^*\cap\theta=F_\theta$. To see $2$ of the proposition, note that $D\setminus \theta=\{(\lusim{\delta}^{(0)}_i)_G\mid i<\lambda\}\subseteq F^*$, it follows that $D\cup E\setminus\theta\subseteq F^*$. Moreover from $(6)$ it follows that for every $k,i,j$, $\lusim{\delta}^{(k)}_{i,j}$ is forced by $\l q', r_\omega\r\in G$ to be below some $\lusim{\delta}^{(0)}_{s,t}$, so $\sup(F_*)=\sup(D)$, hence $\sup(F^*)=\sup(D\cup E)$.
 To see $3$, let $\langle\lambda_\xi\mid \xi<otp(F_*)=:\rho\rangle$ be the increasing enumeration of $F_*$, clearly $|\rho|\leq \lambda$.% and inside $V[F_*]$,  we can order the sequence of name $\l \lusim{\delta}^{(k)}_{i,j}\mid k,j<\omega,i<\lambda\r=\l \lusim{\lambda}_\xi\mid \xi<\rho\r$ accordingly.
 
  Consider the function $R:\rho \rightarrow[\rho]^{<\omega}$ defined by $R(\xi)=\l \l i_1,...,i_n\r,s\r$ such that for some $i,j,k$,  $$\lambda_\xi=(\lusim{\delta}^{(k)}_{i,j})_G, \ \l(\vec{\lusim{\alpha}}^{(k)}_{i,j}(w))_G\mid w\leq |I^{(k)}_{q^{(k)}_{i,j},i,j}|\r=\l\lambda_{i_1},...,\lambda_{i_n}\r\text{ and } (\lusim{\beta}^{(k)}_{i,j})_G=\lambda_s.$$ By the claim, both $\l(\lusim{\delta}^{(k)}_{i,j})_G\mid k,j<\omega, i<\lambda\r\in V[A]$, hence $R\in V[A]$, since $|\rho|\leq\lambda$, then $R\in V[A]\cap V[C_G\cap\theta]=V[F_\theta]\subseteq V[F^*]$. Notice that by $(7)$, $i_1,...,i_n,s<i$. 
 Let us argue first that $F_*\in V[F^*]$, in $V[F^*]$, we inductively define $\langle\beta_i\mid i<\rho\rangle$.
 Clearly $$\{\lambda_i\mid i<\rho\}\cap \theta+1=\{\lambda_i\mid i<\epsilon\}\in V[F_\theta]$$ so we let $\beta_i=\lambda_i$ for $i<\epsilon$.
Assume that $\langle \beta_j\mid j<i\rangle$ is defined, where $i>\epsilon$, in particular $\beta_{i_1},...\beta_{i_n}$ and $\beta_s$ are defined. Let $I=Ind(F_*\setminus D,F_*)\subseteq \rho$, by the claim, $I\in V[A]\cap V[C_G\cap \theta]= V[F_\theta]\subseteq V[F^*]$. Finally, note that $\{q^{(k)}_{i,j}\mid i<\lambda,j<\omega\}\in V[A]\cap V[C_G\cap\theta]=V[F_\theta]\subseteq V[F^*]$ and let $\kappa_{i,j}$ be the measurable on which $mb(T^{(k)}_{q{(k)}_{i,j},i,j})$ splits. Define $$\beta_i=\min(\{x\in (F^*\setminus \{\beta_j\mid j\in I\cap i\})\setminus\beta_{s}+1\mid \ o^{(\kappa_{i,j})}(x)\geq (F^{(k)}_{q^{(k)}_{i,j},i,j})_{I^{(k)}_{q^{(k)}_{i,j},i,j}}(\beta_{i_1},...,\beta_{i_n})\}).$$ This is a legitimate definition in $V[F^*]$ since we worked hard to ensure all the parameters used are there. Let us prove that $\beta_\xi=\lambda_\xi$, inductively assume that $\langle\beta_j\mid j<\xi\rangle=\langle\lambda_j\mid j<\xi\rangle$, we can assume that $\xi>\epsilon$, then $$\{\beta_j\mid j\in I\cap \xi\}=\{\lambda_j\mid j\in I\cap \xi\}=(F_*\setminus D)\cap \lambda_\xi $$ 
 and therefore
 $$(F^*\setminus\{\beta_j\mid j\in I\cap \xi\})\cap(\beta_s,\lambda_\xi)=[(E\cup F_*)\setminus(F_*\setminus D)]\cap(\beta_s,\lambda_\xi)=(E\cup D)\cap(\beta_s,\lambda_\xi).$$
 Assume that $i,k,j$ are such that $\lambda_\xi=(\lusim{\delta}^{(k)}_{i,j})_G$, then by induction hypothesis, $\beta_s=\lambda_s=(\lusim{\beta}^{(k)}_{i,j})_G$ and
  $$\l(\vec{\lusim{\alpha}}^{(k)}_{i,j}(w))_G\mid w\leq |I^{(k)}_{q'_{i,j},i,j}|\r=\l\lambda_{i_1},...,\lambda_{i_n}\r=\l\beta_{i_1},...,\beta_{i_n}\r.$$ By $(3)$ it follows that $$(F^{(k)}_{q^{(k)}_{i,j},i,j})_{I^{(k)}_{q^{(k)}_{i,j},i,j}}(\l\beta_{i_1},...,\beta_{i_n}\r)=o^{(\kappa_{i,j})}((\lusim{\delta}^{(k)}_{i,j})_G)=o^{(\kappa_{i,j})}(\lambda_\xi).$$
 By $(6)$, it follows that in the interval $(\beta_s,\lambda_\xi)$, there are no ordinals $x\in F^*\setminus\{\beta_j\mid j\in I\cap \xi\}$ such that $(F^{(k)}_{q^{(k)}_{i,j},i,j})_{I^{(k)}_{q^{(k)}_{i,j},i,j}}(\l\beta_{i_1},...,\beta_{i_n}\r)\leq o^{(\kappa_{i,j})}(x)$ so $\beta_\xi\geq \lambda_\xi$.  Also $\lambda_\xi\in F^*\setminus\{\beta_j\mid j\in I\cap \xi\}$ and $F^{(k)}_{q'_{i,j},i,j}(\beta_{i_1},...,\beta_{i_k})=o^{(\kappa_{i,j})}(\lambda_\xi)$ hence $\lambda_\xi= \beta_\xi$. Thus $F_*\in V[F^*]$. 
 From this $(3)$ easily follows, indeed, $D\setminus\theta,F_*\setminus E\in V[F^*]$ since their indices inside $F_*$ are subsets of $\theta$, hence $$E\setminus \theta=[(E\cup F_*)\setminus (F_*\setminus E)]\setminus\theta=F^*\setminus[\theta\cup(F_*\setminus E)]\in V[F^*].$$ Also $D\cap\theta,E\cap\theta\in V[F_\theta]\subseteq V[F^*]$ and therefore $D,E\in V[F^*]$ which is what we needed.$\blacksquare$

The following corollary provides a sufficient condition for the main result. It roughly says that given that $\kappa$ changes cofinality in $V[A]$, and given a single $C'\subseteq C_G$ which captures all the initial segments of $A$, we can glue the information needed to capture $A$. 
\begin{lemma}\label{lemmaforsubsetkappa}
Assume $o^{\vec{U}}(\kappa)<\kappa^+$ and $(IH)$. Let $A \in V[G] ,\ A\subseteq\kappa $ and assume that $\exists C^*\subseteq C_G$ such that
\begin{enumerate}
 \item$C^*\in V[A]$ and
 $\forall\alpha<\kappa \  A\cap\alpha\in V[C^*]$.
 \item $cf^{V[A]}(\kappa)<\kappa$.
 \end{enumerate}
Then $ \exists C' \subseteq C_G $ such that $ V[A]=V[C']$.
\end{lemma}
\pr Let $\lambda:=cf^{V[A]}(\kappa)<\kappa$ and $\langle \alpha_i\mid i<\lambda\rangle\in V[A]$ unbounded and cofinal in $\kappa$ witnessing this. By \ref{very short}, there is $C_*\subseteq C_G$ such that $|C_*|\leq\lambda$ and $V[C_*]=V[\langle \alpha_i\mid i<\lambda\rangle]$. Use \ref{countableinformation} to find $C_0\subseteq C_G$ such that $C_0\in V[A]$ and $C_*,C^*\in V[C_0]$. In $V[C_0]$, let $\pi_i:2^{\alpha_i}\leftrightarrow P(\alpha_i)$ be any bijection. Since $A\cap\alpha_i\in V[C_0]$, there is $\delta_i$ such that $$\pi_i(\delta_i)=A\cap\alpha_i.$$ Note that the sequence    $\langle \delta_i\mid i<\lambda\rangle$ might not be inside $V[C_0]$, but it is in $V[A]$. Again by \ref{very short} we can find $C''\subseteq C_G$ such that $|C''|\leq\lambda$ such that $$V[\langle\delta_i\mid i<\lambda\rangle]=V[C''].$$ By Proposition \ref{countableinformation}, we can find some $C'\subseteq C_G$, $C'\in V[A]$,  such that $C_0, C''\in V[C']$. Now in $V[C']$ we can compute $A$ as follows, since $C_0\in V[C']$, also $\l\pi_i\mid i<\lambda\r\in V[C']$, and since $C''\in V[C']$ also $\l\delta_i\mid i<\lambda\r\in V[C']$. It follows that $A=\cup_{i<\lambda}A\cap\alpha_i=\cup_{i<\lambda}\pi_i(\delta_i)\in V[C']$. $\blacksquare$

\subsection{Subsets of $\kappa$ which do not stabilize}

In this section we assume that  $o^{\vec{U}}(\kappa)<\kappa^+$, $A$ does not stabilize and $(IH)$. We do not assume in general that $A\subseteq\kappa$.  However, if $A\in V[G]$ is such that $A\subseteq\kappa$ and does not stabilize, then by \ref{nonstabcofinalitysubsetofkappa}, $cf^{V[A]}(\kappa)<\kappa$. By Lemma \ref{lemmaforsubsetkappa}, to conclude the main result for $A$, it remains to find $C^*\in V[A]$ such that for every $\alpha<\kappa$, $A\cap\alpha\in V[C^*]$. Along this chapter we construct such $C^*$. The naive approach is the following:
Fix a cofinal sequence $\l\alpha_i\mid i<cf^{V[A]}(\kappa)\r\in V[A]$, since for every $i,$ $A\cap\alpha_i$ is bounded, apply \ref{very short} to find $C_i\subseteq C_G$ such that $V[A\cap\alpha_i]=V[C_i]$ and let $C^*=\cup_{i<cf^{V[A]}(\kappa)}C_i$. There are several reasons why $C^*$ in not the desired set: 
\begin{enumerate}
    \item[(I)] The sequence $\l C_i\mid i<cf^{V[A]}(\kappa)\r$ is defined in $V[G]$ and by adding finitely many elements to each $C_i$ we might accumulate an infinite sequence which is not in $V[A]$. 
    \item[(II)] As we have seen in \ref{problemunion}, a union of two sets might lose information, so it is possible that for some $j$,  $C_j\notin V[\cup_{i<cf^{V[A]}(\kappa)}C_i]$. 
\end{enumerate}
For problem (I), we need to ensure that the choice we make is inside $V[A]$, for this we use the definition of a Mathias set, in $V[A]$ we can choose a sequence $\l D_i\mid i<cf^{V[A]}(\kappa)\r$ such that $V[A\cap\alpha_i]=V[D_i]$ and each $D_i$ is a Mathias set. By Proposition  \ref{FiniteNoise}, $D_i\subseteq^* C_G$, so it might be that $D_i\setminus C_G\neq\emptyset$. By fixing problem I,
we have created a new problem: The sequence $D:=\cup_{i<cf^{V[A]}(\kappa)}D_i$ might accumulate infinite noise i.e. $|D\setminus C_G|\geq \omega$. Lemma \ref{Noise} and corollaries \ref{stabilization of star increasing sequences}, \ref{subsets star bound}, show we can remove this noise and stay inside $V[A]$.
\begin{lemma}\label{Noise}

Let $\langle D_i\mid i<\lambda\rangle\in V[A]$ such that $\lambda<\kappa$ and:
\begin{enumerate}
\item $D_i$ is a Mathias set.
    \item $min(D_i)\geq\lambda$.
\end{enumerate} 
Then there is $\langle D^*_i\mid i<\lambda\rangle\in V[A]$ such that:
\begin{enumerate}
    \item $\underset{i<\lambda}{\bigcup}D^*_i$ is Mathias.
    \item $\forall i<\lambda,  D_i=^*D^*_i\subseteq D_i$.
\end{enumerate}
\end{lemma}
\pr By removing finitely many elements from every $D_i$, we can assume that $otp(D_i)$ is a limit ordinal. If every $D_i=\emptyset$, then the claim is trivial.  Otherwise, since $D_i$ is a Mathias set, $\sup(D_i)\in \overline{X}_A$. Denote $D=\underset{i<\lambda}{\bigcup}D_i$ and $\nu^*=\sup(D)>\lambda$. Note that $\nu^*\in \overline{X}_A$, since $\nu^*=\sup(
\sup(D_i)\mid i<\lambda)$ and $\overline{X}_A$ is closed. 

Proceed by induction on $\nu^*$, by Lemma \ref{FiniteNoise}, $D_i\setminus C_G$ is finite. It follows that $|D\setminus C_G|\leq\lambda<\nu^*$. We would like to remove the noise accumulated in $D$ by intersecting it with sets in $\cap\vec{U}(\nu^*)$. Since $\nu^*\in Lim(C_G)$, we can apply \ref{MagidorHousdorf} to $D\setminus C_G$ and find a set $Y^*\in\cap\vec{U}(\nu^*)$ such that $Y^*\cap(D\setminus C_G)=\emptyset$.
Denote $D^*=D\cap Y^*\subseteq C_G$. Note that $D^*\in V[A]$ since $D\in V[A]$ and $Y^*\in V$.

Consider the set $$
Z^{(0)}=\{\nu<\nu^*\mid Y^*\cap\nu\in\cap\vec{U}(\nu)\}$$
to see that $Z^{(0)}\in \cap\vec{U}(\nu^*)$, let $i<o^{\vec{U}}(\nu^*)$, then $j_{U(\nu^*,i)}(Y^*)\cap\nu^*=Y^{*}\in\underset{\xi<i}{\bigcap}U(\nu^*,\xi)$. By coherency, the order of $\nu^*$ in $j_{U(\nu^*,i)}(\vec{U})$ is $i$, which implies that   $$\underset{\xi<i}{\cap}U(\nu^*,\xi)=\cap j(\vec{U})(\nu^*).$$ By definition $\nu^*\in j(Z^{(0)})$ thus $Z^{(0)}\in U(\nu^*,i)$ for every $i<o^{\vec{U}}(\nu^*)$ and $Z^{(0)}\in\bigcap\vec{U}(\nu^*)$. By Proposition \ref{genericproperties}(3), there is $\eta_0<\nu^*$ such that $C_G\cap (\eta_0,\nu^*)\subseteq Z^{(0)}$.  

Consider the sequence of Mathias sets $\l D_i\cap\eta_0\mid i<\lambda\r$,
apply the induction hypothesis to it
and find $\langle D'_i\mid i<\lambda\rangle$ such that
\begin{enumerate}
    \item $\underset{i<\lambda}{\bigcup}D'_i$ is Mathias.
    \item $D_i\cap\eta_0=^*D'_i\subseteq \eta_0$.
\end{enumerate}
Define $$D^*_i=D'_i\uplus (D_i\cap Y^{*}\setminus \eta_0).$$

Let us argue that $\langle D^*_i\mid i<\lambda\r$ is as wanted:
to see condition $(1)$, note that the set
$$\underset{i<\lambda}{\cup}D_i^*=D^*\setminus\eta_0\cup (\underset{i<\lambda}{\cup}D_i')$$ is a Mathias sets as the union of two Mathias sets.

For condition $(2)$, it is clear that $D_i^*\subseteq^* D_i$. Toward a contradiction, assume that there is $i<\lambda$ and $\delta\leq\sup(D_i)$ is minimal such that 
$$|(D_i\cap\delta)\setminus (D^*_i\cap\delta)|\geq\omega.$$
By the definition of $D^*_i$, $\delta>\eta_0$ and $\delta\in Lim(D_i)$.  By the definition of $\eta_0$, $\delta\in C_G\cap(\delta_0,\nu^*)\in Z^{(0)}\cup\{\nu^*\}$ which means that $\delta\cap Y^{*}\in \bigcap\vec{U}(\delta)$.  Since $D_i$ is Mathias, there is $\xi<\delta$ such that $D_i\cap(\xi,\delta)\subseteq Y^{*}$, in particular $$D_i\cap(\xi,\delta)=D_i\cap Y^{*}\cap(\xi,\delta)=D^*_i\cap(\xi,\delta).$$
So $(D_i\cap\delta)\setminus(D^*_i\cap\delta)=(D_i\cap\xi)\setminus(D^*_i\cap\xi)$,
this is a contradiction to the minimality of $\delta$.
$\blacksquare$
 \begin{corollary} \label{stabilization of star increasing sequences}
 Let $\langle D_i\mid i<\theta\rangle\in V[A]$ such that $\theta<\kappa^+$ and:
\begin{enumerate}
\item $D_i$ is a Mathias set.
    \item $D_i\cap\kappa^*=F_{\kappa^*}$ where $V[F_{\kappa^*}]=V[A]\cap V[C_G\cap\theta]$.
    \item$ \l D_i\mid i<\theta\r$ is $\subseteq^*$-increasing.
\end{enumerate} 
Then there is $\langle D^*_i\mid i<\theta\rangle\in V[A]$ such that:
\begin{enumerate}
    \item $\underset{i<\theta}{\bigcup}D^*_i$ is a Mathias set.
    \item $\forall i<\theta,  D_i=^*D^*_i\subseteq D_i$.
    \item $D_i^*\cap \kappa^*=F_{\kappa^*}$
\end{enumerate}
 \end{corollary}
\pr Let $\lambda=cf^{V[A]}(\theta)\leq\kappa$. Since $\kappa$ is singular in $V[A]$,  $\lambda<\kappa$ and let $\l \theta_i\mid i<\lambda\r\in V[A]$ be cofinal in $\theta$. We split each $D_{\theta_i}$ to three intervals:
$$D_{\theta_i}=D_{\theta_i}\cap\kappa^*\uplus D_{\theta_i}\cap(\kappa^*,\lambda)\uplus D_{\theta_i}\setminus\lambda.$$
Denote these sets by $A_i,B_i,C_i$ respectively.
By assumption, $A_i$ is constantly $F_{\kappa^*}$. Apply \ref{Noise} to the sequence $\l C_i\mid i<\lambda\r$ to obtain $\l C^*_i\mid i<\lambda\r\in V[A]$ such that 
$C^*_i=^*C_i$ and $C^*:=\cup_{i<\lambda}C^*_i$ is Mathias.
As for the sequence $\l B_i\mid i<\lambda\r$, either $\lambda\leq\kappa^*$ in which case $B_i=\emptyset$. Otherwise $\lambda>\kappa^*$, and by removing finitely many points from $B_i$, we can assume that $\sup(B_i)\in Lim(B_i)\subseteq \overline{X}_A\setminus\kappa^*\subseteq X_A$ i.e. $\sup(B_i)$ is singular in $V[A]$. Since $\lambda>\kappa^*$ is regular in $V[A]$, it follows that $\lambda\notin X_A$, hence, $\sup(B_i)\leq\max(X_A\cap\lambda):=\mu<\lambda$.  Since $\mu>\aleph_0$ is a strong limit cardinal and $\lambda$ is regular, the sequence $\l B_i\mid i<\lambda\r$ satisfies the assumption of 
 Theorem \ref{Modfinitestab}, hence there is $\lambda'<\lambda$ such that for every $\lambda'\leq\delta<\lambda$, $B_\delta=^*B^*$. 

Note that $F_{\kappa^*}\cup B^*\cup C^*$ is a Mathias set as the union of finitely many of them.
Let $D_i^*:=D_i\cap (F_{\kappa^*}\cup B^*\cup C^*)$.

First, since $\cup_{i<\theta}D^*_i\subseteq F_{\kappa^*}\cup B^*\cup C^*$, and $F_{\kappa^*}\cup B^*\cup C^*$ is Mathias, then also $\cup_{i<\theta}D^*_i$ by the criteria of \ref{FiniteNoise}.
Also $(3)$ follows trivially. To see $(2)$,
It suffices to see that for each interval $$D_i^*\cap\kappa^*=D_i\cap\kappa^*, \ D_i^*\cap (\kappa^*,\lambda)=^*D_i\cap (\kappa^*,\lambda), \ D^*_i\setminus\lambda=^*D_i\setminus\lambda$$
indeed $D_i^*\cap\kappa^*=D_i\cap\kappa^*=F_{\kappa^*}$.
Find $\lambda'\leq\delta<\lambda$ such that $i<\theta_{\delta}$, then $D_i\subseteq^* D_{\theta_\delta}$. In particular, $$D_i\setminus\lambda\subseteq^* D_{\theta_\delta}\setminus\lambda= C_\delta=^*C^*_\delta\subseteq C^*\text{ and }D_i\cap(\kappa^*,\lambda)\subseteq^* D_{\theta_\delta}\cap(\kappa^*,\lambda)= B_\delta=^*B^*.$$
So 
$$D_i\setminus\lambda=^*D_i\cap C^*\setminus\lambda=D^*_i\setminus\lambda\text{ and }D_i\cap(\kappa^*,\lambda)=^*D_i\cap B^*\cap(\kappa^*,\lambda)=D^*_i\cap(\kappa^*,\lambda).$$
Therefore  $D_i=^*D^*_i$.$\blacksquare$
\begin{corollary}\label{subsets star bound}
Let $\langle D_i\mid i<\theta\rangle\in V[A]$ such that $\theta<\kappa^+$ and:
\begin{enumerate}
\item $D_i$ is a Mathias set.
    \item $D_i\cap\kappa^*=F_{\kappa^*}$.
    \item$ \l D_i\mid i<\theta\r$ is $\subseteq^*$-increasing.
\end{enumerate} 
Then in $V[A]$ there is a Mathias set $E\subseteq\sup(\cup_{i<\theta}D_i)$  which is a $\subseteq^*$-bounded for the sequence $\langle D_i\mid i<\theta\rangle$ such that $E\cap\kappa^*=F_{\kappa^*}$.
\end{corollary}
\pr Simply apply \ref{stabilization of star increasing sequences}, to find $\l D^*_i\mid i<\theta\r$ then $E=\cup_{i<\theta}D^*_i$ will be as wanted.$\blacksquare$

As for problem II mentioned in the beginning of this section, the first step will be to take $C^*$, the union of the $C_i$'s. Then we $\subseteq^*$ increase every $C^*\cap\alpha_i\subseteq^* C^{(1)}_i$ so the $C_i\in V[C^{(1)}_i]$. Repeating this process transfinitely, this will eventually stabilize to obtain the desired set. Note also that this definition must take place inside $V[A]$.

The following three propositions formally describe this process, we prove them by induction on $\nu\in X_A$. Recall that under the assumption of this section $\kappa\in X_A$.
\begin{theorem}\label{Making a sequence stae increasing}
Assume that $\nu\in X_A$, $\theta<\nu^+$ and let $\langle D_i\mid i<\theta\rangle\in V[A]$ such that: 
\begin{enumerate}
     \item $D_i\subseteq\theta_i<\nu$ is a Mathias set, $\l \theta_i\mid i<\theta\r$ is non decreasing.
     \item $D_i\cap\kappa^*=F_{\kappa^*}$.
 \end{enumerate} 
 Then there is $\langle D^*_i\mid i<\theta\rangle\in V[A]$ such that \begin{enumerate}
 \item $D^*:=\underset{i<\theta}{\bigcup}D^*_i$ is Mathias.
    \item $D_i\subseteq^*D^*_i\subseteq\theta_i$ and $D_i\in V[D^*_i]$.
    \item $D^*_i\cap \kappa^*=F_{\kappa^*}$.
    \item $\langle D^*_i\mid i<\theta\rangle$ is $\subseteq^*$-increasing.
\end{enumerate}
\end{theorem}

\begin{theorem}\label{maintheoremnonstab}
Assume that $\nu\in X_A$, $\theta<\nu^+$ and let $\langle D_i\mid i<\theta\rangle\in V[A]$ such that: 

\begin{enumerate}
     \item $D_i\subseteq\theta_i<\nu$ is a bounded in $\nu$ Mathias set, $\l \theta_i\mid i<\theta\r$ is non decreasing.
     \item $D_i\cap\kappa^*=F_{\kappa^*}$.
 \end{enumerate} 
 Then there is $\langle D^*_i\mid i<\theta\rangle\in V[A]$ such that \begin{enumerate}
 \item $D^*:=\underset{i<\theta}{\bigcup}D^*_i$ is Mathias.
    \item $\forall i<\theta. D^*_i\in V[D^*]$.
    \item $D_i\subseteq^*D^*_i\subseteq\theta_i$ and $D_i\in V[D^*_i]$.
    \item $D^*_i\cap \kappa^*=F_{\kappa^*}$.
    \item $\langle D^*_i\mid i<\theta\rangle$ is $\subseteq^*$-increasing.
\end{enumerate}
\end{theorem}

\begin{proposition}\label{boundedUnionOfgenerics}
Assume that $\nu\in X_A$, $D,D'\in V[A]$ are such that,
\begin{enumerate}
    \item $D,D'\subseteq\nu$ are Mathias sets.
    \item $D\cap\kappa^*=D'\cap\kappa^*=F_{\kappa^*}$.
\end{enumerate}
Then there is $D^*\in V[A]$ such that
\begin{enumerate}
\item $D^*$ is a Mathias set
    \item $D\cup D'\subseteq D^*\subseteq \sup(D\cup D')$.
    \item $D,D'\in V[D^*]$.
    \item $D^*\cap\kappa^*=F_{\kappa^*}$
    
\end{enumerate} 
\end{proposition}

As mentioned before, the proof of \ref{Making a sequence stae increasing}, \ref{maintheoremnonstab} and \ref{boundedUnionOfgenerics} is by induction on $\nu$. For $\nu\in X_A$ denote:
\begin{enumerate}
    \item $(18)_\nu$ is  Theorem \ref{Making a sequence stae increasing} for $\nu$.
    \item $(19)_{\nu}$ is  Theorem \ref{maintheoremnonstab} for $\nu$.
    \item $(20)_{\nu}$ is  Proposition \ref{boundedUnionOfgenerics} for $\nu$.
    \end{enumerate}
 Clearly, for every  $\nu\leq\kappa^*$, $(18)_{\nu}+(19)_{\nu}+(20)_{\nu}$ holds. Assume that $\nu>\kappa^*$, in particular, $cf^{V[A]}(\nu)<\nu$. Inductively assume that $(18)_{<\nu}+(19)_{<\nu}+(20)_{<\nu}$. The plan is to derive the induction step gradually from the following implications:
 \begin{enumerate}
     \item $(18)_{<\nu}+(19)_{<\nu}+(20)_{<\nu}\Longrightarrow (18)_{\nu}$.
     \item $(18)_{\nu}+(19)_{<\nu}+(20)_{<\nu}\Longrightarrow (19)_{\nu}$.
     \item $(18)_{\nu}+(19)_{\nu}+(20)_{<\nu}\Longrightarrow (20)_{\nu}$.
 \end{enumerate}

\textit{Proof of implication 1 (Theorem \ref{Making a sequence stae increasing}).}  Let us define inductively in $V[A]$ the sequence $\l D^*_i\mid i<\theta\r$, define $D^*_0=D_0$. At successor stage, the sets $D^*_\alpha,D_{\alpha+1}$ are bounded in $\nu$ in $\nu$, apply the induction hypothesis $(20)_{\theta_{\alpha+1}}$ to these sets, to find a Mathias set $D^*_{\alpha+1}$ such that $D^*_\alpha\cup D_{\alpha+1}\subseteq D^*_{\alpha+1}\subseteq\theta_{\alpha+1}$, $D^*_{\alpha+1}\cap\kappa^*=F_{\kappa^*}$ and $D_{\alpha+1}\in V[D^*_{\alpha+1}]$. 

At limit stage $\delta<\theta$, the sequence $\l D^*_i\mid i<\delta\r$ is defined and $\subseteq^*$ increasing. By \ref{subsets star bound}, there is a Mathias set $E^*$ such that $E^*\cap\kappa^*=F_{\kappa^*}$, $E^*\subseteq\sup\{\theta_i\mid i<\delta\}\leq \theta_\delta<\nu$ which is a $\subseteq^*$-bound. Again apply $(20)_{\theta_{\delta}}$ to $E^*,D_{\delta}$ to obtain $D^*_{\delta}\subseteq \theta_{\delta}$. Then $(2),(3),(4)$ are clear. At stage $\theta$, we also need to ensure $(1)$, by \ref{stabilization of star increasing sequences}, we can change the constructed $\l D_i^*\mid i<\theta\r$ to $\l D^{**}_i\mid i<\theta\r$ such that $D^*_i=^*D^{**}_i$, $D^{**}_i\cap\kappa^*=F_{\kappa^*}$ and $(1)$ holds. It suffices to note that $(2),(3),(4)$ still hold if we only change finitely many elements of $D^*_i$.$\blacksquare_{\text{Implication 1}}$

\textit{Proof of implication 2 (Theorem \ref{maintheoremnonstab}).} In the second implication we assume the induction hypothesis and also $(18)_\nu$ which was derived by the first implication.
The crucial difference between \ref{maintheoremnonstab} and \ref{Making a sequence stae increasing} is requirement $(2)$ that $D^*_i\in V[\cup_{j<\theta}D^*_j]$. 

Apply $(18)_\nu$ to the sequence $\l D_i\mid i<\theta\r$ to get $\langle D^0_i\mid i<\theta\rangle$ such that: \begin{enumerate}
\item $\underset{i<\theta}{\bigcup}D^0_i$ is Mathias.
    \item $D_i\subseteq^*D^0_i\subseteq\theta_i$ and $D_i\in V[D^0_i]$.
    \item $D^0_i\cap\kappa^*=F_{\kappa^*}$.
    \item $\langle D^0_i\mid i<\theta\rangle$ is $\subseteq^*$-increasing.
    
\end{enumerate}
Define a matrix of sets $\langle D^\xi_i\mid i<\theta,\xi<\nu^+\rangle$ recursively on the row $\xi<\nu^+$ such that:
\begin{enumerate}
    \item For each $\xi<\nu^+$, $\langle D^\xi_i\mid i<\theta\r$ is $\subseteq^*$- increasing. (Each row is $\subseteq^*$ increasing)
    \item For each $i<\theta$, $\langle D^\xi_i\mid \xi<\nu^+\rangle$ is $\subseteq^*$-increasing. (Each column is $\subseteq^*$ increasing)
    \item $D^\xi_i\subseteq\theta_i$ and $D_i\in V[D^\xi_i]$. (sets in column $i$ are subsets of $\theta_i$)
    \item $D^{(\xi)}:=\underset{j<\theta}{\bigcup}D^\xi_j$ is Mathias. (The union of each row is a Mathias set)
    \item  $D^\xi_i\cap\kappa^*=F_{\kappa^*}$. (All the sets are the same up to $\kappa^*$)
    \item  For every $i<\theta$ and every $\xi<\nu^+$, $D^{(\xi)}\cap\theta_i\subseteq^* D^{(\xi+1)}_i$. (The $i$-th set in a successor row, $\subseteq^*$ includes the union of the previous row up to $\theta_i$)
\end{enumerate}
At successor row, assume $\langle D^{\alpha}_i\mid i<\theta\rangle$ is defined.
For each $i<\theta$ apply $(20)_{\theta_1}$ to $D_i$ and $D^{(\alpha)}\cap\theta_i$ to obtain the sets $E^{(\alpha+1)}_i$ which satisfies $(2),(3),(5),(6)$. Apply $(18)_\nu$ to the sequence $\l E^{(\alpha+1)}_i\mid i<\theta\r$, obtain $E^{(\alpha+1)}_i\subseteq^*D^{(\alpha+1)}_i\subseteq\theta_i$, then also $(1),(4)$ holds without ruining $(2),(3),(5),(6)$. 

For limit $\delta<\nu^+$ the sequences $\langle D^{(\rho)}_i\mid i<\theta\rangle$ are defined for every $\rho<\delta$. For each $i<\theta$, the sequence $\l D^{(\rho)}_i\mid \rho<\delta\r$ is $\subseteq^*$-increasing hence by corollary \ref{subsets star bound}, there is a Mathias $E^{(\delta)}_i\subseteq\theta_i$ which is a $\subseteq^*$-bound, this ensures $(2),(5)$. Apply $(20)_{\theta_i}$ to $E^{(\delta)}_i$ and $D_i$ to obtain $F^{(\delta)}_i$ to ensure $(3)$ and finally apply $(18)_{\nu}$ to the sequence $\l F^{(\delta)}_i\mid i<\theta\r$, obtain the sequence $\l D^{(\delta)}_i\mid i<\theta\r$ which satisfy $(1)-(5)$.

Hence the sequence $\langle D^{(\xi)}_j\mid j<\theta\rangle$ is defined for every $\xi<\nu^+$. For every column $j<\theta$, $\langle D^{(\xi)}_j\mid \xi<\nu^+\rangle$ is a $\subseteq^*$-increasing sequence of subsets of $\theta_j$, thus there is $\xi_j<\nu^+$ from which this sequence stabilizes. Let $\xi^*=\sup(\xi_j\mid j<\theta)<\nu^+$.

Let us prove that $D^{(\xi^*)}_i$ is as wanted. By the construction of the sequence $(1),(3),(4),(5)$ of the theorem follows directly. To see $(2)$, for every $\xi^*\leq\xi'<\nu^+$ and for every $i<\theta$, $D^{(\xi^*)}_i=^*D^{\xi'}_i$. In particular $D^{\xi^*+1}_i=^*D^{\xi^*}_i$. Hence
$$D^{\xi^*}_i\subseteq D^{(\xi^*)}\cap\theta_i\subseteq^* D^{\xi^*+1}_i=^* D^{\xi^*}_i.$$
Hence $D^{\xi^*}_i=^* D^{(\xi^*)}\cap\theta_i\in V[D^{(\xi^*)}]$.$\blacksquare_{\text{Implication 2}}$

\textit{Proof of Implication 3 (Proposition \ref{boundedUnionOfgenerics}).} Assume the induction hypothesis, $(18)_\nu$ and $(19)_\nu$. Let us derive $(20)_\nu$. Let $cf^{V[A]}(\nu)=\lambda<\nu$ and fix a cofinal sequence $\l \nu_i\mid i<\lambda\r\in V[A]$. 
For each $i<\lambda$, apply $(20)_{\nu_i}$ to find 
$$D\cap \nu_i,D'\cap\nu_i\subseteq E_i\subseteq \nu_i$$
such that $D\cap\nu_i,D'\cap\nu_i\in V[E_i]$ and $E_i\cap\kappa^*=F_{\kappa^*}$.
Apply $(19)_\nu$ to the sequence $\l E_i\mid i<\lambda\r$ to find a sequence $\l E^*_i\mid i<\lambda\r$, such that $E_i\subseteq^* E_i^*$, $E_i\in V[E^*]$, where $E^*:=\cup_{i<\lambda}E^*_i$ is a Mathias set. Then
$|D\cup D'\setminus E^*|\leq\lambda$. As in the proof of \ref{lemmaforsubsetkappa}, in the model $V[E^*]$ we have $$\forall i<\lambda.  D\cap\nu_i,D'\cap\nu_i\in V[E^*]$$ so the sequences $\l D\cap \nu_i\mid i<\lambda\r,\l D'\cap\nu_i\mid i<\lambda\r$ can be coded as a single sequence of ordinals $\l \delta_i\mid i<\lambda\r$ (fixing enumerations of $P^{V[E^*]}(\nu_i)$). By \ref{very short}, there is a Mathias set $R\in V[A]$, $|R|\leq\lambda$ such that $V[R]=V[\l \delta_i\mid i<\lambda\r]$. Apply \ref{countableinformation} to $D\cup D'\setminus E^*, R$ and $E^*$ to find $G\in V[A]$ Mathias such that 
$D\cup D'\setminus \lambda, E^*\setminus\lambda\subseteq G $, $G\cap\lambda=F_\lambda$ and $E^*, R\in V[G]$. Let $G_0=F_{\kappa^*}\cup G\cap(\kappa^*,\lambda)$, recall that $\lambda<\nu$, hence we can apply $(20)_\lambda$ to $G_0,(D\cup D')\cap \lambda$ and find $G_1\subseteq \lambda$ such that $(D\cup D')\cap\lambda,G_0\subseteq G_1$, $G_1\cap\kappa^*=F_{\kappa^*}$ and $G_0\in V[G_1]$. 
Finally let $$D^*=F_{\kappa^*}\cup(G_1\cap(\kappa^*,\lambda))\cup(G\setminus\lambda).$$

Clearly, $D^*$ is a Mathias set, $D^*\cap\kappa^*=F_{\kappa^*}$ thus $(1),(4)$ of Theorem \ref{boundedUnionOfgenerics} hold. For $(2)$, $\sup(D^*)=\sup(G)=\sup(D\cup D')$
 $$(D\cup D')\cap\kappa^*=F_{\kappa^*}\subseteq D^*, \ D\cup D'\cap(\kappa^*,\lambda)\subseteq G_1\cap(\kappa^*,\lambda)\subseteq D^*$$ and $$D\cup D'\setminus \lambda\subseteq G\setminus\lambda\subseteq D^*.$$ Hence $D\cup D'\subseteq D^*$.
 Finally to see $(3)$,  
$$D^*\cap \kappa^*=F_{\kappa^*},\ D^*\cap(\kappa^*,\lambda)=G_1\cap(\kappa^*,\lambda), \ D^*\setminus\lambda=G\setminus\lambda.$$
Hence $F_{\kappa^*},G_1\cap (\kappa^*,\lambda),G\setminus\lambda\in V[D^*]$, so $G_1\cap \kappa^*\in V[A]\cap V[C_G\cap\kappa^*]=V[F_{\kappa^*}]\subseteq V[D^*]$, so $G_1\in V[D^*]$. It follows that $G_0\in V[G_1]\subseteq V[D^*]$. By definition of $G_0$, $G\cap(\kappa^*,\lambda)=G_0\setminus\kappa^*\in V[D^*]$, and clearly $G\cap\kappa^*\in V[F_{\kappa^*}]\subseteq V[D^*]$. Therefore $$G\cap\kappa^*,G\cap(\kappa^*,\lambda),G\setminus\lambda\in V[D^*]\text{ and }G\in V[D^*]$$ By definition of $G$, $E^*,R\in V[G]\subseteq V[D^*]$, hence $\l\delta_i\mid i<\lambda\r$ and the coding of $P^{V[E^*]}(\nu_i)$ is in $V[D^*]$ so the sequences $\l D\cap\nu_i\mid i<\lambda\r,\l D'\cap\nu_i\mid i<\lambda\r\in V[D^*]$. Therefore, $D,D'\in V[D^*]$, as wanted.$\blacksquare_{\text{ Implication 3}}$

This concludes the induction for \ref{Making a sequence stae increasing}-\ref{boundedUnionOfgenerics} for every $\nu\in X_A$ and in particular for $\kappa$. Let us conclude the main result for subsets of $\kappa$ which do not stabilize:
\begin{corollary}\label{resaultsubsetkappa}
Assume that $o^{\vec{U}}(\kappa)<\kappa^+$, $(IH)$, $A\subseteq\kappa$, $A\in V[G]$ and $A$ does not stabilize, then there is $C'\subseteq C_G$ such that $V[A]=V[C']$.
\end{corollary}
\pr
By \ref{nonstabcofinalitysubsetofkappa}, $\lambda:=cf^{V[A]}(\kappa)<\kappa$ and let $\l\beta_i\mid i<\lambda\r\in V[A]$ be cofinal. By \ref{very short} there is a sequence of Mathias sets $\langle D'_i\mid i<\lambda\rangle\in V[A]$ such that $V[D'_i]=V[A\cap\beta_i]$ and $D'_i\subseteq\beta_i$ and denote $D_i=D'_i\setminus\kappa^*\cup F_{\kappa^*}$. Then the sequence $\langle D_i\mid i<\lambda\r\in V[A]$ and $A\cap\beta_i\in V[D_i]$.  Use \ref{maintheoremnonstab} to find $\langle D^*_i\mid i<\lambda\rangle$
and set $D^*=\underset{i<\lambda}{\cup}D^*_i$. Then $D^*$ is Mathias and therefore $D^*\subseteq^* C_G$. Let $C^*=C_G\cap D^*$. Hence $C^*=^*D^*$ and $V[C^*]=V[D^*]$.
Finally, for every $\alpha<\kappa$, find $i<\lambda$ such that $\alpha<\beta_i$. By the properties of $D^*$, $D_i\in V[D^*]$, hence, $A\cap\beta_i\in V[D^*]$. Note that $A\cap\alpha=(A\cap\beta_i)\cap\alpha$ and therefore $A\cap\alpha\in V[D^*]=V[C^*]$. Finally, apply \ref{lemmaforsubsetkappa}.$\blacksquare$
\subsection{subsets of $\kappa$ which stabilize}
In this section assume that $o^{\vec{U}}(\kappa)<\kappa^+$, $(IH)$ hence by \ref{VGenericCardinals}, $\zeta_0:=cf^{V[G]}(\kappa)<\kappa$. Let $A\in V[G]$ be a subset of $\kappa$ such that $A$ stabilizes i.e. there is $\lambda<\kappa$ such that $$\forall\alpha<\kappa \ A\cap\alpha\in V[C_G\cap\lambda].$$
Note that if $A\in V[C_G\cap\beta]$ for some $\beta<\kappa$ then we can use $(IH)$, so we also assume that $A$ is fresh with respect to the model $V[C_G\cap\lambda]$. Again we would like to apply Lemma \ref{lemmaforsubsetkappa}, we will use freshness and work a little bit to prove $cf^{V[A]}(\kappa)<\kappa$, while finding $C^*$ is easy:
 
 Increase $\lambda$ if necessary, and assume $\max\{\kappa^*,\zeta_0\}\leq\lambda<\kappa$. By Proposition \ref{CodingBoundedInformation}, find $F_\lambda\subseteq \lambda$ a Mathias set such that $V[F_\lambda]=V[A]\cap V[C_G\cap\lambda]$. Define $C^*=F_\lambda\cap C_G=^* F_\lambda$, then $C^*\in V[A]$ and $$\forall\alpha<\kappa. A\cap\alpha\in V[A]\cap V[C_G\cap \lambda]=V[F_\lambda]=V[C^*].$$
It remains to see that:
\begin{proposition}
 $cf^{V[A]}(\kappa)<\kappa$.
\end{proposition}
\pr By \ref{definition of equivalent subalgebra}, let  $\mathbb{R}\subseteq RO(\Mfor\restriction\lambda)$ for which $V[C^*]=V[H_{C^*}]$ for some $V$-generic filter $H_{C^*}\subseteq \mathbb{R}$ and denote the quotient forcing (definition \ref{definition of quotient}) by $\mathbb{Q}:=(\mathbb{M}[\vec{U}]\restriction\lambda)/H_{C^*}$.
To complete $V[C^*]$ to $V[G]$, it remains to force above $V[C^*]$ with $\mathbb{P}:=\mathbb{Q}\times\mathbb{M}[\vec{U}]\restriction(\lambda,\kappa]$, let $H_{\mathbb{Q}}\times G\restriction(\lambda,\kappa]\subseteq \mathbb{P}$ be $V[C^*]$-generic such that $V[C^*][H_{\mathbb{Q}}\times G\restriction(\lambda,\kappa]]=V[G]$.  Notice that for every $\lambda\leq\alpha<\kappa$ with $o^{\vec{U}}(\alpha)>0$ we have $$|\mathbb{Q}\times\mathbb{M}[\vec{U}]\restriction(\lambda,\alpha]|<\min\{\nu>\alpha\mid o^{\vec{U}}(\nu)=1\}.$$
Let $\lusim{A}$ be a $\mathbb{P}$-name for $A$ and assume that $$\Vdash_{\mathbb{P}}  \lusim{A}\text{ is fresh}.$$
Let $\l c_i\mid i<\zeta_0\rangle\in V[G]$ be a cofinal continuous subsequence of $C_G$ such that $c_0>\lambda$. Fix  $\langle\lusim{c}_i\mid i<\zeta_0\r\in V[C^*]$ a sequence of $\mathbb{P}$-names for $\l c_i\mid i<\zeta_0\rangle$.
Find $p=\l p_0,p_1\r\in H_{\mathbb{Q}}\times G\restriction (\lambda,\kappa)$
such that $$p\Vdash_{\mathbb{P}} \langle \lusim{c}_i\mid i<\zeta_0\r\text{ is a cofinal continuous subsequence of }\lusim{C_G}.$$
For every $i<\zeta_0$ and $q\in\mathbb{Q}/p_0$, consider the set $D_{i,q}$ of all conditions $p_1\leq r\in \Mfor\restriction (\lambda,\kappa)$ such that one of the following holds:
\begin{enumerate}
    \item $ \exists\alpha.\ \l q,r\r \Vdash_{\mathbb{P}}\lusim{c}_i=\alpha\ \wedge\ \exists B. \  \l q,r\r\Vdash_{\mathbb{P}}   \lusim{A}\cap\alpha=B$. Denote this statement by $\phi_i(q,r)$.
    \item For every $r'\geq r$, $\neg \phi_i(q,r')$.
\end{enumerate} 
Then $D_{i,q}$ is clearly dense open. By the strong Prikry property there is $p_1\leq^*p_{i,q}$, $S_{i,q}$ and sets $A^s_{i,q}$ such that for every $t\in mb(S_{i,q})$, $p_{i,q}^{\smallfrown}\l t,\vec{A}^t_{i,q}\r\in D_{i,q}$. Define $$g_{i,q}:mb(S_{i,q})\rightarrow\{0,1\}\text{ by }g_{i,q}(t)=1\leftrightarrow \phi_i(q,p_{i,q}^{\smallfrown}\l t,\vec{A}^t_{i,q}\r)\text{ holds}.$$ Then we can shrink $S_{i,q}$ to $T_{i,q}$ such that $g_{i,q}$ is constant on $mb(T_{i,q})$. 
Now for every $q\in\mathbb{Q}$ such that $g_{i,q}=1$, and every $s\in mb(T_{i,q})$ let $\alpha_i(q,s),A_i(q,s)$ be the values decided by $\l q, p^{\smallfrown}\l s,\vec{B}^s_{i,q}\r\r$ for $\lusim{c}_i,\lusim{A}\cap\lusim{c}_i$ respectively. Let $N_{i,q}=ht(T_{i,q})$, then $$\alpha_i(q,s)\in \{\kappa_1(p),...,\kappa_{l(p)}(p),s(1),...,s(N_{i,q})\}$$ we can extend $T_{i,q}$ if necessary so that $\max(s)\geq\alpha_i(q,s)$. In particular, $A_i(q,s)\subseteq\max(s)$.

Define by recursion $A_i(q,s)$ for $s\in T_{i,q}\setminus mb(T_{q,i})$. Let $s\in \Lev_{N_{i,q}-1}(T_{i,q})$, by ineffability, we can shrink $\succ_{T_{q,i}}(s)$ and find $A_i(q,s)$ such that for every $\alpha\in\succ_{T_{q,i}}(s)$, $A_i(q,s^{\smallfrown}\alpha)= A_i(q,s)\cap\alpha$.
Generally, take $s\in T_{i,q}$ and assume that for every $\alpha$ in $\succ_{T_{i,q}}(s)$, $A_i(q,s^{\smallfrown}\alpha)$ is defined. 
We can find a single $A_i(q,s)$ and shrink $\succ_{T_{i,q}}(s)$ such that $$\forall\alpha\in\succ_{T_{i,q}}(s). \ A_i(q,s^{\smallfrown}\alpha)\cap\alpha= A_i(q,s)\cap\alpha.$$
We abuse notation by denoting the shrinked trees by $T_{i,q}$. Extend $p_{i,q}\leq^* p^*_{i,q}$ find $B^t_{i,q}\subseteq A^t_{i,q}$ such that extensions from $D_{T_{i,q},\vec{B}_{i,q}}$ are pre-dense above $p^*_{i,q}$  and use $\leq^*$ closure of $\Mfor\restriction(\lambda,\kappa]$ to find a single $p^*$ such that for every $q\in\mathbb{Q}/p_0$ and $i<\zeta_0$, $p^*_{i,q}\leq^* p^*$, in particular $p_1\leq p^*$. As usual, shrink all the trees to $p^*$ and let $T_{i,q}$ be the resulting tree.
\begin{claim*}
For every $i<\zeta_0$ and $q\in\mathbb{Q}/p_0$ there is $q'\geq q$ such that $g_{i,q'}\restriction mb(T_{i,q})\equiv 1$. i.e. 
$$\forall t\in mb(T_{i,q}). \ \exists \alpha,B. \ \l q',p^{*\smallfrown}\l t,\vec{B}^t_{i,q'}\r\r\Vdash_{\mathbb{P}} \lusim{c}_i=\alpha \wedge \lusim{A}\cap\alpha=B.$$
\end{claim*}
\pr Let $p_0\leq q_0$, find some $\l q_0,p^*\r\leq \l q,r\r$  and $\alpha$ such that $$\l q,r\r\Vdash_{\mathbb{P}} \lusim{c}_i=\alpha.$$ By assumption on $\lusim{A}$, $\l q,r\r\Vdash_{\mathbb{P}} \lusim{A}\text{ is fresh}$, which implies that there is some $B\in V[C^*]$ and some $\l q',r'\r$ such that $\l q,r\r\leq\l q',r'\r\Vdash B=\lusim{A}\cap\alpha$. Find some $t\in T_{i,q'}$ such that $p^{*\smallfrown}\l t,\vec{B}^t_{i,q'}\r$ and $r'$ are compatible, then a common extension witnesses that  $g_{i,q'}(t)\neq 0$, hence $g_{i,q'}(t)=1$ as wanted.$\blacksquare_{\text{Claim}}$ %$\lusim{c}_i$, $\lambda<\alpha$ and it follows that there is $j\leq l(r)$ such that $\alpha=\kappa_j(r)$. otherwise, we can cut the large sets of $r$ and directly extend $r\leq^* r'$ so that $$\l q,r'\r\Vdash_{\mathbb{P}} \alpha\notin C_{\lusim{G}}\wedge \alpha=\lusim{c}_i\in C_{\lusim{G}}$$ This is a contradiction, since no condition forces contradictory information. Hence we can split $r$ at $\alpha$, then the upper  is sufficiently closed to decide $\lusim{A}\cap\alpha$. Formally, Let $$H_{\alpha}\subseteq\mathbb{Q}\times \Mfor\restriction(\lambda,\alpha)$$  be $V[C']$-generic filter such that $\l q,r\restriction\alpha\r\in H_{\alpha}$. Denote by $(\lusim{A})_{H_{\alpha$ the $\mathbb{P}\restriction(\alpha,\kappa)$-name in $V[C'][H_{\alpha}]$ derived from $\lusim{A}$. 

Move to $V[A]$, let us compare the sets $A_i(q,s)$ with $A$. 

For every $i$ and $q$ such that $g_{i,q}=1$, define $\rho_q^i(k)$ for $k\leq N_{q,i}$.
Let $$\rho^i_q(0)=\min(A\Delta A_i(q,\l\r))+1.$$ Recursively define
$$\rho^i_q(k+1)=\sup(\min(A\Delta A_i(q,\l\delta_1,...,\delta_k\r))+1\mid \l\delta_1,...\delta_k\r\in \Lev_k(T_{i,q})\cap\prod_{j=1}^k\rho^i_q(j)).$$
Finally we let $$\rho^i(k)=\sup\{\rho^i_q(k)\mid q\in\mathbb{Q}\wedge g_{i,q}=1\}.$$
By the claim, for each $i<\zeta_0$, there is $q_i\in H_{\mathbb{Q}}$ such that $g_{i,q_i}=1$, and since $D_{T_{i,q_i},\vec{B}_{q_i,i}}$ is pre-dense, there is some $\vec{c}_i\in mb(T_{i,q})$ such that $\l q_i,p^{*\smallfrown}\l\vec{c}_i,\vec{B}^{\vec{c}_i}_{q_i,i}\r\r\in H_{\mathbb{Q}}\times G\restriction(\lambda,\kappa]$. By assumption on $T_{i,q_i}$, $\max(\vec{c}_i)\geq (\lusim{c}_i)_{H_{\mathbb{Q}}\times G\restriction(\lambda,\kappa]}$ let us argue that for every $k\leq N_i$, $\rho^i(k)> \min\{c_i,\vec{c}_i(k)\}$.

By construction of the tree $T_{i,q_i}$, $A\cap c_i=A_i(q_i,\vec{c}_i)\cap c_i$. 
Since for every $j\leq N_i$, by definition,
$$A_i(q_i,\vec{c}_i\restriction j)\cap \vec{c}_i(j)=A_i(q_i,\vec{c}_i\restriction j+1)\cap\vec{c}_i(j).$$
It follows that for every $j\leq N_i$,
$$A_i(q_i,\vec{c}_i\restriction j)\cap \min\{c_i,\vec{c}_i(j)\}=A\cap\min\{c_i,\vec{c}_i(j)\}.$$
 In particular, 
$A\cap\min\{c_i,\vec{c}_i(0)\}=A_i(q_i,\l\r)\cap \min\{c_i,\vec{c}_i(0)\}$. 

Since $A\cap\rho^i(0)\neq A_i(q_i,\l\r)\cap\rho^i(0)$, it follows that $\min\{c_i,\vec{c}_i(0)\}<\rho^i(0)$. 

Inductively assume that $\min\{c_i,\vec{c}_i(j)\}<\rho^i(j)$ for every $j\leq k$. If $c_i\leq c_i(k)$ then clearly we are done. Otherwise, $\rho_i(j)>\vec{c}_i(j)$, which implies that $$\vec{c}_i\restriction\{1,...,k\}\in \Lev_k(T_{i,q_i})\cap\prod_{j=1}^k\rho^i(j)$$ and since $A_i(q_i,\vec{c}_i\restriction \{1,...,k\})\cap\min\{c_i,\vec{c}_i(k+1)\}=A\cap\min\{c_i, \vec{c}_i(k+1)\}$,
then $$\min\{c_i,\vec{c}_i(k+1)\}<\min(A^i(q_i,\vec{c}_i\restriction \{1,...,k\})\Delta A)\leq \rho^i(k+1).$$
Since $\vec{c}_i(N_i)\geq c_i$, it follows that $\rho^i(N_i)>c_i$.

Next we argue that $\rho^i(k)<\kappa$.
Again by induction on $k$,
$\rho^i(q,0)<\kappa$ since for every $q\in\mathbb{Q}$ with $g_{i,q}=1$, $A\neq A_i(q,\l\r)$,  as $A_i(q,\l\r)\in V[C^*]$ but $A\notin V[C^*]$. Since $|\mathbb{Q}|<\kappa$ and $\kappa$ is regular in $V[C^*]$, it follows that $\rho^i(0)<\kappa$.

Assume that it holds for every $j\leq k$. Toward a contradiction assume that $\rho^i(k+1)=\kappa$. Again, $|\mathbb{Q}|<\kappa$ and $\kappa$ is regular in $V[C^*]$, there must be $q\in\mathbb{Q}$ such that $g_{i,q}=1$ and $\rho^i(q,k+1)=\kappa$. Consider the collection $$\{A_i(q,\l\alpha_1,...,\alpha_k\r)\mid\l\alpha_1,...,\alpha_k\r\in \Lev_k(T_{i,q})\cap\prod_{j=1}^k \rho^i(j)\}\in V[C^*].$$
Then for every $\gamma<\kappa$ pick any distinct $\vec{\alpha}_1,\vec{\alpha}_2\in \Lev_k(T_{i,q})\cap\prod_{j=1}^k \rho^i(j)$ such that $A_i(q,\vec{\alpha}_1)\neq A_i(q,\vec{\alpha}_2)$, but $A_i(q,\vec{\alpha}_1)\cap\gamma=A_i(q,\vec{\alpha}_2)\cap\gamma$.

To see that there are such $\vec{\alpha}_1,\vec{\alpha}_2$, by assumption that $\rho^i(k+1)=\kappa$ there is $\vec{\alpha}_1$ such that $\eta_1:=\min(A\Delta A_i(q,\vec{\alpha}_1))>\gamma$,
 hence $A_i(\vec{\alpha}_1)\cap\gamma=A\cap\gamma$. Let $\vec{\alpha}_2$ be such that $\min(A\Delta A_i(q,\vec{\alpha}_2))>\eta_1$.
 In particular, $A_i(q,\vec{\alpha}_1)\neq A_i(q,\vec{\alpha}_2)$, but  $A_i(q,\vec{\alpha}_1)\cap\gamma=A\cap\gamma=A_i(q,\vec{\alpha}_2)\cap\gamma$.
Since this is defined in $V[C^*]$, where $\kappa$ is still measurable, and the number of pairs $\l\vec{\alpha}_1,\vec{\alpha}_2\r$ is bounded by the induction hypothesis, we can find unboundedly many $\gamma$'s with the same $\vec{\alpha}_1,\vec{\alpha}_2$, which is clearly a contradiction.

So we found a sequence  $\langle\rho^i(N_i)\mid i<\lambda\rangle\in V[A]$ such that $\rho_i(N_i)>c_i$. Hence $cf^{V[A]}(\kappa)\leq\lambda$.$\blacksquare$

As a result of this section we obtain the following:
\begin{corollary}\label{resultforstabsubset}
 Assume $o^{\vec{U}}(\kappa)<\kappa^+$ and $(IH)$. Let $A\in V[G]$, $A\subseteq\kappa$ be such that $A$ stabilizes, then there is $C'\subseteq C_G$ such that $V[A]=V[C']$.
\end{corollary}

\section{The argument for a general set}

Recall the main theorem of this paper is:
\\
\\
 \textbf{Theorem 1.1} \textit{Let $\vec{U}$ be a coherent sequence with maximal measurable $\kappa$, such that $o^{\vec{U}}(\kappa)<\kappa^+$. Assume the inductive hypothesis:}
 $$(IH) \ \ \  For \ every \ \delta<\kappa, \ any \ coherent \ sequence \ \vec{W}\ with\ maximal\ measurable \ \delta\ and \ any\ set \ of \ ordinals$$
 $$A\in V[H]\ for \ H\subseteq\mathbb{M}[\vec{W}], \ there\ is\  \ C\subseteq C_H,\ such\ that\ V[A]=V[C].$$
\textit{ Then for every $V$-generic filter $G\subseteq\Mfor$ and any set of ordinals $A\in V[G]$, there is $C\subseteq C_G$ such that $V[A]=V[C]$.}
 
\vskip 0.2 cm

\begin{remark}
The authors would like to thanks Gunther Fuchs for pointing out that \ref{MainResaultPartwo} does not automatically generalize to every set in $V[G]$. For example, if $V=L$ and $\l c_n\mid n<\omega\r$ are $\omega$ Cohen reals over $L$ (which is equivalent to adding a single real) then certainly $A:=\{c_n\mid n<\omega\}\in V[\l c_n \mid n<\omega\r]$, but the minimal model containing both $L$ and $A$ is $L(A)$ which is a model of $ZF$ rather than $ZFC$.
This situation can also occur in Prikry-type, namely, there are intermediate models of $ZF+\neg AC$ which are intermediate to a Prikry forcing extension. Suppose that $\l c_n\mid n<\omega\r$ is a Prikry sequence over $V$, split $\omega$ to $\omega$-many infinite disjoint sets $\l T_n\mid n<\omega\r$ and let $D_n=\{c_m\mid m\in T_n\}$. Now consider $R_n=\{t\mid t\text{ is a finite change of }D_n\}$, then clearly $\{ R_n\mid n<\omega\}\in V[G]$. Let $\mathcal{G}$ be the group of all permutations on Prikry forcing permuting the $R_n$'s, generated by a permutation of $\omega$. Let $\mathcal{F}$ be the filter generated by the sets $\text{fix}(E):=\{\pi\in \mathcal{G}\mid \forall x\in E\pi(x)=x\}$ for finite sets $E\subseteq \omega$.  Consider the symmetric submodel $N$ of $V[G]$ (see \cite[Chapter 15]{Jech2003}), then $\{R_n\mid n<\omega\}\in N$, $V\subseteq N$ is a model of $ZF$ which fails to satisfy the axiom of choice.
\end{remark}
Let $A$ be a set of ordinals. We prove Theorem \ref{MainResaultParttwo} by induction of $\lambda:=\sup(A)$. If $\lambda\leq\kappa$ then apply \ref{very short}, \ref{resaultsubsetkappa}, \ref{resultforstabsubset}.  Assume that $\lambda>\kappa$, and let us first resolve the induction step for $cf^{V[G]}(\lambda)\leq\kappa$:
\begin{proposition}
Assume $o^{\vec{U}}(\kappa)<\kappa^+$, $(IH)$, and $cf^{V[G]}(\lambda)\leq\kappa$, then there is $C\subseteq C_G$ such that $V[A]=V[C]$.
\end{proposition} 
\pr Since $\kappa$ is singular in $V[G]$ then $cf^{V[G]}(\lambda)<\kappa$. Since $\Mfor$ satisfies $\kappa^+-cc$ we must have that   $\nu:=cf^V(\lambda)\leq\kappa$. Fix $\langle\gamma_i | \ i<\nu\rangle\in V$ cofinal in $\lambda$. Work in $V[A]$, for every $i<\nu$ find $d_i\subseteq \kappa$ such that $V[d_i]=V[A\cap\gamma_i]$. By induction, there exists $C^*\subseteq C_G$ such that $V[\langle d_i\mid i<\nu\rangle]=V[C^*]$ (Note that the $d_i$'s are subsets of $\kappa$ so we can code $\l d_i\mid i<\kappa\r$ as a sequence of pairs in $V$) so , therefore 
\begin{enumerate}
\item $\forall i<\nu \ A\cap\gamma_i\in V[C^*]$
\item $C^*\in V[A]$
\end{enumerate}
Work in $V[C^*]$, for $i<\nu$ fix a bijection $\pi_i:2^{\gamma_i}\leftrightarrow P^{V[C^*]}(\gamma_i)$. Find $\delta_i$ such that $\pi_i(\delta_i)=A\cap\gamma_i$.

By $\kappa^+$-cc of $\Mfor$, there if a function $F:\nu\rightarrow P(\lambda)$ in $V$ such that for every $i<\nu$, $\delta_i\in F(i)$ and $|F(i)|\leq\kappa$. Let $\epsilon_i<\kappa$ be the index of $\delta_i$ inside $F(i)$. Find $C''\subseteq C_G$ such that $V[C'']=V[\langle \epsilon_i\mid i<\nu\rangle]$. Finally we can find $C'\subseteq C_G$ such that $V[C']=V[C^*,C'']$. To see that
$V[A]=V[C']$, clearly, $C^*\in V[A]$ and therefore $\l \pi_i,\delta_i\mid i<\nu\r\in V[A]$. Since $F\in V$, $\l\epsilon_i\mid i<\nu\r\in V[A]$, hence $C''\in V[A]$. It follows that $C'\in V[A]$. For the other direction, $C^*,C''\in V[C']$, so $\l\epsilon_i\mid i<\nu\r\in V[C']$, and since $F\in V$ then $\l\delta_i\mid i<\nu\r\in V[C']$. Since $C^*\in V[C']$ then also $\l\pi_i\mid i<\nu\r\in V[C']$ so $\l\pi_i(\delta_i)\mid i<\nu\r\in V[C']$. It follows that $A=\cup_{i<\nu}\pi_i(\delta_i)\in V[C']$.
$\blacksquare$

\subsection{The Induction step for $cf(\lambda)>\kappa$}

Assume that $o^{\vec{U}}(\kappa)<\kappa^+$ and $(IH)$. The idea for the induction step where $\sup(A)=\lambda$ with $cf(\lambda)>\kappa$ (typical example is $\lambda=\kappa^+$) is the following:
\begin{enumerate}
    \item There is $C^*\subseteq C_G$ such that $C^*\in V[A]$ and for every $\alpha<\lambda. A\cap\alpha\in V[C^*]$.
    \item The quotient forcing $\Mfor/C^*$ (which completes $V[C^*]$ to $V[G]$) is $\kappa^+$-cc (and therefore $cf(\lambda)$-cc) in $V[G]$.
\end{enumerate}
Then we will apply the following  theorem:
\begin{theorem}\label{no fresh subsets}
Let $W\models ZFC$ and $\mathbb{P}\in W$ be a forcing notion. Let $T\subseteq\mathbb{P}$ be any $W$-generic filter and $\theta$ a regular cardinal in $W[T]$. Assume $\mathbb{P}$ is $\theta$-cc in $W[T]$. Then in $W[T]$ there are no fresh subsets with respect to $W$ of cardinals $\lambda$ such that $cf(\lambda)=\theta$.
\end{theorem}
\begin{remark}
Note that it is crucial that $\mathbb{P}$ is $\theta$-cc in the generic extension, otherwise there are trivial examples which contradict this. Namely, the forcing which adds a branch through a Suslin tree is $ccc$, but the branch added is a fresh subset of $\omega_1$.  
\end{remark}
\pr
 Toward a contradiction, assume that $A\in W[T]\setminus W$ is a fresh subset of $\lambda$ and $cf(\lambda)=\theta$. Pick a name $\lusim{A}$ for $A$ and work within $W[T]$. We define recursively a sequence of conditions $\l r_i,s_i\mid i< \theta\r$ and a sequence of ordinals $\l \beta_i\mid i<\theta\r$.  Let $r_0\in T$ be such that $r_0\Vdash \lusim{A}\text{ is fresh}$. There must be $\beta_0<\lambda$ such that $r_0$ does not force $\lusim{A}\cap\beta_0=A\cap\beta_0$. Otherwise, $A=\cup\{B\mid \exists\beta<\lambda. r_0\Vdash \lusim{A}\cap\beta=B\}\in W$, contradicting the fact that $A\notin W$. Hence one can find $B_0\neq A\cap\beta_0$ and $r_0\leq s_0$ such that $s_0\Vdash \lusim{A}\cap\beta_0=B_0$.
 
 Assume $r_i,s_i,\beta_i$ are defined for every $i<j<\theta$. Let $\beta'_j:=\sup\{\beta_i\mid i<j\}<\lambda$, find $r_j\in T$ such that $r_0\leq r_j\Vdash \lusim{A}\cap\beta'_j=A\cap\beta'_j$. Such $r_j$ exists since $A$ is fresh. Argue as before to find $\beta_j<\lambda$, $B_j\neq A\cap\beta_j$ and $s_j\geq r_j$ such that $s_j\Vdash \lusim{A}\cap\beta_j=B_j$. The contradiction is obtained by noticing that $\l s_j\mid j<\theta\r$ is an antichain. Indeed, if $i<j$ and $s_i,s_j\leq s$ then $s_i\leq s$ implies that $s\Vdash\lusim{A}\cap \beta_i=B_i$, also since $r_j\leq s_j\leq s$, then $s\Vdash \lusim{A}\cap \beta_i=A\cap \beta_i\neq B_i$. Hence $s$ forces contradictory information.$\blacksquare$

\begin{corollary}
Assume $(IH)$ and $A\in V[G]$ stabilizes, then there is  $C\subseteq C_G$ such that $V[A]=V[C]$. 
\end{corollary}
\pr Let $\beta<\kappa$ be such that $\forall\alpha<\sup(A).\ A\cap\alpha\in V[G\restriction\beta]$. If $A\in V[G\restriction\beta]$, then we can apply $(IH)$ and we are done. Otherwise, $A\in V[G]$ is fresh with respect to the model $V[G\restriction\beta]$. The forcing completing $V[G\restriction\beta]$ to $V[G]$ is simply $\Mfor\restriction(\beta,\kappa]$ which clearly is $\kappa^+$-cc in $V[G]$ (since $\kappa^+$ is regular in $V[G]$), this is a contradiction to Theorem \ref{no fresh subsets}.$\blacksquare$

Assume that $A$ does not stabilize, since we assumed that $o^{\vec{U}}(\kappa)<\kappa^+$ and by the induction hypothesis on $\sup(A)=\lambda$, we can apply \ref{nonstabcofinality}(2), to conclude that $cf^{V[A]}(\kappa)<\kappa$.

Let $\l \lambda_i\mid i<cf(\lambda)\r$ be cofinal in $\lambda$, then for each $\alpha
<cf(\lambda)$ we choose some $D_\alpha\subseteq C_G$ such that $V[A\cap\lambda_\alpha]=V[D_\alpha]$. In previous results (\cite{partOne},\cite{TomMoti}), $o^{\vec{U}}(\kappa)<\kappa$ and $|C_G|<\kappa$, therefore $2^{|C_G|}<\kappa<cf(\lambda)$, it followed that there is some $D_\alpha$ that repeated cofinaly many times. Here, since $2^{|C_G|}\geq\kappa^+$, we will need as before to somehow accumulate all the information in a $\subseteq^*$-increasing way.
\begin{proposition}\label{SequenceIncreasing}
Assume $o^{\vec{U}}(\kappa)<\kappa^+$, $(IH)$ and $A\in V[G]$ does not stabilize. Let $\l \lambda_i\mid i<cf(\lambda)\r$ be cofinal in $\lambda$ and $\kappa^*<\kappa$ such that for every $\alpha\in C_G\setminus \kappa^*$, $o^{\vec{U}}(\alpha)<\alpha^+$.  Then there is a sequence $\langle D_\alpha\mid \alpha<cf(\lambda)\rangle\in V[A]$ such that:
\begin{enumerate}
    \item $D_\alpha$ is a Mathias set, $D_\alpha\cap\kappa^*=F_{\kappa^*}$, where $V[F_{\kappa^*}]=V[A]\cap V[C_G\cap\kappa^*]$.
    \item $\langle D_\alpha\mid \alpha<cf(\lambda)\rangle$  is $\subseteq^*$-increasing.
    \item $A\cap \lambda_\alpha\in V[D_\alpha]$
\end{enumerate}
\end{proposition}
\pr Work in $V[A]$. For every $\alpha<cf(\lambda)$, by the induction hypothesis, there is a Mathias set $D'_\alpha\subseteq^*C_G$ such that $A\cap\lambda_\alpha\in V[D'_\alpha]$ and $D'_\alpha\cap\kappa^*=F_{\kappa^*}$. Then $(1),(3)$ hold but $(2)$ might fail. 
Let us construct the sequence $\langle D_\alpha\mid \alpha<cf(\lambda)\rangle$ more carefully to ensure condition $(2)$: We go by induction on $\beta<\kappa^+$.
Assume the sequence $\langle D_\alpha\mid \alpha<\beta\rangle$ is defined. If $\beta=\alpha+1$, then use Proposition \ref{boundedUnionOfgenerics} with $D_\alpha$ and $D'_\beta$ to find $D_{\beta+1}$ such that $D_\alpha\subseteq D_\beta$, $D'_\beta\in V[D_\beta]$ and $D_\beta\cap\kappa^*=F_{\kappa^*}$.
If $\beta$ is limit, let $\delta=cf^{V[A]}(\beta)$ and $\l \beta_i\mid i<\delta\r\in V[A]$ be cofinal.
If $\delta>\kappa$, then by \ref{Modfinitestab}, the sequence $\l D_{\beta_\alpha}\mid \alpha<\delta\r$, $=^*$-stabilizes on some Mathias set $E^*_\beta$, in particuar, $E^*_\beta\cap\kappa^*=F_{\kappa^*}$ and since the sequence $\l D_\alpha\mid \alpha<\beta\r$ is $\subseteq^*$-increasing, then it also stabilizes on $E^*_\beta$. Then $E^*_{\beta}$ is a $\subseteq^*$-bound.

If $\delta\leq\kappa$,
since $\kappa$ is singular in $V[A]$, then $\delta<\kappa$. Apply Lemma  \ref{subsets star bound} to the sequence $\langle D_{\beta_\alpha}\mid \alpha<\delta\rangle$, to find a single $E^*_{\beta}\in V[A]$ Mathias which is a $\subseteq^*$-bound and $E^*_\beta\cap\kappa^*=F_{\kappa^*}$. 

In any case, apply Lemma \ref{boundedUnionOfgenerics} to $E^*_\beta,D'_\beta$ and find a Mathias $D_\beta$ such that $E^*_{\beta}\subseteq D_\beta$ and $D'_\beta\in V[D_\beta]$ and $D_\beta\cap\kappa^*=F_{\kappa^*}$. Clearly the sequence $\langle D_\alpha\mid \alpha<cf(\lambda)\rangle$ is as wanted.$\blacksquare$

\begin{corollary}\label{Candiadate}
There is $C^*\subseteq C_G$, such that $C^*\in V[A]$ and for every $\alpha<\lambda$, $A\cap\alpha\in V[C^*]$.
\end{corollary}
\pr Consider the sequence $\langle D_\alpha\mid \alpha<cf(\lambda)\rangle\in V[A]$ from Proposition \ref{SequenceIncreasing}, then use Theorem \ref{Modfinitestab} to find $\alpha^*<cf(\lambda)$ such that for every $\alpha^*\leq \beta<cf(\lambda)$, $D_\beta=^*D_{\alpha^*}$. In particular, $V[D_{\beta}]=V[D_{\alpha^*}]$. Then $C^*=D_{\alpha^*}\cap C_G$ is as wanted.
$\blacksquare$

Let us turn to the proof that the quotient forcing is $\kappa^+$-cc in $V[G]$ (and therefore $cf(\lambda)$-cc). In \cite{partOne} and \cite{TomMoti}, in order to prove $\kappa^+$-cc of the quotient forcing, a concrete description of the quotient was given. Here we will give an abstract argument to avoid this description.
\begin{example}\label{Example1}
 It is tempting to try and discard the name $\lusim{C}^*$ and define $\Mfor/C^*$ to consist of all $p$ such that there is a $V$-generic $H\subseteq\Mfor$, with $p\in H$ and $C^*\subseteq C_H$. Formally, we suggest that $\Mfor/C^*$ is
 $$\Mfor'=\{p\in\Mfor\mid C^*\subseteq \kappa(p)\cup B(p)\}.$$
 Such a forcing is not $\kappa^+$- cc even above $V[C^*]$. Assume that $o^{\vec{U}}(\kappa)=\kappa$, then $cf^{V[G]}(\kappa)=\omega$. We take for example any $C^*=\{c_n\mid n<\omega\}\subseteq C_G$ unbounded in $\kappa$, such that for every $n$, $o^{\vec{U}}(c_n)=0$. Basically, it is a Prikry sequence for the measure $U(\kappa,0)$. Now $V[C^*]\models \kappa^\omega=\kappa^+$ so let $\langle f_i\mid i<\kappa^+\rangle\in V[C^*]$ be an enumeration of all functions from $\omega$ to $\kappa$. we can factor the forcing to first pick $i<\kappa^+$, then the rest of the forcing ensures that $C_G(f_i(n)+1)=c_n$, this means that $f_i$ determines the places of $c_n$'s in the sequence $C_G$. Since no choice of $i\neq j$ can be compatible, the first part is not $\kappa^+$-cc and therefore also the product.
\end{example}
\begin{example}\label{Example2}
Let us consider another possible simplification of $\Mfor/C^*$,  first we enumerate $C^*=\{c^*_\alpha \mid \alpha<\kappa\}$ and find  $\Mfor-$names $\{\lusim{c}'
_\alpha\mid \alpha<\kappa\}$ for it.
$$\Mfor^*=\{q\in\Mfor\mid \text{ for every finite } a\subseteq \kappa \text{ there is } q_a\geq q,
q_a\Vdash \lusim{c}'_\alpha = \check{c}^*_\alpha, \text{ for every } \alpha\in a \}.$$

Let us prove that for suitable choice of names, $\Mfor'$ is not $\kappa^+$-cc
For every $\alpha<\kappa$, let $$X_\alpha=\{\nu<\kappa \mid o^{\vec{U}}(\nu)=\alpha\}.$$
Pick some different $\rho^0,\rho^1 \in X_0$.
The play would be between two conditions $$p^0=\l \rho^0, \l \kappa, \kappa\setminus \rho^0+1\r\r\text{ and }p^1=\l \rho^1, \l \kappa, \kappa\setminus \rho^1+1\r\r.$$
Above $p^0$ we do something simple - for example, let $\lusim{c'}_\alpha$ be a name for the first element of $X_\alpha$ in the generic sequence $C_G$.

Now above $p^1$, let us do something more sophisticated. We will build a $\kappa-$tree with each of its branches corresponding to a direct extension of $p^1$ in $\Mfor/C'$, where $C':=\lusim{C'}_H$ and $H\subseteq\Mfor$ is a $V$-generic filter with $p^0\in H$. These extensions will be incompatible in $\Mfor/C'$.
Start with a description of the first level:
\\ Fix $Y_1 \in U(\kappa,1), $ such that $Y_1 \subseteq X_1$ and $Z_1 =X_1 \setminus Y_1$ has cardinality $\kappa$.
Split $Z_1$ into two disjoint non-empty  sets $Z_{1,0},Z_{1,1}$.

Now, $\lusim{c}_1'$ be a name such that $p^1$ extended by an element of $Y_1$ forces different values from those which $p^0$ forces, for example, let it be the first element of $X_2$ in $C_G$.

For $i=0,1$, for define the name $\lusim{c}_1'$ so that $p^1$ extended by an element of $Z_{1,i}$ forces the same values as $\lusim{c}_1'$ extended by $p^0$.

The idea behind is to ensure that, $p^1{}^\frown {Z_{1 0} \cup Y_1},p^1{}^\frown {Z_{1 1} \cup Y_1}$ will be in $\Mfor/C'$, but only because of $Z_{1 i}$.
Note that, $p^1{}^\frown {Z_{1 0} \cup Y_1},p^1{}^\frown {Z_{1 1} \cup Y_1}$ are incompatible in $\Mfor/C'$ since $Z_{1 0}$ and $Z_{1 1}$ are disjoint.
Continue in a similar fashion to define the rest of the levels, at the $\alpha$-th level we take $Y_\alpha\subseteq X_\alpha$ such that $Z_\alpha:=X_\alpha\setminus Y_\alpha$ has size $\kappa$, and we split $Z_\alpha$ into two disjoint non empty sets $Z_{\alpha,0},Z_{\alpha,1}$. The definition of $\lusim{c'}_\alpha$ is such that $p^1$ extended by elements of $Y_\alpha$ forces $\lusim{c'}_\alpha$ to be the first member of $X_{\alpha+1}$ in $C_G$. While $p^1$  extended by elements of $Z_\alpha$ will force the same value as $p^0$ did.

Note that the construction is completely inside $V$.

Finally, there are  $\kappa^+-$branches of length $\kappa$ in $T$. Let $p^h$ denote an extension of $p^1$ which corresponds to a $\kappa-$branch $h$ i.e. $p^h=\l \rho_1,\l \kappa, \underset{\alpha<\kappa}{\bigcup}Y_\alpha\uplus Z_{\alpha,h(\alpha)}\rangle\rangle$.

Let $h_1, h_2$ be two different branches. Let $\alpha<\kappa$ be the least such that $h_1(\alpha)\neq h_2(\alpha)$.
Then $p^{h_1}$ and $p^{h_2}$ are incompatible in
$\Mfor/C'$.
 This follows from the choice of $\lusim{c}'_\alpha$ and the definitions of conditions at the level $\alpha$. 
 
 Note that every $p^h$ is in $\Mfor'$, since for every finite $a\subseteq\kappa$, we can extend $p^h$ to some $q_a$  using the elements from $Z_{\alpha,h(\alpha)}$.

\end{example}

\begin{proposition}\label{chracterization of quotient}
For every $q\in\Mfor$,
$$q\in \Mfor/C^*\text{ iff there is a }V\text{-generic  }G'\subseteq\Mfor\text{ such that }\lusim{C}_{G'}=C^*.$$
\end{proposition}
\pr Let $q\in \Mfor/C^*=\Mfor/H_{C^*}$, let $G'\subseteq\Mfor/C^*$ be any $V[C^*]$-generic with $q\in G'$, then $G'\subseteq\Mfor$ is a $V$-generic filter and $\pi_*(G')=\pi_*(G)=H_{C^*}$ (for the definition of $\pi_*$ see Definition  \ref{definition of quotient}). To see that $\lusim{C}_{G'}=C^*$, denote $C':=\lusim{C}_{G'}$, toward a contradiction, assume that $s\in C^*\setminus C'$, then there is $$q\leq q'\in G'\text{ such that }q'\Vdash s\notin \lusim{C^*}$$ hence $\pi(q')\leq ||s\notin\lusim{C}||$. It follows that $\pi(q')\bot ||s\in\lusim{C}||\in H_{C^*}$, therefore $\pi(q')\in \pi_*(G')\setminus H_{C^*}$ contradiction. Also if $s\in C'\setminus C^*$, then there is $q\leq q'\in G$ such that $q'\Vdash s\in \lusim{C}$, then $\pi(q')\leq ||s\in\lusim{C}||$, then $\pi(q')\bot ||s\notin\lusim{C}||\in H_{C^*}$. In any case $\pi(q')\in \pi_*(G')\setminus H_{C^*}$ which is again a contradiction. 

For the other direction, if $q\in G'$ for some $G'\subseteq\Mfor$ such that $\lusim{C}_{G'}=C^*$, then $X\cap \pi_*(G')=X\cap \pi_*(G)$, where $X=\{||\alpha\in\lusim{D}||\mid \alpha<\kappa\}$ is the generating set of $\mathbb{P}_{\lusim{C}}$. Since $\pi$ is a projection, $\pi_*(G')$ is a $V$-generic filter for $\mathbb{P}_{\lusim{C}}$ and therefor it is uniquely determined by the intersection with the set of generators $X$. It follows that $\pi_*(G')=\pi_*(G)=H_{C^*}$. Finally, for every $a\in G'$, $\pi(a)\in \pi_*(G)$, thus $a\in\pi^{-1''}H_{C^*}:=\Mfor/H_{C^*}$.$\blacksquare$

\begin{remark}
\begin{enumerate}
\item Example \ref{Example1} produces a much larger forcing than $\Mfor/C^*$ so we can obviously find $q\in \Mfor'$ such that $q\Vdash c^*_\alpha\neq\lusim{c}'_\alpha$ for some $\alpha$.
    \item  In example \ref{Example2}, the conditions $p^{h}$ constructed are not in $\Mfor/C^*$. Otherwise, by the proposition, there is a generic $H$ such that $\{(\lusim{c'}_\alpha)_H\mid \alpha<\kappa\}=C^*$ with $p^{h}\in H$. Since $Y^*:=\underset{\alpha<\kappa}{\bigcup}Y_\alpha\in\cap\vec{U}(\kappa)$, then by Proposition \ref{genericproperties}(3) there is $\xi<\kappa$ such that $C_H\setminus\xi\subseteq Y^*$. It follows that the interpretation $(\lusim{c'}_\alpha)_H$ must be different from the one $p^h$ made, contradiction.
\end{enumerate}
\end{remark}
We will prove that the quotient forcing is $\kappa^+$-cc  for more general Prikry-type forcings which use $P$-point ultrafilters.
\begin{definition}\label{ppoint}
Let $F$ be a uniform $\kappa-$complete filter over a regular uncountable  cardinal $\kappa$.
$F$ is called a \textit{$P-$point filter} iff there is $\pi:\kappa\to \kappa$ such that

\begin{enumerate}
  \item $\pi$ is almost one to one, i.e. there is $X\in F$ such that for every $\alpha<\kappa$, $|\pi^{-1}{}\alpha\cap X|<\kappa$,
  \item For every $\{A_i \mid i<\kappa\}\subseteq F$,
  $\Delta^*_{i<\kappa}A_i=\{\nu<\kappa \mid \forall i<\pi(\nu) (\nu\in A_i)\} \in F$.
\end{enumerate}
\end{definition}
Clearly, every normal filter $F$ is a $P-$point, but there are many non-normal $P-$points as well. For example take a normal filter $U$ and move it to a non-normal by using a permutation on $\kappa$. Also, if $F$ is an ultrafilter, then $\pi$ is just a function representing $\kappa$ in the ultrapower by $F$.

Before proving the main result, we need a generalization of Galvin's theorem (see \cite{Glavin}, or \cite[Proposition 1.4]{GitDensity}):
\begin{proposition}\label{galvinGen}
 Suppose that $2^{<\kappa}=\kappa$ and let $F$ be a $P$-point filter  over $\kappa$ . Let $\langle X_i\mid i<\kappa^+\rangle$ be a sequence of sets such that for every $i<\kappa^+$, $X_i\in F$, and let $\langle Z_i\mid i<\kappa^+\rangle$ be any sequence of subsets of $\kappa$. Then there is $Y\subseteq \kappa^+$ of cardinality $\kappa$, such that 
 \begin{enumerate}
     \item $\bigcap_{i\in Y}X_i\in F$.
     \item there is $\alpha\notin Y$ such that $[Z_{\alpha}]^{<\omega}\subseteq \bigcup_{i\in Y}[Z_i]^{<\omega}$
 \end{enumerate}
\end{proposition}
\pr For every $\vec{\nu}\in[\kappa]^{<\omega}$, $\alpha<\kappa^+$ and $\xi<\kappa$, let
$$H_{\alpha,\xi,\vec{\nu}}=\{i<\kappa^+\mid X_i\cap\xi=X_\alpha\cap\xi \wedge \vec{\nu}\in [Z_i]^{<\omega}\}.$$
\begin{claim*}
There is $\alpha^*<\kappa^+$ such that for every $\xi<\kappa$ and $\vec{\nu}\in[Z_{\alpha^*}]^{<\omega}$, $|H_{\alpha^*,\xi,\vec{\nu}}|=\kappa^+$ 
\end{claim*}
\textit{Proof of claim.} Otherwise, for every $\alpha<\kappa^+$ there is $\xi_\alpha<\kappa$ and $\vec{\nu}_\alpha\in [Z_\alpha]^{<\omega}$ such that $|H_{\alpha,\xi_\alpha,\vec{\nu}_\alpha}|\leq\kappa$.
There is $X\subseteq \kappa^+$, $\vec{\nu}^*\in[\kappa]^{<\omega}$
and $\xi^*<\kappa$, such that $|X|=\kappa^+$ and $$\forall\alpha\in X, \ \vec{\nu}_\alpha=\vec{\nu}^*\wedge \xi_\alpha=\xi^*.$$
Since $\kappa$ is strong limit and $\xi^*<\kappa$, there are less than $\kappa$ many possibilities for $X_\alpha\cap \xi^*$. Hence we can shrink $X$ to $X'\subseteq X$ such that $|X'|=\kappa^+$ and find a single set $E^*\subseteq \xi^*$ such that for every $\alpha\in X'$, $X_\alpha\cap\xi^*=E^*$.
It follows that for every $\alpha\in X'$:
$$H_{\alpha,\xi_\alpha,\vec{\nu}_\alpha}=H_{\alpha,\xi^*,\vec{\nu}^*}=\{i<\kappa^+\mid X_i\cap \xi^*=E^*\wedge \vec{\nu}^*\in [Z_i]^{<\omega}\}.$$
Hence the set $H_{\alpha,\xi_\alpha,\vec{\nu}_\alpha}$ does not depend on $\alpha$, which means it is the same for every $\alpha\in X'$. Denote this set by $H^*$.
To see the contradiction, note that for every $\alpha\in X'$, $\alpha\in H_{\alpha,\xi_\alpha,\vec{\nu}_\alpha}=H^*$, thus $X'\subseteq H^*$, hence $$\kappa^+=|X'|\leq|H^*|\leq \kappa$$
contradiction.$\blacksquare_{\text{Claim}}$

\underline{End of proof of Proposition \ref{galvinGen}:} Let $\alpha^*$ be as in the claim. Let us choose $Y\subseteq \kappa^+$ that witnesses the lemma. First, enumerate $[Z_{\alpha^*}]^{<\omega}$ by $\langle \vec{\nu}_i\mid i<\kappa\rangle$. 
Let $\pi:\kappa\rightarrow\kappa$ be the function in Definitionref{ppoint} guaranteed by 
 $F$ being $P$-point. There is a set $X\in F$ such that for every $\alpha<\kappa$, $X\cap\pi^{-1}{}''\alpha$ is bounded in $\kappa$. So for every $\alpha<\kappa$, we find $\rho_\alpha>\sup(\pi^{-1}{}''[\alpha+1]\cap X)$.

By recursion, define $\beta_i$ for $i<\kappa$. At each step we pick $\alpha^*\neq\beta_i\in H_{\alpha^*,\rho_i+1,\vec{\nu}_i}\setminus\{\beta_j\mid j<i\}$. It is possible find such $\beta_i$, since the cardinality of $H_{\alpha^*,\rho_i+1,\vec{\nu}_i}$ is $\kappa^+$,  and $\{\beta_j\mid j<i\}$ is of size less than $\kappa$.
Let us prove that $Y=\{\beta_i\mid i<\kappa\}$ is as wanted. Indeed, by definition, it is clear that $|Y|=\kappa$. Also, if $\vec{\nu}\in [Z_{\alpha^*}]^{<\omega}$, then $\vec{\nu}=\vec{\nu}_i$ for some $i<\kappa$. By definition, $\beta_i\in H_{\alpha^*,\rho_i+1,\vec{\nu}_i}$, hence $\vec{\nu}\in [Z_{\beta_i}]^{<\omega}$, so $$[Z_{\beta_i}]^{<\omega}\subseteq\bigcup_{x\in Y}[Z_x]^{<\omega}.$$

Finally, we need to prove that $\bigcap_{\gamma\in Y}X_\gamma=\bigcap_{i<\kappa}X_{\beta_i}\in F$. By the $P$-point assumption about $F$, $$X^*:=X\cap X_{\alpha^*}\cap\Delta^*_{i<\kappa}X_{\beta_i}\in F.$$ Thus it suffices to prove that $X^*\subseteq \bigcap_{i<\kappa}X_{\beta_i}$. Let $\zeta\in X^*$, then for every $i<\pi(\zeta)$, $\zeta\in X_{\beta_i}$. For $i\geq\pi(\zeta)$, $\zeta\in \pi^{-1''}i+1\cap X$, and by definition of $\rho_i$, $\zeta<\rho_i$. Recall that $\beta_i\in H_{\alpha^*,\rho_i+1,\vec{\nu}_i}$  $$X_{\alpha^*}\cap(\rho_i+1)=X_{\beta_i}\cap(\rho_i+1)$$ and since $\zeta\in X_{\alpha^*}\cap(\rho_i+1)$, $\zeta\in X_{\beta_i}$. We conclude that $\zeta\in\bigcap_{i<\kappa}X_{\beta_i}$, thus $X^*\subseteq\bigcap_{i<\kappa}X_{\beta_i}$. $\blacksquare$

\begin{theorem}\label{kappapluscc}
Let $\pi:\Mfor\rightarrow \mathbb{P}$ be a projection and $G\subseteq\Mfor$ be $V$-generic and $H=\pi_*(G)$ be the induced generic for $\mathbb{P}$. Then $V[G]\models \Mfor/H$ is $\kappa^+$-cc 
\end{theorem}
\pr Assume otherwise, and let $\langle p_i\mid i<\kappa^+\rangle\in V[G]$ be an anthichain in $\Mfor/H$. Let $\langle\lusim{p}_i\mid i<\kappa^+\rangle$ be a sequence of $\Mfor$-names for them and $r\in G$ such that
$$r\Vdash \langle\lusim{p}_i\mid i<\kappa^+\rangle \text{ is an antichain in } \Mfor/\lusim{H}.$$
Work in $V$, for every $i<\kappa^+$, let $r\leq r_i\in\Mfor$ and $\xi_i\in\Mfor$ be such that $r_i\Vdash \lusim{p}_i=\xi_i$. 
\begin{claim*}
$\forall i<\kappa^+\exists q\geq \xi_i\forall q'\geq q\exists r''\geq r_i \ r''\Vdash q'\in\Mfor/\lusim{H}$
\end{claim*}
\textit{Proof of claim.} Otherwise, there is $i$ such that for every $q\geq \xi_i$, there is $q'\geq q$ such that every $r''\geq r_i$, $r''\not\Vdash q'\in\Mfor/\lusim{H}$. In particular, the set
$$E=\{q\geq \xi_i \mid \forall r''\geq r_i. r''\not\Vdash q\in\Mfor/\lusim{H}\}$$
is dense above $\xi_i$. To obtain a contradiction, let $G'$ be any generic for $\Mfor$ such that $r_i\in G'$. Since $ r_i\geq r$, $r\in G'$ and therefore $\xi_i=(\lusim{p}_i)_{G'}\in\Mfor/\lusim{H}_{G'}$. Denote $H'=\lusim{H}_{G'}$. Then by Proposition \ref{chracterization of quotient}, there is a $V$-generic filter $G''$ for $\Mfor$ such that $\xi_i\in G''$ and $\lusim{H}_{G''}=H'$. By density of $E$, there is $\xi_i\leq q\in E\cap G''$ and in particular, $q\in \Mfor/H'$. Thus, there is $r_i\leq r''\in G'$ such that $r''\Vdash q\in\Mfor/\lusim{H}$, contradicting $q\in E$.$\blacksquare_{\text{Claim}}$

For every $i<\kappa^+$ pick $q_i\geq \xi_i$ such that 
$$(*)_i \ \ \ \ \  \forall q'\geq q_i.\exists r''\geq r_i.r''\Vdash q'\in\Mfor/\lusim{H}.$$
Denote $q_i=\langle t_{i,1},...,t_{i,n_i},\langle \kappa,A(q_i)\rangle\rangle$ and $r_i=\langle s_{i,1},...,s_{i,m_i},\langle\kappa,A(r_i)\rangle\rangle$. Stabilize the sequences $\langle t_{i,1},...,t_{i,n_i}\rangle$ and $\langle s_{i,1},...,s_{i,m_i}\rangle$ i.e. find $X\subset \kappa^+$ such that $|X|=\kappa^+$ and $\vec{t}=\langle t_1,...,t_n\rangle,\vec{s}=\langle s_1,...,s_m\rangle$ such that for every $i\in X$
$$\langle t_{i,1},...,t_{i,n_i}\rangle=\langle t_1,...,t_n\rangle,\text{ and } \langle s_{i,1},...,s_{i,m_i}\rangle=\langle s_1,...,s_m\rangle.$$
This means that for every $i\in X$, $q_i=\vec{t}^{\smallfrown}\langle\kappa,A(q_i)\rangle$ and $r_i=\vec{s}^{\smallfrown}\langle\kappa ,A(r_i)\rangle$.
Let $$A^*(r_i)=\{\nu\in A(r_i)\mid \nu\cap A(r_i)\in\cap\vec{\nu}\}$$ by \ref{BetterSet}, $A^*(r_i)\in\cap\vec{U}(\kappa)$, it follows that for every $\vec{\nu}\in[A^*(r_i)]^{<\omega}$, $r_i^{\smallfrown}\vec{\nu}\in \Mfor$. By Lemma \ref{galvinGen},
there is $Y\subseteq X$ of cardinality $\kappa$, such that 
\begin{enumerate}
    \item $\bigcap_{i\in Y}A(q_i)\in \bigcap_{i<\kappa} U(\kappa,i)$.
    \item There is $\alpha^*\in Y$ such that $[A^*(r_{\alpha^*})]^{<\omega}\subseteq \bigcup_{i\in Y\setminus\{\alpha^*\}}[A^*(r_i)]^{<\omega}$
\end{enumerate}
 Consider the set $A=\bigcap_{i\in Y}A(q_i)$. For every $i\in Y$, $q_i\leq \vec{t}^{\smallfrown}\langle \kappa, A\rangle=:q^*$. Then by $(*)_{\alpha^*}$, there is $r''\geq r_{\alpha^*}$ such that $r''\Vdash q^*\in\Mfor/\lusim{H}$. Hence there is $\vec{s}\leq s''\in\Mfor\restriction \max(\kappa(\vec{s}))$, $k<\omega$, $\vec{\nu}\in[A(r_{\alpha^*})]^{k}$ and $B_1,...,B_k$ such that $$r''=\langle s'',\langle\nu_1,B_1\rangle,...,\langle\nu_k,B_k\rangle,\langle \kappa,A(r'')\rangle\rangle.$$

Since $r''\in\Mfor$, then $\vec{\nu}\in [A^*(r_{\alpha})]^{k}$ and by the property of $\alpha^*$, $\vec{\nu}\in\cup_{j\in Y\setminus\{\alpha^*\}}[A^*(r_j)]^{<\omega}$ and so there is $j\in Y$ such that $\vec{\nu}\in [A^*(r_j)]^{k}$. Since $r_{\alpha^*}$ and $r_j$ have the same lower part, and $\vec{\nu}\in [A^*(r_j)]^{<\omega}$, it follows that $r''$ and $r_j$ are compatible by the condition:
$$r^*=\langle s'', \langle\nu_1, B_1\cap A(r_j)\rangle,...\langle \nu_k, B_k\cap A(r_j)\rangle,\langle\kappa, A(r_j)\cap A(r'')\rangle\rangle.$$
To see the contradiction, note that $r^*\geq r_{\alpha^*},r_j$ and $r$, thus $$r^*\Vdash \lusim{p}_{\alpha^*}=\xi_{\alpha^*},\lusim{p}_j=\xi_j\text{ are incompatible in }\Mfor/\lusim{H}$$
but also $r^*\geq r''$, therefore
$$r^*\Vdash q^*\in\Mfor/\lusim{H}.$$
Since $q^*\geq q_{\alpha^*}\geq\xi_{\alpha^*}$ and $q^*\geq q_j\geq\xi_j$, then $r^*\Vdash \lusim{p}_{\alpha^*},\lusim{p}_j$ are compatible in $\Mfor/\lusim{H}$, contradiction.$\blacksquare$

This suffices to finish the induction step for $cf(\lambda)>\kappa$ and in turn \ref{MainResaultParttwo}.
\begin{corollary}\label{Inductionstephigh}
 Assume that $o^{\vec{U}}(\kappa)<\kappa^+$, $(IH)$ and let $A\in V[G]$ be a set of ordinals such that $cf(\sup(A))>\kappa$. Let $C^*$ be as in \ref{Candiadate}, then $A\in V[C^*]$ and $V[A]=V[C^*]$.
\end{corollary}
\pr  By \ref{Candiadate}, $C^*\subseteq C_G$ is such that $C^*\in V[A]$ and  $\forall\alpha<\lambda. A\cap\alpha\in V[C^*]$. Toward a contradiction assume that $A\notin V[C^*]$, and let $W:=V[C^*]$. 
 The quotient forcing $\Mfor/C^*\in W$  is $\kappa^+$-cc and therefore $cf(\lambda)$-cc in $V[G]=W[G]$ and $A$ is a fresh subsets of $\lambda$ contradicting Theorem \ref{no fresh subsets}.$\blacksquare_{\text{\ref{Inductionstephigh}}}$ $\blacksquare_{\text{\ref{MainResaultParttwo}}}
$ 
\subsection{The Quotient Forcing}
For $\Mfor/C^*$ which is $\kappa^+$-cc in $V[C^*]$, we can use a more abstract and direct argument:
 
 Suppose we have an iteration $P*\lusim{Q}$ of forcing notions. It is a classical result about the iteration that if for a regular cardinal $\lambda$ we have

\begin{enumerate}
  \item $P$ has $\lambda-$cc,
  \item $\Vdash_P \lusim{Q} \text{ has } \lambda$-cc,
  \end{enumerate}
then $P*\lusim{Q}$ satisfies $\lambda-$cc.

Also, if $P$ has $\lambda-$cc, $P*\lusim{Q}$ has $\lambda-$cc, then $\Vdash_P \lusim{Q} \text{ has } \lambda-cc$.
\\Suppose otherwise. Then there are $p\in P$ and a sequence of $P-$names $\l \lusim{q}_\alpha \mid \alpha<\lambda\r$ such that
$$p\Vdash_P  \l \lusim{q}_\alpha \mid \alpha<\lambda\r \text{ is an antichain in } \lusim{Q}.$$
Consider now $\{\l p,\lusim{q}_\alpha\r \mid \alpha<\lambda \}\subseteq P*\lusim{Q}$. By  $\lambda-$cc, there are $\alpha, \beta<\lambda, \alpha\not =\beta$ such that
$\l p,\lusim{q}_\alpha\r$ and $\l p,\lusim{q}_\beta\r$ are compatible.  Hence, there are $\l p', \lusim{q}'\r\geq \l p,\lusim{q}_\alpha\r,\l p,\lusim{q}_\beta\r$. But then
$$p' \Vdash_P \lusim{q}' \text{ is stronger than both } \lusim{q}_\alpha,\lusim{q}_\beta,$$
which is impossible, since $p'$ forces that they are members of an antichain.

However, in \ref{kappapluscc}, we address a different question:

\emph{Suppose that $P*\lusim{Q}$ satisfies $\lambda-$cc. Let $G*H$ be a generic subset of $P*\lusim{Q}$.
Consider the interpretation $Q$ of $\lusim{Q}$ in $V[G*H]$. Does it satisfy $\lambda-$cc?
}

Clearly, this is not true in general. For a simple example, let $P$ be trivial and $Q$ be the forcing for adding a branch to a Suslin tree.
Then, in $V^Q$, $Q$ will not be ccc anymore.

Our attention in Theorem \ref{kappapluscc} is to subforcings and projections of $\Mfor$, however the argument given is more general:

\begin{theorem}\label{kappaplusccgeneral}
Suppose that $ \calP$ is either Prikry or Magidor or Magidor-Radin or Radin or Prikry with a product of $P$-point ultrafilters forcing and $\lusim{Q}$ is a projection of $\calP$. Let $G(\calP)$ be a generic subset of $\calP$.
\\Then, the interpretation of  $\lusim{Q}$  in $V[G(\calP)]$, satisfies $\kappa^+-$cc there.
\end{theorem}

We do not know how to generalize this theorem to wider classes of Prikry type forcing notions.

For example the following may be the first step:

\begin{question}
Is the result valid for a long enough Magidor iteration of Prikry forcings?
\end{question}
The problem is that there is no single complete enough filter here, and so the Galvin theorem (or its generalization) does not seem to apply.
\begin{definition}

Let $F$ be a $\kappa-$complete  uniform filter over a set $X$, for a regular uncountable cardinal $\kappa$.
We say that $F$ has:
\begin{enumerate} 
\item The \emph{Galvin property} iff every family of $\kappa^+$ members of $F$ has a subfamily of cardinality $\kappa$ with intersection in $F$. 
\item The \emph{generalized Galvin property} iff it satisfies the conclusion of \ref{galvinGen}.
\end{enumerate}

\end{definition}

The following question looks natural in this context:
\begin{question}
Characterize filters (or ultrafilters) which satisfy the Galvin property (or the generalized Galvin property).
\end{question}
Construction by U. Abraham and S. Shelah \cite{AbrahamShelah1986} may be relevant here.
They constructed a model in which there is a sequence $\l C_i\mid i<2^{\mu^+}\r$ in $Cub_{\mu^+}$ such that the intersection of any $\mu^+$ clubs in the sequence is of cardinality less than $\mu$.  So the filter $Cub_{\mu^+}$ does not have the Galvin property. However $GCH$ fails there.
Lately, other results related to the Galvin property have been published  \cite{MR3787522}, \cite{MR3604115}, \cite{ghhm}. 
The following questions seem to be open:

\begin{question} Assume $GCH$. Let $\kappa$ be a regular uncountable cardinal.
Is there a $\kappa$-complete filter over $\kappa$ which fails to satisfy the Galvin property?
\end{question}
Let us note that if the ultrafilter is not on $\kappa$, then there is such an ultrafilter, namely, any fine $\kappa$-complete filter $U$ over $P_\kappa(\kappa^+)$ does not satisfy the Galvin property:

For every $\alpha<\kappa^+$, let $X_\alpha=\{Z\in P_\kappa(\kappa^+)\mid \alpha\in Z\}$, then $X_\alpha\in U$ since $U$ is fine but the intersection of any $\kappa$ elements from this sequence of sets is empty. 

A fine normal ultrafilter on $P_\kappa(\lambda)$ is used for the supercompact Prikry forcing (see \cite{Gitik2010} for the definition). Hence, the following question is natural:

\begin{question}
Assume $GCH$ and let $\lambda>\kappa$ be a regular cardinal. Is every quotient forcing of the supercompact Prikry forcing also $\lambda^+$-cc in the generic extension?
\end{question}

One particular interesting case is of filters which extend the closed unbounded filter.
\begin{question} Assume $GCH$. Let $\kappa$ be a regular uncountable cardinal. Is there a $\kappa-$complete filter which extends the closed unbounded filter $Cub_\kappa$  which fails to satisfy the Galvin property?
\end{question}

Our prime interest is on $\kappa-$complete ultrafilters over a measurable cardinal $\kappa$.
\\Note the following:

\begin{proposition}

It is consistent that every $\kappa-$complete(or even $\sigma$-complete) ultrafilter over a measurable cardinal $\kappa$ has the generalized Galvin property.

\end{proposition}
\pr 
This holds in the model $L[U]$, where $U$ is a unique normal measure on $\kappa$.
In this model every $\kappa$-complete ultrafilter is Rudin-Keisler equivalent to a finite power of $U$ (see for example \cite[Lemma 19.21]{Jech2003}). By \ref{product galvin}, it is easy to see that all such ultrafilters satisfy the generalized Galvin property. Note that since in $L[U]$ there is a unique measurable cardinal, every $\sigma$-complete ultrafilter $W$ is actually $\kappa$-complete. Indeed, let $\lambda$ be the completness degree of $W$\footnote{The degree of completness of an ultrafilter $\mathcal{V}$ is the minimal cardinal $\theta$ such that $\mathcal{V}$ is $\theta$-complete.}, it is the critical point of the embedding
$$j_W:L[U]\rightarrow Ult(L[U],W).$$
Since $W$ is $\sigma$-complete, $Ult(L[U],W)$ is well-founded, hence $crit(j_W)$ is measurable in $L[U]$ and $\lambda=\kappa$.
$\blacksquare$

In context of ultrafilters over a measurable cardinal, the following is unclear:

\begin{question}
Is it consistent to have a $\kappa$-complete ultrafilter over $\kappa$ which does not have the Galvin property?
\end{question}

\begin{question}  Is it consistent to have a measurable cardinal $\kappa$ carrying a $\kappa-$complete ultrafilter which extends the closed unbounded filter $Cub_\kappa$ 
(i.e., $Q-$point) which fails to satisfy the Galvin property?
\end{question}

It is possible to produce more examples of ultrafilters (and  filters) with generalized Galvin property.
The simplest example of this kind will be $U\times W$, where $U, W$ are normal ultrafilters over $\kappa$.
We will work in a bit more general setting.

\begin{definition}\label{incresing P-points}
Let $F_1,...,F_n$ be $P$-point filters over $\kappa$, and let $\pi_1,...,\pi_n$ be the witnessing functions for it.
Denote $[\kappa]^{n*}$, the set of all $n$-tuples $\l \alpha_1,...,\alpha_n\r$ such that for every $2\leq i\leq n$, $\alpha_{i-1}<\pi_i(\alpha_i)$.
\end{definition}
Note that if the $F_i$'s are normal, the $\pi_i=id$ and $[\kappa]^{n*}=[\kappa]^n$.
\begin{definition}\label{product of p-points}
Let $F_1,...,F_n$ be $P$-point filters over $\kappa$, and let $\pi_1,...,\pi_n$ be the witnessing functions for it.
Define a filter $\prod_{i=1}^{n*} F_i$ over $[\kappa]^{n*}$ recursively. For $X\subseteq[\kappa]^{n*}$:
\begin{center}
    $ X \in\prod_{i=1}^{n*} F_i\Leftrightarrow \Big\{\alpha<\kappa\mid X_{\alpha}\in \prod_{i=2}^{n*} F_i\Big\}\in F_1.$
    
\end{center}
Where $X_{\alpha}=\{\l\alpha_2,...,\alpha_n\rangle\in[\kappa]^{n-1*} \mid \l\alpha,\alpha_2,...,\alpha_n\r\in X\}$. 
\end{definition}
Again, if the filters are normal, this is simply a product.

\begin{proposition}\label{generate P-point}
Let $F_1,...,F_n$ be $P$-point filters over $\kappa$, and let $\pi_1,...,\pi_n$ be the witnessing functions for it. Then for every $X\in\prod_{i=1}^{n*} F_i$, there are $X_i\in F_i$ such that $\prod_{i=1}^{n*}X_i\subseteq X$.
\end{proposition}
\pr By induction on $n$, for $n=1$, it is clear. Let $X\in\prod_{i=1}^{n*} F_i$. Let 
\begin{center}
$X_1=\Big\{\alpha<\kappa\mid X_{\alpha}\in \prod_{i=2}^{n*} F_i\Big\}\in F_1$
\end{center}
For every $\alpha\in X_1$, find by the induction hypothesis $X_{\alpha,i}\in F_i$ for $2\leq i\leq n$ such that $\prod_{i=2}^{n*}X_{\alpha,i}\subseteq X_{\alpha}$.
Define 
$$X_i=\Delta^*_{\alpha<\kappa}X_{\alpha,i}$$
since $F_i$ is $P$-point, $X_i\in F_i$.
Let us argue that $\prod_{i=1}^{n*} X_i\subseteq X$. Let $\l\alpha_1,...,\alpha_n\r\in \prod_{i=1}^{n*} X_i$, then for every $2\leq i\leq n$, $\alpha_1<\pi(\alpha_i)$, hence $\alpha_i\in X_{\alpha_1,i}$. It follows that $\l \alpha_2,...,\alpha_n\r\in \prod_{i=2}^{n*} X_{\alpha_1,i}\subseteq X_{\alpha_1}$. By definition of $X_{\alpha_1}$, $\langle \alpha_1,\alpha_2...\alpha_n\r\in X$.$\blacksquare$ 
\begin{corollary}\label{product galvin}
Assume that $2^{<\kappa}=\kappa$. Let $F_1,...,F_n$ be $P$-point filters over $\kappa$, and let $\pi_1,...,\pi_n$ be the witnessing functions for it. Then $\prod_{i=1}^{n*} F_i$
 also satisfies the generalized Galvin property of \ref{galvinGen}.
\end{corollary}
\pr 
Let $\langle Y_\alpha\mid \alpha<\kappa^+\r$ and $\l Z_\alpha\mid \alpha<\kappa^+\r$ be as in \ref{galvinGen}. By Proposition \ref{generate P-point}, for every $1\leq i\leq n$, and $\alpha<\kappa^+$, find $X^{(\alpha)}_i\in F_i$ such that $\prod_{i=1}^{*n} X^{(\alpha)}_i\subseteq Y_\alpha$.

For every $\vec{\alpha}=\l \alpha_1,...,\alpha_n\r\in [\kappa]^{n*}$ every $\vec{\nu}\in [\kappa]^{<\omega}$ and every $\xi<\kappa^+$, define
$$H_{\xi,\vec{\alpha},\vec{\nu}}=\Big\{\gamma<\kappa^+\mid \forall 1\leq i\leq n. X^{(\gamma)}_i\cap \alpha_i= X^{(\xi)}_i\cap \alpha_i\text{ and } \vec{\nu}\in [Z_\gamma]^{<\omega}\Big\}.$$
As in \ref{galvinGen}, since there are less than $\kappa^+$ many possibilities for $\langle X^{(\gamma)}_1\cap\alpha_1,X^{(\gamma)}_2\cap\alpha_2,...,X^{(\gamma)}_n\cap\alpha_n\r$, we can find $\alpha^*<\kappa^+$, such that for every $\vec{\alpha}$ and $\vec{\nu}$, $|H_{\alpha^*,\vec{\alpha},\vec{\nu}}|=\kappa^+$.

Enumerate $[Z_{\alpha^*}]^{<\omega}$ by $\langle \vec{\nu}_i\mid i<\kappa\rangle$ and also 
each $F_i$ is $P$-point, so for every $j<\kappa$, there is $\rho^{(j)}_i>\sup(\pi_i^{-1''}[j]\cap B_i)$ for some set $B_i\in F_i$.
Define the sequence $\beta_j$ by induction, $$\beta_j\in H_{\alpha^*,\l \rho^{(j)}_1,...,\rho^{(j)}_n\r,\vec{\nu}_j}\setminus\{\beta_k\mid k<j\}.$$
We claim once again that \begin{center}$X_{\alpha^*}\cap\bigcap_{j<\kappa} X_{\beta_j}\in \prod_{i=1}^{n*} F_i$\end{center}
To see this, define for every $1\leq i\leq n $ $$C_i:=X^{(\alpha^*)}_i\cap \Delta^*_{j<\kappa} X^{(\beta_j)}_i\in F_i.$$
Let $\vec{\alpha}\in\prod_{i=1}^{n*} C_i$, and let $j<\kappa$. For every $1\leq i\leq n$, if $j<\pi(\alpha_i)$ then $\alpha_i\in X^{(\beta_j)}_i$. If $\pi(\alpha_i)\leq j$, then $\alpha_i<\rho^{(j)}_i$, so $\alpha_i\in X^{(\alpha^*)}\cap\rho^{(j)}_i$. Since $\beta_j\in H_{\alpha^*,\l\rho^{(j)}_1,...,\rho^{(j)}_n\r,\vec{\nu}_j}$, $$\alpha_i\in X^{(\alpha^*)}\cap\rho^{(j)}_i=X^{(\beta_j)}\cap\rho^{(j)}_i.$$ Therefore, $\vec{\alpha}\in \prod_{i=1}^{n*}X^{(\beta_j)}_i\subseteq Y_{\beta_j}$. The continuation is as in \ref{galvinGen}.$\blacksquare$

\section{Fresh sets}

Let us conclude this paper with the following result about fresh sets in Magidor generic extensions.  A very close variation of this can be found in \cite{Omer1}.
\begin{theorem}\label{Fresh}
Let $\vec{U}$ be a coherent sequence on $\kappa$.
Let $G\subseteq\Mfor$ be $V$-generic. If $A\in V[G]$ is a fresh set of ordinals with respect to $V$, then $cf^{V[G]}(\sup(A))=\omega$.\footnote{ Clearly, if $\kappa$ changes its cofinality to $\omega$, then any cofinal in $\kappa$ sequence of the order type $\omega$ will be fresh.}
\end{theorem}
Note that we do not restrict the order of $\kappa$ and by taking $o^{\vec{U}}(\kappa)=1$ we obtain the Prikry forcing.

\pr 
By induction on $\kappa$, which is the supermum of $C_G$. Let $A$ be a fresh subset, if $A\in V[C_G\cap\alpha]$ for some $\alpha<\kappa$, by the induction hypothesis we are done. Assume that $\forall\alpha<\kappa. A\notin V[C_G\cap\alpha]$, in particular $\sup(A)\geq \kappa$. Let us start with the difficult part, where $\sup(A)=\kappa$.
\begin{lemma}
If $A\in V[G]$ is fresh subset of $\kappa$ with respect to $V$ such that $\sup(A)=\kappa$, then $cf^{V[G]}(\kappa)=\omega$.
\end{lemma}

\pr Toward a contradiction assume that $\lambda:=cf^{V[G]}(\kappa)>\omega$ and let $\lusim{A}$ be a name for $A$.

First we deal with the case that $\kappa$ is singular in $V[G]$, hence $\omega<\lambda<\kappa$.

Since $\Mfor$ decomposes to the part below $\lambda$ and the part above $\lambda$, we can ensure sufficient closure, by working in $V[C_G\cap\lambda]$, and force with the part of the forcing above $\lambda$. Note that $A$ is fresh also with respect to $V[C_G\cap\lambda]$.

Let $\l c_\alpha\mid \alpha<\lambda\rangle$ be a cofinal continuous subsequence of $C_G$ such that $c_0>\lambda$. Let $\langle\lusim{c'}_\alpha\mid \alpha<\lambda\r$ be a sequence of $\Mfor\restriction(\lambda,\kappa]$-names for it.

Find $p\in G\restriction (\lambda,\kappa)$
such that $$p\Vdash \lusim{A}\text{ is fresh}\wedge \langle \lusim{c'}_\alpha\mid \alpha<\lambda\r\text{ is a cofinal continuous subsequence of }C_G.$$

For every $i<\lambda$, the set  $$D_i=\Big\{q\mid \exists\vec{\alpha} \ \exists B. \ p^{\smallfrown}\vec{\alpha}\leq^* q\wedge \ q\Vdash \lusim{c}_i=\max{\vec{\alpha}}\wedge \lusim{A}\cap\max{\vec{\alpha}}=B\Big\}$$

is dense. To see that, let $q_0\geq p$, find any $q_0\leq q$ and $\vec{\beta}$ such that $$p^{\smallfrown}\vec{\beta}\leq^* q\text{ and }q\Vdash \max(\vec{\beta})=\lusim{c}_i.$$ Above $\max(\vec{\beta})$ there is enough closure to decide $\lusim{A}\cap\max(\vec{\beta})$. Find $q\restriction(\max(\vec{\beta}),\kappa]\leq^*q_{>\max(\vec{\beta})}$ in $\Mfor\restriction(\max(\vec{\beta}),\kappa]$ which decides $\lusim{A}\cap\max(\vec{\beta})$ and $q\restriction\max(\vec{\beta})\leq q_{\leq\max(\vec{\beta})}$ in $\Mfor\restriction(\lambda,\max(\vec{\beta})]$ (not necessarily a direst extension) such that for some $B\subseteq \max(\vec{\beta})$,
$$q^*:=\langle q_{\leq\max(\vec{\beta})},q_{>\max(\vec{\beta})}\r\Vdash\lusim{A}\cap\max(\vec{\beta})=B\wedge \lusim{c}_i=\max(\vec{\beta}).$$
Let $\vec{\alpha}$ be such that $p^{\smallfrown}\vec{\alpha}\leq^* q^*$, then by construction $\max(\vec{\alpha})=\max(\vec{\beta})$ and $q^*$ is as wanted.

By \ref{strongPrkryProperty}, find a condition $p\leq^*p_i$, a $\vec{U}$-fat tree of extensions of $p_i$, $T_i$, and sets $B_i^t$ such that for every $ t\in mb(T_i)$ there is $A_i(t)\subseteq \max(t)$ such that $$p_i{}^{\smallfrown}
\l t, \vec{B}^t_i\r \Vdash\lusim{A}\cap\max(t)=A_i(t)\wedge\lusim{c}_i=\max(t).$$
Since we have sufficient closure in the forcing above $\lambda$, we can find a single $p\leq^* p^*\in G\restriction(\lambda,\kappa]$ such that for every $i<\lambda$, $p_i\leq^* p^*$.

Keep defining by recursion sets $A_i(s)$ for $s\in T_i\setminus mb(T_i)$. Let $s\in \Lev_{ht(T_i)-1}(T_i)$, then we can shrink $\succ_{T_i}(s)$ and find $A_i(s)$ such that for every $\alpha\in\succ_{T_i}(s)$, $A_i(s^{\smallfrown}\alpha)= A_i(s)\cap\alpha$.
 
Generally, take $s\in T_i$ and assume that for every $\alpha$ in $\succ_{T_i}(s)$, $A_i(s^{\smallfrown}\alpha)$ is defined. 
We can find a single $A_i(s)$ and shrink $\succ_{T_i}(s)$ such that 

$$(\star) \ \ \ \forall\alpha\in\succ_{T_i}(s).A_i(s^{\smallfrown}\alpha)\cap\alpha= A_i(s)\cap\alpha.$$

Move to $V[A]$, let us compare the sets $A_i(s)$ with $A$. 
For every $i$, define recursively $\rho^i(k)$ for $k\leq N_i:=ht(T_i)$.
Let $\rho^i(0)=\min(A\Delta A_i(\l\r))+1$. Recursively define
$$\rho_i(k+1)=\sup(\min(A\Delta A_i(\l\delta_1,...,\delta_k\r))+1\mid \delta_1<\rho^i(0),...,<\delta_k<\rho^i(k))).$$
Let $\vec{c}_i\in mb(T_i)$ be such that $p^{*}{}^{\smallfrown}\l\vec{c}_i,\vec{B}^{\vec{c}_i}_i\r\in G$, let us argue that for every $k\leq N_i$, $\rho^i(k)> \vec{c}_i(k)$.
By construction of the tree $T_i$, $c_i=(\lusim{c}_i)_G=\max(\vec{c}_i)$ and $A\cap c_i=A_i(\vec{c}_i)\cap c_i$. 

By $(\star)$, for every $j\leq N_i$,
$$A_i(\vec{c}_i\restriction j)\cap \vec{c}_i(j)=A_i(\vec{c}_i\restriction j+1)\cap\vec{c}_i(j).$$
It follows that for every $j\leq N_i$,
$$A_i(\vec{c}_i\restriction j)\cap \vec{c}_i(j)=A\cap\vec{c}_i(j)$$
 In particular, 
$A\cap\vec{c}_i(0)=A_i(\l\r)\cap \vec{c}_i(0)$. 

Since $A\cap\rho^i(0)\neq A_i{\l\r}\cap\rho^i(0)$, it follows that $\vec{c}_i(0)<\rho^i(0)$. 

Inductively assume that $\vec{c}_i(j)<\rho^i(j)$ for every $j\leq k$. Since $A_i(\vec{c}_i\restriction k+1)\cap\vec{c}_i(k+1)=A\cap \vec{c}_i(k+1)$,
then $$\vec{c}_i(k+1)<\min(A_i(\vec{c}_i\restriction k+1)\Delta A)\leq \rho^i(k+1)$$

Before proving that $cf^{V[G]}(\kappa)=\omega$, let us argue that $\rho^i(k)<\kappa$.
Again by induction on $k$,
$\rho^i(0)<\kappa$ since $A\neq A_i(\l\r)$,  as $A_i(\l\r)\in V[C_G\cap\lambda]$ and $A\notin V[C_G\cap\lambda]$.

Toward a contradiction assume that $\rho^i(k+1)=\kappa$. Back to $V[C_G\cap\lambda]$, consider the collection $$\{A_i(\l\alpha_0,...,\alpha_k\r)\mid\alpha_0<\rho^i(0),...,\alpha_k<\rho^i(k)\}$$
Then for every $\gamma<\kappa$ pick any distinct $\vec{\alpha}_1,\vec{\alpha}_2$ such that $A_i(\vec{\alpha}_1)\neq A_i(\vec{\alpha}_2)$, but $A_i(\vec{\alpha}_1)\cap\gamma=A_i(\vec{\alpha}_2)\cap\gamma$.

To see that there are such $\vec{\alpha}_1,\alpha_2$, if $\rho^i(k+1)=\kappa$ there is $\vec{\alpha}_1$ such that $\eta_1:=\min(A\Delta A_i(\vec{\alpha}_1))>\gamma$,
 hence $A_i(\vec{\alpha}_1)\cap\gamma=A\cap\gamma$. Let $\vec{\alpha}_2$ be such that $\min(A\Delta A_i(\vec{\alpha}_2))>\eta_1$.
 In particular, $A_i(\vec{\alpha}_1)\neq A_i(\vec{\alpha}_2)$, but  $A_i(\vec{\alpha}_1)\cap\gamma=A\cap\gamma=A_i(\vec{\alpha}_2)\cap\gamma$.
Since this is all in $V[C_G\cap\lambda]$, where $\kappa$ is still measurable, we can find unboundedly many $\gamma$'s with the same $\vec{\alpha}_1,\vec{\alpha}_2$, which is clearly a contradiction.

So we found a sequence  $\langle\rho^i(N_i)\mid i<\lambda\rangle\in V[A]$ such that $\rho^i(N_i)>c_i$. Let $Z$ be the closure of $\{\rho^i(N_i)\mid i<\lambda\}$. Since $\lambda>\omega$, there is some limit $\alpha<\lambda$ such that $c_\alpha<\kappa$ is a limit point of $Z$. 

To see the contradiction, note that on one hand, $A\cap c_\alpha\in V[C_G\cap\lambda]$, and therefore the set $Z\cap c_\alpha$, $|Z\cap c_{\alpha}|=\lambda$ is defined in $V[C_G\cap\lambda]$ from $A\cap c_\alpha$, on the other hand, $c_\alpha>\lambda$, thus $c_\alpha$ should stay measurable in $V[C_G\cap\lambda]$, contradiction.

Next we eliminate the case that $\kappa$ is regular in $V[G]$ i.e. $\lambda=\kappa$. Many of the ideas for $\lambda<\kappa$ will also work here.

We no longer work over the model $V[C_G\cap\lambda]$, instead,  simply force over $V$. Let $p\in G$ be such that
$$p\Vdash \lusim{A}\text{ is fresh}$$

By induction we construct a $\leq^*$-increasing sequence of conditions $p_n$ and a tree of trees, i.e. a tree $T_0$, trees $ T_{1,t_0}$ for $t_0\in mb(T_0)$,  and generally  trees $ T_{n+1,t_0,...,t_n}$ where $$t_0\in mb(T_0), t_1\in mb(T_{1,t_0})... t_n\in mb(T_{n, t_0,...,t_{n-1}})$$

First find a condition $p\leq^*p_0$ and take  the tree $T_0$ to be simply the tree with one level which decides $C_{\lusim{G}}(0)$ if it is not already decided, or $T_0=\{\l\r\}$ otherwise. Necessarily,  for each $\alpha\in T_0$, $\min(\kappa(p^{\smallfrown}\alpha))=\alpha$, hence there is enough $\leq^*$-closure to decide $\lusim{A}\cap\alpha$, so we find $p^{\smallfrown}\alpha\leq^* p_\alpha$ and a set $A_0(\alpha)$ such that $p_\alpha\Vdash \lusim{A}\cap\alpha= A_0(\alpha)$. Then $p_0$ is obtained by diagonally intersecting all the sets in $p_\alpha$, and $p_0$ has the following property
$$\forall \alpha\in T_0.\  p_0{}^{\smallfrown}\alpha \Vdash\lusim{A}\cap\alpha=A_0(\alpha)\wedge C_{\lusim{G}}(0)=\alpha$$

For clarity, let us present also the construction of  $p_1$ and $T_{1,t_0}$ for every $t_0\in mb(T_0)$. The proof regarding the construction will be addressed later, in the general definition.

If necessary, find a direct extension of $p_0$ and use ineffability to find a set $A_0(\l\r)\subseteq\kappa$, such that for every $\alpha\in T_0$, $A_0(\l\r)\cap\alpha=A_0(\alpha)\cap\alpha$.

In $V[A]$, define $\eta_0=\min(A\Delta A_0(\l\r))$, since $A\notin V$ and $A_0(\l\r)\in V$, $\eta_0<\kappa$ is well defined. Clearly, for every $V$-generic filter $H$ with $p_0\in H$, $\eta_0> C_H(0)$, since then $p_0^{\smallfrown} C_H(0)\in H$ forces the correct value of $A$. Let $\lusim{\eta}_0$ be a name such that $p_0\Vdash \lusim{\eta}_0=\min(\lusim{A}\Delta A_0(\l\r))$.

Fix $t_0\in mb(T_0)$, consider $p_0^{\smallfrown}t_0$. In the general case we will prove that we can find $p_0^{\smallfrown}t_0\leq^* p_{t_0}$, $T_{1,t_0}$ and sets $Y_1^t$ for $t\in mb(T_{1,t_0})$ such that for every $t_1\in mb(T_{1,t_0})$ there is $A_1(t_0,t_1)\subseteq\max(t_1)$, such that
$$p_{t_0}^{\smallfrown}\l t_1, \vec{Y}^{t_1}_1\r\Vdash  A_1(t_0,t_1)\cap\max(t_1)= \lusim{A}\cap\max(t_1)\wedge \max(t_1)=C_{\lusim{G}}(\lusim{\eta}_0)$$
Note that $$p_0^{\smallfrown}t_0\Vdash\max(t_0)<\lusim{\eta}_0\leq C_{\lusim{G}}(\lusim{\eta}_0)$$
Hence $\max(t_1)>\max(t_0)$.

Find a single $p_0\leq^*p_1$ such that for every $t_0\in mb(T_0)$, $p_{t_0}\leq ^*p_1^{\smallfrown}t_0$.

If necessary, directly extend $p_1$ to get $N_1<\omega$, such that for every $t_0\in mb(T_0)$, $ht(T_{1,t_0})=N_1$.

Define the sets $A_1(t_0,s)$ for every $s\in T_{1,t_0}\setminus mb(T_{1,t_0})$. 
Let $s\in \Lev_{N_1-1}(T_{1,t_0})$, we can shrink $\succ_{T_{1,t_0}}(s)$ and find $A_1(t_0,s)\subseteq\kappa$ such that for every $\alpha\in\succ_{T_{1,t_0}}(s)$, $A_1(t_0,s^{\smallfrown}\alpha)= A_0(s)\cap\alpha$.

Recursively, let $s\in T_{1,t_0}\setminus mb(T_{1,t_0})$ and assume that for every $\alpha$ in $\succ_{T_0}(s)$, $A_1(t_0,s^{\smallfrown}\alpha)$ is defined. 
Find a single $A_1(t_0,s)$ and shrink $\succ_{T_{1,t_0}}(s)$ such that $$\forall\alpha\in\succ_{T_{1,t_0}}(s). \  A_1(t_0,s^{\smallfrown}\alpha)\cap\alpha= A_1(t_0,s)\cap\alpha$$

In $V[A]$, define $\rho^1(k)$ for every $k\leq N_1$. For $k=0$,
$$\rho^1(0)=\sup(\min(A\Delta A_1(t_0,\l\r))\mid t_0\in mb(T_0)\cap[\eta_0]^{<\omega})$$
Recursively, 
$$\rho^1(k+1)=\sup(\min(A\Delta A_1(t_0,\l \alpha_0,...,\alpha_k\r))\mid t_0\in mb(T_0)\cap[\eta_0]^{<\omega}\wedge \alpha_i<\rho^1(i))$$
Note that for every $t_0\in mb(T_0)$ and $s\in T_{1,t_0}$, $A\neq A_1(t_0,s)$, as $A_1(t_0,s)\in V$  and $A\notin V$. Therefore $\rho^1(k)\leq\kappa$ is well defined for every $k\leq N_1$. In the general case we will also prove that $\rho^1(k)<\kappa$.
Finally, define $\eta_1=\rho^1(N_1)$ and let $\lusim{\eta}_1$ be a name such that $p_1$ forces $\lusim{\eta}_1$ is computed by comparing the sets $A_1(t_0,s)$ and $\lusim{A}$, the way we defined it.

Now for the general definition, assume we have defined $p\leq^* p_1\leq^* p_2...\leq^*p_{n}$,   trees $T_{n,t_0,...,t_{n-1}}$ for $t_0\in mb(T_0),t_1\in mb( T_{1,t_0}),...,t_{n-1}\in T_{n-1,t_0,...,t_{n-2}}$, sets    $A_n(t_0,...,t_{n-1},t_n)$ for every $t_n\in mb(T_{n,t_0,...,t_{n-1}})$  and $Y^{t_0}_0,...,Y^{t_{n}}_{n}$, also a name $\lusim{\eta}_{n-1}$ such that,
$$p_{n}^{\smallfrown}\l t_0,\vec{Y}^{t_0}_0\r^{\smallfrown}\l t_1,\vec{Y}^{t_1}_1\r...^{\smallfrown}\l t_n,Y^{t_n}_n\r\Vdash \lusim{A}\cap \max(t_n)=A_n(t_0,...,t_{n-1})\cap\max(t_n)\wedge\max(t_n)= C_{\lusim{G}}(\lusim{\eta}_{n-1})$$

Define recursively the sets $A_n(t_0,...,t_{n-1},s)$ for $s\in T_{n,t_0,...,t_{n-1}}\setminus mb(T_{n,t_0,...,t_{n-1}})$. Assume that $$A_n(t_0,...,t_{n-1},s^{\smallfrown}\alpha)$$ is defined, for every $\alpha\in \succ_{T_{n,t_0,...,t_{n-1}}}(s)$. Directly extend $p_n$ if necessary, shrink  $ \succ_{T_{n,t_0,...,t_{n-1}}}(s)$ and find by ineffability $A_n(t_0,...,t_{n-1},s)$ so that
for every $\alpha\in\succ_{T_{n,t_0,...,t_{n-1}}}(s)$, $$A_n(t_0,...,t_{n-1},s)\cap\alpha=A_n(t_0,...,t_{n-1},s^{\smallfrown}\alpha)\cap\alpha.$$

In $V[A]$, we have defined $\eta_0,...,\eta_{n-1}$, and so we can define $$\rho^n(0)=\sup\Big[\min\Big(A_n(t_0,...,t_{n-1},\l\r)\Delta A\Big)\mid t_0\in mb(T_1)\cap[\eta_0]^{<\omega},...,t_{n-1}\in mb(T_{n-1})\cap[\eta_{n-1}]^{<\omega}\Big]$$
keep defining $\rho^n(k+1)$ recursively as
$$\sup\Big[\min\Big(A_n(t_0,...,t_{n-1},\vec{\alpha})\Delta A\Big)\mid t_0\in mb(T_1)\cap[\eta_0]^{<\omega},...,t_{n-1}\in mb(T_{n-1})\cap[\eta_{n-1}]^{<\omega}, \ \vec{\alpha}\in\prod_{j=1}^k\rho^n(j)\Big]$$
Finally define, $\eta_n=\rho^n(N_n)$.

Again note that $\rho^n(k)\leq\kappa$ is a well define ordinal. Let us prove that $\rho^n(k)<\kappa$. 
\begin{claim*}
For every $k\leq N_n$, $\rho^n(k)<\kappa$.
\end{claim*}
\textit{Proof of claim.} The proof is similar to the case that $\kappa$ is singular in $V[G]$. Toward a contradiction assume that $\rho^n(k)=\kappa$. Back in $V$, consider the collection $$\Big\{A_n(t_0,...,t_{n-1},\vec{\alpha}) \mid t_0\in mb(T_1)\cap[\eta_0]^{<\omega},...,t_{n-1}\in mb(T_{n-1})\cap[\eta_{n-1}]^{<\omega}, \ \vec{\alpha}\in\prod_{j=1}^{k-1}\rho^n(j)\Big\}$$
 Then for every $\gamma<\kappa$ pick any distinct $t_1,...,t_{n-1},\vec{\alpha}$ and $s_0,...,s_{n-1},\vec{\beta}$ such that $$A_n(t_1,...,t_{n-1},\vec{\alpha})\neq A_n(s_0,...,s_{n-1},\vec{\beta}),\text{ but } A_n(t_1,...,t_{n-1},\vec{\alpha})\cap\gamma=A_n(s_0,...,s_{n-1},\vec{\beta})\cap\gamma$$
 To see that there are such $t_1,...,t_{n-1},\vec{\alpha}$ and $s_0,...,s_{n-1},\vec{\beta}$, by assumption that $\rho^n(k)=\kappa$ there are $t_1,...,t_{n-1},\vec{\alpha}$ such that $\xi_1:=\min(A\Delta A_n(t_1,...,t_{n-1},\vec{\alpha}))>\gamma$,
 hence $$A\Delta A_n(t_1,...,t_{n-1},\vec{\alpha})\cap\gamma=A\cap\gamma$$ Find $s_0,...,s_{n-1},\vec{\beta}$, such that $\min(A\Delta A_n(s_0,...,s_{n-1},\vec{\beta}))>\xi_1$.
 In particular, $$A_n(t_1,...,t_{n-1},\vec{\alpha})\neq A_n(s_0,...,s_{n-1},\vec{\beta})\text{ but }A_n(t_1,...,t_{n-1},\vec{\alpha})\cap\gamma=A\cap\gamma=A_n(s_0,...,s_{n-1},\vec{\beta})\cap\gamma$$
Since this is all in $V$, where $\kappa$ is measurable, we can find unboundedly many $\gamma$'s with the same $t_1,...,t_{n-1},\vec{\alpha},s_0,...,s_{n-1},\vec{\beta}$, which is clearly a contradiction.$\blacksquare_{\text{Claim}}$

 Find a name $\lusim{\eta}_n$ such that $p_n$ forces $\lusim{\eta}_n$ is obtained by comparing $\lusim{A}$ with the sets $A_n(t_0,...,t_n,\vec{\alpha})$ as above using $\lusim{\eta}_0,...,\lusim{\eta}_{n-1}$.

Now for the definition of the trees, fix $t_0,...,t_n$ such that
$$t_0\in mb(T_0),t_1\in mb( T_{1,t_0}),...,t_{n}\in mb(T_{n,t_0,...,t_{n-1}})$$ The set $D$ of all conditions $q$ such that for some $\vec{\alpha}\in[\kappa]^{<\omega}$:
\begin{enumerate}
    \item $p_{n}^{\smallfrown}\l t_0,\vec{Y}^{t_0}_0\r^{\smallfrown}...{}^{\smallfrown}\l t_n,\vec{Y}^{t_n}_n\r^{\smallfrown}\vec{\alpha}\leq^*q$.
    \item $q\Vdash \lusim{A}\cap\max(\vec{\alpha})=A(\vec{\alpha})\wedge \max(\vec{\alpha})=C_{\lusim{G}}(\lusim{\eta}_{n})$.
\end{enumerate}
 
is dense above $p_{n}^{\smallfrown}\l t_0,\vec{Y}^{t_0}_0\r^{\smallfrown}...{}^{\smallfrown}\l t_n,\vec{Y}^{t_n}_n\r$. The proof is as for the case $\kappa$ is singular.

By \ref{strongPrkryProperty} and \ref{densetree}, find a condition $p_{n}^{\smallfrown}t_0^{\smallfrown}...^{\smallfrown}t_{n-1}{}^{\smallfrown}t_n\leq^*p_{t_0,...,t_n}$, a $\vec{U}$-fat tree of extensions of $p_{t_0,...,t_n}$, $T_{n+1,t_0,...,t_n}$, and sets $Y^s_{n+1}$, such that for every $ t\in mb(T_{n+1,t_0,...,t_n})$ there is $A_{n+1}(t_0,...,t_n,t)\subseteq \max(t)$ for which $$p_{t_0,...,t_n}{}^{\smallfrown}\l t,\vec{Y}^t_{n+1}\r \Vdash\lusim{A}\cap\max(t)=A_{n+1}(t_0,...,t_n,t)\wedge C_{\lusim{G}}(\lusim{\eta}_n)=\max(t)$$
and the set $D_{T_{n+1,t_0,...,t_n},\vec{Y}_n+1}$ is dense above $p_{t_0,...,t_n}$.

By \ref{amalgamate}, find a single $p_n\leq^* p_{n+1}$ and shrink the trees $T_{i,t_0,...,t_i}$ such that for every $t_0,...,t_n$. $p_{t_0,...,t_n}\leq^* p_{n+1}^{\smallfrown}t_0,...,t_n$

By shrinking even more if necessary, we can assume that  there is $N_{n+1}$,  such that for every $t_0,...,t_n$ $ht(T_{n+1,t_0,...,t_n})=N_{n+1}$. This concludes the recursive definition.

%By construction of the tree $T_i$, $A\cap c_i=A^i_{\vec{c}_i}\cap c_i$. now for every $j\leq n_i$, $$A\cap (\vec{c}_i)_j=A^i_{\l(\vec{c}_i)_0,...,(\vec{c}_i)_{j-1}\r}\cap (\vec{c}_i)_j$$ In particular, 
%$A\cap(\vec{c}_i)_0=A_{\l\r}\cap (\vec{c}_i)_0$. Since $A\cap\rho^i_0\neq A_{\l\r}\cap\rho^i_0$, it follows that $(\vec{c}_i)_0<\rho^i_0$. 
%Assume that $(\vec{c}_i)_j<\rho^i_j$ for every $j\leq k$. Since $A^i_{\l (\vec{c}_i)_0,...,(\vec{c}_i)_k\r}\cap(\vec{c}_i)_{k+1}=A\cap (\vec{c}_i)_{k+1}$,
%then $$(\vec{c}_i)_{k+1}<\min(A^i_{\l (\vec{c}_i)_0,...,(\vec{c}_i)_k\r}\Delta A)\leq \rho^i_{k+1}$$

By $\sigma$-completeness, there is $p_{\omega}$ such that $p_n\leq^* p_\omega$. By density, there is such $p_{\omega}\in G$.

In $V[A]$ we have the sequence $\l \eta_n\mid n<\omega\r$.
Clearly, $C_G(\eta_n)\geq \eta_n$
and as we have seen, $C_G(0)<\eta_0$. Let us prove that $\eta_{n+1}>C_G(\eta_n)$.

\begin{claim*}
For every $0<n<\omega$, $\eta_{n}>C_G(\eta_{n-1})$.
\end{claim*}
\textit{Proof of claim.}
Find $\vec{c}_0\in mb(T_0)$,$\vec{c}_1\in mb(T_{1,\vec{c}_0})$,...$\vec{c}_n\in mb(T_{n,\vec{c}_0,...,\vec{c}_{n-1}})$ such that $$p_\omega^{\smallfrown} \l \vec{c}_0,\vec{Y}^{\vec{c}_0}_0\r^{\smallfrown}...^{\smallfrown}\l\vec{c}_n,\vec{Y}^{\vec{c}_n}_n\r\in G$$
It follows that $\vec{c}_n(N_n)=C_G(\eta_{n-1})$ and that $A\cap C_G(\eta_{n-1})=A_{n+1}(\vec{c}_0,...\vec{c}_n)\cap C_G(\eta_{n-1})$.
Since for every $j\leq N_n$, by definition,
$$A_{n+1}(\vec{c}_0,...\vec{c}_{n-1},\vec{c}_n\restriction j)\cap \vec{c}_n(j)=A_{n+1}(\vec{c}_0,...\vec{c}_{n-1},\vec{c}_n\restriction j+1)\cap\vec{c}_n(j)$$
It follows that for every $j\leq N_n$,
$$A_{n+1}(\vec{c}_0,...\vec{c}_{n-1},,\vec{c}_n\restriction j)\cap \vec{c}_n(j)=A\cap\vec{c}_n(j)$$
 In particular, 
$A\cap(\vec{c}_n)(0)=A_{n+1}(\vec{c}_0,...\vec{c}_{n-1},\l\r)\cap (\vec{c}_n)(0)$. 

Let us argue that for every $k\leq N_n$, $\rho^n(k)> \vec{c}_n(k)$

Since by definition $A\cap\rho^n(0)\neq A_{n}(\vec{c}_0,...\vec{c}_{n-1},\l\r)\cap\rho^n(0)$, it follows that $\vec{c}_n(0)<\rho^n(0)$. Inductively assume that $\vec{c}_n(j)<\rho^n(j)$ for every $j\leq k$. Since $A_{n}(\vec{c}_0,...,\vec{c}_{n-1},\vec{c}_n\restriction k+1)\cap\vec{c}_n(k+1)=A\cap \vec{c}_n(k+1)$,
then $$\vec{c}_n(k+1)<\min(A_n(\vec{c}_0,...,\vec{c}_{n-1},\vec{c}_n\restriction k+1)\Delta A)\leq \rho^n(k+1)$$
Hence $C_G(\eta_{n-1})=\vec{c}_n(N_n)<\rho^n(N_n)=\eta_n$.$\blacksquare_{\text{Claim}}$

We conclude that
$$C_G(0)<\eta_0\leq C_G(\eta_0)<\eta_1\leq C_G(\eta_1)...$$
Let $\kappa^*=\sup_{n<\omega}\eta_n$, then $\kappa^*\in Lim(C_G)$ and therefore regular in $V$. Also, by assumption, $cf^{V[G]}(\kappa)>\omega$, hence $\kappa^*<\kappa$. By freshness, $A\cap\kappa^*\in V$.

This means that in $V$ we can construct the sequence $\l \eta_n\mid n<\omega\r$ which is a contradiction. This concludes the proof for sets with supremum $\kappa$.$\blacksquare_{\text{Lemma}}$

%$$A(t_1)=A(t_2)\leftrightarrow t_1\restriction I_i=t_2\restriction I_2$$

%We also make sure that the functions $t\mapsto \xi_t$ where $\xi_t$ is such that $\succ_{T_i}(t)\in U(\kappa,\xi_t)$ have important coordinates, $J_{i,k}$. Note that $J_{i,k}\subseteq I_i$, since if $k\in I_i$, then $J_{i,k}\subseteq I_i$.

%We claim that for every $k\leq n_i$ and every $s_1,s_2\in \Lev_k(T_i)$, if $s_1\restriction I_i=s_2\restriction I_i$ then $A^i_{s_1}=A^_i_{s_2}$.
%By induction of $k$, for $k=n$ this is true. Assume it holds for $k$ and prove for $k-1$, let $s_1,s_2\in \Lev_{k-1}(T_i)$, if $s_1\restriction I_i=s_2\restriction I_i$, then we split into two cases. First, if $k\notin I_i$, then for any $\alpha_1\in \succ_{T_i}(s_1)$ and any $\alpha_1<\alpha_2\in\succ_{T_i}(s_2)$, $s_1^{\smallfrown}\alpha_1\restriction I_i=s_2^{\smallfrown} \alpha_2\restriction I_i$, thus $A^i_{s_1^{\smallfrown}\alpha_1}=A^i_{s_2^{\smallfrown}\alpha_2}$.
%By the definition, if follows that
%$$A^i_{s_1^{\smallfrown}\alpha_1}\cap\alpha_1=A^i_{s_1}\cap\alpha_1\wedge A^i_{s_2^{\smallfrown}\alpha_2}\cap\alpha_2=A^i_{s_2}\cap\alpha_2$$
%since $\alpha_1<\alpha_2$ and it follows that $A^i_{s_1}\cap\alpha_1=A^i_{s_2}\cap\alpha_1$. This is true for unbounded $\alpha_1$ so $A^i_{s_1}=A^i_{s_2}$.

%If $k\in I_i$, then $J_{i,k}\subseteq I_i$, and therefore $\xi_{s_1}=\xi_{s_2}$, which means that $\succ_{T_i}(s_1)\cap\succ_{T_i}(s_2)\in U(\kappa,\xi_{s_1})$. For every $\alpha\in \succ_{T_i}(s_1)\cap\succ_{T_i}(s_2)$, $A^i_{s_1\smallfrown\alpha}=A^i_{s_2^{\smallfrown}\alpha}$ and we conclude that $A^i_{s_1}=A^i_{s_2}$.

Now for the remaining cases of Theorem \ref{Fresh}:
 \begin{lemma}\label{freshabovek}
If $A\in V[G]$ is a fresh set of ordinals with respect to $V$,
such that $\sup(A)>\kappa$, then $cf^{V[G]}(\sup(A))=\omega$.
\end{lemma} 
\pr  Let $\mu:=cf^{V}(\sup(A))$, by Theorem \ref{no fresh subsets} $\mu\leq\kappa$. There is a fresh set $X\subseteq \mu$ such that $V[A]=V[X]$. To see this, pick in $V$ a cofinal sequence $\langle \eta_i\mid i<\mu\r$ in $\sup(A)$. Then by $\kappa^+$-c.c, there is $F\in V$, such that 
\begin{enumerate}
    \item $Dom(F)=\mu$.
    \item For every $i<\mu$, $|F(i)|=\kappa$.
    \item $A\cap\eta_i\in F(i)$.
\end{enumerate} 
For each $i<\mu$, find in $V$, an enumeration $\langle x^i_j\mid j<\kappa\rangle$ of $F(i)$, such that for every $W\in F(i)$, $\{j<\kappa\mid x^i_j=W\}$ is unbounded in $\kappa$.

Move to $V[A]$, inductively define $\l\gamma_i\mid i<\mu\r$ increasing such that $x^i_{\gamma_i}=A\cap\eta_i$.

Set $\gamma_0=\min(j\mid x^0_j=A\cap\eta_0)$.
Assume that $\gamma_i$ was defined for every $i\leq k<\mu$, define $\gamma_{k+1}=\min(j>\gamma_k\mid x^{k+1}_j=A\cap\eta_{k+1})$. 
Note that at limit stage $\delta$, the sequence $\langle\gamma_i\mid i<\delta\rangle$ is definable using only the enumeration and $A\cap\eta_\delta$ which is all available in $V$. hence $\gamma_{\delta}'=\sup(\gamma_i\mid i<\delta)<\kappa$ and we define $\gamma_\delta=\min(j>\gamma_\delta'\mid x^{\delta}_j=A\cap\eta_\delta)$.

Let $X=\{\gamma_i\mid i<\mu\}\subseteq \kappa$.
Since $\langle \gamma_i\mid i<\mu\rangle$ is increasing, $cf^{V[G]}(\sup(X))=cf^{V[G]}(\mu)$, $V[A]=V[X]$ and $X$ is fresh. It follows by the proof for subsets of $\kappa$ that $cf^{V[G]}(\mu)=\omega$, hence $cf^{V[G]}(\sup(A))=\omega$.
$\blacksquare_{\text{Lemma }\ref{freshabovek}}$ $\blacksquare_{\text{Theorem }\ref{Fresh}}$
\section{ Open problems }

Here are some related open problems:

Distinguishing from the case where $o^{\vec{U}}(\kappa)<\kappa$, we do not have here a classification of subforcings of $\Mfor$. 

\begin{question} Classify subforcings of $\Mfor$.
\end{question}

For $o^{\vec{U}}(\kappa)<\kappa^+$, using Theorem \ref{MainResaultParttwo}, it suffices to consider models of the form $V[C']$ for some $C'\subseteq C_G$, and try to classify the forcings which generate these models.

Our conjecture, at least for $o^{\vec{U}}(\kappa)=\kappa$ is the following:
\begin{conjecture}
Let $G\subseteq \Mfor$ be a $V$-generic filter, where  $\forall\alpha\leq \kappa.o^{\vec{U}}(\alpha)\leq\alpha$.
If $V\subseteq M\subseteq V[G]$ is a transitive $ZFC$ model, then either it is a finite iteration of Magidor-like forcings as in \cite{partOne}, or there is a tree $T\subseteq[\kappa]^{<\omega}$ in $V$ such that $ht(T)=\omega$ and for every $t\in T$ and every $\alpha\in \succ_T(t)$, there is a name $\lusim{{\Mfor}^*}_{t^{\smallfrown}\alpha}$ for a Magidor-like forcing,  such that if $H$ is $V$-generic filter for the forcing adding a branch through the tree $T$ along with the forcings $\lusim{{\Mfor}^*}_{t^{\smallfrown}\alpha}$ corresponding to the branch, then $M=V[H]$. 
\end{conjecture}

\begin{question}\label{question3} Suppose that $o^{\vec{U}}(\kappa)=\kappa^+$. Is still every set of ordinals in the extension equivalent to a subsequence of a generic sequence?
\end{question}

Note that the situation here is more involved since $\kappa$ stays regular in $V[G]$ and it is no longer possible to separate the measures.  

\begin{question} The same as \ref{question3}, but with $o^{\vec{U}}(\kappa)>\kappa^+$.\end{question}

\begin{question} What can we say about other Prikry type forcing notions ?\end{question} 

In \cite{TomMoti}, an example of a non-normal ultrafilter is given which adds a Cohen function to $\kappa$. So in general, not every intermediate model of Prikry type extensions is a Prikry type extension.

The following questions were stated in Section $5$:

In attempt to generalize \ref{kappaplusccgeneral} to a wider class of forcings, the simplest would probably be to deal with a long enough Magidor iteration of Prikry forcings and to analyze its subforcings.

\begin{question}
Is the result of Theorem \ref{kappaplusccgeneral} valid for a long enough Magidor iteration of the Prikry forcings?
\end{question}

\begin{question}
Characterize filters (or ultrafilters) which satisfy the Galvin property (or the generalized Galvin property).
\end{question}

\begin{question} Assume $GCH$. Let $\kappa$ be a regular uncountable cardinal.
Is there a $\kappa$-complete filter on $\kappa$ which fails to satisfy the Galvin property?
\end{question}

\begin{question}
Assume $GCH$. Let $\kappa$ be a regular uncountable cardinal. Is there a $\kappa-$complete filter which extends the closed unbounded filter $Cub_\kappa$  and fails to satisfy the Galvin property?
\end{question}

\begin{question}
Is it consistent to have a $\kappa$-complete ultrafilter over $\kappa$ which does not have the Galvin property?
\end{question}

\begin{question}  Is it consistent to have a measurable cardinal $\kappa$ carrying a $\kappa-$complete ultrafilter which extends the closed unbounded filter $Cub_\kappa$ 
(i.e., $Q-$point) and fails to satisfy the Galvin property?
\end{question}

In section $5$ we have seen that a fine $\kappa$-complete ultrafilter over $P_\kappa(\lambda)$ does not satisfy the Galvin property. Indeed, if $U$ is a fine normal measure on $P_\kappa(\lambda)$ then supercompact Prikry forcing is not $\kappa^+$-cc, however, under $GCH$ this forcing is $\lambda^+$-cc 

\begin{question}
Assume $GCH$ and let $\lambda>\kappa$ be a regular cardinal. Is every quotient forcing of the supercompact Prikry forcing also $\lambda^+$-cc in the generic extension?
\end{question}
\subsection{Acknowledgements}
The authors would like to thank the referee for his careful examination of the paper. Also they would like to thank the participants of the Tel-Aviv University Set Theory Seminar and CUNY Set Theory Seminar, especially to Menachem Magidor and Gunther Fuchs for their interesting insights and observations. 

\bibliographystyle{amsplain}
\bibliography{ref}
\end{document}